\documentclass{amsart}

\usepackage{amssymb,amsmath,amscd,amsthm,xspace,color,tocvsec2, array, 
	tikz-cd, enumitem}
    \usepackage{comment}
\usepackage[all,cmtip]{xy}
\usepackage[mathscr]{euscript}
\usepackage[breaklinks=true, hidelinks]{hyperref} 
\usepackage[ampersand]{easylist}
\usepackage{stmaryrd}
\usepackage[myheadings]{fullpage}
\usepackage{csquotes}
\usepackage{mathtools}
\usepackage[skip=5pt, indent=10pt]{parskip}

\newtheorem{Thm}{Theorem}[subsection]
\newtheorem{Prop}[Thm]{Proposition}
\newtheorem{Cor}[Thm]{Corollary}
\newtheorem{Lem}[Thm]{Lemma}
\theoremstyle{definition}
\newtheorem{Ex}[Thm]{Example}
\newtheorem{defi}[Thm]{Definition}
\newtheorem{Rem}[Thm]{Remark}

\newcommand{\fl}{\mathfrak{l}}
\newcommand{\uHom}{\underline{\on{Hom}}}
\newcommand{\comod}{\text{-Comod}}
\newcommand{\fh}{\mathfrak{h}}
\newcommand{\fm}{\mathfrak{m}}
\newcommand{\fb}{\mathfrak{b}}
\newcommand{\halpha}{\check{\alpha}}
\newcommand{\fk}{\mathfrak{k}}

\newcommand{\RepF}{\Rep(\mathsf{F})}
\newcommand{\sY}{\mathscr{Y}}
\newcommand{\diglet}{\mathrm{F}}

\newcommand{\OOF}{\OO_{\mathsf{F}}}
\newcommand{\OOG}{\OO_{\mathsf{G}}}
\newcommand{\RepG}{\Rep(\mathsf{G})}
\newcommand{\Rep}{\on{Rep}}

\newcommand{\slf}{\mathfrak{sl}}

\newcommand{\lmod}{\on{-mod}}

\newcommand{\Oblv}{\on{Oblv}}

\newcommand{\coker}{\on{coker}\hspace{.5mm}}
\newcommand{\gl}{\mathfrak{gl}}

\newcommand{\fa}{\mathfrak{a}}

\newcommand{\Fun}{\Ring}
\newcommand{\Promod}{\on{Pro}(\Ring\mod)}
\newcommand{\sC}{\mathscr{C}}
\newcommand{\Proc}{\on{Pro}(\sC)}
\newcommand{\Pro}{\on{Pro}}
\newcommand{\sss}{\subsubsection{}}
\newcommand{\fg}{\mathfrak{g}}
\newcommand{\Spec}{\on{Spec}\hspace{.5mm}}
\newcommand{\ft}{\mathfrak{t}}

\newcommand{\fr}{^\on{(1)}}
\newcommand{\R}{\Ring}
\newcommand{\D}{\sD}
\newcommand{\sD}{\mathscr{D}}
\newcommand{\Dh}{\mathscr{D}^\hbar}
\newcommand{\sF}{\mathsf{F}}
\newcommand{\Sing}{\mathsf{S}}
\newcommand{\Ting}{\mathsf{T}}

\newcommand{\QCoh}{\on{QCoh}}
\newcommand{\sM}{\mathscr{M}}
\newcommand{\OO}{\mathscr{O}}

\newcommand{\sE}{\mathscr{E}}
\newcommand{\sT}{\mathscr{T}}

\newcommand{\Hom}{\on{Hom}}
\newcommand{\sN}{\mathscr{N}}

\renewcommand{\mod}{\text{-}\on{Mod}}
\newcommand{\Sym}{\on{Sym}\hspace{.5mm}}
\newcommand{\on}{\mathrm}
\renewcommand{\and}{\quad \on{and} \quad}
\newcommand{\Ring}{\mathsf{R}}
\newcommand{\g}{\mathfrak{g}}
\newcommand{\hg}{\hat{\mathfrak{g}}}

\newcommand{\gr}{\operatorname{gr}}
\newcommand{\h}{\mathfrak{h}}
\newcommand{\hh}{\hat{\mathfrak{h}}}
\newcommand{\ffr}{\mathrm{ffr}}
\newcommand{\sP}{\mathscr{P}}
\newcommand{\HC}{\mathrm{HC}} 
\newcommand{\F}{\mathbb{F}}
\newcommand{\Z}{\mathbb{Z}}
\newcommand{\im}{\operatorname{im}}
\newcommand{\Wak}{\on{Wak}}
\newcommand{\lf}{\mathfrak{l}}
\newcommand{\A}{\mathscr{A}}
\newcommand{\Jet}{\mathscr{J}}
\newcommand{\Loop}{\mathscr{L}}
\newcommand{\K}{\mathbb{K}}
\newcommand{\ff}{\mathfrak{f}}
\newcommand{\frp}{\mathfrak{p}}
\newcommand{\AS}{\mathsf{AS}}
\newcommand{\CDO}{\mathrm{CDO}}
\newcommand{\Am}{\Aut_-(\pD)}

\newcommand{\pD}{\sD^\times}
\newcommand{\frz}{\mathfrak{z}}
\newcommand{\bone}{\mathbf{1}}

\newcommand{\hl}{\hat{\mathfrak{l}}}
\newcommand{\bk}{\mathbf{k}}
\newcommand{\Iw}{\mathsf{Iw}}
\newcommand{\WakDelta}{\mathsf{W}\Delta}
\newcommand{\sx}{\mathsf{x}}
\newcommand{\End}{\operatorname{End}}
\newcommand{\Ffr}{\mathrm{{Ffr}}}

\newcommand{\wotimes}{\widehat{\otimes}_{\bk[\hbar]}}
\newcommand{\cotimes}{\hspace{.5mm}\widetilde{\otimes}\hspace{.5mm}}

\newcommand{\ckappa}{\check{\kappa}}
\newcommand{\cft}{\check{\fh}}

\newcommand{\Conn}{\on{Conn}}
\newcommand{\Op}{\on{Op}}
\newcommand{\wsB}{\widetilde{\mathscr{B}}}
\newcommand{\sB}{\mathscr{B}}
\newcommand{\Aut}{\on{Aut}}

\newcommand{\ASk}{\ckappa^p - \ckappa}
\newcommand{\cM}{\Pi}
\newcommand{\fZ}{\mathfrak{Z}}
\newcommand{\bK}{\mathbb{K}}
\renewcommand{\bK}{\mathbb{K}} 

\newcommand{\sO}{\mathscr{O}}

\numberwithin{equation}{section}
\title{Affine Harish-Chandra center in positive characteristic}
\author{Gurbir Dhillon and Ivan Losev}

\begin{document}
\begin{abstract} { Let $G$ be a split reductive group defined over a field of characteristic bigger than the Coxeter numbers of its simple factors. We identify the Harish--Chandra centers of the associated Kac--Moody vertex algebras and enveloping algebras at noncritical levels  with algebras of functions on moduli spaces of connections for the Langlands dual group $\check{G}$.

Namely, given a noncritical level $\kappa$ for $G$,  consider the associated Kac--Moody vertex algebra $V_\kappa(\fg)$ defined on a formal disc $\sD$, and its Harish--Chandra center of arc group invariants $$V_\kappa(\fg)^{\Jet G} \subset V_\kappa(\fg).$$
Write $\ckappa$ for the dual level for $\check{G}$, $\ckappa^p - \ckappa$ for its image under the Artin--Schreier map, and $\Op_{\check{G}}(\sD\fr)_{\ASk}$ for the moduli space of $(\ckappa^p - \ckappa)$-opers on the Frobenius twisted disc $\sD\fr$. We establish a canonical isomorphism 
$$\Spec V_\kappa(\fg)^{\Jet G} \simeq \Op_{\check{G}}(\sD\fr)_{\ASk}.$$
Similarly, we canonically identify the Harish--Chandra center of the filtered complete enveloping algebra at level $\kappa$ with the filtered complete algebra of functions on the moduli space of $(\ASk)$-opers on the Frobenius twisted punctured disc ${\pD}\fr$:
$$\Spec \widehat{U}_{\kappa}(\hat{\fg})^{\Loop G} \simeq \Op_{\check{G}}({\pD}\fr)_{\ASk}.$$

The presence of these Harish--Chandra centers is  a genuinely new phenomenon for loop groups at noncritical level in positive characteristic: these centers  are nontrivial, unlike for loop groups at noncritical level in characteristic zero, and in particular are not the reductions mod $p$ of the characteristic zero Harish--Chandra centers, unlike for  finite dimensional reductive groups or loop groups at critical level. 
}
\end{abstract}
\maketitle
\setcounter{tocdepth}{1}
 \tableofcontents

\section{Introduction}

\subsection{Overview}

\sss {In this paper we determine the Harish--Chandra center of an affine Kac--Moody algebra in positive characteristic, and show that its spectrum is canonically identified with a moduli space of opers for the Langlands dual group. After providing some contextual discussion, we state our main result in Theorem \ref{t:thmintro}.  
}

\sss \label{sss:intro1} Let $F$ be an algebraic group over a field $\bK$ with Lie algebra $\ff$. Let $U(\ff)$ denote the universal enveloping algebra
of $\ff$. For many basic representation theoretic questions about $\ff$ one needs to understand the center $Z(U(\ff))$ of  $U(\ff)$. Note that if $F$ is connected, and $\bK$ has characteristic $0$, then the center $Z(U(\ff))$ coincides with the subalgebra $U(\ff)^F$ of $F$-invariants in $U(\ff)^F$. 

Assume now that the characteristic of $\bK$ is $p>0$. Then  $U(\ff)^F$ is still a subalgebra of $Z(U(\ff))$, usually called the {\it Harish-Chandra center} (shortly, HC center).
One of the reasons why this central subalgerba is important is that it acts by endomorphisms on modules over the Harish-Chandra pair 
$(\ff, \underline{F})$ for every algebraic subgroup $\underline{F}\subset F$. 

The subalgebra $U(\ff)^F\subset Z(U(\ff))$ can be shown to be proper if $F$ is not abelian. More precisely, there is a subalgebra in $Z(U(\ff))$ called the {\it p-center} or the {\it Frobenius center}. Denote it by $Z_{Fr}(U(\ff))$ (with ``Fr'' from Frobenius). It is a polynomial algebra in $\dim F$-variables such that $U(\ff)$ is a free module over $Z_{Fr}(U(\ff))$ of rank $p^{\dim F}$. The subalgebra $Z_{Fr}(U(\ff))$ is $F$-stable, but the action of $F$ on it is nontrivial if $F$ is non-abelian, hence the p-center cannot be contained in the HC center. We will elaborate on the p-centers in the main body of the paper. 

\sss
For a connected reductive group one can determine the HC center. In this case, we will denote the group by $G$ instead of $F$. Let $h$ denote the maximum of the Coxeter numbers for the simple factors of the Lie algebra $\g$. To simplify the exposition we assume that $p>h$ (one can relax this restriction, see, e.g., \cite[Sec. 6.4]{Jantzen_lectures} or \cite[Sec. 2.2]{Riche_Kostant}). Let $\fh\subset \g$ be a Cartan subalgebra, 
and $W$ be the Weyl group. Consider the $\rho$-shifted action 
$(w,x)\mapsto w\cdot x$ of $W$ on $\fh^*$. Then we have the Harish-Chandra isomorphism: $U(\g)^G\xrightarrow{\sim}\bK[\fh^*]^{(W,\cdot)}$. This explains the name ``Harish-Chandra center''. Below we will need the following notation: let $r=\dim \fh$ and let $d_1,\ldots,d_r$ be the degrees of free homogeneous generators of $\F[\h^*]^W$, for example, for $G=\operatorname{SL}_n$, we have $d_i=i+1$ for $i=1,\ldots,n-1$.

In fact, one can also describe the entire center $Z(U(\g))$ in this case. Namely, by a theorem of Veldkamp, \cite{Veldkamp}, 
$$Z(U(\g))=U(\g)^G\otimes_{[Z_{Fr}U(\g)]^G}Z_{Fr}(U(\g)).$$ 



\sss  
In this paper we want to determine the Harish-Chandra center for \underline{affine} Kac--Moody algebras in characteristic $p$. Let us recall necessary definitions. Let $G$ be a connected reductive group.
We assume that $\operatorname{char}\bK$ is $0$ or is bigger than $h$. To an element $\kappa\in \Sym^2(\g^*)^G$ one can assign a
central extension $\hat{\g}_\kappa$ of $\g\otimes \bK(\!(t)\!)$ by $\bK$.
We will use the notation $\mathbf{1}$ for the unit element in this central $\bK$.
Our convention is that the corresponding cocycle corresponds not to $\kappa$ itself but to $\kappa-\frac{1}{2}\kappa_\g$, where we write $\kappa_\g$ for the Killing form. We refer to $\kappa$ as the (shifted) level, $\kappa=0$ is referred to as the {\it critical level}. The loop group $\Loop G$ acts on $\hat{\g}_\kappa$ and hence on its universal enveloping algebra. It is customary to replace $U(\hat{\g}_\kappa)$ with its quotient $U_\kappa(\hg):=U(\hg_\kappa)/(\mathbf{1}-1)$, where $1$ is the unit in $U(\hg_\kappa)$.

\sss\label{SSS_FF_Thm_intro}
However, the algebra $U_\kappa(\hg)^{\Loop G}$ is not expected to be particularly interesting -- it is likely just the scalars, cf. 
\cite{YZ} for a related result. To get something interesting one replaces $U_\kappa(\hg)$ with its completion $\widetilde{U}_\kappa(\hg)$, a topological algebra whose category of modules with discrete topology consists precisely of smooth representations of $U_\kappa(\hg)$. For this completion the description of the HC center is classically known when $\operatorname{char}\bK=0$ (and the HC center coincides with the entire center, \cite[Proposition 4.3.8]{Frenkel_loop}). 
For now we state a very rough result, and give a more precise description later. To simplify the statement assume $G$ is simple. 
Then if $\kappa\neq 0$, then the center consists of scalars. 
A famous theorem of Feigin and Frenkel, \cite{FF91}, asserts that for $\kappa=0$ the center is a suitable completion of the algebra of polynomials in infinitely many variables.  More precisely, we have $r$ (where $r=\operatorname{rk}\g$) families of generators, $S_{i,j}$, where $j\in \Z$. For now, let us just mention that the PBW degree of $S_{i,j}$ equals $d_i$ for all $j$. 

\sss The argument of the proof of the Feigin-Frenkel theorem is based on studying the center of a related and somewhat easier algebra: the affine vertex algebra $V_\kappa(\g)$. As a space
(and a $\hat{\g}_\kappa$-module) it is 
$$\operatorname{Ind}_{\g[[t]]\oplus \bK \mathbf{1}}^{\hg_\kappa}\bK,$$
where $\g[[t]]$ acts on $\bK$ by zero and $\mathbf{1}\in \bK$ acts
by $1$. This space carries an additional structure: that of a vertex algebra so one can talk about the center, which, in fact, coincides with the $\g[[t]]$-invariants. Therefore, it makes sense to talk about the Harish-Chandra center: the arc group $\Jet G$ acts on $V_\kappa(\g)$ and we can consider the invariants for this action. This is a commutative vertex subalgebra of $V_\kappa(\g)$. A version of the Feigin-Frenkel theorem on the level of vertex algebras 
can be stated as follows: for $\kappa=0$, the center of $V_\kappa(\g)$ is the polynomial algebra in the elements $S_{i,j}$ for $j\leqslant -d_i$ (no completion is needed). It turns out that one can deduce the description of $Z(\widetilde{U}_\kappa(\hg))$
from here, see, e.g., \cite[Sec. 4.3]{Frenkel_loop}. 

\sss The descriptions of the centers for $\widetilde{U}_0(\hg)$ and $V_0(\g)$ can be made much more precise and, in particular, coordinate-free. It was shown by Feigin and Frenkel that both are identified with the algebras of polynomial functions on spaces of {\it opers} for the Langlands dual group $\check{G}$ on suitable 1-dimensional schemes $X$ over $\bK$. An oper is a connection on a principal bundle over a curve with some additional structure. 
The center of $V_0(\g)$ (resp., $\widetilde{U}_0(\hg)$) is the algebra of functions on the $\check{G}$-opers on $\D:=\operatorname{Spec}(\bK[[t]])$ (resp., $\D^\times:=\operatorname{Spec}(\bK(\!(t)\!))$).

\sss We now assume that $\operatorname{char}\bK=p>0$.
The center of $V_0(\g)$ was studied in \cite{ATV}. In 
\cite[Theorem 1.3]{ATV} the authors proved that a direct analog of the Feigin-Frenkel theorem for the HC center holds (in the case when $\g$ is exceptional, they had to require that $p$ is very large, in fact, the techniques of the present paper allow to remove that restriction). Moreover, they proved that a direct analog of the Veldkamp's theorem holds for the entire center of $V_0(\g)$, roughly, the center is generated by the HC center and the p-center (the latter still makes sense
in the affine setting, see \cite{AW} or Section \ref{S_ffr_p_cent} below). 

\sss \label{sss:introlast} The focus of this paper is that case when $\kappa$ is non-degenerate. Here we have several new phenomena (very different from both the characteristic $p$ story for semisimple Lie algebras and the characteristic $0$ 
story for affine Lie algebras). We continue to assume that $p$ is bigger than the Coxeter number of any simple factor of $\g$. We show that the HC center of $V_\kappa(\g)$ is the polynomial algebra in the elements $\underline{S}_{i,j}$ with $j\leqslant -d_i$, where the PBW degree of $\underline{S}_{i,j}$ is now $pd_i$. Informally, this center is ``$p$ times smaller'' than the Feigin-Frenkel center. The situation with the completed universal enveloping algebra is more complicated. Note that $\widetilde{U}_\kappa(\hg)$ carries the PBW filtration but the filtration is not exhaustive. Let $\widehat{U}_\kappa(\hg)$ denote the union of filtered pieces. We show that $\widehat{U}_\kappa(\hg)^{\Loop G}$
is a suitable completion of the polynomial algebra 
in the variables $\underline{S}_{i,j}$ now for $j\in \Z$. We expect that  $\widehat{U}_\kappa(\hg)^{\Loop G}$ is dense (with respect to the inverse limit topology) in 
$\widetilde{U}_\kappa(\hg)^{\Loop G}$.

 We still expect an analog of the Veldkamp's theorem to hold for non-integral $\kappa$.

{ 
\sss Let us now state our results more precisely, and in particular coordinate invariantly. To do so, we first review some relevant concepts and notation. 

\sss Let us suppose our Kac--Moody level $\kappa \in \Sym^2(\fg^*)^G$ is nondegenerate, and let us write $\check{\fg}$ for the Lie algebra of $\check{G}$. Then to $\kappa$ we may canonically associate a dual level   $\ckappa \in \Sym^2(\check{\fg}^*)^{\check{G}}$, defined as follows. 

Note that if we write $\ft$ for the abstract Cartan of $\fg$, i.e., the quotient of any Borel subalgebra by its nilpotent radical, $\kappa$ induces a nondegenerate invariant bilinear form $\kappa_\ft$ on $\ft$. We next recall that the abstract Cartan $\check{\ft}$ of $\check{\fg}$ is simply the dual vector space of $\ft$, i.e., $$\check{\ft} \simeq \ft^*.$$With this, $\ckappa$ is characterized by the property that $\kappa_\ft$ and $\ckappa_{\check{\ft}}$ are dual bilinear forms in the sense of linear algebra. 

\sss We next note that $\Sym^2(\check{\fg}^*)^{\check{G}}$ has a naturally defined $\mathbb{F}_p$-form, as $\check{\fg}$ does. In particular, viewing $\Sym^2(\check{\fg}^*)^{\check{G}}$ as a vector group over $\K$ we have a canonical Artin--Schreier map
$$\on{AS}: \Sym^2(\check{\fg}^*)^{\check{G}} \rightarrow \Sym^2(\check{\fg}^*)^{\check{G}}, \quad \quad \ckappa \mapsto \ckappa^p - \ckappa.$$

\sss Finally, if we write $\sD\fr$ for the Frobenius twist of $\sD$, and similarly ${\pD}\fr$ for the Frobenius twist of $\pD$, we may form the moduli spaces of $(\ckappa^p - \ckappa)$-opers on $\sD\fr$ and ${\pD}\fr$, which we denote respectively by 
$$\Op_{\check{G}}(\sD\fr)_{\ckappa^p - \ckappa} \quad \text{and} \quad \Op_{\check{G}}({\pD}\fr)_{\ckappa^p - \ckappa};$$
we refer the reader to Sections \ref{ss:kappaopers} and \ref{sss:operspd} for the precise definition of these moduli spaces.

\sss Let us write $\Aut(\sD)$ for the automorphism group ind-scheme of $\sD$, and similarly for $\Aut(\pD)$. Our main theorem then reads as follows

\begin{Thm} \label{t:thmintro} There are canonical $\Aut(\sD)$ and $\Aut(\pD)$-equivariant isomorphisms 
$$\Spec V_\kappa(\fg)^{\Jet G} \simeq \Op_{\check{G}}(\sD\fr)_{\ckappa^p - \ckappa} \quad \text{and} \quad \Spec \widehat{U}_\kappa(\hat{\fg})^{\Loop G} \simeq \Op_{\check{G}}({\pD}\fr)_{\ckappa^p - \ckappa}.$$
That is, $V_\kappa(\fg)^{\Jet G}$ is canonically identified with the algebra of polynomial functions on the space of $(\ckappa^p - \ckappa)$-opers on $\sD\fr$, and similarly $\widehat{U}_\kappa(\hat{\fg})^{\Loop G}$ is canonically identified with the filtered complete algebra of functions on the space of $(\ckappa^p - \ckappa)$-opers on ${\pD}\fr$. 
\end{Thm}

In particular, the natural action of $\Aut(\sD)$ on $V_\kappa(\fg)^{\Jet G}$ factors through the relative Frobenius homomorphism 
$$\on{Fr}: \Aut(\sD) \rightarrow \Aut(\sD)\fr \simeq \Aut(\sD\fr),$$and similarly for the action of $\Aut(\pD)$ on $\widehat{U}_\kappa(\hat{\fg})^{\Loop G}$. 

}

\sss We now briefly discuss motivations for our work. 
The Feigin-Frenkel theorem is one of the cornerstones of the categorical geometric Langlands program. Similarly, we expect that our description of the HC center will play an important role in the modular quantum categorical geometric Langlands program, which is currently in its infancy. The second, somewhat related, motivation comes from a more classical representation theoretic problem: studying a category $\sO$. Roughly, this is the category of modules over the Harish-Chandra pair $(\hg, \mathsf{Iw})$, where $\mathsf{Iw}$ 
is the Iwahori subgroup in the Kac-Moody group corresponding to $\hg$. One reason to care about this category is that it is supposed to be governed by a version of the double affine Hecke category (somewhat imprecisely, suitable blocks of the category $\sO$ categorify the polynomial representation of the double affine Hecke algebra).

\subsection{Harish-Chandra isomorphism, revisited}\label{SS_HC_iso_intro}
To explain our approach (and also a proof of the Feigin-Frenkel theorem) we start with sketching a proof of the Harish-Chandra isomorphism $U(\g)^G\xrightarrow{\sim} \bK[\h^*]^{(W,\cdot)}$. It is quite similar to proofs typically found in basic Lie theory textbooks, but deviates in some aspects (for example, we do not use Verma modules or category $\sO$). 

In this section we assume that $\operatorname{char}\bK=0$ or is bigger than $h$, the maximum of the Coxeter numbers of the simple factors of $\g$. We warn the reader that some of the steps outlined below are nearly obvious but we emphasize them because the proof of our main result essentially follows the pattern we outline and the analogs of these steps in our setting are very far from being obvious. 

\sss\label{SSS_HC_iso_map}
Choose a parabolic subgroup $P\subset G$ with Levi decomposition $P=L\ltimes N$. Let $N^-$ denote the unipotent radical of the opposite parabolic $P^-$. 
This yields the triangular decomposition $\g=\mathfrak{n}^-\oplus \mathfrak{l}\oplus \mathfrak{n}$. We notice that $[U(\g)/U(\g)\mathfrak{n}]^{Z(L)}\xrightarrow{\sim}U(\fl)$, hence $[U(\g)/U(\g)\mathfrak{n}]^P\xrightarrow{\sim} U(\lf)^L$.
Then we have the following algebra homomorphism
$$U(\g)^G\hookrightarrow U(\g)^P\rightarrow [U(\g)/U(\g)\mathfrak{n}]^P\xrightarrow{\sim} U(\lf)^L$$
that we denote by $\HC_L$. Now choose a Borel subgroup $B$
with Levi decomposition $H\ltimes \tilde{N}$, where $H\subset L$
and $N\subset \tilde{N}$. Then we also have the maps $\HC_H:U(\g)^G\rightarrow U(\h)^H$ and $\HC^L_H: U(\lf)^L\rightarrow U(\h)$. Tracking the construction one arrives at the transitivity property
\begin{equation}\label{eq:HC_transitivity_baby}
\HC_H=\HC^L_H\circ \HC_L.
\end{equation}

\sss\label{SSS_HC_target} The target for $\HC_H$ is very easy to describe: $H$ acts trivially on $U(\h)$, so $U(\h)^H=\bK[\h^*]$.

\sss\label{SSS_HC_graded} 
The homomorphism $\HC_L$ preserves the PBW filtrations. 
The inclusions $\gr[U(\g)^G]\hookrightarrow S(\g)^G$ and 
$\gr[U(\lf)^L]\hookrightarrow S(\lf)^L$ intertwine 
the associated graded map $\gr \mathsf{HC}_L$ with the Chevalley restriction map $\bK[\g^*]^G\rightarrow \bK[\lf^*]^L$. The latter is injective, in particular, $\HC_L$ is injective. 

Thanks to the injectivity,
$U(\g)^G$ is identified with its image under $\HC_H$. Now let $P_1,\ldots, P_r$ be the minimal parabolics containing  $B$, and let $L_1,\ldots,L_r$ be their Levi subgroups containing  $H$. Thanks to  
(\ref{eq:HC_transitivity_baby}), we have 
\begin{equation}\label{eq:HC_containment_basic}
\operatorname{im}\HC_H\subset \bigcap_{i=1}^r \operatorname{im}\HC^{L_i}_{H}.
\end{equation}

\sss\label{SSS_HC_sl2}
The determination of $\operatorname{im}\HC^{L_i}_{H}$ reduces to the case of $\g=\mathfrak{sl}_2$. This case is easy and one sees that 
$\operatorname{im}\HC^{L_i}_H=\bK[\h^*]^{(s_i,\cdot)}$, where $s_i$
is the simple reflection corresponding to $L_i$.

\sss\label{SSS_HC_intersection_baby}
Clearly, $\bigcap_{i=1}^r \bK[\h^*]^{(s_i,\cdot)}=\bK[\h^*]^{(W,\cdot)}$ -- the simple reflections $s_i$ generate the group $W$. So, 
(\ref{eq:HC_containment_basic}) implies $\operatorname{im}\HC_H\subset \bK[\h^*]^{(W,\cdot)}$.

\sss\label{SSS_HC_big_baby}
It just remains to observe that $U(\g)^G$ is big enough: 
$\gr U(\g)^G\xrightarrow{\sim} S(\g)^G$ thanks to the symmetrization map (in positive characteristic we symmetrize free homogeneous generators of $S(\g)^G$, they sit in degrees less than $p$
thanks to  $p>h$). We note that a suitable version of the symmetrization map exists under milder restrictions on $p$, cf. \cite{Friedlander_Parshall}. We also thank Lewis Topley for helpful comments on the draft of this paper.

\subsection{Sketch of proof for $\kappa=0$}
In this section we will explain steps of a proof of the Feigin-Frenkel theorem, \S\ref{SSS_FF_Thm_intro}, and its positive characteristic version. Our proof works both in zero and positive characteristic (bigger than $h$) and deviates from the Feigin-Frenkel proof 
in several aspects, while still using many of the ideas of that proof. We follow the steps outlined in Section 
\ref{SS_HC_iso_intro} and work with vertex algebras instead of the completed universal enveloping algebras.
In the end of the section we explain how to pass from the former to the latter.

\sss\label{SSS_HC_vertex} We have vertex algebra homomorphisms
$\HC_{0,L}:V_0(\g)^{\Jet G}\rightarrow V_0(\lf)^{\Jet L}$, where, recall, $\Jet G$ denotes the arc group of $G$. The construction is, however, significantly more complicated than in the semisimple case and is based on the so called {\it free field realization} map, see 
\cite[Section 5]{Frenkel_loop} for the characteristic $0$ case and \cite{AW} for the characteristic $p$ case. We will revisit this construction in Section \ref{S_CDO_ffr} using chiral differential operators.
With some work one can show that the direct analog of 
(\ref{eq:HC_transitivity_baby}) holds. 

\sss The target of $\HC_{0,H}$ is still easy to determine as the action of $\Jet H$ on $V_0(\h)$ is trivial. This commutative algebra is identified with 
the $\bK[\Jet \h^*]$, the algebra of functions on the arc space $\Jet\h^*$ of $\h^*$.

\sss\label{SSS_HC_PBW_injectivity} The homomorphism $\HC_H$ preserves the PBW filtrations. Similarly to \S\ref{SSS_HC_graded},
$\HC_{0,H}$ is injective. Further,  a direct analog of (\ref{eq:HC_containment_basic})
holds with the same argument.

\sss Again, the determination of $\operatorname{im}\HC_{0,H}^{L_i}$ boils down to the case $\g=\slf_2$. Here $V_0(\g)^{\Jet G}$ is the polynomial algebra in the Sugawara modes $S_j, j\leqslant -2,$ and it is not that difficult to compute their images under $\HC_{0,H}^{L_i}$.

\sss\label{SSS:opers_image_intro} Now our task is to determine $\bigcap_{i=1}^r\operatorname{im}\HC_{0,H}^{L_i}$ getting an affine analog of \S\ref{SSS_HC_intersection_baby}.
The desired answer is that this is the image of the algebra $\bK[\Op_{\check{G}}(\D)]$ (i.e., the algebra of regular functions on the scheme of $\check{G}$-opers on the disc $\D$) under the Miura map. Proving this claim takes a significant portion of \cite[Secs. 7,8.2]{Frenkel_loop} and uses several ingredients that are either quite complicated and do not carry easily to positive characteristic (screening operators) or just fail in characteristic $p$ (the claim that the Miura map is the pullback under the projection from a principal $\check{N}$-bundle over $\Op_{\check{G}}(\D)$). Fortunately, we found a  more elementary proof identifying the intersection with the image of the Miura map.   

\sss\label{SSS:center_affine_big_intro}
The remaining step in our plan is to show that 
$V_0(\g)^{\Jet G}$ is big enough so that the embedding $V_0(\g)^{\Jet G}\hookrightarrow \bK[\Op_{\check{G}}(\D)]$ is an isomorphism. The argument in the case of semisimple Lie algebras involved the symmetrization map, it completely fails in the affine setting as the symmetrization map that does not make sense in the affine setting: we have the symmetrization map $S(\hat{\g})\rightarrow U(\hat{\g})$ but it neither induces a map on the level of vertex algebras nor extends to the completions of interest. Instead, the proof that the center is big enough in the characteristic $0$ case requires an identification of the Verma $\hat{\g}_0$-module with zero highest weight
and a suitable version of the Wakimoto module, see 
\cite[Secs. 6.3, 8.1]{Frenkel_loop}. We do not expect that this identification continues to hold in characteristic $p$. Instead, to show that the center is sufficiently large, we can use an easier isomorphism between Verma and Wakimoto modules with sufficiently generic highest weights.  

\sss\label{SSS_from_V_to_U_intro} The isomorphism $V_0(\g)^{\Jet G}\hookrightarrow \bK[\Op_{\check{G}}(\D)]$ then yields 
a homomorphism $$\bK[\Op_{\check{G}}(\D^\times)]\rightarrow \widetilde{U}_0(\hg)^{\Loop G}.$$ That this is an isomorphism requires an argument involving the associated graded construction, \cite[Sec. 4.3]{Frenkel_loop}. A crucial ingredient is the fact that $\gr[V_0(\g)^{\Jet G}]=\bK[\Jet\g^*]^{\Jet G}$.

\subsection{Approach and results of this paper} \label{s:maintheorem}
Now suppose $\operatorname{char}\bK=p>h$ and $\kappa$ is nondegenerate. Our goal here is to sketch of description of $V_\kappa(\g)^{\Jet G}$. 

\sss Completely analogously to \S\ref{SSS_HC_vertex}  we  get vertex algebra homomorphisms
$\HC_{\kappa,L}:V_\kappa(\g)^{\Jet G}\rightarrow V_\kappa(\lf)^{\Jet L}$. They satisfy the analog of the transitivity property, (\ref{eq:HC_transitivity_baby}).

\sss The first new challenge comes when determining 
$V_\kappa(\h)^{\Jet H}$: the action of $\Jet H$ on $V_\kappa(\h)$ is no longer trivial. We prove that there is an embedding 
$\Psi:\bK[\Jet\h^*]\hookrightarrow \bK[\Jet\h^*]$ such that $V_\kappa(\h)^{\Jet H}$ is identified with its image. The linear generators of $\bK[\Jet\h^*]$ are sent to elements of PBW degree $p$.

\sss Claims in \S\ref{SSS_HC_PBW_injectivity} continue to hold with the same argument. From here we already can handle the case when $\kappa$ is integral. From now on we assume that $\g$ is simple and $\kappa$ is not integral. 

\sss Now we need to determine $V_\kappa(\g)^{\Jet G}$ for $\g=\mathfrak{sl}_2$
and its image under $\HC_{\kappa,H}$. This is significantly more challenging than in the case $\kappa=0$ and requires several new ingredients. What we find is that $V_\kappa(\g)^{\Jet G}$
is the polynomial algebra in certain elements of PBW degree $2p$. However, we  have the following remarkable equality
\begin{equation}\label{eq:images_sl2_coincidence}
\operatorname{im}\HC_{\kappa,H}=\Psi(\operatorname{im}\HC_{0,H}),
\end{equation} an equality of subalgebras in $V_\kappa(\h)^{\Jet H}\xrightarrow{\sim}\Psi(\bK[\Jet \h^*])$. 

\sss Using (\ref{eq:images_sl2_coincidence}) and 
\S\ref{SSS:opers_image_intro} we show that 
$\bigcap_{i=1}^r \operatorname{im}\HC^{L_i}_{\kappa,H}$ coincides with $\Psi(\bK[\Op_{\check{G}}(\D)])$ where we identify the algebra of the functions on the scheme of opers with its image under the Miura map. 

\sss Similarly to \S\ref{SSS:center_affine_big_intro} we get an identification  $V_\kappa(\g)^{\Jet G} \xrightarrow{\sim}\Psi(\bK[\Op_{\check{G}}(\D)])$. We will see that this implies the first isomorphism in Theorem \ref{t:thmintro}.

\sss Similarly to \S\ref{SSS_from_V_to_U_intro} we get a homomorphism 
\begin{equation}\label{eq:intro_homom_HC_center_U}
\Psi(\bK[\Op_{\check{G}}(\D^\times)])\rightarrow \widetilde{U}_\kappa(\hg)^{\Loop G}
\end{equation}
However, since the inclusion $\gr[V_0(\g)^{\Jet G}]\hookrightarrow\bK[\Jet\g^*]^{\Jet G}$ is not surjective in this case, the argument in \cite[Sec. 4.3]{Frenkel_loop} no longer goes through and we do not know how to prove that (\ref{eq:intro_homom_HC_center_U}) is an isomorphism, although we expect that it is (and know how to prove this in the case when $G$ is a torus). We do however prove that the restriction of (\ref{eq:intro_homom_HC_center_U}) to the subalgebras of elements of finite PBW degree is an isomorphism, which is done by studying the analog of $\HC_{\kappa,H}$ for the universal enveloping algebras. This proves the second isomorphism in Theorem \ref{t:thmintro}.


\subsection{Organization of the paper}
\label{s:orgpaper}
Sections \ref{S_group_ind_schemes}-\ref{S_ffr_p_cent} are preparatory: while many results there are new, they are mostly variants of results already in the literature. In Section \ref{S_group_ind_schemes} we discuss basics on group ind-schemes (such as affine Kac-Moody groups and Virasoro groups) and their Lie algebras including restricted structures.In Section \ref{S_vertex} we recall basics of vertex algebras over rings.
In Section \ref{S_CDO_ffr} we revisit the free field realization map and its parabolic analog based on the study of chiral differential operators. Finally, in Section 
\ref{S_ffr_p_cent} we revisit the $p$-center for affine Kac-Moody algebras and the behavior of the free field realization maps on the $p$-centers providing a more general and conceptual version of results from \cite{AW}.

In Section \ref{S_affine_HC} we introduce the Harish-Chandra map $\HC_{\kappa,H}$, show that it is injective, and establish its transitivity property. We then study the Harish-Chandra center in the case of the torus, and use this result to settle the case of general $G$ when $\kappa$ is integral.

In Section \ref{S_sl2} we handle the case when $G=\operatorname{SL}_2$ determining the Harish-Chandra center of $V_\kappa(\mathfrak{sl}_2)$ and its image under $\HC_{\kappa,H}$. The main ingredient is the compatibility with Poisson structures. 

In Section \ref{S_Verma_Wak} we introduce and study Verma and Wakimoto modules as well as modules intermediate between them. Our main result is that certain maps between such modules are isomorphisms for generic enough highest weights. The proof is mostly based on a reduction to the classical level (i.e., taking the associated graded under the PBW filtration).  

In Section \ref{S_Verma_endom} we introduce   a certain algebra, denoted by $\,^I\!\A$, realized as a subalgebra in the polynomials in infinitely many variables. Then we construct an action of $\,^I\!\A$ on a universal Verma module by endomorphisms. The constructions in previous sections give an action on a certain localization of the universal Verma. To show that the original Verma is preserved we utilize some constructions with ``generic'' twisted equivariant modules over reductive groups.

In Section \ref{S_IA_description} we show that $\,^I\!\A$ coincides with the image  of the algebra of functions on the space of opers with regular singularities under the Miura map. Using this description we finally describe the Harish-Chandra center of $V_\kappa(\g)$ in Section \ref{S_V_inv} identifying it with the algebra of functions on the opers on $\D$. This description depends on the coordinate, $t$, on $\D$, and in Section \ref{S_coord_indep} we show that the Harish-Chandra center is canonically identified with the algebra of polynomial functions on the space of opers on the Frobenius twisted disc $\D^{(1)}$ (isomorphic to $\D$ but not in a coordinate-independent way).

\subsection{Acknowledgements} It is a pleasure to thank Pramod Achar, Tomoyuki Arakawa, Alexander Braverman, Roman Bezrukavnikov, Pavel Etingof, Boris Feigin, Edward Frenkel, Dennis Gaitsgory, Ian Grojnowski, Shrawan Kumar,  Sam Raskin, Simon Riche, Lewis Topley,  Weiqiang Wang, Juan Villareal, Zhiwei Yun, and David Yang for helpful conversations and correspondences which have influenced our thinking on this subject.

Gemini was used in working out the proof of Lemma \ref{Lem:KM_cocycle_group}, and Claude was used to correct typos in Section \ref{S_coord_indep}. All edits were checked by the authors and the suggested proof of Lemma \ref{Lem:KM_cocycle_group} was rewritten by the authors.

  G.D. was partially supported by NSF grants Nos.	2103387 and 2502740. I.L. was partially supported by NSF grants No. 2001139 and 2501558.

\section{Group ind-schemes and their restricted Lie algebras}\label{S_group_ind_schemes}

\subsection{Introduction}
\label{sss:introdisc}
 \sss We would like to discuss Lie algebras of group ind-schemes, possibly of infinite type, and their restricted Lie algebra structure in the case when the group is defined over a field of positive characteristic, as well as some basic properties thereof. 

The particular Lie algebras of interest to us are the affine Kac--Moody algebra and the Virasoro algebra. These are {Lie algebras of} central extensions of the loop group and the {group of} automorphisms of a punctured disc, respectively. Both the relevant central extensions of groups are in turn pulled back, via certain natural actions on a Tate vector space $V$,  from the Tate extension $$1 \rightarrow \mathbb{G}_m \rightarrow GL(V)^\flat \rightarrow GL(V) \rightarrow 1.$$It is therefore natural to include in our discussion the case of $GL(V)^\flat$ itself. As $GL(V)^\flat$ is not a group ind-scheme of {\em countable type}, i.e., presentable as a countable filtered colimit of affine schemes along closed embeddings, this inclusion requires a little bit of care. 

\sss The basic idea is then the following. We will see that $GL(V)^\flat$ is in fact an ind-affine ind-scheme, albeit with a presentation of uncountable type. Therefore, in describing its algebra of functions, and coherent sheaves on it, we will need to work with those genuinely as pro-vector spaces; for countable ind-schemes, we recall that one can alternatively work with the perhaps more familiar language of topological vector spaces, complete and separated, admitting a countable basis for their topology. With this in mind, we first recall some basic assertions regarding pro-categories. 

\subsection{Recollections on pro-completions}

\sss Let us review some relevant facts about the pro-completion of a category $\sC$; as the pro-completion of $\sC$ is the opposite of the ind-completion of the opposite category $\sC^{op}$, one possible reference for the following material is Chapters 6 and 8 of \cite{KS}.

\sss Let $\sC$ be a category. Let us denote the pro-completion of $\sC$ by $\Proc$. Explicitly, objects of $\Proc$ are formal cofiltered limits of objects of $\sC$. We will denote such an object by $\varprojlim_\alpha M_\alpha$, where $\alpha$ runs over the essentially small cofiltered diagram category. Morphisms between pro-objects are computed as 
$$\Hom_{\Proc}(\varprojlim_\alpha M_\alpha, \varprojlim_\beta N_\beta) := \varprojlim_{\beta} \varinjlim_{\alpha} \Hom_{\sC}(M_\alpha, N_\beta).$$
In particular, we have a tautological fully faithful embedding $\sC \rightarrow \Proc$, consisting of one-point diagrams. 

It will be convenient in what follows to recall the following standard fact, cf. \cite[Corollary 6.1.14]{KS}. Given objects $M^\wedge$ and $N^\wedge$ in $\Proc$, and a morphism $\phi: M^\wedge \rightarrow N^\wedge$, one can always choose a cofiltered category $I$ and presentations $M^\wedge \simeq \varprojlim M_i, N^\wedge \simeq \varprojlim N_i$, such that the morphism $\phi$ is induced by a system of compatible maps $\phi_i: M_i \rightarrow N_i$, for $i \in I$. 

\sss We next recall that if $\sC$ is abelian, then $\Pro(\sC)$ is again abelian, and the tautological functor $\sC \rightarrow \Proc$ is exact, cf.   \cite[Theorem 8.6.5]{KS}. 

Explicitly, given a map $\phi: M^\wedge \rightarrow N^\wedge$ as above, if we compute kernels and cokernels level-wise, i.e.,
$$0 \rightarrow \ker \phi_i \rightarrow M_i \rightarrow N_i \rightarrow \coker \phi_i \rightarrow 0;$$
then the inverse limits $\varprojlim \ker \phi_i$ and $\varprojlim \coker \phi_i$ satisfy the required universal properties of the kernel and cokernel of $\phi$, respectively, and the first isomorphism theorem in $\sC$ then readily implies it in $\Pro(\sC)$. In particular, any exact sequence $$0 \rightarrow A^\wedge \rightarrow B^\wedge \rightarrow C^\wedge \rightarrow 0$$
may be presented as a formal inverse limit of exact sequences
\begin{equation} \label{e:ses}0 \rightarrow A_i \rightarrow B_i \rightarrow C_i \rightarrow 0,\end{equation}
for some cofiltered index category $I$, cf.   \cite[Proposition 8.6.6(a)]{KS}. 

\sss We next recall that if $\sC$ is (symmetric) monoidal, then $\Pro(\sC)$ inherits a (symmetric) monoidal structure. Explicitly, at the level of binary products, if we denote the product on $\sC$ by $\sC \times \sC \rightarrow \sC,$ $(M, N) \mapsto M \otimes N$, then for $M^\wedge$ and $N^\wedge$ as above, we have that $$M^\wedge \otimes N^\wedge \simeq \varprojlim M_\alpha \otimes N_\beta.$$We recall two basic consequences. First, it follows that if $\sC$ has monoidal unit $\mathbf{1}$, then the image of $\mathbf{1}$ along $\sC \rightarrow \Proc$ is again the monoidal unit.

Second, if $\sC$ is abelian and monoidal, and the functor $- \otimes - $ is (left) right exact in the first (or second) variable, then it follows that the same exactness holds for $\Proc$, e.g., by using \eqref{e:ses}. 

\sss The relevant consequences of the above discussion may then be summarized as follows. 

\begin{Cor}Given a commutative ring $\Ring$, the category $\Promod$ is naturally a symmetric monoidal $\Ring$-linear abelian category, for which the tensor product is right exact in each factor. If $\Ring$ is a field, then the tensor product on $\Promod$ is moreover exact in each factor. 
\end{Cor} 

\subsection{Restricted Lie algebras from symmetric monoidal categories}
\label{ss:restrI}

\sss\label{SSS_restricted_Lie} For the reader's convenience, we recall the definition of a restricted Lie algebra. Let $\Ring$ denote a commutative $\mathbb{F}_p$-algebra. 

\begin{defi}\label{Defi:restricted} A {\em restricted Lie algebra} $\ff$ over $\Ring$ is a Lie algebra $\ff$ $$\ff \wedge \ff \rightarrow \ff, \quad \quad X \wedge Y \mapsto [X,Y] =: \on{ad}_X(Y), \quad \quad X,Y \in \ff,$$
equipped with a further map of sets, the {\em restricted power map},
$$(-)^{[p]}: \ff \rightarrow \ff, \quad \quad X \mapsto X^{[p]},$$
satisfying the following conditions. 

\begin{enumerate}
\item For any $X \in \ff$, one has the equality of $\Ring$-linear endomorphisms of $\ff$
$$\on{ad}_X^p = \on{ad}_{X^{[p]}}.$$
    \item For any $r \in \Ring$, and $X \in \ff$, one has $(r\cdot X)^{[p]} = r^p \cdot X^{[p]}$.
    \item For any $X, Y \in \ff$, we have the equality 
    $$(X+Y)^{[p]} = X^{[p]} + \sum_{i=1}^{p-1} \sigma_i(X,Y) + Y^{[p]},$$
    where for fixed $i$, $1 \leqslant i \leqslant p-1$, the map $i \cdot \sigma_i$ is the coefficient of $t^{i-1}$ in $$(\on{ad}_{t \cdot X + Y})^{p-1}(X) \in \ff \otimes_{\Ring} \Ring[t].$$
\end{enumerate}
\end{defi}

\sss In what follows, we would like to discuss the restricted Lie algebra associated to a group ind-scheme, possibly of infinite type, as well as some basic properties thereof. 

To do so, it suffices to more generally associate a Lie algebra to any Hopf algebra in a symmetric monoidal category $\sC$, with the expected properties, which moreover inherits a restricted structure if $\sC$ is $\F_p$-linear.

\sss Let $\sC$ be an additive monoidal category. In particular,  the functor $- \otimes -: \sC \times \sC \rightarrow \sC$ is bilinear on Hom spaces, i.e., induces maps 
$$\Hom(M_1, N_1) \underset{\mathbb{Z}} \otimes \Hom(M_2, N_2) \rightarrow \Hom(M_1 \otimes N_1, M_2 \otimes N_2).$$

\sss Let $A$ be an algebra in $\sC$, and denote its multiplication map by $$\mu: A \otimes A \rightarrow A.$$ Let us denote by $A\text{-Bimod}$ the category of $A$-bimodules in $\sC$. For an $A$-bimodule $M$, with corresponding action maps $$\alpha_\ell: A \otimes M \rightarrow M, \quad \quad \alpha_r: M \otimes A \rightarrow M,$$we recall that a {\it derivation} $d: A \rightarrow M$ is a map in $\sC$ satisfying the Leibnitz rule
\begin{equation} \label{e:leibrule}d \circ \mu = \alpha_\ell  \circ (\on{id}_A \otimes d)  + \alpha_r \circ (d \otimes \on{id}_A) \in \Hom_{\sC}(A \otimes A, M).\end{equation}

\begin{Rem} \label{r:splitderiv} Given a bimodule $M$, one can form the associated square zero extension $\epsilon: A \oplus M \rightarrow A$. Then it is straightforward to check that a derivation $A \rightarrow M$ is the same data as a map of algebras $\sigma: A \rightarrow A \oplus M$ splitting $\epsilon$. \end{Rem}

\sss Given a bimodule $M$ and derivations $d_1, d_2: A \rightarrow M$, it is straightforward to see from \eqref{e:leibrule} that their sum $d_1 + d_2$ is a derivation, as is $-d_1$. Let us denote by 
$$\on{Der}(A, M) \subset \Hom_{\sC}(A, M)$$
the {subgroup} of derivations. {If $\sC$ is linear over a commutative ring $\Ring$, then $\on{Der}(A, M)$ is an $\Ring$-submodule}. We note in addition that, given a map $\phi: M_1 \rightarrow M_2$ of $A$-modules, the map $\Hom_{\sC}(A, M_1) \rightarrow \Hom_{\sC}(A, M_2)$ restricts to a map 
$$\on{Der}(A, M_1) \rightarrow \on{Der}(A, M_2).$$

\begin{Rem} \label{r:kahlerdiff}Suppose that $A\text{-Bimod}$ admits cokernels. Then the functor of derivations $\on{Der}(A, -)$ is corepresentable, i.e., there exists an $A$-bimodule $\Omega_A$ and a natural isomorphism of functors $$\on{Der}(A, -) \simeq \Hom_{A\text{-bimod}}(\Omega_A, -).$$Namely, $\Omega_A$ can be presented as the cokernel of the alternating cyclic sum 
$$\mu \otimes \on{id} \otimes \on{id} - \on{id} \otimes \mu \otimes \on{id} + \on{id} \otimes \on{id} \otimes \mu: A \otimes A \otimes A \otimes A \rightarrow A \otimes A \otimes A,$$where we view $A^{\otimes 4}$ and $A^{\otimes 3}$ as $A$-bimodules via left multiplication on their leftmost tensor factors, and right multiplication on their rightmost tensor factors.  
\end{Rem}

\sss\label{sss:restrassoc} Recall that, given an associative algebra $E$, we may view it as a Lie algebra via the commutator $[e_1, e_2] = e_1 e_2 - e_2 e_1$. Moreover, if $E$ is $\F_p$-linear, then we may further view it as a restricted Lie algebra, via setting $e^{[p]} := e^p,$ where the latter denotes the multiplication of $e$ with itself $p$ times. 

In particular, if $\sC$ is $\F_p$-linear, then for any object $M$ of $\sC$, its endomorphisms $\Hom_{\sC}(M, M)$  naturally form a restricted Lie algebra. 

\begin{Lem} If $\sC$ is $\F_p$-linear, and $A$ is an algebra in $\sC$, then its {space of} $A$-valued derivations 
$$\on{Der}(A,A) \subset \Hom_{\sC}(A, A)$$
is a restricted Lie subalgebra.     
\end{Lem}

\begin{proof} We must show that, given two derivations $d, D$, their commutator $d D - D d$ is again a derivation, and that $d^p$ is again a derivation. To see these, note that a derivation $\delta: A \rightarrow A$ is a map satisfying the identity 
\begin{equation} \label{e_selfderiv} \delta \circ \mu = \mu \circ (\delta \otimes \on{id} + \on{id} \otimes \delta).\end{equation}
Using this, the verification of the asserted identities is straightforward. 
\end{proof}

\sss\label{SSS_restricted_definition} Suppose now that $\sC$ is symmetric monoidal, so that we may speak of a commutative Hopf algebra $\OOF$ in $\sC$. 
Write $\RepF$ for its category of left comodules over $\OOF$. Recall that $\RepF$ is naturally symmetric monoidal, compatibly with the forgetful functor $\on{Oblv}: \RepF \rightarrow \sC$.

In particular, given an algebra object $A$ in $\RepF$, we may form its restricted Lie algebra of derivations. This comes equipped with a map of restricted Lie algebras
$$\on{Oblv}: \on{Der}_{\RepF}(A,A) \hookrightarrow \on{Der}_{\sC}(A,A),$$
and explicitly consists of the subspace of derivations $d: A \rightarrow A$ in $\sC$ which are maps of $\OOF$-comodules. 

\begin{defi} For $\OOF$ viewed as a comodule over itself via the right translation coaction, we set 
$$\ff := \on{Der}_{\RepF}(\OOF, \OOF),$$equipped with its natural structure of restricted Lie algebra.    {We call $\ff$ the {\it Lie algebra} of $\OOF$}. 
\end{defi}

\sss The relationship between $\ff$ and the tangent space at the identity is as follows. Write 1 for the monoidal unit of $\sC$ and  $e: \OOF \rightarrow 1$ for the augmentation map of algebras. Passing to the corresponding categories of bimodules, $e$ gives rise to induction and restriction functors 
$$e^*: \OOF\text{-Bimod} \rightleftarrows 1\text{-Bimod} \simeq \sC: e_*,$$
so that in particular we may form the $\OOF$-bimodule $e_* 1$.

\begin{Lem} \label{l_LieAlg=TanSpace}There is a canonical isomorphism of abelian groups\label{l:restrid} 
$$\ff \simeq \on{Der}(\OOF, e_* 1).$$
\end{Lem}

\begin{proof} We will construct maps in either direction and check they are mutually inverse. Given an element $X$ of $\ff$, we may consider the composition 
$$\OOF \xrightarrow{X} \OOF \xrightarrow{e} 1;$$this assignment yields a map $\ff \rightarrow \on{Der}(\OOF, e_* 1)$. Conversely, given an element $Z$ in $\on{Der}(\OOF, e_* 1)$, one can produce an element of $\ff$ as follows. Write $\Delta: \OOF \rightarrow \OOF \otimes \OOF$ for the comultiplication map, and attach to $Z$ the composition 
$$\OOF \xrightarrow{\Delta} \OOF \otimes \OOF \xrightarrow{ Z \otimes \on{id}} e_* 1 \otimes \OOF \simeq \OOF.$$As $\Delta$ is a map of algebras, and $Z \otimes \on{id}$ is a derivation of $\OOF \otimes \OOF$-bimodules, the composition is an $\OOF$-derivation, and each map is a morphism of $\OOF$-comodules, this yields the desired map $\on{Der}(\OOF, e_* 1) \rightarrow \ff$. 

It is straightforward to see that the composition $$\on{Der}(\OOF, e_* 1) \rightarrow \ff \rightarrow \on{Der}(\OOF, e_* 1)$$is the identity, namely 
$$e \circ (Z \otimes \on{id}) \circ \Delta = Z \circ (\on{id} \otimes e) \circ \Delta = Z \circ \on{id} = Z.$$
To see that the composition $\ff \rightarrow \on{Der}(\OOF, e_* 1) \rightarrow \ff$ is the identity, for an element $X$ of $\ff$ one uses the commutative diagram 
$$\xymatrix{ \OOF \ar[d]_{\Delta} \ar[r]^{X} & \OOF \ar[d]^{ \Delta} \ar[r]^{\on{id}} & \OOF \ar[d]^{\on{id}} \\ \OOF \otimes \OOF \ar[r]^{X \otimes \on{id}} & \OOF \otimes \OOF \ar[r]^{e \otimes \on{id}} & \OOF.}
$$

\end{proof}

\sss Let us now discuss the action of the Lie algebra $\ff$ on objects of $\RepF$. 

{
\begin{Lem}\label{l:derivaction} For any object $A$ in $\RepF$, there is a canonical map of restricted Lie algebras
$$\ff \rightarrow \on{Hom}_{\sC}(A,A), \quad \quad X \mapsto \mathscr{L}_X.$$Moreover, this assignment is functorial in $A$, i.e., given a map $\phi: A \rightarrow B$ in  $\RepF$, and any $X \in \ff$, the following diagram commutes
$$\xymatrix{A \ar[d]_\phi \ar[r]^{\mathscr{L}_X} & A \ar[d]^\phi \\ B \ar[r]^{\mathscr{L}_X}  & B.}$$    
\end{Lem}
}
In what follows, we call $\mathscr{L}_X$ the {\it Lie derivative}. 

\begin{proof} Consider the coaction map 
$$\on{coact}_A: A \rightarrow  \OOF \otimes A.$$
Note that this is a map of $\OOF$-comodules, if we equip $A$ with its initial coaction, and ${\OOF} \otimes A$ with the tensor product of the left translation coaction on ${\OOF}$ and trivial coaction on $A$. We will henceforth refer to these as the `left' coactions on $A$ and $\OOF \otimes A$, respectively. 

In addition, consider the commuting trivial `right' coaction of $\OOF$ on $A$, and the commuting `right' coaction of $\OOF$ on $\OOF \otimes A$ given by the tensor product of the right translation coaction on $\OOF$ on itself and $\on{coact}_A$ and. Then, with respect to these additional `right' coactions, $\on{coact}_A$ is in fact a map of $\OOF \otimes \OOF$-comodules. 

Moreover, we claim that 
\begin{itemize}
\item[(*)]
$\on{coact}_A$ exhibits $A$ as the invariants of $\OOF \otimes A$ with respect to the `right' coaction. \end{itemize}
Recall that by the invariant in this context one means  the equalizer of the `right' coaction and the trivial coaction $\OOF \otimes A \rightrightarrows \OOF \otimes \OOF \otimes A$. Note that the equalizers exist in $\mathscr{C}$.

To prove (*),  consider the automorphism $\theta$ of $\OOF \otimes A$ in $\sC$ given by the composition
$$\OOF \otimes A \xrightarrow{\on{id} \otimes \on{coact}_A} \OOF \otimes \OOF \otimes A \xrightarrow{\mu \otimes \on{id}} \OOF \otimes A.$$
By construction, $\theta$ is an automorphism of $\OOF$-comodules with respect to previously defined `left' coaction, and transforms the `right' coaction on the left hand copy of $\OOF \otimes A$
 with the coaction of $\OOF$ on the right hand copy of $\OOF \otimes A$ given by  the tensor 
product of the right translation coaction on $\OOF$ and the trivial coaction on $A$. Therefore, having applied $\theta$, (*) follows from the observation that for any object $M$ of $\sC$, if one considers the tensor product of the left translation and trivial coaction of $\OOF$ on $\OOF \otimes M$, its invariants are given by the split monomorphism $$M \simeq 1 \otimes M \xrightarrow{1 \otimes \on{id}} \OOF \otimes M.$$

With these preparations, for $X \in \ff$, i.e., a right invariant derivation $X: \OOF \rightarrow \OOF$, consider the derivation
$$X \otimes \on{id}: \OOF \otimes A \rightarrow \OOF \otimes A.$$
By our assumption of right invariance, it follows $X \otimes \on{id}$ commutes with the `right' coaction of $\OOF$ on $\OOF \otimes A$, and in particular restricts to a derivation of its invariants, which is the desired map 
$\mathscr{L}_X: A \rightarrow A.$ The claimed functoriality of this assignment follows from the construction. 
\end{proof}

Let us record a formula for the Lie derivative $\mathscr{L}_X$.

\begin{Lem}\label{l:Liederiv} For a commutative algebra $A$ in $\RepF$, and an element $X \in \ff$, {the derivation} $\mathscr{L}_X$ agrees with the composition 
$$A \xrightarrow{\on{coact}} \OOF \otimes A \xrightarrow{(e \circ X) \otimes \on{id}} 1 \otimes A \simeq A.$$    
\end{Lem}

\begin{proof} This follows from the construction of $\mathscr{L}_X$ and the recollection that the composition $$A \xrightarrow{\on{coact}} \OOF \otimes A \xrightarrow{e \otimes \on{id}} 1 \otimes A \simeq A$$is the identity.  
\end{proof}

{
Let us deduce from the preceding lemma the Leibnitz rule for Lie derivatives. As temporary notation, for an object $A$ of $\Rep(F)$, and $X \in \ff$, let us write $\mathscr{L}_X^A: A \rightarrow A$ for the corresponding endomorphism of $A$, i.e., what was previously denoted simply by $\mathscr{L}_X: A \rightarrow A$.

\begin{Cor} The following are true. 
\begin{enumerate}
 
\item Given a pair of objects $A_1$ and $A_2$ of $\Rep(F)$, for any $X \in \ff$ we have the equality $$\mathscr{L}^{A_1 \otimes A_2}_X = \mathscr{L}^{A_1}_X \otimes \on{id}_{A_2} + \on{id}_{A_1} \otimes \mathscr{L}^{A_2}_X.$$

\item If $A$ is an algebra object of $\Rep(F)$, the map $\ff \rightarrow \Hom_{\sC}(A,A)$ factors as 
$$\ff \rightarrow \on{Der}_\sC(A, A) \hookrightarrow \End_\sC(A,A).$$
I.e., the Lie algebra of $\OOF$ acts by derivations of $A$. 

\end{enumerate} 
\end{Cor}

\begin{proof} For ease of notation, we will use the identification $\ff \simeq \on{Der}(\OOF, e_* 1)$ of Lemma \ref{l_LieAlg=TanSpace}, so that for $X \in \ff$ and an object $A$ of $\Rep(F)$ the corresponding Lie derivative $\mathscr{L}_X$ by Lemma \ref{l:Liederiv} is given by 
$$A \xrightarrow{\on{coact}} \OOF \otimes A \xrightarrow{X \otimes \on{id}} 1 \otimes A \simeq A.$$

For (1), note that $\mathscr{L}^{A_1 \otimes A_2}_X$ is given by the composition
$$A_1 \otimes A_2 \xrightarrow{\on{coact} \otimes \on{coact}} \OOF \otimes A_1 \otimes \OOF \otimes A_2 \simeq \OOF \otimes \OOF \otimes A_1 \otimes A_2 \xrightarrow{\mu \otimes \on{id} \otimes \on{id}} \OOF \otimes A_1 \otimes A_2 \xrightarrow{X \otimes \on{id} \otimes \on{id}} 1 \otimes A_1 \otimes A_2 \simeq A_1 \otimes A_2.$$
On the other hand, $\mathscr{L}^{A_1}_X \otimes \on{id}_{A_2}$ 
is given by the composition 
$$A_1 \otimes A_2 \xrightarrow{\on{coact} \otimes \on{coact}} \OOF \otimes A_1 \otimes \OOF \otimes A_2 \xrightarrow{X \otimes \on{id} \otimes \on{e} \otimes \on{id}} 1 \otimes A_1 \otimes 1 \otimes A_2 \simeq A_1 \otimes A_2,$$
and similarly for $\on{id}_{A_1} \otimes \mathscr{L}^{A_2}_X$.

Therefore, it suffices to identify the composition 
\begin{equation} \label{e:deriv1} \OOF \otimes \OOF \xrightarrow{\mu} \OOF \xrightarrow{X} 1\end{equation}
with the composition
\begin{equation} \label{e:deriv2}\OOF \otimes \OOF \xrightarrow{X \otimes e + e \otimes X} 1 \otimes 1 \simeq 1.\end{equation}
To see this, we note both are derivations $\OOF \otimes \OOF \rightarrow (e \times e)_* 1.$ It is therefore enough to see they agree when restricted to the subalgebras $\OOF \otimes 1$ and $1 \otimes \OOF$. To see this final assertion, it suffices to note that any self derivation $\theta: 1 \rightarrow 1$ necessarily vanishes, as the derivation property implies $\theta = 2 \theta$, i.e., $\theta = 0$. 

Finally, (2) follows from applying (1) to the multiplication map $A \otimes A \rightarrow A$, cf. Equation \eqref{e_selfderiv}.     
\end{proof}

}

Let us also recover from the two preceding lemmas and corollary the functoriality of Lie algebras with respect to maps of Hopf algebras. 

\begin{Lem}\label{l:bert} The assignment $\OOG\mapsto \mathfrak{g}$ naturally enhances to a functor from the category of commutative Hopf algebras in $\sC$ to the category of restricted Lie algebras. 
\end{Lem}

\begin{proof} We must show that a map $\phi: \OOG \rightarrow \OOF$ of commutative Hopf algebras induces a map of restricted Lie algebras
$$d\phi: \ff \rightarrow \fg,$$
compatible with composition of maps, i.e., $d(\phi \circ \psi) = d\psi \circ d\phi$, $d(\on{id}) = \on{id}$.

Let us regard $\OOG$ as an $\OOF \otimes \OOG$-comodule, via the restriction of the left translation action to $F$ and the right translation action. Lemma \ref{l:derivaction}, applied to $\sC = \RepG$, then yields the desired map 
$$\ff \rightarrow \on{Der}_{\RepG}(\OOG, \OOG) = \fg.$$The compatibility with composition and the differential of the identity map follow from the construction.  \end{proof}

\sss As a particular case of the above lemma, we may deduce the following. For a Hopf algebra $\OOF$ in $\sC$, let us write $\on{Aut}(\OOF)$ for its group of Hopf algebra automorphisms. 

\begin{Cor}\label{c:ernie} There is a natural action of $\on{Aut}(\OOF)$ on $\ff$ by restricted Lie algebra automorphisms.     
\end{Cor}

\begin{Rem} Here is a further special case, namely the adjoint action by {\em inner automorphisms}. First, note that the set of maps of algebras $\mathscr{F} := \Hom_{\on{Alg}(\sC)}(\OOF, 1)$ naturally forms a group, where the underlying binary product of two elements $g, h$ is defined as the composition 
$$\OOF \xrightarrow{\Delta} \OOF \otimes \OOF \xrightarrow{g \otimes h} 1 \otimes 1 \simeq 1.$$
With respect to this, the augmentation $e: \OOF \rightarrow 1$ is the identity element, and the inverse of $g$ is given by $g \circ \tau$, where $\tau$ denotes the antipode map of $\OOF$.

Given such a $g \in \mathscr{F}$, we may associate to it two maps in $\sC$
$$\ell_g: \OOF \xrightarrow{\Delta} \OOF \otimes \OOF \xrightarrow{g \otimes \on{id}} 1 \otimes \OOF \simeq \OOF$$
$$r_g: \OOF \xrightarrow{\Delta} \OOF \otimes \OOF \xrightarrow{\on{id} \otimes g} \OOF \otimes 1 \simeq \OOF.$$
One may check that these give rise to an action of $\mathscr{F} \times \mathscr{F}$ on $\OOF$ by algebra automorphisms, where $(g,h)$ acts by $\ell_g \circ r_h = r_h \circ \ell_g.$ In addition, one may check that one obtains an action of $\mathscr{F}$ on $\OOF$ by Hopf algebra automorphisms, where $g$ acts by conjugation, i.e., $\ell_g \circ r_{g^{-1}}.$

The corollary therefore implies we have an action of $\mathscr{F}$ on $\ff$ by restricted Lie algebra automorphisms. 
\end{Rem}

\sss Finally, let us note the compatibility of the previous constructions with change of the symmetric monoidal category. That is, consider a symmetric monoidal functor 
$$\mu: \sC_1 \rightarrow \sC_2.$$
Let $\OOF$ be a Hopf algebra in $\sC_1$, and $\mu(\OOF)$ the corresponding Hopf algebra in $\sC_2$. Let us write $\on{Lie}(\OOF)$ and $\on{Lie}(\mu(\OOF))$ for the corresponding restricted Lie algebras.

\begin{Lem} \label{l:changecat}There is a canonical map of restricted Lie algebras $$\mu: \on{Lie}(\OOF) \rightarrow \on{Lie}(\mu(\OOF)),$$
equivariant for the map of groups $\on{Aut}(\OOF) \rightarrow \on{Aut}(\mu(\OOF))$ acting by restricted Lie algebra automorphisms. 
\end{Lem}

\begin{proof} The existence of the map $\on{Lie}(\OOF) \rightarrow \on{Lie}(\mu(\OOF))$ follows from noting that $\mu$ sends derivations of $\OOF$ to derivations of $\mu(\OOF)$, and $p^{\on{th}}$ powers of endomorphisms of $\OOF$ to $p^{\on{th}}$ powers of endomorphisms of $\mu(\OOF)$, i.e., functoriality induces a map of restricted Lie algebras $$\mu: \on{Der}(\OOF) \rightarrow \on{Der}(\mu(\OOF)),$$which moreover preserves the subspace of invariant derivations. The compatibility with the actions by automorphisms again follows by functoriality. 
\end{proof}

\begin{Rem} \label{Re_allgoodchar0} Note that, in the results and constructions of Section \ref{ss:restrI}, if we pass from symmetric monoidal categories enriched over $\mathbb{F}_p$-modules to more general rings, everything save for the restricted structure still makes sense. 
\end{Rem}

\subsection{The Tate extension as a restricted Lie algebra}

\sss In this section, we will apply the generalities of Section \ref{ss:restrI} to obtain the restricted Lie algebra structure on the Tate extension of the continuous endomorphisms of a Tate vector space; the  formula appears in a canonical formulation in Section \ref{sss:Tateformula} and more explicitly in Corollary \ref{c:Tateformula}. 

\sss We should emphasize that the formula is straightforward to guess directly. However, to deduce from it similar formulas for the restricted Lie algebras of the Kac--Moody and Virasoro groups, as we do in Section \ref{ss:restrIII} below, we need to obtain the formula from the corresponding Tate extension of the group of automorphisms of a Tate vector space. This necessity accounts for the length of this section.

\sss  We begin by recalling some preliminary material on Tate modules; some useful references are Section 2 of \cite{BBE} and Sections 3-5 of \cite{DrInf}.

 Let $\Ring$ be a commutative ring. Let $P$ be a projective $\Ring$-module, not necessarily finitely generated. Let us view $P$ as an object of $\Pro(\Ring\mod)$ via the tautological embedding $\Ring\mod \rightarrow \Pro(\Ring\mod)$.

Recall that its dual $P^* := \Hom_{\Ring\mod}(P, \Ring)$ canonically lifts to an object of $\Pro(\Ring\mod)$, in a more interesting way, as follows. We may write $P$ as the direct limit of its finitely generated submodules $P = \varinjlim M_i$, and the expression $$P^* \simeq \varprojlim \Hom(M_i, \Ring)$$provides the desired lift.

\sss  We recall three properties of this construction. First, given a collection of projective modules $P_i, i \in I$, then there is a canonical isomorphism in $\Pro(\Ring\mod)$
$$(\oplus_i P_i)^* = \Pi_i (P_i^*).$$    
This follows from noting that finitely generated submodules of the form $\oplus_i M_i$, where each $M_i$ is a finitely generated submodule of $P_i$, and all but finitely many are zero, are cofinal among all finitely generated $M \subset P$.

Second, given a free module $P \simeq \oplus_i \Ring$, by cofinality we have similarly that $P^* \simeq  \varprojlim_J \Ring^{J},$ where $J$ runs over all finite subsets of $I$.  Note the latter is simply the product $\Ring^{\Pi I}$ taken in $\Pro(\Ring\mod)$; in what follows, we always use $\Ring^{\Pi I}$ to denote this product, rather than image of the product in $\Ring\mod$. 

Therefore, for a general $P$, if we write it as a summand of $\Ring^{\oplus I}$, then $P^*$ is the corresponding summand of $\Ring^{\Pi I}$.

Finally, note that one has a canonical isomorphism of $\Ring$-modules
$$\Hom_{\Pro(\Ring\mod)}(P^*, \Ring) \simeq P;$$
this follows by taking summands from the case of a free module $P \simeq \Ring^{\oplus I}$.

\sss We recall that an object $V$ of $\Pro(\Ring\mod)$ is called an {\em elementary Tate module} if it is  isomorphic to a direct sum $P \oplus Q^*$, for some projective $\Ring$-modules $P$ and $Q$. 

We recall that a {\em Tate module} $V$ is an object which can be realized as a  direct summand of an elementary Tate module, equivalently of $\Ring^{\oplus I} \times \Ring^{\Pi J}$, for some sets $I$ and $J$.

\begin{Rem} In fact,  Drinfeld proved that for any Tate module $\Ring$-module $V$, there exists a Nisnevich cover $\Ring \rightarrow \Ring'$ for which $V  \otimes_\Ring \Ring'$ is elementary Tate; see Theorem 3.4 of \cite{DrInf} and Proposition 2.12 of \cite{BBE}. \end{Rem}

The basic example of a Tate $\Ring$-module, which explains its relevance to affine Lie algebras, is the following. 

\begin{Ex} The ring of Laurent series $\Ring(\!(t)\!)$, with its usual topology, is a Tate module in $\Pro(\Ring\mod)$. E.g, one may use the splitting 
  $$\Ring(\!(t)\!) \simeq \Ring[[t]] \oplus t^{-1} \Ring[t^{-1}];$$the latter summand is free, and the former summand is the dual of the free module $\Ring(\!(t)\!)dt/\Ring[[t]]dt \simeq t^{-1}\Ring[t^{-1}]dt,$ via the residue pairing.   
\end{Ex}

\sss \label{sss:tateduality}Recall that Tate modules carry a natural involutive duality. Concretely, for an elementary Tate module $P \oplus Q^*$, its dual is the elementary Tate module $P^* \oplus Q$, and one passes to summands for the general case.

\sss Given a commutative algebra $S$ in $\on{Pro}(\Ring\mod)$, we have the associated induction functor
$$\Pro(\Ring\mod) \rightarrow S\mod(\Pro(\Ring\mod)), \quad \quad  M \mapsto M \otimes_{\Ring} S. $$

Given Tate $\Ring$-modules $V$ and $W$, consider the set-valued functor on commutative algebras in $\on{Pro}(\Ring\mod)$ given by
$$\Hom(V,W)(S) := \Hom_{S\mod}( V \otimes_{\Ring} S, W \otimes_{\Ring}  S).$$
Similarly, we may consider the subfunctor of isomorphisms
$$\on{Isom}(V,W)(S) := \on{Isom}_{S\mod}(V \otimes_{\Ring} S, W \otimes_{\Ring} S),$$where the right hand side denotes the set of invertible homomorphisms between $V \otimes_{\Ring} S$ and $W \otimes_{\Ring} S$.

\sss \label{SSS:ind_schemes} We would like to describe in what sense the previously introduced functors are representable. To do so, let us recall some standard definitions pertaining to ind-schemes. Let us write $\on{Sch}_\Ring$ for the category of schemes over $\Ring$.  We refer to a general object of its ind-completion $\on{IndSch}_\Ring^{ns} := \on{Ind}(\on{Sch}_\Ring)$, i.e., a filtered colimit of $\Ring$-schemes along arbitrary maps, as a {\em non-strict ind-scheme}, we emphasize that the filtered diagrams, here and below, need not be countable. If it is moreover isomorphic to a filtered colimit of  $\Ring$-schemes under closed embeddings, we refer to it as a {\em strict ind-scheme}, and denote by $$\on{IndSch}_\Ring^{s} \subset \on{IndSch}_\Ring^{ns}$$the corresponding full subcategory. 

 By a {\em non-strict ind-affine ind-scheme}, we mean an object of $\on{IndSch}_{\Ring}^{ns}$ isomorphic to a filtered colimit of affine $\Ring$-schemes under arbitrary maps; if the maps may be taken to be closed embeddings, we refer to it as a {\em strict ind-affine ind-scheme}. Let us denote the full subcategories by $$\on{IndAffSch}^{s}_\Ring \subset \on{IndAffSch}^{ns}_\Ring \subset \on{IndSch}_\Ring^{ns}.$$

  We will follow the standard abuse of notation that strict (ind-affine) ind-schemes will be referred to simply as (ind-affine) ind-schemes, unless there is a risk of confusion.

\sss Note that there is a tautological equivalence between non-strict ind-affine ind-schemes and the opposite of the pro-category of (commutative) $\Ring$-algebras, i.e.,   
$$\on{IndAffSch}_\Ring^{ns} \simeq \on{Pro}(\on{CommAlg}(\Ring\mod))^\on{op}.$$Explicitly, this assigns to a filtered diagram $(Z_i)_{i \in I}$ of affine schemes the cofiltered diagram of $\Ring$-algebras obtained by passing to functions and their pullbacks $(\OO_{Z_i})_{i \in I^{\on{op}}}$. In particular, this exchanges ind-affine ind-schemes with cofiltered diagrams of $\Ring$-algebras under surjective transition maps. 

The following basic lemma will be of use to us. 
\begin{Lem} \label{l:limcolim}The category $\on{IndAffSch}_{\Ring}^{ns}$ is complete, i.e., every small limit of representable contravariant functors is again representable, and cocomplete, i.e., every small colimit of corepresentable covariant functors is again corepresentable.  
\end{Lem}

\begin{proof} We equivalently must show that $\on{Pro}(\on{CommAlg}(\Ring\mod))$ is complete and cocomplete.

We first recall that $\on{CommAlg}(\Ring\mod)$ is complete and cocomplete. Indeed, for completeness, the formation of limits commutes with the forgetful functor to $\Ring\mod$. For cocompleteness, recall that if a category contains finite coproducts and arbitrary filtered colimits, it contains arbitrary coproducts, and if a category contains  coequalizers and arbitrary coproducts  then it is cocomplete. In particular, it suffices to show we have all finite coproducts, filtered colimits, and coequalizers.

In our situation, finite coproducts are given by tensor product over $\Ring$, the formation of filtered colimits commutes with the forgetful functor to $\Ring\mod$, and the coequalizer of a diagram $$f: A \rightarrow B \leftarrow A:g$$is given by the quotient of $B$ by the ideal generated by the elements $f(a) - g(a), a \in A$.

From here, we recall some generalities. Let $\sC$ be a category. First, if $\sC$ is closed under taking limits of finite diagrams, then $\Pro(\sC)$ is complete, i.e., closed under taking limits of arbitrary diagrams, cf. Proposition 6.1.18(iii) of \cite{KS} for the dual statement for ind-completions. In addition, if $\sC$ is cocomplete, then so is $\Pro(\sC)$, see Proposition 11.1 of \cite{Isak}. In particular, if $\sC$ is complete and cocomplete, the same holds for $\Pro(\sC)$. Applying this to $\sC = \on{CommAlg}(\Ring\mod)$, we are done. \end{proof}

Finally, note that there is a tautological full embedding 
$$\on{Pro}(\on{CommAlg}(\Ring\mod)) \hookrightarrow \on{CommAlg}(\on{Pro}(\Ring\mod)).$$
In particular, we may speak of a set-valued functor on commutative algebras $S$ in $\on{Pro}(\Ring\mod)$ being representable by a strict or non-strict ind-affine ind-scheme $Z$.

\sss With those preparations in hand, we may state the following. 

\begin{Lem}\label{l:homtate} For any two Tate vector spaces $V,W$, the functors $\Hom(V,W)$ and $\on{Isom}(V,W)$ are representable by ind-affine ind-schemes.  
\end{Lem}

\begin{proof} We will show that the case of $\on{Isom}(V,W)$ follows from that of $\Hom(V,W)$. We will in turn reduce the latter, in a series of steps, to the case of $V$ and $W$ being the sum of a free module and the dual of a free module.

{\em Step 1.} We first consider the case of $\Hom(V,W)$, beginning with $V = \Ring$. For a free discrete module $W \simeq \oplus_I \Ring$, we have \begin{equation} \label{tatest}\Hom(\Ring,W) \simeq \varinjlim_J \Pi_J \mathbb{A}^1,\end{equation}where $J$ runs over all finite subsets of $I$, and to an inclusion $J \subset J'$ one associates the standard linear closed embedding $\Pi_J \mathbb{A}^1 \hookrightarrow \Pi_{J'} \mathbb{A}^1$; in particular \eqref{tatest} exhibits $\Hom(\Ring, W)$ as an ind-affine ind-scheme, namely ind-$\mathbb{A}^\infty$.
For the dual of a free discrete module $W \simeq \Pi_I \Ring$ we have $\Hom(\Ring,W)$ is the affine scheme $\Pi_I \mathbb{A}^1$, i.e., pro-$\mathbb{A}^\infty$. 

{\em Step 2.} We next consider the case of $\Hom(V,W)$, for $V$ and $W$ of the form $$V \simeq \Ring^{\oplus I} \times \Ring^{\Pi J}, \quad \quad W \simeq \Ring^{\oplus K} \times \Ring^{\Pi L},$$for sets $I,J,K,L$.  Note we have a tautological isomorphism of functors
$$\Hom(\oplus_I \Ring \oplus \Pi_J \Ring, \oplus_K \Ring \oplus \Pi_L \Ring) \simeq \Pi_I (\oplus_K \Ring) \times \Pi_{I \times L} \Ring \times \oplus_{J \times K} \Ring \times \oplus_J \Pi_L \Ring.$$For the first term, note that ind-affine ind-schemes are closed under arbitrary products, e.g. by the argument of Proposition 11.1 of \cite{Isak}, hence $\Pi_I (\oplus_K \Ring)$ is again ind-affine. We saw above that $\Pi_{I \times L} \Ring$ and $\oplus_{J \times K} \Ring$ are ind-affine, and for the final term we have $\oplus_J \Pi_L \Ring$ is the colimit of the finite products $\Pi_{J' \times L} \Ring$, where $J'$ runs over all finite subsets of $J$, under the standard linear closed embeddings.

{\em Step 3.} We now consider the case of $\Hom(V,W)$, for general Tate modules $V$ and $W$.  Choose modules $V' \simeq \oplus_I \Ring \oplus \Pi_J \Ring$, $W' \simeq \oplus_K \Ring \oplus \Pi_L \Ring$ and embeddings $\iota_V: V \rightarrow V'$, $\iota_W: W \rightarrow W'$ with splittings $\pi_V: V' \rightarrow V$, $\pi_W: W' \rightarrow W$.  It follows that $\Hom(V,W)$ is the equalizer of $\Hom(V', W') \rightrightarrows \Hom(V',W'),$ where the two maps send $\phi \in \Hom(V', W')$ to $\phi$ and $\iota_W \circ \pi_W \circ \phi \circ \iota_V \circ \pi_V$, respectively. 

As an ind-affine ind-scheme is in particular an inductive limit of separated schemes along closed embeddings, it follows that $\Hom(V,W)$ is a closed subfunctor of $\Hom(V', W')$, and in particular again an ind-affine ind-scheme, as desired. 

{\em Step 4.} We finally consider the case of $\on{Isom}(V,W)$. Note that $\on{Isom}(V,W)$ is the subfunctor of $\Hom(V,W) \times \Hom(W,V)$ consisting of pairs $(\phi, \psi)$ satisfying $\phi \circ \psi = \on{id}_W$, $\psi \circ \phi = \on{id}_V$. Again by the separatedness of the appearing ind-schemes, it follows $\on{Isom}(V,W)$ is a closed subfunctor, and hence again ind-affine.  \end{proof}

\sss As a particular case of Lemma \ref{l:homtate}, it follows that, for a given Tate  module $V$, its automorphism group 
$$GL(V) := \on{Isom}(V,V)$$
is a group ind-affine ind-scheme. In particular, it corresponds to a commutative Hopf algebra in $\on{Pro}(\Ring\mod)$. Note this commutative Hopf algebra has an underlying commutative algebra which is moreover {\em strict}, i.e., representable as an inverse limit of discrete $\Ring$-algebras along surjective maps of algebras. 

\sss Explicitly, if we write $GL(V)$ as a filtered colimit of affine $\Ring$-schemes along closed embeddings $GL(V) \simeq \varinjlim_i \Spec(T_i)$, we have an isomorphism of algebras $$\OO_{GL(V)} \simeq \varprojlim T_i \in \Pro(\Ring\mod)$$where the latter formal inverse limit of commutative algebras under algebra maps carries its natural commutative algebra structure. Namely, if we denote the multiplication and projection maps by 
$$\mu: \varprojlim T_i \otimes \varprojlim T_i \rightarrow \varprojlim T_i, \quad \mu_i: T_i \otimes T_i \rightarrow T_i, \quad \pi_i: \varprojlim T_i \rightarrow T_i, \quad i \in I,$$
for any $i$ the composite 
$$\varprojlim T_i \otimes \varprojlim T_i \xrightarrow{\mu} \varprojlim T_i \xrightarrow{\pi_i} T_i$$factors as
$$\varprojlim T_i \otimes \varprojlim T_i \xrightarrow{\pi_i \otimes \pi_i} T_i \otimes T_i \xrightarrow{\mu_i} T_i.$$

\sss So, if $\Ring$ is an $\mathbb{F}_p$-algebra, we obtain from the discussion of Section \ref{ss:restrI}, see, in particular, Section \ref{SSS_restricted_definition}, its restricted Lie algebra $\mathfrak{gl}(V)$, which we now explicitly identify.

\begin{Prop} \label{p:endtate}Suppose that $\Ring$ is an $\mathbb{F}_p$-algebra. Then there is a canonical isomorphism of restricted Lie algebras
$$\mathfrak{gl}(V) \simeq \Hom_{\Pro(\Ring\mod)}(V,V),$$
where the latter carries its tautological structure of restricted Lie algebra as in Section \ref{sss:restrassoc}.     
\end{Prop}

\begin{proof} We begin with the identification as abelian groups. As in Lemma \ref{l:restrid}, we have 
\begin{equation} \label{e:glV}\mathfrak{gl}(V) \simeq \on{Der}_{\on{Pro}(\Ring\mod)}(\OO_{GL(V)}, e_* \Ring).\end{equation}
As in Remark \ref{r:splitderiv}, \eqref{e:glV} may identified with maps of algebras $\OO_{GL(V)} \rightarrow \Ring[\epsilon]$ in $\Pro(\Ring\mod)$, where $\epsilon$ is square zero, lifting $e: \OO_{GL(V)} \rightarrow \Ring$, i.e., automorphisms of $V \otimes \Ring[\epsilon]$ which reduce to $\on{id}_V$ modulo $\epsilon$. However, the latter are simply endomorphisms of the form $\on{id}_V + \epsilon X$, for any element $X$ in $\on{Hom}_{\Pro(\Ring\mod)}(V,V)$, as desired.

Let us verify the obtained identification is one of restricted Lie algebras. Consider the natural action of $GL(V)$ on $V$, and in particular the induced map of restricted Lie algebras
$\mathfrak{gl}(V) \rightarrow \on{Der}(\OO_V).$ For $X \in \mathfrak{gl}(V)$, consider the associated map 
\begin{equation} \label{e:acttate}\on{Spec} \Ring[\epsilon] \times V \rightarrow GL(V) \times V \rightarrow V.\end{equation}
Note that $V$ is representable by the algebra $\Sym V^*$, i.e., the free commutative algebra in $\Pro(\Ring\mod)$ associated to its dual Tate module $V^*$, cf. Section \ref{sss:tateduality}, so in particular an $S$-point of $V$ is the same data as a map $V^* \rightarrow S$ in $\Pro(\Ring\mod)$.

With this, the map \eqref{e:acttate} is given on $S$-points by the formula
$$(\eta, v) \mapsto v + \eta \cdot X(v),$$where $\eta: \Ring[\epsilon] \rightarrow S$, $v: V^* \rightarrow S$, and $\eta \cdot X(v)$ denotes the composition
$$V^* \simeq \Ring \otimes V^* \xrightarrow{-\epsilon \otimes X^*} \Ring[\epsilon] \otimes V^* \xrightarrow{\eta \otimes v} S \otimes S \rightarrow S,$$and $X^*: V^* \rightarrow V^*$ denotes the dual of $X: V \rightarrow V$. 

Consider within all vector fields on $V$ the linear ones
$$\Hom_{\Pro(\Ring\mod)}(V^*, V^*) \hookrightarrow \Hom_{\Pro(\Ring\mod)}(V^*, \Sym V^*) \simeq \on{Der}(\Sym V^*, \Sym V^*),$$and note that the above composition is in fact a homomorphism of restricted Lie algebras. The identity \eqref{e:acttate} implies that the vector field $\mathscr{L}_X$ associated to $X$ agrees with the linear vector field associated to $-X^*$, whence it follows that the identification $\mathfrak{gl}(V) \simeq \on{Hom}_{\Pro(\Ring\mod)}(V,V)$ is one of restricted Lie algebras, as desired. 
\end{proof}

\sss \label{sss:Tateformula}Let us now consider the Tate extension of $GL(V)$, a group ind-affine ind-scheme  which we denote by $GL(V)^\flat$: 
\begin{equation} 1 \rightarrow \mathbb{G}_m \rightarrow GL(V)^\flat \rightarrow GL(V) \rightarrow 1. \label{e:Tategrp}\end{equation}
A detailed construction of the underlying functor on $\Ring$-algebras appears in Sections 2.10 and 2.11 of \cite{BBE}, cf. also Section 5 of \cite{DrInf}.

Let us recall some salient aspects. First, given $V = V_1 \oplus V_2$, it is known that the pullback of the Tate extension for $GL(V)$ along the tautological map  
$$GL(V_1) \hookrightarrow GL(V), \quad \quad \phi \mapsto \phi \oplus \on{id}_{V_2},$$
yields the Tate extension for $GL(V_1)$, cf. Section 2.10(iii) of \cite{BBE}. Therefore, we may restrict our discussion to the automorphisms of an elementary Tate module. 

Recall that a submodule $L$ of $V$ is called a {\em c-lattice} if it may be realized as the summand $Q^*$ with respect to some elementary decomposition $V = P \oplus Q^*$, i.e., is a summand of $V$ such that $V/L$ is a projective $\Ring$-module. Given two $c$-lattices $L_1, L_2$, their sum $L_1 + L_2$ is again a $c$-lattice, and the quotients $(L_1 + L_2)/L_i$ are finite rank projective $\Ring$-modules, for $i = 1,2$. In particular,  we may form the relative determinant line bundle$$\on{rel.det}(L_1, L_2) := \det(L_1 + L_2/L_1)^\vee \otimes_{\Ring} \det(L_1 + L_2 / L_2).$$ 

With these preliminaries in hand, we may recall the functor of points description of $GL(V)^\flat$. Namely, for a commutative algebra $S$ in $\Ring\mod$, an $S$-point of $GL(V)^\flat$ consists of an $S$-point of $GL(V)$, i.e., an automorphism of the corresponding $S$-Tate module $$\phi: V \otimes_{\Ring} S \simeq V \otimes_{\Ring} S,$$along with a trivialization of the relative determinant line $$\tau: S \simeq \on{rel.det}( Q^* \otimes_{\Ring} S, \phi(Q^* \otimes_{\Ring} S)).$$The multiplicativity of relative determinant lines endows the set of $S$-points with a group structure. For a map $S \rightarrow T$, one sends a pair $(\phi, \tau)$ as above to their base changes $(\phi \underset{S} \otimes T, \tau \underset{S} \otimes T)$, and this is a map of groups.  

The forgetful map $(\phi, \tau) \mapsto \phi$ yields the desired sequence \eqref{e:Tategrp}. 

\begin{Lem}\label{l:homtateext} $GL(V)^\flat$ is representable by an ind-affine ind-scheme.     
\end{Lem}

\begin{proof} It is known that for any $\Ring$-scheme $X$ and a map to $X \rightarrow GL(V)$, 
the pullback $$X^\flat := GL(V)^\flat \underset{GL(V)} \times X,$$equipped with its natural action of $\mathbb{G}_{m}$, is in fact a $\mathbb{G}_m$-torsor over $X$, cf. Section 5.4.3 of \cite{DrInf}. In particular, if we write $GL(V)$ as a filtered colimit of affine schemes under closed embeddings $(Z_i)_{i \in I}$, cf. Lemma \ref{l:homtate}, it follows that $GL(V)^\flat$ is the filtered colimit of affine schemes under closed embeddings $(Z_i^\flat)_{i \in I}.$     
\end{proof}

In particular, as we did for $GL(V)$, we may again pass from $GL(V)^\flat$ to the corresponding commutative Hopf algebra in $\Pro(\Ring\mod)$ to obtain a canonical restricted structure on the Lie algebra of $GL(V)^\flat$.   Below, we will refer to this as the restricted structure of Lemma \ref{l:homtateext}.

\sss To make this obtained restricted Lie algebra more explicit, we start reviewing a canonical description of the underlying  Tate extension of (non-restricted) Lie algebras
$$0 \rightarrow \Ring \rightarrow \mathfrak{gl}(V)^\flat \rightarrow \mathfrak{gl}(V) \rightarrow 0$$cf. Section 2.13 of \cite{BBE}. Again, it is enough, by the compatibility of Tate extensions and direct sums, to discuss the case of an elementary Tate module.

We first recall some definitions. An element $X \in \gl(V)$ is called {\em discrete} if it has an open kernel. Equivalently, it factors through a usual $\Ring$-module, i.e., admits a factorization $V \rightarrow I \rightarrow V$, where $I$ is an object of $\Ring\mod \hookrightarrow \Pro(\Ring\mod)$. It is straightforward to see that the set of discrete endomorphisms, which we denote by $\gl_d(V) \hookrightarrow \gl(V)$, is an ideal of $\gl(V)$, stable under the adjoint action of $GL(V)$. Moreover, as it is an ideal, $\gl_d(V)$ is in particular closed under the operation of taking $p^{\on{th}}$ powers, and therefore inherits the structure of a restricted Lie algebra. 

A subobject of a Tate module $L \hookrightarrow V$ is called {\em bounded} if for any map from $V$ to a usual $\Ring$-module the image of $L$ is finitely generated.

An element $X \in \gl(V)$ is called {\it bounded} if it has bounded image. It is straightforward to see that the set of bounded endomorphisms, which we denote by $\gl_c(V) \hookrightarrow \gl(V)$, is an ideal of $\gl(V)$, stable under the adjoint action of $GL(V)$. 

In particular, the intersection $\gl_f(V) := \gl_c(V) \cap \gl_d(V)$ is the $\on{Ad}$-invariant restricted ideal of maps whose images are finitely generated $\Ring$-modules. This admits a trace homomorphism 
$$\on{tr}: \gl_f(V) \rightarrow \Ring, \quad \phi \mapsto \on{tr}(\phi);$$
{we recall an elementary definition of the latter in the next remark.}

\begin{Rem} Concretely, if we exhibit $V$ as a summand of $\Ring^{\oplus I} \times \Ring^{\Pi J}$, the trace of $\phi$ agrees with the trace of $\phi \oplus 0 \in \gl_f(\Ring^{\oplus I} \times \Ring^{\Pi J}).$ {To define} the latter, given $\psi \in \gl_f(\Ring^{\oplus I} \times \Ring^{\Pi J})$, for any finite subsets $I' \subset I, J' \subset J$, consider the composition $\psi_{I', J'}$ given by
$$\Ring^{I' \sqcup J'} \hookrightarrow \Ring^{\oplus I} \times \Ring^{\Pi J} \xrightarrow{\psi} \Ring^{\oplus I} \times \Ring^{\Pi J} \twoheadrightarrow \Ring^{I' \sqcup J'}.$$
Then $\on{tr}(\psi_{I', J'})$ is the usual trace of the endomorphism of a finitely generated {free} module. This is independent of $I', J'$, for all sufficiently large finite subsets $I', J'$, and for such subsets we have $\on{tr}(\psi) = \on{tr}(\psi_{I', J'})$. \end{Rem}

Using e.g. the elementary definition recalled in the remark, it is straightforward to see that $\on{tr}: \gl_f(V) \rightarrow \Ring$ is a map of restricted Lie algebras, where we view $\Ring$ as the restricted Lie algebra of $GL(\Ring)$, i.e., set $1^{[p]} = 1.$

\sss \label{sss:humgum}
 Note that any element $\phi$ of $\gl(V)$ may be written (non-uniquely) as the sum of a bounded element and a discrete element. Indeed, for $V = \Ring^{\oplus I} \oplus R^{\Pi J}$, consider the projections $p_-$, $p_+$ to $\Ring^{\oplus I}$ and $\Ring^{\Pi J}$, respectively. Then $\phi =  p_- \phi + p_+  \phi$, and the two summands are bounded and discrete, respectively. For a general $V$, choose a split inclusion $$\iota: V \leftrightarrows \Ring^{\oplus I} \oplus R^{\Pi J}: \pi.$$Then we have $$\phi =        \pi  \iota  \phi  \pi   \iota = \pi  (p_-\iota\phi \pi + p_+ \iota\phi \pi)\iota = \pi p_-\iota\phi\pi \iota + \pi p_+ \iota \phi \pi;$$the two summands are again discrete and bounded, respectively. 
 
 By the preceding paragraph, we therefore have a short exact sequence of restricted Lie algebras
$$0 \rightarrow \gl_f(V) \rightarrow \gl_d(V) \rightarrow \gl(V)/\gl_c(V) \rightarrow 0,$$where the appearing maps are the tautological inclusion and projection. By taking the pushout of this along $\on{tr}: \gl_f(V) \rightarrow \Ring$, we obtain a central extension of restricted Lie algebras
\begin{equation} 0 \rightarrow \Ring \rightarrow (\gl(V)/\gl_c(V))^\flat \rightarrow \gl(V)/\gl_c(V) \rightarrow 0.\label{e:tateasymp}\end{equation} 
{ Finally, pulling this back along the projection $\gl(V) \rightarrow \gl(V)/\gl_c(V)$ yields the Tate extension $\gl(V)^\flat$, with an {\em a priori} different structure of restricted Lie algebra from that of Lemma \ref{l:homtateext}. We will refer below to this newly obtained restricted structure as the restricted structure of Section \ref{sss:humgum}.

\sss  We will show momentarily in Lemma \ref{Lem:Tate_restricted_comparison} below that the two obtained restricted structures on the Tate extension coincide. Before doing so, it will be useful to record a concrete cocycle presentation of the restricted structure of Section \ref{sss:humgum}.

\sss In this section we will give an explicit description of the  restricted Lie algebra from Section \ref{sss:humgum} for an elementary Tate module with fixed decomposition $V = P \oplus Q^* =: V^- \oplus V^+$. Consider the associated decomposition $$\Hom_{\Pro(\Ring\mod)}(V,V) \simeq \Hom_{\Pro(\Ring\mod)}(V, V^-) \oplus \Hom_{\Pro(\Ring\mod)}(V, V^+),$$and for an element $\phi \in \gl(V)$ consider its corresponding expression $\phi = \phi_- + \phi_+$, where $\phi_{\pm}: V \rightarrow V^{\pm}.$ 

\begin{Cor}\label{c:Tateformula} The Tate extension of $\mathfrak{gl}(V)$ for $V = V^+ \oplus V^-$ may be written as follows. As an $\Ring$-module, we have 
$$\gl(V)^\flat \simeq \gl(V) \oplus \Ring \cdot 1.$$
As a Lie algebra, $1$ is central, and for $\phi, \psi \in \gl(V)$ we have 
$$[\phi, \psi]_{\gl(V)^\flat} =  [\phi, \psi]_{\gl(V)}  + \on{tr}( [\phi, \psi]_- - [\phi_-, \psi_-]) \cdot 1.$$For the restricted structure, we have $1^{[p]} = 1$, and for $\phi \in \gl(V)$ we have 
\begin{equation} \label{e:restrtate}\phi^{[p]} = \phi^p  + \on{tr}( (\phi_-)^p - (\phi^p)_-) \cdot 1.\end{equation}
\end{Cor}

\begin{proof} By the construction of Section \ref{sss:humgum}, we have a commutative diagram of short exact sequences of restricted Lie algebras 
$$\xymatrix{0 \ar[r] & \gl_f(V)  \ar[d]_{\on{tr}}\ar[r] & \gl_d(V) \ar[r] \ar[d]^{\delta} & \gl(V)/\gl_c(V) \ar[r] \ar[d]^{\on{id}} & 0 \\ 0 \ar[r] & \Ring \ar[r] & (\gl(V)/\gl_c(V))^\flat \ar[r] & \gl(V)/\gl_c(V) \ar[r] & 0 \\ 0 \ar[r] & \Ring \ar[r] \ar[u]^{\on{id}} & \gl(V)^\flat \ar[r] \ar[u] & \gl(V) \ar[r] \ar[u]_q & 0, }$$
where $q$ is the tautological quotient map. In particular, given an element $\phi \in \gl(V)$, a lift of $q(\phi)$ to $(\gl(V)/\gl_c(V))^\flat$ is given by $\delta(\phi_-)$, which provides the desired splitting of $\gl(V)^\flat$ as an extension of $\Ring$-modules. The cocycle expressions for the commutator and restricted structures on $\gl(V)^\flat$ then follow by construction, and the identity $1^{[p]} = 1$ comes from the identification of restricted Lie algebras $\Ring \simeq \on{Lie}(\mathbb{G}_m)$.. 
\end{proof}

\begin{Rem} For later use in calculations, we recall that one may rewrite the cocycle for the commutator as  $$ \on{tr}( [\phi, \psi]_- - [\phi_-, \psi_-])) = \on{tr}( \phi_- \psi_+ - \psi_- \phi_+)  .$$
\end{Rem}

\sss We now verify that the two exhibited restricted structures on $\gl(V)^\flat$ coincide.

\begin{Lem}\label{Lem:Tate_restricted_comparison} The restricted Lie algebra structure on $\gl(V)^\flat$ {viewed as the Lie algebra of $GL(V)^\flat$, see} Lemma \ref{l:homtateext}, agrees with its restricted structure constructed in Section \ref{sss:humgum}.

\end{Lem}

\begin{proof} By the compatibility of Tate extensions and direct sums of Tate modules, we may reduce to the case of $V \simeq V^- \oplus V^+$, where $V^- := \Ring^{\oplus I}, V^+ := \Ring^{\Pi J}$.

Consider the object $V^\infty$ of $\on{Ind}(\Ring\mod) \hookrightarrow \on{Ind}(\Pro(\Ring\mod))$ given by the formal colimit
$$V^\infty := \varinjlim V/L,$$
where $L$ runs over all bounded submodules of $V$; explicitly a cofinal set of $L$'s is provided by $R^{\oplus I'} \oplus V^+ \hookrightarrow  V$, where $I'$ runs over all finite subsets of $I$. Consider the associated space of endomorphisms 
$$\Hom(V^\infty, V^\infty) := \varprojlim_i \varinjlim_j \Hom(V/L_i, V/L_j);$$
where each $\Hom(V/L_i, V/L_j)$ is the strict affine ind-scheme according to  Lemma \ref{l:homtate}. As the category  $\Pro(\on{CommAlg}(\Ring\mod))$ contains all limits and  colimits, cf. Lemma \ref{l:limcolim}, it follows that $\Hom(V^\infty,V^\infty)$ is again representable by an object of $\on{IndAffSch}^{ns}_\Ring$, though it is no longer strict. As $\on{IndAffSch}^{ns}_{\Ring}$ is in particular closed under taking equalizers, it follows that the invertible asymptotic endomorphisms 
$$GL(V^\infty) \hookrightarrow \Hom(V^\infty, V^\infty) \times \Hom(V^\infty, V^\infty)$$are again representable by a  Hopf algebra in $\Pro(\Ring\mod)$. By an argument similar to that of Proposition \ref{p:endtate}, its restricted Lie algebra is given by 
$$\gl(V^\infty) := \varprojlim \varinjlim \Hom_{\Pro(\Ring\mod)}(V/L_i, V/L_j),$$
{ where the restricted structure is again given by sending an element $\phi \in \gl(V^\infty)$ to its $p$-fold iterate $\phi^p$. }

As $V^\infty$ is not a Tate module, note the previously discussed Tate extension does not literally apply to $GL(V^\infty)$. Instead, consider the Tate extension $GL^\flat(V^\infty)$ of $GL(V^\infty)$ as defined in \cite[Section 2.10]{BBE}:
\begin{equation} \label{e:asymptate}1 \rightarrow \mathbb{G}_m \rightarrow GL^\flat(V^\infty) \rightarrow GL(V^\infty) \rightarrow 1.\end{equation}
Let us recall three relevant properties from {\em loc. cit.} of this extension. 

\begin{enumerate}
\item The pullback of $GL^\flat(V^\infty)$ along the natural map $GL(V) \rightarrow GL(V^\infty)$ yields the Tate extension $GL(V)^\flat$ as defined previously.

\item The pullback of $GL^\flat(V^\infty)$ along $GL(V^-) \rightarrow GL(V^\infty)$ is canonically split. 

\item Write $GL(V^-)_f$ for the subgroup of endomorphisms which are the identity on a summand of $V^-$ of the form $R^{\oplus I'}$ for $I' \subset I$ a subset with finite complement. More intrinsically, this is the subgroup of automorphisms $\alpha$ of $V^-$ such that there exists a direct sum decomposition $V^- = Q \oplus Q'$, where $Q'$ is a finitely generated projective module, and $\alpha$ restricts to the identity on $Q$. Note that the composition $$GL(V^-)_f \rightarrow GL(V^-) \rightarrow GL(V^\infty)$$ is the trivial homomorphism, so this affords another splitting of the pullback of the Tate extension to $GL(V^-)_f$. Then the obtained homomorphism $GL(V^-)_f \rightarrow \mathbb{G}_m$ from the pair of splittings  is simply the determinant. 
\end{enumerate}

Let us pass to the restricted Lie algebras of these commutative Hopf algebras in $\Pro(\Ring\mod)$. In particular, let us consider from Equation \eqref{e:asymptate} the associated Tate extension of restricted Lie algebras
$$0 \rightarrow \Ring \rightarrow \gl^\flat(V^\infty) \rightarrow \gl(V^\infty) \rightarrow 0.$$
We deduce from the above three statements respectively that:
\begin{enumerate}
\item[(1$'$)] The pullback of $\gl^\flat(V^\infty)$ along the natural map $\gl(V) \rightarrow \gl(V^\infty)$ yields the Tate extension $\gl(V)^\flat$, with its restricted structure as given by Lemma \ref{l:homtateext}. 

\item[(2$'$)] The pullback of $\gl^\flat(V^\infty)$ along $\gl(V^-) \rightarrow \gl(V^\infty)$ is canonically split. 

\item[(3$'$)] Write $\gl(V^-)_f$ for the restricted subalgebra of endomorphisms of $V^-$ which are identically zero on a summand of $V^-$ of the form $R^{\oplus I'}$ for $I' \subset I$ a subset with finite complement. The composition $$\gl(V^-)_f \rightarrow \gl(V^-) \rightarrow \gl(V^\infty)$$is the trivial homomorphism of restricted Lie algebras, so this affords another splitting of the pullback of the Tate extension to $\gl(V^-)_f$. Then the obtained homomorphism $\gl(V^-)_f \rightarrow \Ring$ of restricted Lie algebras from the pair of splittings is simply the trace. 
\end{enumerate}
Using these three statements, let us now finish the proof of the lemma. First, we note that the tautological map $\gl(V) \rightarrow \gl(V^\infty)$ factors through an isomorphism of restricted Lie algebras $\gl(V)/\gl(V)_c \simeq \gl(V^\infty)$. Second, we note that the tautological map $\gl(V^-) \rightarrow \gl(V^\infty)$ is surjective, and that we have a commutative diagram with exact rows
$$\xymatrix{0 \ar[r] & \gl(V)_f \ar[r] & \gl(V)_d \ar[r] & \gl(V)/\gl(V)_c \ar[r] & 0 \\ 0 \ar[r] & \gl(V^-)_f \ar[u] \ar[r] & \gl(V^-) \ar[r] \ar[u] & \gl(V^\infty) \ar[r] \ar[u]^\sim & 0.}$$We may push out both extensions of restricted Lie algebras along the trace maps $\on{tr}: \gl(V)_f \rightarrow \Ring$ and $\on{tr}: \gl(V^-)_f \rightarrow \Ring$ to obtain a commutative diagram of restricted Lie algebras with exact rows 
$$\xymatrix{0 \ar[r] & \Ring \ar[r] & \widetilde{\gl(V)/\gl(V)_c} \ar[r] & \gl(V)/\gl(V)_c \ar[r] & 0 \\ 0 \ar[r] & \Ring \ar[u]^{\on{id}} \ar[r] & \widetilde{\gl(V^\infty)} \ar[r] \ar[u]^\iota & \gl(V^\infty) \ar[r] \ar[u]^\sim & 0.}$$
By definition, the extension $\widetilde{\gl(V)/\gl(V)_c}$ is the Tate extension $(\gl(V)/\gl(V)_c)^\flat$ with the restricted structure of Section \ref{sss:humgum}. By points (2$'$) and (3$'$) above, the extension $\widetilde{\gl(V^\infty)}$ is the Tate extension $\gl^\flat(V^\infty)$ of $\gl(V^\infty)$ with its single restricted Lie algebra structure considered above. As the middle vertical arrow $\iota$ is tautologically an isomorphism of restricted Lie algebras, we are now done by (1$'$). \end{proof} }

{ In what follows, it will also be convenient to have an explicit formula for the adjoint action of $GL(V)$ on $\gl(V)^\flat$. We record the answer here.

\begin{Cor} \label{c:tategroupcocycle}As in Corollary \ref{c:Tateformula}, consider an elementary Tate module $V = V^+ \oplus V^-$ and  the associated splitting of the underlying $\Ring$-module of the Tate extension: $$\gl(V)^\flat \simeq \gl(V) \oplus \Ring \cdot 1.$$Then, with respect to this splitting, for any $\Ring$-algebra $\Sing$, the adjoint action of $GL(V)$ on $\gl(V)^\flat$ is given on $\Sing$-points by the formulas 
\begin{equation} \label{e:tategroupadj}\on{Ad}_g(1) = 1 \quad \text{and} \quad \on{Ad}_g(\phi) = g\phi g^{-1}  + \on{tr}( (g \circ \phi \circ g^{-1})_- - g \circ \phi_- \circ g^{-1}), \quad \quad g \in GL(V)(\Sing),  \quad \phi \in \gl(V)(\Sing).\end{equation}
  
\end{Cor} 

\begin{proof} Note that we have {\em a priori} two different natural actions of $GL(V)$ on $\gl(V)^\flat$, given as follows.

\begin{enumerate} \item The adjoint action is obtained via thinking of $\gl(V)^\flat$ as the Lie algebra of the central extension $GL(V)^\flat$. 
\item The second action is obtained via Section \ref{sss:humgum}. Namely, $\gl_f(V), \gl_d(V)$, and $\gl_c(V)$ are all $GL(V)$-invariant submodules of $\gl(V)$, and the trace map $$\on{tr}: \gl_f(V) \rightarrow \Ring$$is $GL(V)$-equivariant, where $\Ring$ is equipped with the trivial action of $GL(V)$, so the construction of Section \ref{sss:humgum} yields an action of $GL(V)$ on $\gl(V)^\flat$. 
\end{enumerate}

The proof of Lemma \ref{Lem:Tate_restricted_comparison} shows that these two actions coincide. Moreover, the argument of Corollary \ref{c:Tateformula} yields the formula  \eqref{e:tategroupadj} for the action (2), but this agrees with the adjoint action (1) as noted above, so we are done.  \end{proof}
}

\subsection{Kac--Moody and Virasoro I: the restricted structure}
\label{ss:restrIII}
\sss Let us now deduce the desired formulae for the restricted Lie algebra structures on the Kac--Moody and Virasoro Lie algebras. 

\sss \label{SSS_KM_groups} We begin with Kac--Moody. Let $F$ be a linear algebraic group over $\Ring$ with Lie algebra $\ff$. Recall that the loop group $F(\!(t)\!)$ is the group affine ind-scheme over $\Ring$ with $S$-points given by $F(S(\!(t)\!))$, for any $\Ring$-algebra $S$. If we view $\OO_{F(\!(t)\!)}$ as a commutative Hopf algebra in $\Pro(\Ring\mod)$, recall that there is a canonical identification of its Lie algebra with $\ff(\!(t)\!) := \ff \otimes_{\Ring} \Ring(\!(t)\!).$

\begin{Lem} \label{l:restrloop}The restricted Lie algebra structure on $\ff(\!(t)\!)$ is given by the formula
$$(X \otimes f)^{[p]} = X^{[p]} \otimes f^p, \quad \quad X \in \ff, f \in \Ring(\!(t)\!).$$    
\end{Lem}

\begin{proof} Choose a closed embedding of $F$ into the endomorphisms of a projective $\Ring$-module $V$ of finite rank, and consider the induced action of $F(\!(t)\!)$ on $V(\!(t)\!)$. This induces an {injective} homomorphism of restricted Lie algebras $\ff(\!(t)\!) \rightarrow \gl(V(\!(t)\!))$. The formula now follows from Proposition \ref{p:endtate}.
\end{proof}

\sss\label{SSS_KM_definition}
Suppose in addition we fix a representation $\rho: F \rightarrow GL(W)$, where $W$ is a projective $\Ring$-module of finite rank. Consider the induced map 
$\rho:  F(\!(t)\!) \rightarrow GL(W(\!(t)\!)),$ and denote the pullback of the Tate extension by $F(\!(t)\!)^\flat$. This is a strict group ind-scheme of countable type, i.e., presentable as a countable filtered colimit of affine schemes along closed embeddings. 
In particular, we obtain its corresponding Lie algebra $\ff(\!(t)\!)^\flat$ that fits into the exact sequence
$$0 \rightarrow \Ring \rightarrow \ff(\!(t)\!)^\flat \rightarrow \ff(\!(t)\!) \rightarrow 0.$$
\begin{Prop} \label{p:kacmoodyrestr}We have a splitting of $\Ring$-modules $\ff(\!(t)\!)^\flat \simeq \ff(\!(t)\!) \oplus \Ring \cdot \mathsf{1}$, with respect to which the commutator is given by the formula 
$$[X \otimes f, Y \otimes g]_{\ff(\!(t)\!)^\flat} = [X,Y] \otimes fg - \on{tr}_W(XY) \cdot \on{Res} f dg \cdot \mathsf{1},$$
and the restricted structure satisfies the identities $\mathsf{1}^{[p]} = \mathsf{1}$, $(X \otimes t^i)^{[p]} = X^{[p]} \otimes t^{ip}$, for any $i \in \mathbb{Z}$. 
\end{Prop}
Before giving the proof, we make some orienting remarks. 

\begin{Rem}Note first that, by continuity, and compatibility of restricted structures and sums, cf. Definition \ref{Defi:restricted}(3), the above formulas determine the restricted Lie algebra structure uniquely.     
\end{Rem}

\begin{Rem} \label{r:traceforms} Note in addition that, by inspection, the obtained extension $\ff(\!(t)\!)^\flat$ only depends, as a restricted Lie algebra, on the trace form $$\on{tr}_W \in \Sym^2(\ff^*)^F: \quad \quad X \otimes Y \mapsto \on{tr}_W(XY).$$
We flag for the reader that, unlike in this preliminary discussion, beginning in Section \ref{SSS_affine_VA_defn} below we will incorporate a critical shift in our discussion of trace forms; we hope this does not cause confusion. 
\end{Rem}

\begin{Rem}\label{r:centext}If $\ff$ is simple or one dimensional, as long as the characteristic of $\Ring$ is not too small relative to $\ff$,  note that $\Sym^2(\ff^*)^F$ is a free rank one $\Ring$-module.  It follows that for two representations $V$ and $W$ with $\on{tr}_V$ and $\on{tr}_{W}$ both units in $\Ring$, we have an identification of the underlying Lie algebras of the two central extensions, 
$$ \ff(\!(t)\!)^\flat_V \simeq \ff(\!(t)\!)^\flat_W, \quad \quad X \otimes f \mapsto X \otimes f, 
  \quad \quad \mathsf{1} \mapsto  \frac{\on{tr}_{V}}{\on{tr}_{W}} \cdot \mathsf{1}.$$
Moreover, as $\on{tr}_V$ and $\on{tr}_W$ may be recovered from their restriction to the image of the cocharacter lattice of a maximal torus, it follows that the ratio $\frac{\on{tr}_V}{\on{tr}_W}$
lies in $\mathbb{F}_p^\times$,  whence the above identification is moreover one of restricted Lie algebras. 

The discussion of the preceding paragraph applies, {\em mutatis mutandis}, to the Lie algebra of any split reductive group, where one considers each simple factor and the center as above; below we will only consider cases where the trace forms are nondegenerate. Finally, note that given a restricted Lie algebra $\mathfrak{a}$ and a finite set of restricted central extensions $\widetilde{\mathfrak{a}}_i, i \in I$, of $\mathfrak{a}$ by $\Ring$, one has a natural restricted central extension of $\mathfrak{a}$ by $\Ring^{\oplus I}$, 
$$0 \rightarrow \Ring^{\oplus I} \rightarrow \widetilde{\mathfrak{a}}_I \rightarrow \mathfrak{a} \rightarrow 0,$$such that the pushout along each projection $\pi_i: \Ring^{\oplus I} \rightarrow \Ring, i \in I$, rescovers $\widetilde{\mathfrak{a}}_i$. So, in particular, in the reductive case, we may equally well discuss a Kac--Moody central extension where we include different central elements for different simple and central factors.  
\end{Rem}

\begin{proof}[Proof of Proposition \ref{p:kacmoodyrestr}.] We will apply Corollary \ref{c:Tateformula} with respect to the splitting $$W(\!(t)\!) = W \otimes_{\Ring} \Ring[[t]] \oplus W \otimes_{\Ring} t^{-1} \Ring[t^{-1}].$$ For the restricted structure, we apply the identity
 \eqref{e:restrtate}, i.e., 
 $$\phi^{[p]} = \phi^p  + \on{tr}( (\phi_-)^p - (\phi^p)_-) \cdot 1$$
 to $\phi = X \otimes t^i \in \ff(\!(t)\!)$.

 We saw in Lemma \ref{l:restrloop} that $\phi^{[p]} = X^{[p]} \otimes t^{ip}$.  We claim that the cocycle term $ \on{tr}( (\phi_-)^p - (\phi^p)_-) \cdot 1$ vanishes. Indeed, this follows from considering the homogeneity of $(\phi_-)^p$ and $(\phi^p)_-$ with respect to the loop rotation grading on $$W \otimes_{\Ring} \Ring[t^{\pm 1}].$$Namely, for $i \neq 0$, it is clear that $\on{tr}( (\phi_-)^p - (\phi^p)_-)$ vanishes, as both appearing operators have degree $ip \neq 0$, hence have vanishing trace. For $i = 0$,  the vanishing follows from noting that $(\phi_-)^n = (\phi^n)_-$, for any $n \in \mathbb{Z}$, and in particular for $n = p$. Therefore, we have $(\phi_-)^p - (\phi^p)_-$ vanishes, and in particular has vanishing trace.    
\end{proof}

\sss\label{SSS_Virasoro}
Let us now discuss the Virasoro case, beginning with some recollections on the Witt group and Lie algebra. Recall, cf. {Section \ref{SSS:ind_schemes}}, that a non-strict ind-affine ind-scheme over $\Ring$ is an object in the opposite category of $\on{Alg}(\Pro(\Ring\mod))$. Let us write $\D^{\times}$ for the formal punctured disc, i.e., the ind-affine ind-scheme corresponding to $\Ring(\!(t)\!).$

Consider the functor $\on{End}(\D^\times)$ from $\Ring$-algebras to sets given by the formula 
$$\on{End}(\D^\times)(S) := \Hom_{\on{Alg}(\Pro(S\mod))}(S(\!(t)\!), S(\!(t)\!)).$$Explicitly, for fixed $S$, such a map is entirely determined by the image of $t$, which must be a Laurent series $\phi(t) = \sum a_i t^i  \in S(\!(t)\!)$ with $a_i$ nilpotent for $i \leqslant 0$. In particular $\on{End}(\D^\times)$ is a strict ind-affine ind-scheme of countable type. We may further consider the automorphisms of the formal disc
$$\on{Aut}(\D^\times) \hookrightarrow \on{End}(\D^\times) \times \on{End}(\D^\times),$$
consisting of pairs $(\phi(t), \psi(t))$ satisfying $\phi(\psi(t)) = \psi(\phi(t)) = t.$ It follows that $\on{Aut}(\D^\times)$ is a group ind-affine ind-scheme of countable type, and in particular we form its commutative Hopf algebra of functions in $\Pro(\Ring\mod)$, and its restricted Lie algebra $\mathfrak{Witt}$. 

\begin{Rem} One may also check that the first projection $\on{Aut}(\D^{\times}) \rightarrow \on{End}(\D^{\times})$ exhibits it as the open subfunctor consisting of series $\phi(t) = \sum a_i t^i$ as before which moreover satisfy the condition that $a_1$ is a unit in $S$.  
\end{Rem}

Consider  $\Ring(\!(t)\!) \partial_t$, the restricted Lie algebra of vector fields on the formal punctured disc, i.e., derivations of $\Ring(\!(t)\!)$ in $\Pro(\Ring\mod)$. Explicitly, the commutator is the usual Lie bracket of vector fields, and the restricted structure sends a derivation to its $p^{th}$ power, cf. Section \ref{sss:restrassoc}. 

\begin{Prop} \label{p:wittrestr} There is a canonical isomorphism of restricted Lie algebras between $\mathfrak{Witt}$ and $\Ring(\!(t)\!)\partial_t$. In particular, if for $i \in \mathbb{Z}$ we write $L'_i := - t^{i+1} \partial_t$, these satisfy the identities $$L_i'^{[p]} = \begin{cases} L'_{ip}, & (i,p) \neq 1, \\ 0 & (i,p) = 1. \end{cases}$$
\end{Prop}

\begin{proof} The identification $\mathfrak{Witt} \simeq \Ring(\!(t)\!)\partial_t$ as abelian groups follows by considering $\Ring[\epsilon]$ points of $\on{Aut}(\D^\times)$ extending the identity,  cf. the beginning of the proof of Proposition \ref{p:endtate}.

If we write $\D^{\times}$ for `$\Spec \Ring(\!(t)\!)$', i.e.,  the object of $(\on{CommAlg}(\Pro(\Ring\mod)))^{op}$ corresponding to $\Ring(\!(t)\!)$, we have a tautological action of $\on{Aut}(\D^{\times})$ on $\D^{\times}$. Moreover, for $X \in \Ring(\!(t)\!)\partial_t$, the associated derivation $\mathscr{L}_X$ of $\Ring(\!(t)\!)$ agrees with its standard action, cf. Lemma \ref{l:Liederiv}, from which the first claim of the proposition follows. The second statement concerning $L_i'^{[p]}$ then follows by applying $L'_i$, viewed as an endomorphism of $\Ring(\!(t)\!)$, $p$ times to  $t$. 
\end{proof}

\sss\label{SSS_Virasoro_pcenter} Let us now discuss the central extension. Consider the tautological linear action of $\on{Aut}(\D^\times)$ on the Tate $\Ring$-module $\Ring(\!(t)\!)$ by change of coordinates. Viewing this as a homomorphism $\on{Aut}(\D^\times) \rightarrow GL(\Ring(\!(t)\!)),$
we may pull back the Tate extension to obtain the Virasoro group ind-scheme {$\on{Aut}(\D^{\times})^\flat$ that fits} into the following short exact sequence of group ind-schemes
$$1 \rightarrow \mathbb{G}_m \rightarrow \on{Aut}(\D^{\times})^\flat \rightarrow \on{Aut}(\D^\times) \rightarrow 1;$$
and in particular its restricted Lie algebra
$$0 \rightarrow \Ring \cdot 1 \rightarrow \mathfrak{Vir} \rightarrow \mathfrak{Witt} \rightarrow 0.$$

\begin{Prop} We have a splitting of $\Ring$-modules $\mathfrak{Vir} \simeq \mathfrak{Witt} \oplus \Ring \cdot 1$, with respect to which the commutator is determined by the formula 
$$[L'_m , L'_n] = (m-n) L'_{m+n}  -  \delta_{m, -n} \cdot \binom{m+1}{3}  \cdot 1,$$
where $\delta_{m, -n}$ denotes the Kronecker delta function, and the restricted structure is determined by the identities $$1^{[p]} = 1, \quad \quad L_i'^{[p]} = \begin{cases} L'_{ip}, & (i,p) \neq 1, \\ 0 & (i,p) = 1. \end{cases}$$. \end{Prop}

\begin{proof} This follows from applying Corollary \ref{c:Tateformula} with respect to the splitting $\Ring(\!(t)\!) = \Ring[[t]] \oplus t^{-1} \Ring[t^{-1}].$ Let us describe the calculation of the cocycle for the commutator; the argument for the restricted structure proceeds similarly to that of Proposition \ref{p:kacmoodyrestr}. 

 Note the cocycle modifying the Lie bracket vanishes for $[L'_m, L'_n]$ unless $m = -n$ for degree reasons. For $m > 0$, note that $(L'_m)_- (L'_{-m})_+$ vanishes, hence we have 
 $$[L'_m, L'_{-m}] = 2m \cdot L'_0 +  \on{tr}( (L'_{-m})_- (L'_m)_+).$$
For degree reasons, the trace of $(L'_{-m})_- (L'_m)_+$ agrees with the trace of its restriction to the span of $t^{-1}, \ldots, t^{-m-1}$, where $(L'_{-m})_- (L'_m)_+$ agrees with the operator 
$$(t^{-m+1} \partial_t) \circ (t^{m+1} \partial_t) = (t \partial_t)^2 + m  \cdot t \partial_t.$$In particular, its trace is given by 
$$\sum_{i = 1}^m (i^2 - i \cdot m) = - \binom{m+1}{3},$$
as desired. 
\end{proof}

We would like to point out that this restricted Lie algebra structure on $\mathfrak{Vir}$ has been previously constructed in \cite{JLM}.

{ \sss{} \label{sss:autd} For future use, let us recall two relevant subgroups of $\Aut(\pD)$. 

First, if in the discussion of Section \ref{SSS_Virasoro} we replace $\Sing(\!(t)\!)$ with $\Sing[[t]]$, we obtain the automorphisms of the formal disc $\Aut(\sD)$. Explicitly, an $\Sing$-point of $\Aut(\sD)$ is a series 
$$\phi(t) = \underset{i \geqslant 0} \Sigma \hspace{.5mm} a_it^i \in \Sing[[t]],$$
with (i) $a_0$ nilpotent, (ii) $a_1$ a unit, and (iii) $a_i, i > 1$, arbitrary, and the group law is given by composition of series. In particular, $\Aut(\D)$ is a group ind-scheme, with its underlying scheme the product of (i) the formal completion of $\mathbb{A}^1$ at 0, (ii) $\mathbb{G}_m$, and (iii) a countable product of copies of $\mathbb{A}^1$, respectively. 

Within $\Aut(\sD)$, we the subfunctor of automorphisms which preserve the center of the disc, i.e., series as above with $a_0 = 0$. This is now a group scheme, which we denote by $\Aut^\bullet(\sD)$. To summarize, we have closed embeddings of of group ind-schemes
$$\Aut^\bullet(\sD) \hookrightarrow \Aut(\sD) \hookrightarrow \Aut(\pD)$$
with restricted Lie algebras respectively given as in Proposition \ref{p:wittrestr} by 
$$t\Ring[[t]]\partial_t \hookrightarrow \Ring[[t]]\partial_t \hookrightarrow \Ring(\!(t)\!)\partial_t. $$}

\subsection{Kac--Moody and Virasoro II: functor of points}
\label{s:kmspace}

\sss Let us retain our assumption that $\Ring$ is an $\F_p$-algebra. So far, the methods used in this section have produced the Kac--Moody and Virasoro algebras as abstract restricted Lie algebras over $\Ring$, equipped with adjoint actions of the abstract groups of $\Ring$-points of the loop group and $\on{Aut}(\D^\times)$, respectively. 

On the other hand, it will be useful in our analysis later that the Kac--Moody and Virasoro algebras are naturally Tate modules over $\Ring$, and the previous actions of $\Ring$-points naturally extend to actions of the corresponding group ind-schemes. 

In this final section, we explain how to recover the latter {\em a priori} richer structures of algebro-geometric nature from the previous construction on rational points.

\sss\label{SSS_spaces} 
Recall that an $\Ring$-{\em space} $\sY$ is an arbitrary functor from $\Ring$-algebras to sets. This category is naturally symmetric monoidal with respect to the categorical product of functors, i.e., $$\sY_1 \times \sY_2 (\Sing) := \sY_1(\Sing) \times \sY_2(\Sing).$$

Using this, one can speak of a group $\Ring$-space $\mathscr{G}$; explicitly, this consists of the data of a group structure on each $\mathscr{G}(\Sing)$, such that for each $\Sing \rightarrow \Ting$ the map $\mathscr{G}(\Sing) \rightarrow \mathscr{G}(\Ting)$ is a group homomorphism.

Similarly, recall that one can speak of an $\Ring$-module structure on an $\Ring$-space $\sM$, this is the datum of `scaling' and `addition' maps 
$$\mathbb{A}^1_\Ring \times \sM \rightarrow \sM \quad \quad \sM \times \sM \rightarrow \sM$$
satisfying the usual identities. Explicitly, an $\Ring$-module structure on $\sM$ is the data of an $\Sing$-module structure on $\sM(\Sing)$, for any $\Ring$-algebra $\Sing$, such that for any map $\Sing \rightarrow \Ting$ the map $\sM(\Sing) \rightarrow \sM(\Ting)$ is a map of $\Sing$-modules. 

Now assume that $\Ring$ is an $\F_p$-algebra. Given an $\Ring$-module $\sM$, we may moreover speak of a restricted Lie algebra structure on $\sM$; again this consists of $\Sing$-linear restricted Lie algebra structures on each $\sM(\Sing)$, such that for any map of $\Ring$-algebras $\Sing \rightarrow \Ting$ the map $\sM(\Sing) \rightarrow \sM(\Ting)$ is an $\Sing$-linear map of restricted Lie algebras.  

\sss \label{sss:adjointactions} With these generalities in hand, we return to the problem of reconstructing the topologies on the Kac--Moody and Virasoro algebras, as well as their adjoint actions by the corresponding group ind-schemes. 

Given a commutative $\Ring$-algebra $\Sing$, note we have a tautological symmetric monoidal functor
$$\Pro(\Ring\lmod) \rightarrow \Pro(\Sing\lmod), \quad \quad \varprojlim M_i \mapsto \varprojlim (M_i \otimes_\Ring \Sing).$$
In particular, given a commutative Hopf algebra $\OOF$ in $\Pro(\Ring\lmod)$,  we have the corresponding base changed commutative Hopf algebra $\OOF \otimes_\Ring \Sing$ in $\Pro(\Sing\lmod)$. As in Lemma \ref{l:changecat}, this induces  homomorphisms of abstract groups 
$$\sF(\Ring) \rightarrow \sF(\Sing)$$
and homomorphism of restricted Lie algebras
$$\ff(\Ring) \rightarrow \ff(\Sing),$$
which intertwines the adjoint actions of $\sF(\Ring)$ and $\sF(\Sing)$. It is straightforward to see that these assignments are compatible with maps of $\Ring$-algebras $\Sing \rightarrow \Ting$, i.e., define a group $\Ring$-space $\sF$ and and restricted Lie algebra $\Ring$-space $\ff$, along with an action of $\sF$ on $\ff$ by Lie algebra automorphisms. 

Moreover, it is straightforward to see that the group $\Ring$-space $\sF$ is simply the group ind-scheme corresponding to $\OOF$, and in particular is representable by an ind-scheme.  

In particular, if $\ff$ is also representable by an ind-scheme over $\Ring$, then we deduce that $\ff$ is naturally a restricted Lie algebra in ind-schemes, and $\sF$ acts on the ind-scheme $\ff$ by restricted Lie algebra automorphisms; in particular the action of $\Ring$-points is naturally continuous with respect to the pro-topology on the $\Ring$-points of $\ff$. Finally, we note that the representability of $\ff$ holds in particular for the Kac--Moody and Virasoro Lie algebras, as they are representable by Tate vector spaces.

To proceed, we further observe that if $\sF$ is a central extension 
 of group $\Ring$-spaces$$1 \rightarrow Z \rightarrow \sF \rightarrow \overline{\sF} \rightarrow 1,$$
then tautologically the adjoint action of $\sF(\Sing)$ on $\ff(\Sing)$ factors through $\overline{\sF}(\Sing)$ for every $\Ring$-algebra $\Sing$, i.e., the adjoint action of $\sF$ on $\ff$ factors through 
an action of the group $\Ring$-space $\overline{\sF}$.

The relevant corollaries of the preceding discussion are then the following. 

\begin{Cor} \label{c:kmadj} The Kac--Moody Lie algebra, i.e., the restricted Lie algebra in ind-schemes over $\Ring$, which assigns to a $\Ring$-algebra $\Sing$ the restricted Lie algebra $\ff(\!(t)\!)^\flat$ over $\Sing$ as in Remark \ref{r:centext}, carries a canonical action of the group ind-scheme $F(\!(t)\!)$ by restricted Lie algebra automorphisms. 
\end{Cor}

\begin{Cor} The Virasoro Lie algebra, i.e., the restricted Lie algebra in ind-schemes over $\Ring$, which assigns to a $\Ring$-algebra $\Sing$ the restricted Lie algebra $\mathfrak{Vir}(\Sing)$ as in Section \ref{SSS_Virasoro_pcenter}, carries a canonical action of the group ind-scheme $\on{Aut}(\D^\times)$ by restricted Lie algebra automorphisms.
\end{Cor}

{
\sss It will be useful below to have an explicit formula, in the setting of Corollary \ref{c:kmadj}, for the adjoint action of $F(\!(t)\!)$ on $\ff(\!(t)\!)^\flat$. I.e., we wish to see that we really do recover the standard formula for the Kac--Moody cocycle from the Tate formalism. Let us check this now.

As notation, let us denote the adjoint actions of $F(\!(t)\!)$ on $\ff(\!(t)\!)$ and $\ff(\!(t)\!)^\flat$ respectively by 
$$\on{Ad}: F(\!(t)\!) \times \ff(\!(t)\!) \rightarrow \ff(\!(t)\!) \quad \text{and} \quad \on{Ad}^\flat: F(\!(t)\!) \times \ff(\!(t)\!)^\flat \rightarrow \ff(\!(t)\!)^\flat.$$
Similarly, as preparation for stating the form of the Kac--Moody cocycle, note that for any $\Sing$-point $g$ of $GL_W(\!(t)\!)$, i.e., an element  $g \in GL_W(\Sing(\!(t)\!))$, we may consider the element 
$$ g^{-1}dg \in \gl_W \otimes_{\Ring} \Sing(\!(t)\!)dt.$$
Moreover, note that we have a natural composition map 
$$- \circ - : (\gl_W \underset{\Ring}\otimes \Sing(\!(t)\!)) \underset{\Sing(\!(t)\!)} \otimes (\gl_W \underset{\Ring}\otimes \Sing(\!(t)\!)dt ) \simeq (\gl_W \otimes_{\Ring} \gl_W) \otimes_{\Ring} \Sing(\!(t)\!)dt \rightarrow \gl_W \underset{\Ring}\otimes \Sing(\!(t)\!)dt$$
and in particular may form the composition 
$$(\gl_W \underset{\Ring}\otimes \Sing(\!(t)\!)) \underset{\Sing(\!(t)\!)} \otimes (\gl_W \underset{\Ring}\otimes \Sing(\!(t)\!)dt ) \xrightarrow{ - \circ - } \gl_W \underset{\Ring}\otimes \Sing(\!(t)\!)dt \xrightarrow{\on{tr}_W \otimes \on{id}} \Ring \otimes_{\Ring} \Sing(\!(t)\!)dt \simeq \Sing(\!(t)\!)dt \xrightarrow{\on{Res}} \Sing.$$

\begin{Lem}\label{Lem:KM_cocycle_group}
Consider the Kac--Moody Lie algebra $\ff(\!(t)\!)^\flat$, with its splitting of underlying $\Ring$-modules in $\Ring$-spaces 
$$\ff(\!(t)\!)^\flat \simeq \ff(\!(t)\!) \oplus \Ring \cdot 1$$
as in Proposition \ref{p:kacmoodyrestr}. Then with respect to this splitting, the adjoint action of $F(\!(t)\!)$ on $\ff(\!(t)\!)^\flat$ is given on $\Sing$-points by the formulas
\begin{equation}\label{e:kmformula} \on{Ad}_g^\flat(1) = 1 \quad \text{and} \quad \on{Ad}_g^\flat(X)  = \on{Ad}_g(X) + \on{Res} (\on{tr}_W( X \circ g^{-1}dg)) \cdot 1, \quad \quad X \in \ff(\Sing(\!(t)\!)), \quad g \in F(\Sing(\!(t)\!)).  \end{equation}

\end{Lem}

\begin{proof} Applying Corollary \ref{c:tategroupcocycle} to the elementary Tate module 
$$W \otimes_{\Ring}  \Sing(\!(t)\!) \simeq (W \otimes_{\Ring} \Sing[[t]]) \oplus (W \otimes_{\Ring} t^{-1} \Sing[t^{-1}]) =: W(\!(t)\!)^+ \oplus W(\!(t)\!)^-,$$
we have that 
$$\on{Ad}^\flat_g(X) = \on{Ad}_g(X)  + \on{tr}_{W \otimes_{\Ring} \Sing(\!(t)\!)}( (g \circ X \circ g^{-1})_- - g \circ X_- \circ g^{-1}) \cdot 1,$$
i.e., we must show the equality 
\begin{equation} \label{e:opensesame} \on{tr}_{W \otimes_{\Ring}  \Sing(\!(t)\!)}( (g \circ X \circ g^{-1})_- - g \circ X_- \circ g^{-1}) \overset{?}= \on{Res} (\on{tr}_W( X \circ g^{-1}dg)).\end{equation}
We will check this by a direct computation. Namely, let us consider the Laurent expansions 
$$X = \underset{i} \Sigma \hspace{.5mm} X_i \otimes t^i, \quad g = \underset{j} \Sigma \hspace{.5mm} g_j \otimes t^j, \quad g^{-1} = \underset{k} \Sigma \hspace{.5mm} (g^{-1})_k \otimes t^k, \quad \quad i, j, k \in \Z, \quad X_i, g_j, (g^{-1})_k \in \gl_W.$$
If we write $\delta_{< 0}: \Z \rightarrow \Z$ for the indicator function of the negative integers, i.e., 
$$\delta_{< 0}(i) = \begin{cases} 1 & \text{if } i < 0, \\ 0 & \text{otherwise}. \end{cases},$$
    by considering the subspaces of $W(\!(t)\!)$ of the form $W \otimes t^i, i \in \Z$, we obtain that 
\begin{align*} \on{tr}_{W \otimes_{\Ring}  \Sing(\!(t)\!)}( (g \circ X \circ g^{-1})_- - g \circ X_- \circ g^{-1}) & = \underset{i, j, k} \Sigma \hspace{.5mm}  \on{tr}_W(g_{i - k} \circ X_{k - j} \circ (g^{-1})_{j - i}) \cdot ( \delta_{< 0}(i) - \delta_{< 0}(k)). \intertext{If we make the change of variables $a = j-i$, $b = k - j$, so that $k = i - a - b$, we may rewrite this as} &= \underset{i, a, b} \Sigma \hspace{.5mm}  \on{tr}_W(g_{-a-b} \circ X_{b} \circ (g^{-1})_{a}) \cdot ( \delta_{< 0}(i) - \delta_{< 0}(i - a- b)), \intertext{where the sum runs over $i, a, b \in \Z$. Summing over $i$, note that for fixed $a,b$ we have $$\underset{i} \Sigma \hspace{.5mm} \delta_{< 0}(i) - \delta_{<0}(i - a - b) =    - a - b,$$we may continue as} &= \underset{a,b} \Sigma \hspace{.5mm} (-a-b) \cdot \on{tr}_W(g_{-a-b} \circ X_b \circ (g^{-1})_a) \intertext{which using the cyclicity of the trace we rewrite as} &= \underset{a,b} \Sigma \hspace{.5mm} (-a-b) \cdot \on{tr}_W( X_b \circ (g^{-1})_a \circ g_{-a-b})\intertext{Noting that explicitly $dg = \underset{j} \Sigma \hspace{.5mm} j \cdot g_j \otimes t^{j-1}dt$, and $g^{-1} dg =  \underset{i,j} \Sigma \hspace{.5mm} j \cdot (g^{-1})_i \circ g_j \otimes t^{i + j - 1}dt$, it is now straightforward to recognize this as} &= \on{Res}(\on{tr}_W( X \circ g^{-1} dg)), \end{align*}
as desired.    
\end{proof}

}

\section{Basics on vertex algebras}\label{S_vertex}
The goal of this section is to recall some basic definitions, constructions and results concerning vertex algebras over general rings.  

\subsection{Vertex algebras}\label{SS_vertex}
Let $\Ring$ be a commutative associative unital ring.

\sss
The following definition of vertex algebras over $\Ring$, as far as we know, is due to Borcherds and Ryba \cite{BR}. 

\begin{defi}\label{defi:vertex_ring} A vertex algebra $V$ over $\Ring$ is an $\Ring$-module equipped with a map of $\Ring$-modules
$$V \otimes_{\Ring} V \rightarrow V(\!(z)\!), \quad \quad a \otimes b \mapsto Y(a,z)b := \sum_n a_{(n)}b z^{-n-1},$$
(known as the {\it state-field correspondence}),
and a vacuum vector
$|\varnothing\rangle\in V$
satisfying the following identities.

\begin{enumerate}

   \item For any $a \in V$, we have that
$Y(|\varnothing \rangle,z) = \operatorname{id}$, i.e., that $$| \varnothing \rangle_{(n)} = \begin{cases} \operatorname{id}_V, & \text{if } n = -1, \text{ and} \\ 0, & \operatorname{otherwise}.\end{cases}$$

\item For any $a\in V$, we have, in addition, that 
$$Y(a,z) | \varnothing \rangle \in a + z V[[z]],$$
i.e. that $a_{(n)} | \varnothing \rangle = 0$ for $n \geqslant 0$, and $a_{(-1)} | \varnothing \rangle = a$.

\item (Locality) For any $a, b$ in $V$, there exists an integer $N \gg 0$ such that for any $c$ in $V$ we have
$$(z - w)^N \left(Y(a,z)Y(b,w)c- Y(b,w)Y(a,z) c\right)=0$$
in $V[[z^{\pm 1},w^{\pm 1}]]$ (the set of series in $z,w$ that are infinite in both directions).

\item (Translation) For $n \geqslant 0$, define the operator $T^{(n)}:V\rightarrow V$ by
$$T^{(n)} a := a_{(-n-1)} | \varnothing \rangle.$$
Then one has the identity
$$T^{(n)}Y(a,z) = \partial_z^{(n)}Y(a,z),$$
i.e.,  for any $j \in \mathbb{Z}$ one has the equality of endomorphisms of $V$
\begin{equation*}
(T^{(n)}a)_{(j)}  = (-1)^n {j\choose n} a_{(j-n)}. 
\end{equation*}
\end{enumerate}
\end{defi}

One defines a homomorphism of vertex algebras in an obvious way. 


\begin{defi}\label{defi:vertex_commutative}
We say that $V$ is {\it commutative}, if $a_{(n)}=0$ for $n\geqslant 0$.
\end{defi}

\sss
We note that the locality axiom, (3), is equivalent to the following formula, cf. \cite[(2.3-7)]{Frenkel_loop}:
\begin{equation}\label{eq:Y_bracket}
[Y(a,z),Y(b,w)]=\sum_{k\geqslant 0}Y(a_{(\ell)}b,w)\partial_w^{(\ell)}\delta(z-w),
\end{equation}
where $\delta(z-w)=\sum_{m\in \Z} z^mw^{-m-1}$. Equivalently, 
\begin{equation}\label{eq:VA_modes_brackets}
[a_{(m)},b_{(n)}]=\sum_{\ell\geqslant 0}{m\choose \ell}(a_{(\ell)}b)_{(m+n-\ell)}.
\end{equation}

We will also need the following identity that holds for $n<0$, cf. \cite[(2.3-4)]{Frenkel_loop}:
\begin{equation}\label{eq:Y_product}
Y(a_{(n)}b,z)=:[\partial_z^{(-n-1)}Y(a,z)]Y(b,z):,
\end{equation}
where $:-:$ indicates the normally ordered product. 

Later on in this section we will give some examples. 

\sss\label{SSS_vertex_filtration}
We will need to deal with filtrations on vertex algebras. 

\begin{defi}
Let $V$ be a vertex algebra over $\Ring$. By a {\it filtration} on $V$ we mean an exhaustive ascending $\Ring$-module filtration $V=\bigcup_{i\geqslant 0}V_{\leqslant i}$ satisfying the following conditions:
\begin{enumerate}
\item $| \varnothing \rangle\in V_{\leqslant 0}$,
\item and for $a\in V_{\leqslant i}, b\in V_{\leqslant j}$, one has $a_{(n)}b\in V_{\leqslant i+j}$ for all $n\in \Z$. 
\end{enumerate}
If in (2), we additionally have $a_{(n)}b\in V_{\leqslant i+j-1}$
for $n\geqslant 0$, then we say that the filtered vertex algebra $V$ is {\it almost commutative}. 
\end{defi}

\sss
We also have two  notions of a grading on a vertex algebra, $V=\bigoplus_{i\in \Z}V_i$.  An {\it energy grading} is defined as in  \cite[Section 2.2.2]{Frenkel_loop}: we require that $|0\rangle\in V_0$, while for $a\in V_i, b\in V_j$ we have $a_{(n)}b\in V_{i+j+1-n}$. We will also consider {\it naive gradings}, where under the above assumptions we require $a_{(n)}b\in V_{i+j}$. 
For a filtered vertex algebra, it makes sense to speak about its associated graded. This is a vertex algebra which carries a naive grading. The filtered vertex algebra is almost commutative if and only if its associated graded vertex algebra is commutative.

\sss If $V$ is a filtered vertex algebra over $\Ring$, then we can form its Rees algebra $R_\hbar(V)$, a naively graded vertex algebra over $\Ring[\hbar]$. As usual, we have natural isomorphisms
$R_\hbar(V)/(\hbar-1)\cong V, R_\hbar(V)/(\hbar)\cong \gr V$. Below we will often write $V^\hbar$ for $R_\hbar(V)$ and $V^0$ for $\gr V$.

\sss
We note that there is an obvious notion of a vertex algebra homomorphism, so vertex algebras form a category. Thanks to this it makes sense to speak about vertex subalgebras and vertex algebra ideals. 

Also, for two vertex algebras $V^1,V^2$ we can form their tensor product $V^1 \otimes_{\Ring} V^2$ of $\Ring$-modules, and this tensor product carries a natural vertex algebra structure (with $Y(a^1\otimes a^2, z)(b^1\otimes b^2)=Y(a^1,z)b^1\otimes Y(a^2,z)b^2$). 

Finally, we mention that there is an obvious base change construction: if $V$ is a vertex algebra over $\Ring$, and $\tilde{\Ring}$ is a commutative $\Ring$-algebra, then $\tilde{\Ring}\otimes_{\Ring}V$ is naturally a vertex algebra over $\tilde{\Ring}$.

\sss\label{SSS_Poisson_vertex}
One can define the notion of a Poisson vertex algebra over $\Ring$ similarly to \cite[Section 3.3]{ATV}. By definition, a Poisson vertex algebra is a tuple $(V^0,|\varnothing\rangle, Y_+,Y_-)$ satisfying certain axioms. First, we require that $(V^0,\varnothing,Y_+)$ is a commutative vertex algebra, in particular, giving operators $T^{(n)}$
for all $n\geqslant 0$. The operation $Y_-: V^0\rightarrow z^{-1}\operatorname{End}(V^0)[z^{-1}]$ is a vertex algebra analog of the Poisson bracket. It is required to satisfy the following axioms:
\begin{enumerate}
\item $Y_-(T^{(n)}a,z)=\partial^{(n)}Y_-(a,z)$.
\item $Y_-(a,z)b=\sum_{n=0}^\infty z^n T^{(n)}Y_-(b,-z)a$.
\item If we write $Y_-(a,z)$ as $\sum_{n\geqslant 0}a_{(n)}z^{-n-1}$, then $$[a_{(n)},b_{(m)}]=\sum_{\ell=0}^\infty {m\choose \ell}(a_{(\ell)}b)_{(m+n-\ell)}$$ for all $m,n\geqslant 0$.
\item $Y_-(a,z)(b_{(-1)}c)=(Y_-(a,z)b)_{(-1)}c+(Y_-(a,z)c)_{(-1)}b$ that should be thought of as a vertex algebra analog of the Leibniz identity.
\end{enumerate}

The main example for us is as follows. Let $(V,|\varnothing\rangle, Y)$
be an almost commutative filtered vertex algebra. We take $V^0:=\operatorname{gr}V$. For $Y_+$ (resp., $Y_-$) we take the top degree terms of $\sum_{n\geqslant 0}a_{(-n-1)}z^n$ (resp., 
$\sum_{n<0}a_{(-n-1)}z^n$). Note that $Y_+$ has degree $0$ and $Y_-$
has degree $-1$.

\subsubsection{}\label{SSS_VA_modules} Our final topic in this section is modules over vertex algebras. Let $V$ be a vertex algebra over $\Ring$. By a $V$-module structure on an $\Ring$-module $M$ we mean a map $V\rightarrow \operatorname{End}_{\Ring}(M)[[z^{\pm 1}]], a\mapsto Y_M(a,z)$ satisfying the following conditions:
\begin{enumerate}
\item $Y_M(a,z)m\in M(\!(z)\!)$ for all $a\in V, m\in M$.
\item $Y(|\varnothing\rangle,z)=\operatorname{id}_M$.
\item $Y_M(T^{(n)}a,z)=\partial_z^{(n)}Y_M(a,z)$ for all $n>0$ and $a\in V$.
\item $[Y_M(a,z),Y_M(b,w)]=\sum_{k\geqslant 0}Y_M(a_{(\ell)}b,w)\partial_w^{(\ell)}\delta(z-w)$, for all $a,b\in V$.
\end{enumerate}
For example, $V$ becomes a module over itself. 

Similarly to \cite[(2.3-4)]{Frenkel_loop}, we have 
\begin{equation}\label{eq:module_normal_product}
Y_M(a_{(n)}b,z)=:[\partial_z^{(-n-1)}Y_M(a,z)]Y_M(b,z):.
\end{equation}


\subsection{Universal enveloping algebras}\label{SS_universal_enveloping}
In this section we will review the construction of the universal enveloping algebra associated to a vertex algebra $V$. The case when $\Ring$ is a characteristic $0$ field can be found, e.g., in \cite[Section 3.2]{Frenkel_loop}, and the general case is easily adapted from there. 

\sss
First, we construct the Lie algebra  $F_V$ of ``Fourier modes''. Consider the $\Ring[t^{\pm 1}]$-module $V[t^{\pm 1}]$. For $n\in \Z_{\geqslant 0}$, define the endomorphism $\partial^{(n)}$ of $V[t^{\pm 1}]$ by
$$\partial^{(n)}=\sum_{j=0}^n T^{(j)}\partial_t^{(n-j)},$$
where $\partial_t^{(n-j)}$ is the divided power of the derivative with respect to $t$.

By the {\it Lie algebra of Fourier modes} one means 
$$F_V:=V[t^{\pm 1}]/\left(\sum_{n\geqslant 0}\operatorname{im}\partial^{(n)}\right)$$
We write $a_{[n]}$ for the image of $at^n$ in $F_V$. Then $F_V$ is the $\Ring$-module with generators $a_{[n]}$ and relations 
$(T^{(n)}a)_{[j]}=(-1)^n{j\choose n}a_{[j-n]}$.

The bracket on $F_V$ is uniquely recovered from (cf. \cite[(3.2-2)]{Frenkel_loop}) 
$$[a_{[m]},b_{[n]}]=\sum_{\ell\geqslant 0}{m\choose \ell}(a_{(\ell)} b)_{[m+n-\ell]}.$$
The proof that $F_V$ is a Lie algebra with respect to this bracket repeats that of 
\cite[Proposition 3.2.1]{Frenkel_loop}.


\sss\label{SSS:Utilde}
Now assume that $V=\bigoplus_{i\in \Z}V_i$ is an energy graded vertex algebra. We form the completed universal enveloping algebra $\widetilde{U}(F_V)$, the inverse limit 
\begin{equation}\label{eq:inverse_limit}
\varprojlim_{m\rightarrow \infty} U(F_V)/I_m,
\end{equation}
where $I_m$ is the left ideal in $U(F_V)$ spanned by the element $a_{[n]}$, where $a\in V_i$ with $n\geqslant m+i$. Then $\widetilde{U}(F_V)$  is an associative algebra. 

The (completed) universal enveloping algebra $\widetilde{U}(V)$ is  defined as the quotient of $\widetilde{U}(F_V)$ by relations (\ref{eq:univ_envel_add_relations})
below. To state this relation we introduce the notation $Y[a,z]:=\sum_{n\in \Z}a_{[n]}z^{-n-1}\in F_V[[z^{-1},z]]$. The relations are as follows:
\begin{equation}\label{eq:univ_envel_add_relations}Y[a_{(-1)}b,z]=:Y[a,z]Y[b,z]:.
\end{equation}

We note that the assignment $V\mapsto \widetilde{U}(V)$ is a functor from the category of energy graded vertex algebras to the category of topological associative algebras. Namely, let $\varphi:V^1\rightarrow V^2$ be a vertex algebra homomorphism. Then there is a unique topological algebra homomorphism
$\widetilde{U}(\varphi): \widetilde{U}(V^1)\rightarrow \widetilde{U}(V^2)$ such that
\begin{equation}\label{eq:U_homom}
\widetilde{U}(\varphi)(Y[a,z])=Y[\varphi(a),z], \forall a\in V^1.
\end{equation}

\subsubsection{} In the setting of Section \ref{SSS:Utilde}, a module over $V$ is the same thing as a discrete module over the topological algebra $\tilde{U}(V)$. Indeed, every discrete $\tilde{U}(V)$-module is a $V$-module: $Y_M(a,z)$
is the action of $\sum_{n\in \Z}a_{[n]}z^{-n-1}$. Conversely, let $M$ be a $V$-module. Thanks to
(2)-(4) in Section \ref{SSS_VA_modules}, $M$ becomes a $U(F_V)$-module. Thanks to (1) there, $M$ is a discrete module over (\ref{eq:inverse_limit}). And using (\ref{eq:module_normal_product}), we see that the action of (\ref{eq:inverse_limit}) on $M$ factors through $\tilde{U}(V)$.

\subsubsection{} Now we discuss the compatibility between the universal enveloping algebras and tensor products. 
\begin{Lem}\label{Lem:univ_env_tensor}
For two vertex algebras $V^1,V^2$, we have a natural isomorphism of topological algebras (where in the source we have the completed tensor product):
$$\widetilde{U}(V^1)\, \widetilde{\otimes}_{\Ring} \,\widetilde{U}(V^2) 
\xrightarrow{\sim} \widetilde{U}(V^1
\otimes_{\Ring} V^2).$$
\end{Lem}
\begin{proof}
First, note that the vertex algebra homomorphisms $\iota_i:V^i\rightarrow V^1\otimes_{\Ring}V^2$ (e.g. $\iota_1(a^1)=a^1\otimes |\varnothing\rangle$) give rise to topological algebra homomorphisms $\tilde{U}(\iota_i): \tilde{U}(V^i)\rightarrow \tilde{U}(V^1\otimes_\Ring V^2)$ with commuting images. So they extend to a homomorphism $\tilde{U}(V^1)\otimes_\Ring \tilde{U}(V^2)\rightarrow \tilde{U}(V^1\otimes_\Ring V^2)$ that is uniquely characterized by the condition that it sends $a^i_{[n]}$ to $a^i_{[n]}$ for all $i=1,2, n\in \Z, a^i\in V^i$. This homomorphism extends to  $\widetilde{U}(V^1)\, \widetilde{\otimes}_{\Ring} \,\widetilde{U}(V^2)$ by continuity. 

On the other hand, we have a base of neighborhoods of $0$ consisting of left ideals $I_m\subset \widetilde{U}(V^1)\, \widetilde{\otimes_{\Ring}} \,\widetilde{U}(V^2)$ such that 
$$\widetilde{U}(V^1)\, \widetilde{\otimes_{\Ring}} \,\widetilde{U}(V^2)=\varprojlim \widetilde{U}(V^1)\, \widetilde{\otimes_{\Ring}} \,\widetilde{U}(V^2)/I_m.$$
Namely, for $I_m$ we take the closure of the left ideal generated by elements of energy degree $\geqslant m$.

Each space $\widetilde{U}(V^1)\, \widetilde{\otimes_{\Ring}} \,\widetilde{U}(V^2)/I_m$ carries commuting $V^i$-module structures for $i=1,2$. They combine into a $V^1\otimes_\Ring V^2$-module structure, and hence give a $\widetilde{U}(V^1\otimes_\Ring V^2)$-module structure. This gives a topological module  homomorphism of action on $1$:
\begin{equation}\label{eq:map_from_tensor}
\widetilde{U}(V^1\otimes_\Ring V^2)\rightarrow 
\widetilde{U}(V^1) \widetilde{\otimes}_{\Ring} \,\widetilde{U}(V^2).
\end{equation}

It is easy to see that the action of $\tilde{U}(V^1\otimes_\Ring V^2)$ on  $\widetilde{U}(V^1)\widetilde{\otimes}_{\Ring} \,\widetilde{U}(V^2)$ commutes with the right multiplication action of $\widetilde{U}(V^1)\, \widetilde{\otimes_{\Ring}} \,\widetilde{U}(V^2)$. It follows that (\ref{eq:map_from_tensor}) is an algebra homomorphism. 
This homomorphism is uniquely characterized by the property that $a^i_{[n]}$ is sent to $a^i_{[n]}$ for all $n\in \Z$ and $a^i\in V^i$. So we get mutually inverse topological algebra homomorphisms proving the isomorphism in the statement of the lemma. 
\end{proof}

\sss\label{SSS:filtered_enveloping}
Now we discuss a variant of $\widetilde{U}(V)$ in the presence of filtrations. Suppose that $V$ is, in addition, filtered, $V=\bigcup_{j}V_{\leqslant j}$, where the filtration is compatible with the energy grading, which means that $V_{\leqslant j}$ is graded for all $j$. This gives rise to an exhaustive Lie algebra filtration on $F_V$ (with $F_{V,\leqslant j}$ being the span of $a_{[n]}$ for $n\in \Z$ and $a\in V_{\leqslant j}$) and hence to an ascending associative algebra  filtration $\widetilde{U}(V)_{\leqslant i}\subset \widetilde{U}(V)$. We note that the latter is non-exhaustive, and we define $\widehat{U}(V)$ as $\bigcup_i \widetilde{U}(V)_{\leqslant i}$.  This is a topological algebra but it is not complete. Instead, it is {\it filtered complete} in the sense of the following definition. 

\begin{defi}\label{defi:filtered_complete}
Let $A$ be a topological associative algebra equipped with an exhaustive ascending algebra filtation $A=\bigcup_{i\geqslant 0}A_{\leqslant i}$. We say that $A$ is {\it filtered complete} if
\begin{itemize}
\item 
the topology on each $A_{\leqslant i}$ is complete and separated, 
\item the inclusion $A_{\leqslant i}\hookrightarrow A_{\leqslant i+1}$ is closed for any $i$, and $A$ has the topology of the union of its filtered pieces.
\end{itemize}
Similarly, if $A$ is equipped with a grading, $A=\bigoplus_{i\geqslant 0}A_i$, then we say that $A$ is {\it graded complete} if the topology on each $A_i$ is complete and separated and the topology on $A$ is that of the direct sum. 
\end{defi}

The naive grading on $\gr V$ induces a grading on $\widehat{S}(\gr V):=\widehat{U}(\gr V)$ and we have a natural epimorphism $\widehat{S}(\gr V)\twoheadrightarrow \gr \widehat{U}(V)$. The algebra $\widehat{S}(\gr V)$ is graded complete. Note  that if the filtration on $V$ is almost commutative, then so is the filtration on $\widehat{U}(V)$. 

The grading on $R_\hbar(V)$ gives rise to a grading on $\widehat{U}^\hbar(V):=\widehat{U}(R_\hbar(V))$, and 
$\widehat{U}_\hbar(V)$ is identified with the Rees algebra of $\widehat{U}(V)$.

\sss\label{SSS_from_univ_to_vertex} Suppose that $\varphi:V^1\rightarrow V^2$ be a filtration preserving homomorphism between filtered vertex algebras. Set $\Phi:=\widehat{U}(\varphi):\widehat{U}(V^1)\rightarrow \widehat{U}(V^2)$. Here we explain how to recover $V^i$ from $\widehat{U}(V^i)$ and $\varphi$ from $\Phi$. One can identify $V^i$ with the quotient of $\widehat{U}(V^i)$ by the left ideal generated by the images of the elements $a_{[n]}$ for $a\in V^i$ with $n\geqslant 0$. The homomorphism $\Phi$ sends the image of $a_{[n]}$ to the image of $\varphi(a)_{[n]}$ for all $n$ and hence induces a  $\widehat{U}(V^1)$-module homomorphism $V^1\rightarrow V^2$. It is easy to see that the induced homomorphism coincides with $\varphi$.

\subsection{Arc and loop spaces and groups}\label{SS_arcs_loops}
Let $X$ be a finite type affine scheme over $\Ring$. 
\sss
Recall that we have another affine $\Ring$-scheme $\Jet X$, the {\it jet} or {\it arc space} of $X$, representing the functor $\mathsf{S}\mapsto X(\mathsf{S}[[t]])$. This is an affine scheme of infinite type. This scheme, by construction, is the inverse limit of the finite type schemes $\Jet_nX$, the $n$th jet schemes. 

Note that $\Ring[\Jet X]$ carries a natural structure of a commutative vertex algebra. In order to see this, we first equip $\Ring[\Jet \mathbb{A}^n]$ with operators $T^{(\ell)}$ for $\ell\geqslant 0$. Namely, let $x^1,\ldots,x^n$ be the coordinates on $\mathbb{A}^n$. Then, by construction, $\Jet\mathbb{A}^n$ is an infinite dimensional affine space with coordinates $x^i_j$ with $i=1,\ldots,n$, and $j< 0$. Then there is a unique collection of operators $T^{(\ell)}: \Ring[\Jet\mathbb{A}^n]\rightarrow \Ring[\Jet\mathbb{A}^n]$ subject to the following conditions
\begin{itemize}
\item $T^{(\ell)}x^i_j=(-1)^\ell {j\choose \ell}x^i_{j-\ell}, \forall i=1,\ldots,n,j\in \Z_{<0},\ell\in \Z_{\geqslant 0}$,
\item $T^{(\ell)}(fg)=\sum_{j=0}^k (T^{(j)}f)(T^{(\ell-j)}g), \forall f,g\in \Ring[\Jet \mathbb{A}^n],\ell\in \Z_{\geqslant 0}$.
\end{itemize}

Now suppose that $X$ is given in $\mathbb{A}^n$ by equations $F^1,\ldots,F^\ell$. Then $\Jet X$ is given in $\Jet\mathbb{A}^n$ by equations $T^{(\ell)}F(x^1_{-1},\ldots,x^n_{-1})$. In particular, the operators $T^{(\ell)}$ descend to 
$\Ring[\Jet X]$. Moreover, they are independent of the choice of an embedding $X\hookrightarrow  \mathbb{A}^n$. 

Now we can define a vertex algebra structure on $\Ring[\Jet X]$: we set $$Y(a,z)b=(\sum_{j=0}^\infty T^{(j)}az^j)b.$$

\begin{Rem}\label{Rem:jet_vertex_coordinate_free}
Here is a coordinate independent way to think about 
the vertex algebra structure on $\Ring[\Jet X]$. 
For $a\in \Ring[\Jet X]$, we define $a_{-z}\in \Ring[\Jet X][[z]]$ as follows. Let $\tilde{\Ring}$
be an $\Ring$-algebra. Note that a regular function on a scheme, $\mathscr{X}$, is the same as a morphism $\mathscr{X}\rightarrow \mathbb{A}^1$, equivalently a functorial (in $\tilde{\Ring}$) map $\mathscr{X}(\tilde{\Ring})\rightarrow \tilde{\Ring}$. An $\tilde{\Ring}$-point of $\Jet X$ is an $\tilde{\Ring}[[t]]$-point of $X$, denote it by $x(t)$. So to specify $a_{-z}$ we need to define $a(x(t))\in \tilde{\Ring}[[z]]$. Expanding $x(t-z)$, a formal power series in $t-z$, in the powers of $z$, we can view $x(t-z)$ as an $\tilde{\Ring}[[z]]$-point of $\Jet X$. 
We define $a_{-z}$ by $x(t)\mapsto a(x(t-z))\in \tilde{R}[[z]]$. Equivalently, $a_{-z}=(\sum_{j=0}^\infty z^j T^{(j)})a$, so $Y(a,z)$ is the multiplication by $a_{-z}$.
\end{Rem}
The assignments $X\mapsto \Jet X,\Jet_nX$ are easily seen to be functorial, and $\Jet(X\times Y)$ is naturally identified with $\Jet X\times \Jet Y$, the same for $\Jet_n$. 

 Remark \ref{Rem:jet_vertex_coordinate_free} implies that $X\rightarrow \Ring[\Jet X]$ can be viewed as a functor from the category of finite type affine schemes to the category of commutative vertex algebras. 

\sss We now discuss arc groups.
\begin{Rem}\label{Rem:arc_group_vertex}
If $F$ is a group scheme over $\Ring$, then 
$\Jet F$ and $\Jet_nF$ are also group schemes over $\Ring$. It is straightforward to check that the coproduct 
$\Ring[\Jet F]\rightarrow \Ring[\Jet F]\underset{\Ring}\otimes\Ring[\Jet F]$ is a vertex algebra homomorphism. It follows that $\Ring[\Jet F]$ is a Hopf algebra object in the category of vertex algebra. 
\end{Rem}

We will mostly care about the situation when $F$ is a split connected reductive group $G$ over $\Ring$.  We have the following lemma. 

\begin{Lem}\label{Lem:jet_smoothness_flatness}
The following claims hold:
\begin{enumerate}
\item The scheme $\Jet_n G$ is smooth over $\operatorname{Spec}(\Ring)$.
\item The scheme $\Jet G$ is flat over $\operatorname{Spec}(\Ring)$.
\end{enumerate}
\end{Lem}
\begin{proof} 
We note that $G$ is smooth over $\operatorname{Spec}(\Ring)$ (to see this, one can, for example, cover $G$ by translates of the open Bruhat cell, which is smooth as the product of an affine space and an algebraic torus, note that the open Bruhat cell is defined over the integers).
Since $G$ over $\Ring$ satisfies infinitesimal lifting property, so does $\Jet_n G$ over $G$, and (1) follows. (2) follows from (1) because $\Ring[\Jet G]$ is a filtered colimit of the algebras $\Ring[\Jet_n G]$.      
\end{proof}

\sss
Now we proceed to loop spaces and loop groups that already made an appearance in Section \ref{SSS_KM_groups}. 
As before, let $X$ be a finite type affine scheme. We write $\Loop X$ for the loop space of $X$. This is a strict ind-affine ind-scheme of countable type, so it is given by its algebra of functions, that is a complete and separated topological algebra with a countable basis of open neighborhoods of zero that are ideals. 

\begin{Lem}\label{Lem:loop_vs_univ_envel}
The algebra $\Ring[\Loop X]$ is naturally identified with the universal enveloping algebra $\widetilde{U}(\Ring[\Jet X])$. 
\end{Lem}
\begin{proof}
One first checks this for $X=\mathbb{A}^n$ and then shows that $\Ring[\Loop X]$ and  $\widetilde{U}(\Ring[\Jet X])$ are quotients of 
$\Ring[\Loop \mathbb{A}^n]=\widetilde{U}(\Ring[\Jet \mathbb{A}^n])$ by the same relations. 
\end{proof}

Of course, if $X=F$ is a group scheme over $\Ring$, then $\Loop F$ is the same $F(\!(t)\!)$ from Section \ref{SSS_KM_groups}. 

\sss\label{SSS_Poisson_structure_uniqueness}
Assume now that $\Ring[X]$ is a Poisson algebra. The axioms for $Y_-$ listed in Section \ref{SSS_Poisson_vertex} imply that there is at most one structure of a Poisson vertex algebra on $\Ring[J X]$ such that
for $a,b\in \Ring[X]\subset \Ring[JX]$ we have $a_{(0)}b=\{a,b\}$. In fact, one can show that this structure always exists, but we will not need this. 

\subsection{Affine vertex algebras}\label{SS_affine_vertex}
We now proceed to the most important class of vertex algebras appearing in this paper. 

\sss \label{SSS_affine_VA_defn}
Assume $2$ is invertible in $\Ring$.
Let $\ff$ be a Lie algebra over $\Ring$ that is a finitely generated projective $\Ring$-module. Fix a representation $W$ of $\ff$, a finitely generated projective $\Ring$-module. Let $\beta_W$ denote the trace form for this representation (in particular, $\beta_{\ff}$ is the Killing form).  We define the loop algebra $\mathfrak{f}(\!(t)\!)$ and its central extension $\hat{\mathfrak{f}}=\mathfrak{f}((t))\oplus \Ring \mathbf{1}\oplus \Ring \mathbf{1}'$, where $\mathbf{1},\mathbf{1}'$ are central and for $x_1 f_1,x_2 f_2\in \mathfrak{f}(\!(t)\!)$ (with $x_i\in \mathfrak{f}, f_i\in \Ring(\!(t)\!)$) their bracket in $\hat{\mathfrak{f}}$ is given by (cf. Remark \ref{r:centext}):
\begin{equation}
[x_1 f_1,x_2f_2]:=[x_1,x_2]f_1f_2-\operatorname{Res}_{t=0}(f_1df_2)\left(\beta_W(x_1,x_2)\mathbf{1}-\frac{1}{2}\beta_{\ff}(x_1,x_2)\mathbf{1}'\right). 
\end{equation}
The algebra $\hat{\mathfrak{f}}$ is the {\it affine Lie algebra} associated with $\mathfrak{f}$. 

\begin{Rem}\label{Rem:general_setting}
We will be interested in the class that appeared in Section \ref{SSS_KM_definition}: when $\ff$ is the Lie algebra of an algebraic group $F$ and $W$ is a representation of $F$.  
\end{Rem}

Pick $\kappa\in \Ring$.
From $\hat{\mathfrak{f}}$ we can produce the {\it affine vertex algebra} to be denoted by 
$V_\kappa(\mathfrak{f})$. Namely, consider the $\mathfrak{f}[[t]]\oplus \Ring \mathbf{1}\oplus \Ring \mathbf{1}'$-module $\Ring_\kappa$, where $\mathfrak{f}[[t]]$ acts by $0$, $\mathbf{1}$ acts by $\kappa$, and $\mathbf{1}'$ acts by $1$. Then the induced module $V_\kappa(\mathfrak{f})$
comes with a natural energy $\Z_{\leqslant 0}$-graded vertex algebra structure, cf. \cite[Theorem 2.2.2]{Frenkel_loop} and \cite[Section 3.4]{ATV}.

We can also consider the universal version $V_{\bone}(\ff)$, the induction of $\Ring$ from $\ff[[t]]\oplus \Ring \mathbf{1}'$, where $\mathbf{1}'$ acts by $1$. This is a vertex algebra over $\Ring[\mathbf{1}]$. Specializing $\mathbf{1}$ to $\kappa$, we recover $V_\kappa(\ff)$. 

\sss
The standard PBW filtration on $U(\hat{\mathfrak{f}})$ gives rise to a filtration on $V_\kappa(\mathfrak{f})$ and this is easily seen to be a vertex algebra filtration in the sense of Section \ref{SSS_vertex_filtration}. The associated graded vertex algebra is identified with $\Ring[\Jet \ff^*]$. 

In the same way, we get a filtration on $V_{\bone}(\ff)$ (note that $\mathbf{1}$ is in filtration degree $1$). The algebra $\gr V_{\bone}(\ff)$ can be interpreted as the algebra of regular functions on the space $\mathbb{A}^1\times \Jet \ff^*$, the space of ``$\lambda$-connections'' $\{\lambda\partial+ \ff^*[[t]]dt\}$, where $\mathbf{1}$ viewed as a function sends a point $\lambda\partial+ \ff^*[[t]]dt$ to $\lambda$.

\sss
Now we turn to the universal enveloping algebra of $V_\kappa(\ff)$.
Consider the algebra $U_\kappa(\mathfrak{f}):=U(\hat{\mathfrak{f}})/(\mathbf{1}-\kappa, \mathbf{1}'-1)$, its completion $\widetilde{U}_\kappa(\hat{\mathfrak{f}})$,
and the filtered part 
$\widehat{U}_\kappa(\hat{\mathfrak{f}})$ of the latter defined in Section \ref{SSS:filtered_enveloping}.
We note that the filtration on each filtered piece is complete and separated. 

Similarly to \cite[Lemma 3.2.2]{Frenkel_loop}, we see that there is a natural topological algebra isomorphism 
$$\widetilde{U}(V_\kappa(\mathfrak{f}))\xrightarrow{\sim} \widetilde{U}_\kappa(\hat{\mathfrak{f}}).$$
This isomorphism is compatible with the filtrations yielding a filtered algebra isomorphism
$$\widehat{U}(V_\kappa(\mathfrak{f}))\xrightarrow{\sim} \widehat{U}_\kappa(\hat{\mathfrak{f}}).$$

\sss Now suppose that $F$ is an algebraic group scheme over $\Ring$ with Lie algebra $\mathfrak{f}$. We write $F(\Ring(\!(t)\!))$ for the group of $\Ring(\!(t)\!)$-points of $F$. Suppose that $\beta$ is $F$-invariant. According to Lemma \ref{Lem:KM_cocycle_group} the adjoint action of the group  
$F(\Ring(\!(t))\!)$ on $\hat{\mathfrak{f}}$ is via
\begin{equation}\label{eq:loop_adjoint_action}
g(t).x(t)=\operatorname{Ad}(g(t))x(t)+ \operatorname{Res}_{t=0}(\beta_W\mathbf{1}-\frac{1}{2}\beta_{\ff}\mathbf{1}')(g(t)^{-1} dg(t),x(t)).
\end{equation}
Note that this action extends by continuity to actions on $\widetilde{U}_\kappa(\hat{\mathfrak{f}})$
and $\widehat{U}_\kappa(\hat{\mathfrak{f}})$. 

\subsection{Group actions and invariants}
Here we discuss an appropriate notion of an arc group action on a vertex algebra, the invariants for the action, and the interaction of taking the invariants and taking the universal enveloping algebra. 

\sss
Let $\Ring$ be a commutative ring, $V$ be a vertex algebra over $\Ring$ and $F$ be a finite type affine algebraic group scheme over $\Ring$. 
As we mentioned in Remark \ref{Rem:arc_group_vertex},
$\Ring[\Jet F]$ is a Hopf algebra object in the category of vertex algebras.

\begin{defi}\label{defi:vertex_action}
By a {\it vertex  $\Jet F$-action} on $V$ we mean a vertex algebra homomorphism $\alpha:V\rightarrow V\underset{\Ring}\otimes \Ring[\Jet F]$ subject to the usual coassociativity and the counit axioms of coaction homomorphisms. 
\end{defi}

Now we reformulate this definition in terms of the action map (rather than the co-action). 

\begin{defi}\label{defi:vertex_action_equiv}
\begin{enumerate}
\item 
We say that the action of $\Jet F$ on an $\Ring$-module $V$ is {\it rational} if every element of $V$ lies in a submodule $V'$, where the action of $\Jet F$ factors through a rational representation of the algebraic group $\Jet_n F$. 
\item Let $V$ be a vertex algebra. We say that a rational action of $\Jet F$ on $V$ is a {\it vertex action} if $|\varnothing\rangle$ is invariant, and for any point $g(t)\in \Jet F_{\tilde{\Ring}}$ (for an $\Ring$-algebra $\tilde{\Ring}$) and any $a\in V$, we have 
\begin{equation}\label{eq:vertex_action_condition}
Y(g(t).a,z)=g(t-z)Y(a,z)g(t-z)^{-1},
\end{equation} an equality in $\left(\operatorname{End}_{\tilde{\Ring}} (V\otimes_\Ring \tilde{\Ring})\right)[[z^{\pm 1}]]$. 
\end{enumerate}
\end{defi}

To get the action $g.v$ for $\tilde{\Ring}$-points $g,v$ of $\Jet G, V$, we evaluate $\alpha(v)$ on $g$, then the description of $Y$ for the vertex algebra $\Ring[\Jet F]$ in Remark \ref{Rem:jet_vertex_coordinate_free} shows that the condition that $\alpha$ is a vertex algebra homomorphism is equivalent to the invariance of $|\varnothing\rangle$ and (\ref{eq:vertex_action_condition}).

\sss
Recall that, over a characteristic $0$ field, a rational representation of a connected algebraic group $F$ in an algebra is by automorphisms if and only if the corresponding action of the Lie algebra $\mathfrak{f}$ is by derivations. We will need an analog of this claim in the vertex algebra setting. 

\begin{Lem}\label{Lem:vertex_action_Lie_algebra}
Suppose that $F$ is connected, and $\Ring$ is a characteristic $0$ field. Suppose $V$ is a vertex algebra equipped with a rational representation of $\Jet F$. Then the following conditions are equivalent:
\begin{enumerate}
\item The action of $\Jet F$ in $V$ is a vertex action,
\item for all $x\in \mathfrak{f},m\in \Z_{\geqslant 0}, b\in V, n\in \Z$, we have 
\begin{equation}\label{eq:vertex_action_derivation}
[xt^m,b_{(n)}]=\sum_{\ell\geqslant 0}{m\choose \ell}(xt^\ell.b)_{(m+n-\ell)}.
\end{equation}
\end{enumerate}
\end{Lem}
\begin{proof}
We note that it is enough to check (\ref{eq:vertex_action_condition}) for the cases when $g$ is constant, and $g(0)=1$ separately. (\ref{eq:vertex_action_condition}) for all constant $g$ is equivalent to (\ref{eq:vertex_action_derivation}) for all $x$ and $m=0$. Now note that the elements $g$ with $g(0)=1$ are uniquely represented in the form $\operatorname{exp}(\xi)$ for $\xi\in t\mathfrak{f}[[t]]$. Now it is straightforward to check that (\ref{eq:vertex_action_condition}) for all $g$ satisfying $g(0)=1$ is equivalent to
(\ref{eq:vertex_action_derivation}) for all $x\in \ff$ and all $m>0$. This completes the proof. 
\end{proof}

\begin{Rem}\label{Rem:inner_Lie_algebra_action}
Suppose there is a vertex algebra homomorphism $\varphi:V_\kappa(\mathfrak{f})\rightarrow V$. Then we get an action of $\mathfrak{f}[[t]]$ on $V$ via $xt^m.b:=\varphi(xt^{-1})_{(m)}b$. Then (\ref{eq:vertex_action_derivation}) follows from 
(\ref{eq:VA_modes_brackets}).
\end{Rem}

\sss 
Now we discuss an example of this situation.
Suppose $\Ring$ is a domain whose field of fractions has characteristic $0$.
Suppose that $F$ is a smooth connected finite type affine group scheme over $\Ring$.
By the construction of $V_\kappa(\mathfrak{f})$ as the induced module, we have an action of $\Jet F$ on the $\Ring$-module $V_\kappa(\mathfrak{f})$. Note that this action is rational: observe that $V_\kappa(\mathfrak{f})$ is $\Z_{\leqslant 0}$-graded via the energy grading. The action of $\Jet F$ preserves $\bigoplus_{j\geqslant j_0}V_\kappa(\mathfrak{f})^j$ for all $j_0$, where we write $V_\kappa(\mathfrak{f})^j$ for the $j$-th graded component. Moreover,  the action on this $\Ring$-submodule factors through a rational action of $\Jet_n F$ for $n$ sufficiently large. 

\begin{Lem}\label{Lem:vertex_algebra_action_example} The action of $\Jet F$ on $V_\kappa(\mathfrak{f})$ is a vertex  action. The same holds for all base changes. 
\end{Lem}
\begin{proof}
We need to show that the coaction map 
$V_\kappa(\mathfrak{f})\rightarrow V_\kappa(\mathfrak{f})\otimes \Ring[\Jet F]$ is a vertex algebra homomorphism. This amounts to checking some collection of polynomial equations. We split the check into three cases:

{\it Case 1}: $\Ring$ is a characteristic $0$ field.
We note that the action of $\mathfrak{f}[[t]]$ on $V_\kappa(\mathfrak{f})$ arises as in Remark \ref{Rem:inner_Lie_algebra_action}. Our claim follows from
Lemma \ref{Lem:vertex_action_Lie_algebra}.

{\it Case 2}: General case. We note $V_\kappa(\mathfrak{f})$ is flat over $\Ring$ -- this follows from the definition. By Lemma \ref{Lem:jet_smoothness_flatness}, $\Ring[\Jet \Ring]$ is also flat over $\Ring$. So this case follows from the case of $\operatorname{Frac}(\Ring)$.  

{\it Case 3}: We then can change the base to any $\Ring$-algebra.
\end{proof}

Similarly, the action of $\Jet F$ on $V_{\bone}(\ff)$ is a vertex action.

\sss
Now we get back to the general situation where we have action of $\Jet F$ on a vertex algebra $V$ with co-action map $\alpha:V\rightarrow V \underset{\Ring}\otimes \Ring[\Jet F]$. An element 
$v\in V$ is called {\it $\Jet F$-invariant} if $\alpha(v)=v\otimes 1$. 
It is an easy exercise that the invariants form a vertex subalgebra in $V$ to be denoted by $V^{\Jet F}$.



\begin{Lem}\label{Lem:action_on_modes}
The action of $\Jet F$ on $V$ gives rise to 
an action of the loop group $\Loop F$ on the completed universal enveloping algebra $\tilde{U}(V)$. Moreover, the image of the natural homomorphism 
$\tilde{U}(V^{\Jet F})\rightarrow \tilde{U}(V)$ lies in the subalgebra of invariants $\tilde{U}(V)^{\Loop F}$.
\end{Lem}
\begin{proof} 
By  functoriality of $\widetilde{U}$, Section \ref{SSS:Utilde},  the co-action map $\alpha: V\rightarrow V\otimes_\Ring\Ring[\Jet F]$ gives rise to 
$$\widetilde{U}(\alpha):  \widetilde{U}(V)\rightarrow \widetilde{U}(V\otimes_{\Ring}\Ring[\Jet F])=
\widetilde{U}(V)\widetilde{\otimes}_{\Ring}\widetilde{U}(\Ring[\Jet F])).$$
Here $\widetilde{\otimes}_\Ring$ denotes the completed tensor product of complete topological $\Ring$-modules.
By Lemma \ref{Lem:loop_vs_univ_envel}, 
$\widetilde{U}(\Ring[\Jet F]))=\Ring[\Loop F]$. 
Using functoriality of $\widetilde{U}$ combined with Lemma \ref{Lem:univ_env_tensor}, we see that $\widetilde{U}(\alpha)$ is a co-action map for the topological Hopf algebra $\Ring[\Loop F]$. Let $\iota$ denote the inclusion map 
$V^{\Jet F}\hookrightarrow V$. It is easy to see that every element $x$ in the image of $\widetilde{U}(\iota)$ satisfies $\widetilde{U}(x)=x\otimes 1$, which implies the claim of the lemma. 
\end{proof}

\sss Now suppose that $V$ is a filtered vertex algebra. Equip $\Ring[\Jet F]$ with the trivial filtration, where all elements are of degree $0$.
The action of $\Jet F$ on $V$ preserves the filtration if and only if the co-action map $\alpha$ is a homomorphism of filtered vertex algebras (where we take the trivial filtration on $\Ring[\Jet F]$). 
In this case the action of $\Loop F$ on 
$\widetilde{U}(V)$ preserves $\widehat{U}(V)$.

For example, if $V=V_\kappa(\ff)$, then the action of $\Jet F$ preserves the filtration. 
We recover the actions of $\Loop F$ on $\widetilde{U}_\kappa(\hat{\ff}), \widehat{U}_\kappa(\hat{\ff})$ that on the level of $\Ring$-points were discussed in Section 
\ref{SS_affine_vertex}.

\subsection{Chiral differential operators} \label{ss:defcdo}
In this section we introduce another important class of vertex algebras: the algebras of chiral differential operators. 

In the entire section we suppose that $\Ring$ is a Noetherian ring satisfying the following conditions:
\begin{itemize}
\item $\Ring$ is a domain such that $2\in \Ring$ is invertible,
\item and $\K:=\operatorname{Frac}(\Ring)$ is a characteristic $0$ field. 
\end{itemize}
We also assume that $F$ is a smooth connected finite type affine group scheme over $\operatorname{Spec}(\Ring)$.

Let $W$ have the same meaning as in Remark \ref{Rem:general_setting}.
Form the Lie algebra $\hat{\mathfrak{f}}$, see Section \ref{SS_affine_vertex}.

\subsubsection{Construction}\label{SSS_CDO_construction}
%
Note that $\Jet\mathfrak{f} = \mathfrak{f}[[t]]$ acts on $\Ring[\Jet F]$ (in fact, in two different ways: via right -- for the action from the left -- and left invariant vector fields; we are now interested in the former action). So we can form the induced module
\[
\CDO_\kappa(F): = \operatorname{Ind}_{\mathfrak{f}[[t]]\oplus \Ring\mathbf{1}\oplus \Ring\mathbf{1}'}^{\hat{\mathfrak{f}}} \Ring[\Jet F].
\]
where $\mathbf{1}$ acts by $\kappa$ and $\mathbf{1}'$ acts by $1$. By definition, we have an $\Ring$-linear identification 
\begin{equation}\label{eq:CDO_CDO_decomposition}
U(t^{-1}\mathfrak{f}[t^{-1}])\otimes \Ring[\Jet F]\xrightarrow{\sim} \CDO_\kappa(F).
\end{equation}
In particular, together with Lemma \ref{Lem:jet_smoothness_flatness} this implies that $\CDO_\kappa(F)$ is flat as a module over $\Ring$. (\ref{Lem:jet_smoothness_flatness}) equips $\CDO_\kappa(F)$ with a default energy grading.

We need to equip $\CDO_\kappa(F)$ with a vertex algebra structure. 
Define the vacuum element in  $\CDO_\kappa(F)$ as the image of $1\otimes 1$ under (\ref{eq:CDO_CDO_decomposition}). 


We will explain how to define the fields $Y(xt^{-1},z)$ for $x\in \mathfrak{f}$ and $Y(f,z)$ for $f \in \Ring[\Jet F]$, here we slightly abuse the notation and write $xt^{-1}$
for the image of $(xt^{-1})\otimes 1$ under (\ref{eq:CDO_CDO_decomposition}). To define $Y(xt^{-1},z)$ is easy: $$Y(xt^{-1},z)=\sum_{i\in \Z}[xt^i]z^{-(i+1)},$$ 
these fields make sense for any smooth $\hat{\mathfrak{f}}$-module, in particular, $\CDO_\kappa(F)$. To define $Y(f,z)$ is more tricky. First, recall from Section \ref{SS_arcs_loops} that $\Ring[\Jet F]$ is a commutative vertex algebra. We want it to be a vertex subalgebra, which defines $Y(f,z)$ on $\Ring[\Jet F]$. To extend $Y(f,z)$ to $\CDO_\beta(F)$, we note that thanks to (\ref{eq:CDO_CDO_decomposition}) it's enough to specify the commutator
\[
[Y(xt^{-1} z), Y(f,w)].
\]
By formula  (\ref{eq:Y_bracket}), in a vertex algebra, 
\begin{equation}\label{eq:bracket_CDO}
[Y(xt^{-1} z), Y(f,w)]=
\sum_{\ell \geqslant 0}  Y([xt^{\ell}].f,w) \partial_w^{(\ell)}\delta(z-w).
\end{equation}
 
Thanks to (\ref{eq:bracket_CDO}) we can define $Y(f,z)$ on $\CDO_\kappa(\ff)$: we know how to define this series of operators on $\R[\Jet F]$ and (\ref{eq:bracket_CDO}) allows to commute the modes of $Y(f,z)$ past elements of the form $xt^{k}, k<0$. 

\begin{Lem}\label{Lem:CDO_vertex}
There is a unique vertex algebra structure on $\CDO_\kappa(F)$ satisfying the properties listed above in this section.
\end{Lem}
\begin{proof}
Using (\ref{eq:CDO_CDO_decomposition}) and (\ref{eq:Y_product}) allows to uniquely extend the assignment $a\mapsto Y(a,z)$ to any element $a\in \CDO_\kappa(F)$.  
It follows from \cite[Section 3.3]{Arkhipov_Gaitsgory} that $\CDO_\kappa(F)$ becomes a vertex algebra 
after changing the base to $\K$. 

Since $\CDO_\kappa(F)$ is flat over $\Ring$ (as was remarked after (\ref{eq:CDO_CDO_decomposition})), it follows that the axioms in Definition \ref{defi:vertex_ring} hold, so $\CDO_\kappa(F)$ is indeed a vertex algebra.
\end{proof}


\begin{Rem}\label{Rem:affine_CDO_embedding}
Note that $\Ring \subset \Ring[\Jet F]$ (the constants) is an $\mathfrak{f}[[t]]$-submodule. This gives rise to an embedding $\iota_l:V_\kappa(\mathfrak{f})\hookrightarrow \CDO_\kappa(F)$ of $\hat{\mathfrak{f}}_\kappa$-modules. The construction of the vertex algebra structure on $\CDO_\kappa(F)$ implies that this is an embedding of vertex algebras. 
\end{Rem}

Finally, note that for any $\Ring$-algebra $\tilde{\Ring}$ we can consider the base change $\tilde{\Ring}\otimes_\Ring\CDO_\kappa(F)$,this is is a vertex algebra over $\tilde{\Ring}$.

\subsubsection{Localization and filtration}\label{SSS_CDO_local_filtr}
Here we introduce a certain localization of $\CDO_\kappa(F)$, equip it with a (naive) filtration, and study the associated graded. We use the same assumptions on $\Ring$ as in the beginning of the section.

Let $F^\diamond$ be an open affine (say, principal) subset of $F$. Then $\Jet \mathfrak{f}$ still acts on $\Ring[\Jet F^\diamond]$, so we can form the induced module
\begin{equation}\label{eq:induced_localization}
\operatorname{Ind}_{\mathfrak{f}[[t]]\oplus \Ring\mathbf{1}\oplus \Ring\mathbf{1}'}^{\hat{\mathfrak{f}}}
\Ring[\Jet F^\diamond].
\end{equation}

\begin{Lem}\label{Lem:va_CDO_localization}
This induced module carries a natural vertex algebra structure. 
\end{Lem} 
\begin{proof}
In the case when $\Ring$ is a characteristic $0$ field, we can apply the strong reconstruction theorem, see \cite[Section 3.6]{BZF}. We take the elements of form $(xt)^{-1}\otimes 1, 1\otimes f, 1\otimes g^{-1}$, for $x\in \mathfrak{f}, f\in \F[F]$ and $g\in \F[F]$ invertible on $F^\diamond$. It is easy to see that (\ref{eq:bracket_CDO}) still holds on (\ref{eq:induced_localization}), one basically needs to check this when applying both sides to elements $1\otimes g^{-1}$.

Next, we can define an operator $Y(g^{-1},z)$ on 
(\ref{eq:induced_localization}) using a complete analog of (\ref{eq:bracket_CDO}) with $Y(f,w)$ replaced with $Y(g^{-1},w)$ to describe the commutation of $Y(g^{-1},w)$ past the elements $xt^{-k}$. This specifies $Y(g^{-1},z)$ uniquely but does not guarantee that (\ref{eq:bracket_CDO}) holds in the entire generality. From the construction of $Y(g^{-1},w)$ we see that 
$Y(g,w)Y(g^{-1},w)=\operatorname{id}$. We already know (\ref{eq:bracket_CDO}) for $f=g$, and (\ref{eq:bracket_CDO}) for $f=g^{-1}$ then follows. This shows in particular that the fields $Y(xt^{-1},z), Y(f,z), Y(g^{-1},z)$ are mutually local and allows us to apply the strong reconstruction theorem.

And in the case when $\Ring$ is a domain with characteristic $0$ fraction field, we embed $\CDO_\kappa(F^\diamond)$ into $\operatorname{Frac}(\Ring)\otimes_{\Ring}\CDO_\kappa(F^\diamond)$ and see that the image is closed under $(a,b)\mapsto a_{(n)}b$ for all $n$.
\end{proof}

The resulting vertex algebra will be denoted by $\CDO_\kappa(F^\diamond)$. Note that the natural $\hat{\mathfrak{f}}$-module embedding $\CDO_\kappa(F)\hookrightarrow\CDO_\kappa(F^\diamond)$ is a vertex algebra embedding. The target should be thought of as a localization of $\CDO_\kappa(F)$.

The vertex algebra $\CDO_\kappa(F^\diamond)$ is filtered by the degree of the differential operator. This  is the filtration on the induced module coming from the PBW filtration on $U(\hat{\mathfrak{f}}_\kappa)$ (with $\Ring[\Jet F^\diamond]$ in degree $0$). The resulting filtered vertex algebra is easily seen to be almost commutative. Let $\CDO_\kappa(F^\diamond)_{\leqslant m}$ denote the degree $\leqslant m$ piece. 

Note that we have vertex algebra homomorphisms $$\operatorname{gr}V_\kappa(\ff), \Ring[\Jet F^\diamond]\hookrightarrow \gr \CDO_\kappa(F^\diamond).$$ 
By considering the associated graded homomorphism of (\ref{eq:CDO_CDO_decomposition}), we easily see that 
that the induced homomorphism of commutative vertex algebras $\operatorname{gr}V_\kappa(\ff)\otimes \Ring[\Jet F^\diamond]\rightarrow \operatorname{gr}\CDO_\kappa(F^\diamond)$ is an isomorphism. The source vertex algebra is nothing else but $\Ring[\Jet T^*F^\diamond]$. So we get 
\begin{equation}\label{eq:gr_CDO}
\Ring[\Jet(T^*F^\diamond)]\xrightarrow{\sim} \gr \CDO_\kappa(F^\diamond). 
\end{equation}
Note that we can trivialize $T^*F$ using left-invariant vector fields, then $\Jet(T^*F^\diamond)\cong \Jet F^\diamond\times \Jet \ff^*$.

The results and constructions above in this section remain valid after base change to any $\Ring$-algebra $\tilde{\Ring}$.

\sss\label{SSS_CDO_universal}
Similarly we can consider the universal version 
$\CDO_{\bone}(F)$, it is induced from the $\ff[[t]]\oplus \Ring\mathbf{1}$-module $\Ring$, and is a vertex algebra over $\Ring[\mathbf{1}]$. It admits a vertex algebra embedding $\iota_l$ from $V_{\bone}(\ff)$. It also admits a filtration induced from the PBW filtration on $U(\hat{\ff})$. 
And similarly to Section \ref{SSS_CDO_local_filtr} we can consider the localization $\CDO_{\bone}(F^\diamond)$.
The associated graded $\gr \CDO_{\bone}(F^\diamond)$ is the algebra of regular functions on $\Jet F^\diamond\times \{\lambda\partial+\ff^*[[t]]dt\}$. 






\subsubsection{The case of graded unipotent groups}
Consider the case when $F:=N$ is a unipotent group over $\Ring$. We are interested in the structure of $\CDO(N):=\CDO_0(N)$ -- we want to identify it with the vertex algebra known as the $\beta\gamma$-system (one could also call it the Weyl vertex algebra).

We are going to make a simplifying assumption:
\begin{itemize}
\item[(*)] $\mathbb{G}_m$ acts on $N$ by group automorphisms such that the induced grading on $\Ring[N]$ is negative (meaning $\Ring[N]=\Ring\oplus \bigoplus_{i<0}\Ring[N]^i$).
\end{itemize}

The main claim of this section, Proposition \ref{CDO_prop_nilpotent} should hold without (*) (and this is known over characteristic $0$ fields) but the proof simplifies if we impose (*) and all groups we care about satisfy (*).

For example, suppose $G$ is a split connected reductive group over $\Ring$.  If $P \subset G$ is a standard parabolic subgroup, we can take $N:= \operatorname{Rad}_u(P)$ and take the $\mathbb{G}_m$-action corresponding to the co-character $2\rho^\vee$. 

Choose free homogeneous generators $y^1,\ldots,y^\ell$ of $\Ring[N]$. Let $y_1,\ldots,y_\ell\in \mathfrak{n}$ denote the dual basis to $d_1y^1,\ldots, d_1y^\ell$.
Then $\Ring[\Jet N]=\Ring[y^i_n]$, where $i=1,\ldots,\ell$ and $n\leqslant 0$ (before the generators of the functions on jets were labelled by negative integers, but here it is convenient to shift the numbering by $1$). The grading on $\Ring[N]$ gives rise to a grading on $\CDO(N)$
so that (\ref{eq:CDO_CDO_decomposition}) is graded. This is an energy vertex algebra grading.  
The subalgebra  $\Ring[\Jet N]$ is negatively graded (with $\operatorname{deg}T^{(i)}=-i$), while $V(\mathfrak{n})$ is positively graded. 

Let $\partial_1,\ldots, \partial_\ell$ be the constant vector fields of $N\cong \mathbb{A}^\ell$ corresponding to the free generators $y^1,\ldots,y^\ell$ of $\Ring[N]$. Let $\,^Ry_1,\ldots,\,^Ry_\ell$ be the right-invariant vector fields corresponding to $y_1,\ldots,y_\ell$, equivalently, the images of $y_1,\ldots,y_\ell$ under the homomorphism $\mathfrak{n}\rightarrow \mathsf{Vect}(N)$ corresponding to the action of $N$ on itself from the left. Then we can find unique elements $f_{ij}\in \Ring[N]$ satisfying 
$\partial_i=\sum_{j=1}^\ell f_{ij}\,^Ry_j$ for all $i=1,\ldots,\ell$. Consider the following elements of $\CDO(N)$:
\begin{equation}\label{eq:CDO_const_vect_vertex}
\partial_{i,-1}=\sum_{j=1}^\ell (f_{ij})_{(-1)}\iota_l(y_it^{-1}).
\end{equation}
Here we view $f_{ij}$'s as elements of $\Ring[\Jet N]$ via the inclusion $\Ring[N]\hookrightarrow \Ring[\Jet N]$. 

\begin{Prop}\label{CDO_prop_nilpotent}
The following equalities hold:
\begin{enumerate}
    \item $[Y(y^i_0, z), Y(y^j_0, w)] = 0, \forall i,j$.
    \item $[Y(\partial_{i,-1}, z), Y(y^j_0, w)] = \delta_{ij} \delta(z-w), \forall i,j$.
    \item $[Y(\partial_{i,-1}, z), Y(\partial_{j,-1}, w)] = 0,  \forall i,j$. 
\end{enumerate}
\end{Prop}
Over a characteristic $0$ field, this is standard. We give a general proof for reader's convenience. 

\begin{proof}
{\it Part (1)}. This follows because $\Ring[\Jet N]$ is a commutative vertex algebra.

{\it Part (2)}. We have
    \[
    [Y(\partial_{i,-1}, z), Y(y^j_0, w)] = 
  \sum_{k \geqslant 0} Y([\partial_i]_{(k)} y^j_0, w) \partial_w^{(k)}\delta(z-w)  
    \]
    Let's compute $[\partial_i]_{(k)} y^j_0$, i.e., the coefficient of $z^{-n-1}$ in 
    $Y(\partial_{i,-1},z)y^j_0$. 
    Thanks to (\ref{eq:CDO_const_vect_vertex}),  we have
    \begin{equation}\label{eq:CDO_Weyl_relation1}
    Y(\partial_{i,-1},z)y^j_0=\sum_{h=1}^\ell :Y(f_{ih,0},z)Y(\iota_l(y_it^{-1}),z):y^j_0
    \end{equation}
Here we write $f_{ih,0}$ for $f_{ih}\in \Ring[N]$ viewed as an element of $\Ring[\Jet N]$ under the natural inclusion of $\Ring[N]$. We will write $f_{ih,(k)}$ for the Fourier modes of $Y(f_{ih,0},z)$.
    
    Note that $f_{ih,(k)}y^j_0=0$ for $k\geqslant 0$ because $\Ring[\Jet N]$ is a commutative vertex algebra. It turns out that $\iota_l(y_{h}t^{-1})_{(k)}y^j_0=0$ for $k>0$. In order to see this, we will need the following two observations:
    \begin{itemize}
    \item[(a)] $y_0^j\in \Ring[N]\subset \Ring[\Jet N]$,
    \item[(b)] and the endomorphism $\iota_l(y_{h}t^{-1})_{(k)}$ of $\Ring[\Jet N]$ is the action of the element $y_ht^k\in \mathfrak{n}[[t]]$ on the module $\Ring[\Jet N]$. The action of $\Jet N$ on $\Ring[N]$ factors through $\Jet N\twoheadrightarrow N$, hence elements of $t\mathfrak{n}[[t]]$ annihilate $\Ring[N]$.
    \end{itemize}
    
    Thanks to these observations, the right hand side of (\ref{eq:CDO_Weyl_relation1}) simplifies to 
    $$\sum_{h=1}^\ell f_{ih,(-1)}\iota_l(y_{h}t^{-1})_{(0)} y^j_0= \partial_{i,(0)}y^j_0=\delta_{ij}$$ finishing the proof. 
    
{\it Part (3)}.      We need to show that the elements 
    $\nu:=\partial_{i,(k)}\partial_{j,-1}$ vanish for all $k\geqslant 0$. Thanks to part 2 of the proposition, $[Y(\partial_{i,-1},z),Y(\partial_{j,-1},w)]$ commutes with $Y(y^\alpha_0,u)$ for all $\alpha=1,\ldots,\ell$. It follows from (\ref{eq:Y_bracket}) that 
    \begin{equation}\label{eq:CDO_nu_condition}
    \nu_{(\ell)}y^\alpha_0=0, \forall \ell\geqslant 0.
    \end{equation}

    Thanks to the filtration on $\CDO(N)$ by the order of differential operator (that turns $\CDO(N)$ into an almost commutative vertex algebra), see Section \ref{SSS_CDO_local_filtr}, 
    $\nu$ lies in filtration degree $1$. 
    We claim that (\ref{eq:CDO_nu_condition}) implies that $\nu\in \Ring[\Jet N]$. Once we show this, we complete the proof that $\nu=0$: by the construction $\nu$ has positive energy degree, while elements of $\Ring[\Jet N]$ have non-positive degree.

    So, we need to show  (\ref{eq:CDO_nu_condition}) implies $\nu\in \Ring[\Jet N]$. Since $\nu$ is in filtration degree $1$, it is the sum of an element of $\Ring[\Jet N]$ and of a linear combination of elements of the form $y_it^{k}\otimes g_{ik}=\iota_l(y_it^{-1})_{(k)}g_{ik}\in \CDO_\kappa(F)$, where $k<0, g_{ik}\in \Ring[\Jet N]$. Note that one can express $\,^R y_i$ as a linear combination of $\partial_{j}$'s with coefficients from $\Ring[N]$. Using this and induction on $k$, one can express $[\iota_l(y_it^{-1})]_{(k)}g_{ik}$ as a linear combination of an element in $\Ring[\Jet N]$ and elements $(\partial_{i',-1})_{(k')}g_{i',k'}'$ with $g_{i',k'}\in \Ring[\Jet N]$. 
    We have 
    \begin{align*}
&Y((\partial_{i',-1})_{(k')}g_{i',k'}',z)y^\alpha_{0}=
    :[\partial_z^{(-k'-1)}Y(\partial_{i',-1},z)]Y(g_{i',k'},z):y^\alpha_{0}=(-1)^{k'}g_{i',k'}z^{-k'}+h.d.t.,
    \end{align*}
    where ``h.d.t'' means terms with $z$ in degree larger than $-k'$.  We use (\ref{eq:CDO_nu_condition}) to conclude that all $g_{i',k'}$ are equal to $0$. The claim that $\nu\in \Ring[\Jet N]$ follows.  
\end{proof}

\subsubsection{Embedding $V_{-\kappa}(\mathfrak{f})\hookrightarrow \CDO_\kappa(F)$}\label{SSS_iota_r}
Our goal here is to produce an embedding
\[
\iota_r: V_{-\kappa}(\ff) \to \CDO_\kappa(F)
\]
whose image commutes with that of $\iota_l$:
$$[Y(\iota_l(xt^{-1}),z),Y(\iota_r(yt^{-1}),w)]=0.$$ 

We construct $\iota_r$ as the sum $\iota_1+\iota_0$, where $\iota_1$ can be described as the naive expected map and $\iota_0$ is a correction term taking values in $\Ring[\Jet F]$.

First, we explain the construction to $\iota_1$, which is analogous to (\ref{eq:CDO_const_vect_vertex}).
Pick a basis $\{y_i\}$ for $\mathfrak{f}$. To construct $\iota_1$, we express the left invariant vector field $\,^Ly_i$ (=the image of $y_i$ under the homomorphism associated to the $F$-action on itself from the right) as $\sum f_{ij} \,^Ry_j$ and set
\[
\iota_1(y_it^{-1}) = \sum (f_{ij})_{(-1)} \iota_l(y_jt^{-1}), \text{ for }i=1,\ldots,\ell.
\]

To define $\iota_0$, we need a natural map $\mathfrak{f} \to \Ring[\Jet F]$. Let $T$ be the  derivation of $\Ring[\Jet F]$, which is a part of the vertex algebra structure. Consider the composition of $T$ with the natural embedding $\Ring[F]\hookrightarrow \Ring[\Jet F]$. This is a derivation of $\Ring[F]\rightarrow \Ring[\Jet F]$, and hence it factors through a $\Ring$-linear map $\eta: \Omega^1(F) \rightarrow \Ring[\Jet F]$. 

We can trivialize $TF \cong F\times \mathfrak{f}$ using left-invariant vector fields. This gives rise to a trivialization $T^*F\cong F\times \mathfrak{f}^*$. Hence $\eta$ restricts to $\mathfrak{f}^* \to \Ring[\Jet F]$ whose image consists of $F$-invariants for the action of $F$ from the left. Next, we can view any element of $\Sym^2(\ff^*)$ as a map $\mathfrak{f}\rightarrow \mathfrak{f}^*$. Hence we get a map $\mathfrak{f}\rightarrow \Omega^1(F)$
$$\mathfrak{f}t^{-1}\xrightarrow{t\cdot}\mathfrak{f}
\xrightarrow{(*)} \mathfrak{f}^*\xrightarrow{\eta}\Omega^1(F)\rightarrow \Ring[\Jet F],$$
where (*) stands for the form $-\frac{1}{2}\beta_{\ff}-\kappa\beta_W$, where $\beta_{\ff}$ is the Killing form and $\beta_W$ is the trace form of the representation $W$ used to define the central extension of $\Loop \ff$, see Section \ref{SSS_KM_definition}.

If we replace $\Ring$ with $\K$, the following result is a part of \cite[Theorem 3.7]{Arkhipov_Gaitsgory}.
Once we know this over $\K$, the claim over $\Ring$ follows, cf. the proof of Lemma \ref{Lem:CDO_vertex}.

\begin{Prop}\label{Prop:commuting_embedding}
The map $\iota_r:=\iota_0+\iota_1:  \mathfrak{f}t^{-1}\rightarrow  \CDO_\beta(F)$ gives a vertex algebra embedding $V_{-\kappa}(\mathfrak{f}) \to \CDO_\kappa(F)$ whose image commutes with that of $V_\kappa(\mathfrak{f})$.
\end{Prop}

\begin{Rem}\label{Rem:CDO_righthanded}
The homomorphism $\iota_r$ turns $\CDO_\kappa(F)$ into a $U_{-\kappa}(\hat{\mathfrak{f}})$-module. Note that $\iota_0(xt^{-1})_{(n)}$ vanishes on $\Ring[\Jet F]$ for $n\geqslant 0$, and $\iota_1(xt^{-1})_{(n)}$ is the action of $xt^n\in \mathfrak{f}[[t]]$ from the right, cf. \cite[Theorem 3.7(c)]{Arkhipov_Gaitsgory}. So we get a $U_{-\kappa}(\hat{\mathfrak{f}})$-module homomorphism 
\begin{equation}\label{eq:CDO_alternative}\operatorname{Ind}_{\mathfrak{f}[[t]]\oplus \Ring \mathbf{1}\oplus \Ring \mathbf{1}'}^{\hat{\mathfrak{f}}}\Ring[\Jet F]\rightarrow \CDO_\kappa(F).\end{equation}
This homomorphism is an isomorphism. One can see this by passing to associated graded algebras. Also, one can equip the source with a vertex algebra structure as in Section \ref{SSS_CDO_construction}, and (\ref{eq:CDO_alternative}) becomes an isomorphism
of vertex algebras. Although this is straightforward to verify directly using \cite[Theorem 3.7(c)]{Arkhipov_Gaitsgory}, perhaps surprisingly we are unaware of a reference which spells this out explicitly. However, we emphasize this is certainly known to experts in characteristic zero, and for the convenience of the reader we now sketch one way to deduce it, albeit rather circuitously, from written references. Namely, we recall that under the equivalence between $\CDO_\kappa(F)$-modules and $\kappa$-twisted D-modules on $\Loop F$, cf. \cite[Appendix A]{Arkhipov_Gaitsgory}, $\CDO_\kappa(F)$ corresponds to the ($\kappa$-twisted) dualizing sheaf $\omega_{\Jet F, \kappa}$ on $\Jet F$. Moreover, the vertex algebra structure on the $\CDO_\kappa(F)$ corresponds to the natural factorization algebra structure on  $\omega_{\Jet F, \kappa}$, cf. \cite[Section 3.9]{Chiral_Algebras}. Finally, we note that the inversion map $\Loop F \simeq \Loop F$ exchanges $\kappa$-twisted and $-\kappa$-twisted D-modules on $\Loop F$, cf. \cite[Example 1.2]{Quantum_Langlands}, and in particular the factorization algebras $\omega_{\Jet F, \kappa}$ and $\omega_{\Jet F, -\kappa}$.
. 
Finally, notice that the resulting isomorphism $\Ring[\Jet(T^*F)]\xrightarrow{\sim} \gr \CDO_\kappa(F)$ is the same as (\ref{eq:gr_CDO}).
\end{Rem}

We can extend $\iota_r$ to a homomorphism $V_{\bone}(\ff)\rightarrow \CDO_{\bone}(F)$, where the target vertex algebra was defined in Section \ref{SSS_CDO_universal}. It is $\Ring[\mathbf{1}]$-semilinear with respect to the automorphism that sends $\mathbf{1}$ to $-\mathbf{1}$.

\subsubsection{Arc group actions}\label{SSS_arc_actions_CDO}
The goal of this section is to establish two commuting vertex actions of $\Jet F$ on $\CDO_\kappa(F)$
and to identify the images of $\iota_L$ and $\iota_R$ with the subalgebras of invariants. 

We note that $\Jet F$ acts on $\Ring[\Jet F]$ from the right and the action commutes with the action of $\mathfrak{f}[[t]]$ from the left.
So we get a $\Jet F$-action on 
$$\CDO_\kappa(F) = \operatorname{Ind}_{\mathfrak{f}[[t]]\oplus \Ring\mathbf{1}\oplus \Ring\mathbf{1}'}^{\hat{\mathfrak{f}}} \Ring[\Jet F].$$
The invariants of this $\Jet F$-action is exactly the image of $\iota_\ell$. 

Similarly, using the action of $\Jet F$ on $\Ring[\Jet F]$ from the left, we get another action of $\Jet F$
on $\CDO_\kappa(F)$ using the construction in 
Remark \ref{Rem:CDO_righthanded}. Its invariants is the image of $\iota_r$.

\begin{Lem}\label{Lem:arc_actions_CDO}
The following claims hold:
\begin{enumerate}
\item The two actions of $\Jet F$ commute.
\item They are vertex actions. 
\end{enumerate}
\end{Lem}
\begin{proof}
To prove both statements we can base change to $\mathbb{K}$, cf. the proof of Lemma \ref{Lem:CDO_vertex}. To prove (1) we observe that the corresponding actions of $\mathfrak{f}[[t]]$ commute with each other, this is because they are parts of the commuting (see Proposition \ref{Prop:commuting_embedding}) actions of the vertex subalgebras $V_\kappa(\mathfrak{f})$ and $V_{-\kappa}(\mathfrak{f})$. The proof of (2) repeats the proof of Lemma \ref{Lem:vertex_algebra_action_example}.
\end{proof}

Note that both $\Jet F$-actions preserves the filtration on $\CDO_\kappa(F)$ by the order of differential operator. The isomorphism $\Ring[\Jet(T^*F)]\xrightarrow{\sim} \gr \CDO_\kappa(F)$ is $\Jet F\times \Jet F$-equivariant for the action of $\Jet F\times \Jet F$ on $\Jet(T^*F)$ induced from the action of $F\times F$ on $T^*F$. 

Similarly, $\Jet F\times \Jet F$ acts on $\gr \CDO_{\bone}(F)$. The identification $\gr \CDO_{\bone}(F)\cong \Jet F\times \{\lambda\partial+\ff^*[[t]]dt\}$ from Section \ref{SSS_CDO_universal} is $\Jet F\times \Jet F$-invariant (where the right action of $\Jet F$ on 
$\{\lambda\partial+\ff^*[[t]]dt\}$ is trivial, while the left action is by gauge transformations).

\section{Free field realizations}\label{S_CDO_ffr}
As in the approach of Feigin and Frenkel, see \cite{Frenkel_loop}, 
a crucial role in what we do is played by the free field realization of the affine vertex algebra $V_\kappa(\g)$ as well as parabolic analogs of the free field realization. These parabolic analogs are vertex algebra homomorphisms. In these section we construct them and establish some of their properties. 

Here $\Ring$ is a commutative ring such that $2$ is invertible in $\Ring$, $G$ is a split connected reductive group over $\Ring$, $P$ is its parabolic subgroup, $N\subset P$ is the unipotent radical, and $L$ is a Levi subgroup of $P$. We write $V_\kappa(\lf)$ for the affine vertex algebra for $\lf$ and the restriction of $\kappa$ to $\lf$. 

The homomorphism in question is $V_\kappa(\g)\rightarrow \CDO(N)\otimes V_\kappa(\lf)$. It will be called the {\it parabolic free field realization map} and denote it by $\ffr_P$. When $P$ is a Borel subgroup, we drop ``parabolic'' from the name.

\subsection{Decomposition of $\CDO_\kappa(G^\diamond)$}\label{SS:CDO_decomp}
We write $N^-$ for the unipotent radical of the parabolic $P^-$ opposite to $P$. Set 
\begin{equation}\label{eq:CDO_open_subset}
G^\diamond:=N L N^-=(PN^-=NP^-)\subset G,
\end{equation}
this is a principal open subset. It carries an action of $P$ from the left and the action of $P^-$ from the right.


Recall from Section \ref{SSS_CDO_local_filtr} that we have the localization $\CDO_{\kappa}(G^\diamond)$, the base change to $\Ring$ from $\mathbb{Z}[1/2]$. Note that $\Ring[\Jet G^\diamond] $ decomposes as $\Ring[\Jet P]\otimes \Ring[\Jet N^-]$. 

First, we need analogs of subalgebras $D(P),D(N^-)\subset D(G^\diamond)$ of usual differential operators in $\CDO(G^\diamond)$. Consider the following vertex subalgebras in $\CDO(G^\diamond)$:

1) The subalgebra generated by $\Ring[\Jet N^-]$ and the elements $\iota_r(xt^{-1})$ for all $x\in\mathfrak{n}^-$. As an $\Ring$-submodule, it coincides with $\operatorname{Ind}_{\mathfrak{n^-}[[t]]}^{\mathfrak{n^-}(\!(t)\!)}\Ring[\Jet N^-]$, and as a vertex algebra, it is $\CDO(N^-)$.
 
2) The subalgebra generated by $\Ring[\Jet P]$ and
the elements $\iota_l(xt^{-1})$ for $x\in \mathfrak{p}$. As an $\Ring$-submodule it coincides with
$$\operatorname{Ind}_{\frp[[t]]\oplus \Ring \mathbf{1}\oplus \Ring \mathbf{1}'}^{\hat{\frp}'_{\kappa}}\Ring[\Jet P].$$
The meaning of the prime symbol is that the form we consider still comes from the functional $\kappa\operatorname{tr}_W-\frac{1}{2}\operatorname{tr}_{\g}$ (as opposed to
$\kappa\operatorname{tr}_W-\frac{1}{2}\operatorname{tr}_{\mathfrak{p}}$).
As a vertex algebra the induced module is $\CDO'_{\kappa}(P)$.

From the construction it is easy to see that 
the vertex subalgebras 
$$\CDO(N^-),\CDO'_{\kappa}(P)\subset \CDO_\kappa(G^\diamond)$$ commute.

So we get a vertex algebra homomorphism:
\begin{equation}\label{eq:CDO_tensor_decomposition}
\CDO'_{\kappa}(P)\otimes \CDO(N^-) \to \CDO_\kappa(G^\diamond).
\end{equation}

\begin{Lem}\label{Lem:CDO_tensor_decomposition}
(\ref{eq:CDO_tensor_decomposition}) is an isomorphism.
\end{Lem}
\begin{proof}
Note that (\ref{eq:CDO_tensor_decomposition}) is preserves the filtrations by the order of differential operator. The associated graded homomorphism 
coincides with the isomorphism 
$\Ring[\Jet(T^*P)]\otimes \Ring[\Jet(T^*N^-)]\xrightarrow{\sim} \Ring[\Jet(T^*G^\diamond)]$. The latter comes from the decomposition $\Jet(T^*G^\diamond)\cong \Jet(T^*P)\times \Jet(T^*N^-)$, which, in its turn, is induced by $T^*G^\diamond\cong T^*P\times T^*N^-$.
\end{proof}

\begin{Rem}\label{eq:CDO_tensor_decomp_equiv}
We'll need equivariance properties of (\ref{eq:CDO_tensor_decomposition}).  First, note that $\Jet P$ acts from the left, and (\ref{eq:CDO_tensor_decomposition}) is equivariant by the construction. Also, $\Jet P^-$ acts from the right: $\Jet L$ acts diagonally, while $\Jet N^-$ acts on the second factor only. Then (\ref{eq:CDO_tensor_decomposition}) is $\Jet P^-$-equivariant.
\end{Rem}

For similar reasons, we have the decomposition
\begin{equation}\label{eq:decomp_for_parabolic}
\CDO'_\kappa(P)=\CDO(N)\otimes \CDO_\kappa(L).
\end{equation}
Note that the second factor is indeed associated to the functional $\kappa\operatorname{tr}_W-\frac{1}{2}\operatorname{tr}_{\lf}$.

Note that we also have the direct analogs of (\ref{Lem:CDO_tensor_decomposition}) and (\ref{eq:decomp_for_parabolic})
for $\CDO_{\bone}(G^\diamond)$: we can define the vertex subalgebra 
$\CDO'_{\bone}(P)$ whose specialization to $\mathbf{1}=\kappa$ is 
$\CDO'_\kappa(P)$.

\subsection{Parabolic Free Field Realization Map}\label{SS_CDO_parabolic_free_field}
We now construct a parabolic free field realization homomorphism
\begin{equation}\label{eq:CDO_parabolic_ffr}
\mathsf{ffr}_{P}:V_{\kappa}(\g) \to \CDO(N) \otimes V_{\kappa}(\lf)
\end{equation}
using (\ref{eq:CDO_tensor_decomposition}) and its equivariance properties from Remark \ref{eq:CDO_tensor_decomp_equiv}. We will write $\ffr_{P,\kappa}$ if we want to indicate the dependence on $\kappa$.

\sss\label{SSS:ffr_construction} Namely, consider the  inclusion $\iota_l:V_{\kappa}(\g) \to \CDO_\kappa(G)$ and compose it with the localization inclusion $$\CDO_\kappa(G) \to \CDO_{\kappa}(G^\diamond)\xrightarrow{\sim}\CDO'_\kappa(P)\otimes \CDO(N^-)\xrightarrow{\sim} \CDO(N)\otimes \CDO_{\kappa}(L)\otimes \CDO(N^-).$$ Consider the composition $V_\kappa(\g)\hookrightarrow \CDO(N)\otimes \CDO_\kappa(L)\otimes \CDO(N^-)$.
Its image is contained in the $\Jet P^-$-invariants. 
Note that the action of $\Jet N^-$ on $\CDO(N)\otimes \CDO_\kappa(L)$ is trivial. So we get an inclusion:
\[
\theta: V_{\kappa}(\g) \hookrightarrow \CDO(N) \otimes \CDO_\kappa(L)\otimes V(\mathfrak{n}^-).
\]
Now note that we have a vertex algebra epimorphism:
\[
\pi:V(\mathfrak{n}^-) \twoheadrightarrow \Ring
\]
induced by $\mathfrak{n}^-\twoheadrightarrow \{0\}$. 

So we get a vertex algebra homomorphism:
\begin{equation}\label{eq:ffr_prelim}
[\operatorname{id}\otimes \operatorname{id}\otimes \pi]\circ \theta:
V_\kappa(\g)\rightarrow \CDO(N)\otimes \CDO_\kappa(L).
\end{equation}
Note that $\Jet L$ acts on $V(\mathfrak{n}^-)$ by a vertex algebra action (restricted from the action on $V_\kappa(\g)$). The homomorphism $\pi$ is $\Jet L$-invariant. It follows that (\ref{eq:ffr_prelim}) is $\Jet L$-invariant,
and hence the image is contained in 
$\CDO(N)\otimes \CDO_\kappa(L)^{\Jet L}=\CDO(N)\otimes V_\kappa(\lf)$. 

\begin{defi}\label{defi:ffrP}
For $\ffr_P$ we take the composition $(\operatorname{id}_N\otimes \operatorname{id}_L\otimes \pi)\circ \theta$ (where we write $\operatorname{id}_N$ and  
$\operatorname{id}_L$ for the identity endomorphisms of $\CDO(N)$ and $\CDO_\kappa(L)$) viewed as a homomorphism to $\CDO(N)\otimes V_\kappa(L)$. 
\end{defi}

\begin{Rem}\label{Rem:ffr_equivariance}
By the construction,
$\ffr_P$ is $\Jet P$-equivariant. In particular, a choice of a maximal torus $H\subset L$ gives rise to gradings on $V_\kappa(\g),\CDO(N)\otimes V_\kappa(\lf)$ by the root lattice, and $\ffr_P$ is graded for this choice of grading. Also, $\ffr_P$ is graded for the energy grading by the construction.
\end{Rem}

\sss\label{SSS_ffrp_previous}
Assume that $\Ring$ is a characteristic $0$ field. A construction of $\ffr_P$ was given in \cite[Section 10.1-3]{FG} using the quantum Drinfeld-Sokolov reduction (the reference ostensibly only discusses the case of the Borel but applies to general parabolics {\em mutatis mutandis}). We claim that our homomorphism $\ffr_P$ coincides with that one. Namely, in that construction one takes 
the $\Jet N^-$-invariant inclusion 
$$V_\kappa(\g)\hookrightarrow 
\CDO(N)\otimes \CDO_{\kappa}(L)\otimes \CDO(N^-),$$
applies the quantum Drinfeld-Sokolov reduction, to be denoted by $\mathsf{Red}_0$, for the action of $\Jet N^-$ to get 
a vertex algebra homomorphism $V_\kappa(\g)\hookrightarrow 
\CDO(N)\otimes \CDO_{\kappa}(L)$. To show the coincidence of the two constructions of $\ffr_P$ one needs to check that the composition $V(\mathfrak{n}^-) \rightarrow \mathsf{Red}_0(\CDO(N^-))\rightarrow \Ring$ coincides with $\pi$. This is true because both are energy and root lattice graded vertex algebra homomorphisms.

\sss\label{SSS_ffrp_classical}
Next, $\ffr_P$ is a filtered vertex algebra homomorphism. 
Now we discuss the associated graded of $\ffr_P$ to be denoted by $\ffr_P^0$, as we will see just below, it is independent of $\kappa$. 
Consider the morphism $\mu: T^*P\rightarrow \g^*$ that is the restriction of the moment map $T^*(G/N^-)\rightarrow \g^*$ to the open subset $T^*P\subset T^*(G/N^-)$. It is $L$-invariant for the action of $L$ from the right. Then we can consider the induced morphism $\underline{\mu}: T^*N\times \mathfrak{l}^*\rightarrow \g^*$. It is easy to see that $\ffr^0_P=\Jet\underline{\mu}$.

We will write $\ffr_P^\hbar$ for the homomorphism between the Rees algebras $V_\kappa^\hbar(\g)\rightarrow \CDO^\hbar(N)\otimes_{\Ring[\hbar]}V^\hbar_\kappa(\lf)$ induced by $\ffr^P$.

\sss\label{SSS_Cartan_images_ffr}
We will need a formula for the images of the elements of the form $xt^{-1}$ for $x\in \fl$. Let $\Phi_P^+$ denote the subset of all  roots whose root vectors lie in $\mathfrak{n}$. We note that there is an $L$-equivariant isomorphism $\mathfrak{n}\xrightarrow{\sim} N$ whose differential at $0\in \mathfrak{n}$ is the identity.
One way to construct it is as follows. Note that once we choose a parabolic subgroup $P'\subset G$ with Levi decomposition $P'=L'\ltimes N'$ such that $N'\subset N$ and $L\subset L'$, then the product map $N'\times (L'\cap N)\rightarrow N$ is an isomorphism. So in our construction of an isomorphism $\mathfrak{n}\xrightarrow{\sim} N$ reduces to the case when $P$ is a maximal parabolic. Since $p>h$, the identification of $\mathfrak{n}$ with $N$ using the Campbell-Hausdorff formula works. 

Choose root vectors $\partial_\alpha\in \mathfrak{n}, \alpha\in \Phi^+_P$ and dual basis vectors $y^\alpha\in \mathfrak{n}^*\subset \Ring[N]$. So we get elements $y^\alpha_0,\partial_{\alpha,-1}\in \CDO(N)$. We can write the operator $x_{\mathfrak{n}}$ of the action of $x$ on $\mathfrak{n}$ as 
$$\sum_{\alpha,\beta\in \Phi_P^+}c_\alpha^\beta \partial_\beta\otimes y^\alpha, c_\alpha^\beta\in \Ring.$$

\begin{Lem}\label{Lem:ffr_Cartan_images}
We have 
\begin{equation}\label{eq:Cartan_image_ffrp}
\ffr_P(xt^{-1})=-\sum_{\alpha\in \Phi^P_+}c_\beta^\alpha y^\beta_0\partial_{\alpha,-1}\otimes |\varnothing\rangle+ |\varnothing\rangle\otimes (xt^{-1}).
\end{equation}
\end{Lem}
\begin{proof}
We can assume that $x$ is an eigenvector for $\h$.
The analogous formula holds for $\ffr_P^0$ by Section \ref{SSS_ffrp_classical}. So the difference between the two sides of (\ref{eq:Cartan_image_ffrp}) lies in 
$\Ring[\Jet N]\otimes |\varnothing\rangle\subset 
\CDO(N)\otimes V_\kappa(\h)$ (the filtration degree zero component). Since $\ffr_P$ is graded both for the root lattice grading and the energy grading, Remark \ref{Rem:ffr_equivariance}, 
the degree of the difference for the root lattice grading is in the span of roots of $\fl$, and the energy degree is $-1$.
It follows that the difference is zero. 
\end{proof}



\sss\label{SSS:ffr_deformed} We also note that $\ffr_P$ can be extended to $\ffr_{P,\bone}:V_{\bone}(\g)\rightarrow \CDO(N)\otimes V_{\bone}(\lf)$ (so that $\ffr_P$ is obtained from this by setting $\mathbf{1}:=\kappa$). This homomorphism is also filtered, let $\ffr_{P,\bone}^0$ denote the associated graded homomorphism.  We also have the intermediate maps: $\theta_{\bone}, \pi_{\bone}$. They are filtered and their associated graded maps are going to be denoted by  
$\theta_{\bone}^0,  \pi_{\bone}^0$. We note that 
\begin{equation}\label{eq:ffr0_formula}
\ffr_{P,\bone}^0=(\operatorname{id}\otimes \operatorname{id}\otimes \pi_{\bone}^0)\circ \theta_{\bone}^0.
\end{equation}

\subsection{Transitivity} \label{ss:diglets}
Let $B$ be a Borel subgroup in $P$ such that $B_L = P \cap L$ is a Borel subgroup in $L$. Choose a maximal torus $H \subset B_L$ and consider the opposite Borel $B^-$. Let $\tilde{N} = R_u(B)$, $\tilde{N}^- = R_u(B^-)$. We set $\tilde{N}_L = \tilde{N}\cap L, \tilde{N}^-_L = \tilde{N}^-\cap L$ so that we have variety isomorphisms $N\times \tilde{N}_L\xrightarrow{\sim} \tilde{N}, 
\tilde{N}^-_L\times N^-\xrightarrow{\sim} \tilde{N}^-$ given by the multiplication maps $(g_1,g_2)\mapsto g_1g_2$.

We have the following vertex algebra homomorphisms:
\[
\mathsf{ffr}_{B}: V_{\kappa}(\g) \to \CDO(\tilde{N}) \otimes V_{\kappa}(\h),
\]
\[
\mathsf{ffr}_{P}: V_{\kappa}(\g) \to \CDO(N) \otimes V_{\kappa}(\lf),
\]
\[
\mathsf{ffr}_{B_L}: V_{\kappa'}(\lf) \to \CDO(N_L) \otimes V_{\kappa}(\h).
\]

Also note that we have an isomorphism $\iota:\CDO(\tilde{N})\xrightarrow{\sim}\CDO(N) \otimes \CDO(N_L)$ thanks to $N\times \tilde{N}_L\xrightarrow{\sim} \tilde{N}$ (cf. (\ref{eq:CDO_tensor_decomposition}) and Lemma 
\ref{Lem:CDO_tensor_decomposition}). The following claim is what we mean by the transitivity.

\begin{Prop}\label{Prop:ffr_transitivity}
The following diagram is commutative:
\[
\begin{tikzcd}
V_\kappa(\g) \arrow{r}{\mathsf{ffr}_P} \arrow{d}{\mathsf{ffr}_B}
& \CDO(N) \otimes V_{\kappa}(\lf) \arrow{d}{\operatorname{id}_N\otimes \mathsf{ffr}_{B_L}} \\
\CDO(\tilde{N}) \otimes V_{\kappa}(\h) \arrow{r}{\iota\otimes \operatorname{id}_\h} & \CDO(N) \otimes \CDO(N_L)\otimes V_{\kappa}(\h),  
\end{tikzcd}
\]
where we write $\operatorname{id}_N, \operatorname{id}_{\h}$ for the identity endomorphisms of $\CDO(N)$ and $V_\kappa(\h)$, respectively.
\end{Prop}

\begin{proof}

Consider the inclusions:
\[
B\tilde{N}^-=NL^\diamond N^-\subset G^\diamond=PN^-\subset  G
\]
They give rise to localization embeddings of vertex algebras (cf. Section \ref{SSS_CDO_local_filtr}):
\[
\CDO_\kappa(G) \xrightarrow{\lambda} \CDO_\kappa(NLN^-) \xrightarrow{\lambda'} \textup{CDO}_\kappa(NN_LHN_L^-N^-).
\]
Also let $\lambda_L$ denote the localization homomorphism $\CDO_\kappa(L)\rightarrow 
\CDO(N_LHN_L^-)$ so that under the identifications that come from the product maps
\begin{align*}
&\CDO_\kappa(NLN^-)\xrightarrow{\sim} \CDO(N)\otimes \CDO_\kappa(L)\otimes \CDO(N^-), \\
&\CDO(NN_LHN_L^-N^-)\xrightarrow{\sim}
\CDO(N)\otimes \CDO(N^-_L)\otimes\CDO_\kappa(H)\otimes \CDO(N^-_L)\otimes \CDO(N^-),
\end{align*}
we have 
\begin{equation}\label{eq:lambda_formula}
\lambda'=\operatorname{id}_N\otimes\lambda_L\otimes \operatorname{id}_{N^-}.
\end{equation}

Also recall intermediate homomorphisms from Section \ref{SSS:ffr_construction}
\begin{align*}
&\theta_P:
V_{\kappa}(\g) \hookrightarrow \CDO(N) \otimes \CDO_\kappa(L)\otimes V(\mathfrak{n}^-),\\
&\theta_B:V_{\kappa}(\g) \hookrightarrow \CDO(\tilde{N}) \otimes \CDO_\kappa(H)\otimes V(\widetilde{\mathfrak{n}}^-),\\
&\theta_{B_L}: V_{\kappa}(\lf) \hookrightarrow \CDO(N_L) \otimes \CDO_\kappa(H)\otimes V(\mathfrak{n}_L^-),\\
&\pi_N:V(\mathfrak{n}^-)\twoheadrightarrow \Ring,\\
&\pi_{\widetilde{N}}: V(\widetilde{\mathfrak{n}}^-)\twoheadrightarrow \Ring,\\
&\pi_{N_L}:V(\mathfrak{n}_L^-)\twoheadrightarrow \Ring.
\end{align*}
We have 
\begin{equation}\label{eq:vertex_subalg}
\begin{split}
&\CDO(\widetilde{N})\otimes \CDO_\kappa(H)\otimes V(\widetilde{\mathfrak{n}}^-)=\CDO(\widetilde{N}H\widetilde{N}^-)^{\Jet \widetilde{N}^-}\subset\\& 
\CDO(\widetilde{N}H\widetilde{N}^-)^{\Jet N^-}=\CDO(N)\otimes \CDO_\kappa(N_LHN_L^-)\otimes V(\mathfrak{n}^-).
\end{split}
\end{equation}
We claim that under this inclusion, we have
\begin{equation}\label{eq:theta_equality}
\theta_B= \lambda'\circ\theta_P,
\end{equation} 
an equality of maps 
$$V_\kappa(\g)\rightarrow \CDO(N)\otimes \CDO_\kappa(N_LHN_L^-)\otimes V(\mathfrak{n}^-).$$
Indeed, composing with the inclusion of the target into $\CDO(NL^\diamond N^-)$ both sides of 
(\ref{eq:theta_equality}) become the composition 
$V_\kappa(\g)\hookrightarrow \CDO_\kappa(G)\hookrightarrow \CDO_\kappa(NL^\diamond N^-)$.

Note that both vertex algebras in (\ref{eq:vertex_subalg}) have vertex algebra ideals generated by elements of the form $|\varnothing\rangle\otimes |\varnothing\rangle\otimes xt^{-1}$ with $x\in \mathfrak{n}^-$. Quotienting these ideals gives us the inclusion
\begin{equation}\label{eq:vertex_subalg2}
\begin{split}
&\CDO(\widetilde{N})\otimes \CDO_\kappa(H)\otimes V(\mathfrak{n}^-_L)=\CDO(\widetilde{N}H\widetilde{N}^-)^{\Jet \widetilde{N}^-}\subset\\& 
\CDO(\widetilde{N}H\widetilde{N}^-)^{\Jet N^-}=\CDO(N)\otimes \CDO_\kappa(N_LHN_L^-).
\end{split}
\end{equation}
We write $\pi$ for the quotient map $$\CDO(N)\otimes \CDO_\kappa(N_LHN_L^-)\otimes V(\mathfrak{n}^-)\rightarrow \CDO(N)\otimes \CDO_\kappa(N_LHN_L^-).$$  
Observe that the subalgebra $\CDO(\widetilde{N})\otimes \CDO_\kappa(H)\otimes V(\mathfrak{n}^-_L)$ is that of invariants for the right action of $\Jet N^-_L$ on $\CDO(N)\otimes \CDO_\kappa(N_LHN_L^-)$. 

We obtain $\ffr_B$ as the composition 
\begin{align*}&(\operatorname{id}_{\tilde{N}}\otimes\operatorname{id}_H\otimes \pi_{N_L})\circ \pi\circ \theta_B=[(\ref{eq:theta_equality})]=\\&
(\operatorname{id}_{\tilde{N}}\otimes\operatorname{id}_H\otimes \pi_{N_L})\circ \pi\circ \lambda'\circ\theta_P=[(\ref{eq:lambda_formula})]
\\&(\operatorname{id}_{\tilde{N}}\otimes\operatorname{id}_H\otimes \pi_{N_L})\circ (\operatorname{id}_N\otimes \lambda_L)\circ \ffr_P=\\
&\langle\operatorname{id}_N\otimes ([\operatorname{id}_{N_L}\otimes\operatorname{id}_H\otimes \pi_{N_L}]\circ \lambda_L)\rangle\circ \ffr_P.
\end{align*}
To finish the proof notice that the restriction of $[\operatorname{id}_{N_L}\otimes\operatorname{id}_H\otimes \pi_{N_L}]\circ \lambda_L$ to 
$V_\kappa(\mathfrak{l})\subset \CDO_\kappa(L)$ is nothing else but $\ffr_{B_L}$ by the construction in Section \ref{SSS:ffr_construction}.
\end{proof}

\subsection{Homomorphism between universal enveloping algebras}
In this section we will study the filitered associative algebra homomorphism $\Ffr_P:=\widehat{U}(\ffr_P)$.  

\sss
Recall, Section \ref{SSS_CDO_construction}, that $\CDO(N)$ comes with a default energy grading.
First, we are going to examine the structure of the algebras $\widehat{U}(\CDO(N))\subset \widetilde{U}(\CDO(N))$, see Sections \ref{SSS:Utilde} and \ref{SSS:ffr_construction} for the constructions of these algebras.
Recall the elements $y_0^i,\partial_{j,-1}\in \CDO(N)$, see Proposition  \ref{CDO_prop_nilpotent}. They give rise to elements 
$y^i_n, \partial_{j,m}\in \widehat{U}(\CDO(N))$ with $n,m\in \Z$. The commutation relation 
between these elements is $[\partial_{j,m},y^{i}_n]=\delta_{i,j}\delta_{n+m,0}$, the other brackets are zero. The energy degree of both $y^i_n$ and $\delta_{j,n}$ is $n$, while the filtration degree of the elements $y^i_n$ is $0$, and for the elements $\partial_{j,m}$ it is $1$. The algebra $\widetilde{U}(\CDO(N))$ is the completed Weyl algebra, cf. \cite[Section 5.3.3]{Frenkel_loop}. We will denote it by $\widetilde{D}(\Loop N)$.
The algebra $\widehat{D}(\Loop N):=\widehat{U}(\CDO(N))\subset \widetilde{D}(\Loop N)$ consists of all infinite sums of ordered monomials in the elements 
$y^i_n,\partial_{j,m}$ written in the increasing order of energy degree (in degree $0$ we write $y$'s before $\partial$'s) subject to the following conditions:
\begin{itemize}
\item The degree in the elements $\partial$ is bounded,
\item For any $N\in \Z$, there are only finitely many monomials that do not contain $y^i_n,\partial_{i,n}$ with $n>N$.
\end{itemize}

\sss\label{SSS_filt_compl_tensor}
Note that the algebra $\widehat{U}(\CDO(N)\otimes V_\kappa(\lf))$ (where the tensor product is taken over $\Ring$) is the filtered completed tensor product $\widehat{D}(\Loop N)\widehat{\otimes} \widehat{U}_\kappa(\hat{\lf})$ defined as follows. Note that
$\widehat{D}(\Loop N)\otimes \widehat{U}_\kappa(\hat{\lf})$ carries a tensor product filtration $(\widehat{D}(\Loop N)\otimes \widehat{U}_\kappa(\hat{\lf}))_{\leqslant n}$ that is an exhaustive algebra filtration. Then $\widehat{D}(\Loop N)\widehat{\otimes} \widehat{U}_\kappa(\hat{\lf})$ is the union of the completions of the filtered pieces $(\widehat{D}(\Loop N)\otimes \widehat{U}_\kappa(\hat{\lf}))_{\leqslant n}$ with respect to their topology. 
We deduce $\widehat{U}(\CDO(N)\otimes V_\kappa(\lf))=\widehat{D}(\Loop N)\widehat{\otimes} \widehat{U}_\kappa(\hat{\lf})$ from the PBW theorem: these two algebras consist of the same infinite sums.

\sss\label{SSS_ffr_univ_enveloping} We  write $\Ffr_P$ for the homomorphism 
$\widehat{U}_\kappa(\hat{\g})\rightarrow \widehat{D}(\Loop N)\widehat{\otimes} \widehat{U}_\kappa(\hat{\lf})$. Thanks to Remark 
\ref{Rem:ffr_equivariance}, this homomorphism is $\Loop P$-equivariant. Thanks to Proposition \ref{Prop:ffr_transitivity} we have the following commutative diagram:

\begin{equation}\label{eq:ffr_transitivity_algebra}
\begin{tikzcd}
\widehat{U}_\kappa(\hat{\g}) \arrow{r}{\mathsf{Ffr}_P} \arrow{d}{\mathsf{Ffr}_B}
& \widehat{D}(\Loop N) \widehat{\otimes} \widehat{U}_{\kappa}(\hat{\lf}) \arrow{d}{\operatorname{id}\otimes \mathsf{Ffr}_{B_L}} \\
\widehat{D}(\Loop\tilde{N}) \widehat{\otimes} \widehat{U}_{\kappa}(\hat{\h}) \arrow{r}{\sim} & \widehat{D}(\Loop N) \widehat{\otimes} \widehat{D}(\Loop N_L)\widehat{\otimes} \widehat{U}_{\kappa}(\hat{\h}).  
\end{tikzcd}
\end{equation}

\sss\label{SSS_ffr_un} As in the vertex algebra setting, Section \ref{SSS:ffr_deformed}, we have the universal version of $\ffr_P$, a homomorphism 
$$\Ffr_{P,\bone}:\widehat{U}_{\bone}(\hat{\g})\rightarrow \widehat{D}(\Loop N)\widehat{\otimes} \widehat{U}_{\bone}(\hat{\lf}).$$
We can also consider the associated graded homomorphism 
$\Ffr_{P,\bone}^0$ and the Rees algebra homomorphism $\Ffr_{P,\bone}^\hbar$.

%
%
\section{Free field realizations and $p$-centers}\label{S_ffr_p_cent}

The goal of this section is to describe the compatibility between the free field realization and suitable defined $p$-centers of the algebras on both sides. A closely related result has been obtained in \cite{AW} using formulas for the free field realization maps. We instead present a conceptual approach. 

\subsection{$p$-centers for restricted Lie algebras}

\sss We first begin with some generalities on $p$-centers. To do so, we need to recall the absolute and relative Frobenius maps. 

\sss We begin with the absolute Frobenius. Recall that for a commutative $\F_p$-algebra $\Sing$, the absolute Frobenius map $\mathfrak{Fr}_\Sing: s \mapsto s^p$ is a map of algebras, and this is natural in $\Sing$, i.e., for any map $\phi: \Sing \rightarrow \Ting$ of algebras, we have $\mathfrak{Fr}_\Ting \circ \phi = \phi \circ \mathfrak{Fr}_\Sing$.

As a consequence, given an $\F_p$-space $\sY$ (see Section \ref{SSS_spaces}), for any $\Sing$ the tautological maps
$$\sY(\mathfrak{Fr}_\Sing): \sY(\Sing) \xrightarrow{} \sY(\Sing)$$ yield a natural transformation from $\sY$ to itself, i.e., a map of spaces
$$\mathfrak{Fr}_{\sY}: \sY \rightarrow \sY;$$we refer to this as the {\it absolute Frobenius} for $\sY$. Moreover, this is a monoidal natural transformation from the identity map of $\F_p$-spaces to itself. I.e., for any map $\phi: \sY \rightarrow \mathscr{Z}$ of spaces, we again have $\on{Fr}_{\mathscr{Z}} \circ \phi = \phi \circ \mathfrak{Fr}_{\mathscr{Y}}$, and for any pair of spaces $\sY_1, \sY_2$, we have the identity 
$$\mathfrak{Fr}_{\sY_1 \times \sY_2} = \mathfrak{Fr}_{\sY_1} \times \mathfrak{Fr}_{\sY_2}.$$

\sss Let us turn to the relative Frobenius. Let $\Ring$ be a fixed $\F_p$-algebra. Given an $\Ring$-space $\sY$, we may view it as an $\F_p$-space equipped with a map $\sY \rightarrow \text{Spec }\Ring.$ In particular, the tautological commutative diagram 
$$\xymatrix{\sY \ar[d] \ar[r]^{\mathfrak{Fr}_{\sY}} & \sY \ar[d] \\ \text{Spec }\Ring \ar[r]^{\mathfrak{Fr}_\Ring}&  \text{Spec }\Ring }$$
yields the relative Frobenius map of $\Ring$-spaces
$$\on{Fr}_{\sY}: \sY \rightarrow \sY^{(1)} := \sY  \underset{\text{Spec }\Ring} \times \text{Spec }\Ring.$$
This is again natural and monoidal, i.e., given a map $\phi: \sY \rightarrow \mathscr{Z}$, in evident notation we have $\on{Fr}_{\mathscr{Z}} \circ \phi = \phi^{(1)} \circ \on{Fr}_{\sY}$, and with respect to the natural identification 
$$(\sY_1 \times \sY_2)^{(1)} \simeq \sY_1^{(1)} \times \sY_2^{(2)}$$
we have $\on{Fr}_{\sY_1 \times \sY_2} = \on{Fr}_{\sY_1} \times \on{Fr}_{\sY_2}.$ {In other words}, the assignment $\sY \mapsto \sY^{(1)}$ is a symmetric monoidal functor $(-)^{(1)}$ from $\Ring$-spaces to itself, and the map $\on{Fr}: \sY \rightarrow \sY^{(1)}$ carries the datum of a natural transformation between monoidal functors $\on{id} \rightarrow (-)^{(1)}$. 

In what follows, we omit the subscript $\sY$ when discussing the relative Frobenius $\on{Fr}: \sY \rightarrow \sY^{(1)}$, and similarly for the absolute Frobenius $\mathfrak{Fr}: \sY \rightarrow \sY.$ 

\sss A useful consequence of the preceding discussion is the following. 
\begin{Cor} Given a group $\Ring$-space $\mathscr{G}$, its relative Frobenius twist $\mathscr{G}^{(1)}$ inherits a natural structure of group $\Ring$-space, and $\on{Fr}: \mathscr{G} \rightarrow \mathscr{G}^{(1)}$ is a map of group $\Ring$-spaces. 
\end{Cor}

\begin{proof}That $\mathscr{G}^{(1)}$ is again naturally a group $\Ring$-space follows from the monoidality of $(-)^{(1)}$, and that $\on{Fr}$ is a map of group $\Ring$-spaces follows from the monoidality of the natural transformation $\on{Fr}$.  \end{proof}

A similar argument shows the following. 

\begin{Cor} Given a restricted Lie algebra $\ff$ in $\Ring$-spaces, its relative Frobenius twist $\ff\fr$ inherits a natural structure of restricted Lie algebra in $\Ring$-spaces, and $\on{Fr}: \ff \rightarrow \ff^{(1)}$ is a map of restricted Lie algebra $\Ring$-spaces. 
\end{Cor}

\sss\label{SSS:p_center} To discuss $p$-centers, the following setup will be convenient. 

Given two restricted Lie algebras $\ff_1$ and $\ff_2$ in $\Ring$-spaces, we have a natural associated $\Ring$-space $\uHom_{\on{ResLie}_\Ring}(\ff_1, \ff_2)$. Namely, to a $\Ring$-algebra $\Sing$, we define $\uHom_{\on{ResLie}_\Ring}(\ff_1, \ff_2)(\Sing)$ to be the set of all maps of restricted Lie algebra $\Sing$-spaces $$\phi: \ff_1 \underset{\Spec \Ring} \times \Spec \Sing \rightarrow \ff_2 \underset{\Spec \Ring} \times \Spec \Sing.$$ 
Given a map $\Sing \rightarrow \Ting$, the map $\uHom_{\on{ResLie}_\Ring}(\ff_1, \ff_2)(\Sing) \rightarrow \uHom_{\on{ResLie}_\Ring}(\ff_1, \ff_2)(\Ting)$ is given by base change, i.e., assigns to $\phi$ as above the composition \begin{align*}\ff_1 \underset{\Spec \Ring} \times \Spec \Ting \simeq  (\ff_1 \underset{\Spec \Ring} \times \Spec \Sing) \underset{\Spec \Sing} \times \Spec \Ting \\ \xrightarrow{\phi \times \on{id}}  (\ff_2 \underset{\Spec \Ring} \times \Spec \Sing) \underset{\Spec \Sing} \times \Spec \Ting \simeq \ff_2 \underset{\Spec \Ring} \times \Spec \Ting.\end{align*}

Given two (unrestricted) Lie algebras $\fl_1$ and $\fl_2$ in $\Ring$-spaces, we have a similarly defined $\Ring$-space $\uHom_{\on{Lie}_\Ring}(\fl_1, \fl_2)$, where one omits the restricted condition on the maps in the definition of the $\Sing$-points. 

Similarly, given two associative algebras $\mathfrak{A}_1$ and $\mathfrak{A}_2$ in $\Ring$-spaces, we have an analogously defined $\Ring$-space $\uHom_{\on{Alg}_\Ring}(\mathfrak{A}_1, \mathfrak{A}_2)$, and given two modules $\mathfrak{M}_1$ and $\mathfrak{M}_2$ in $\Ring$-spaces, we have $\uHom_{\on{Mod}_\Ring}(\mathfrak{M}_1, \mathfrak{M}_2)$.

\sss Let us define $\on{ResLie}_{\Ring}$ to be the category enriched in $\Ring$-spaces with objects given by restricted Lie algebra $\Ring$-spaces $\ff$, and morphisms from $\ff_1$ to $\ff_2$ given by $\uHom_{\on{ResLie}_{\Ring}}(\ff_1, \ff_2)$, with the evident composition rules. 

We define the enriched categories $\on{Lie}_\Ring, \on{Alg}_\Ring,$ and $\on{Mod}_\Ring$ similarly. 
\sss By construction, we have evident forgetful functors:
$$\xymatrix{ \on{ResLie}_\Ring \ar[rd]_{\Oblv^{\on{ResLie}}_{\on{Lie}}}  & & \on{Alg}_\Ring \ar[ld]^{\Oblv^{\on{Alg}}_{\on{Lie}}}  \\  & \on{Lie}_\Ring \ar[d]_{\Oblv^{\on{Lie}}_{\on{Mod}}} \\ &  \on{Mod}_\Ring;}$$
let us denote the two appearing compositions by $$\on{Oblv}^{\on{ResLie}}_{\on{Mod}} \simeq \on{Oblv}^{\on{Lie}}_{\on{Mod}} \circ \on{Oblv}^{\on{ResLie}}_{\on{Lie}} \quad \text{and} \quad \on{Oblv}^{\on{Alg}}_{\on{Mod}} \simeq \on{Oblv}^{\on{Lie}}_{\on{Mod}} \circ \on{Oblv}^{\on{Alg}}_{\on{Lie}}.$$When the source and target are clear from context, we will  sometimes denote all of these functors simply by $\Oblv$.  
\begin{Lem} The forgetful functor $\on{Oblv}^{\on{Alg}}_{\on{Lie}}: \on{Alg}_\Ring \rightarrow \on{Lie}_\Ring$ admits a left adjoint, i.e., 
$$U(-): \on{Lie}_\Ring \rightleftarrows \on{Alg}_\Ring: \on{Oblv}^{\on{Alg}}_{\on{Lie}}.$$    
\end{Lem}

\begin{proof} Given an $\Ring$-algebra $\Sing$, and Lie algebra $\Ring$-space $\fl$, recall that $\fl(\Sing)$ is a Lie algebra over $\Sing$, and in particular we may form its $\Sing$-linear enveloping algebra $U(\fl(\Sing))$. We set $U(\fl)(S) := U(\fl(\Sing)).$ 
Given a map $\Sing \rightarrow \Ting$, the composition 
$$\fl(\Sing) \rightarrow \fl(\Ting) \rightarrow U(\fl)(\Ting)$$
is tautologically a $\Sing$-linear map of Lie algebras, whence induces a unique map $U(\fl)(\Sing) \rightarrow U(\fl)(\Ting)$, and it is straightforward to see this defines an algebra $U(\fl)$ in $\Ring$-spaces with the required universal property. 
\end{proof}

\sss In particular, given restricted Lie algebras in $\Ring$-spaces, we can form their enveloping algebras, i.e., the functor
\begin{equation} \label{route101}\on{ResLie}_\Ring \xrightarrow{\on{Oblv}} \on{Lie}_\Ring \xrightarrow{U(-)} \on{Alg}_\Ring.\end{equation}
On the other hand, to describe the $p$-center, we would like to also consider the Frobenius twist of the underlying $\Ring$-modules. Namely, we have a functor 
$$\on{Fr}^*: \on{Mod}_\Ring \rightarrow \on{Mod}_\Ring,$$ defined as follows. For a module $\sM$ in $\Ring$-spaces, we have 
$$(\on{Fr}^*\sM)(\Sing) := \Sing \underset{\Sing}\otimes \sM(\Sing),$$where we are inducing along the absolute Frobenius $\mathfrak{Fr}: \Sing \rightarrow \Sing$. Given a map $\phi: \Sing \rightarrow \Ting$, with associated map $\sM(\phi): \sM(\Sing) \rightarrow \sM(\Ting)$, we set $$(\on{Fr}^* \sM)(\phi) := \phi \otimes \sM(\phi): \Sing \underset{\Sing} \otimes \sM(\Sing) \rightarrow \Ting \underset{\Ting} \otimes \sM(\Ting).$$ 
For functoriality, note that we have a tautological map 
$$\uHom_{\on{Mod}_\Ring}(\sM_1, \sM_2) \rightarrow \uHom_{\on{Mod}_\Ring}(\on{Fr}^* \sM_1, \on{Fr}^* \sM_2)$$
which on $\Sing$-points sends a map $\phi$ to $\on{id} \otimes \phi$.

In particular, we may also consider the functor
\begin{equation} \label{e:route102} \on{ResLie}_\Ring \xrightarrow{\on{Oblv}} \on{Mod}_\Ring  \xrightarrow{\on{Fr}^*} \on{Mod}_\Ring .\end{equation}
Formation of the $p$-center may then be stated as follows. 
\begin{Lem} There is a canonical natural transformation from \eqref{e:route102} to $\on{Oblv}^{\on{Alg}}_{\on{Mod}} \circ \eqref{route101}$, i.e., 
$$\fZ_p: \on{Fr}^* \circ \on{Oblv}^{\on{ResLie}}_{\on{Mod}} \rightarrow \on{Oblv}^{\on{Alg}}_{\on{Mod}} \circ U(-) \circ \on{Oblv}^{\on{ResLie}}_{\on{Lie}};$$
\end{Lem}

\begin{proof} Given a restricted Lie algebra $\Ring$-space $\ff$, and an $\Ring$-algebra $\Sing$, we have a natural map 
\begin{equation} \label{e:pcent}\Sing \underset{\Sing}\otimes \ff(\Sing) \rightarrow U(\ff)(\Sing), \quad \quad 1 \otimes X \mapsto X^p - X^{[p]},\end{equation}
and it is straightforward to see this is a natural transformation.
\end{proof}

\begin{Cor} \label{c:funct}Given a restricted Lie algebra $\ff$ in $\Ring$-spaces, the natural action of $\uHom_{\on{ResLie}_\Ring}(\ff, \ff)$ on the $p$-center factors through the natural map 
$$\uHom_{\on{ResLie}_\Ring}(\ff, \ff) \rightarrow \uHom_{\on{Mod}_\Ring}(\on{Fr}^* \ff, \on{Fr}^* \ff).$$
\end{Cor}

\begin{Rem}\label{r:pcencen}Given an algebra $\mathfrak{A}$ in $\Ring$-spaces, we can define its center $Z(\mathfrak{A})$, with $\Sing$-points the elements $z \in \mathfrak{A}(\Sing)$ which are central in $\mathfrak{A}(\Ting)$ for every $\Sing \rightarrow \Ting$. With this, $\fZ_p$ factors for a fixed restricted Lie algebra $\ff$ as $$\on{Fr}^* \ff \rightarrow Z(U(\ff)) \hookrightarrow U(\ff).$$Explicitly, on $\Ring$-points, this is simply the standard generators of the $p$-center of $U(\ff(\Ring))$, cf. Equation \eqref{e:pcent}.
\end{Rem}

\sss We would like to deduce from Corollary \ref{c:funct} the assertion that, if a group $\Ring$-space $\mathscr{G}$ acts on a restricted Lie algebra $\Ring$-space $\ff$, then the action of $\mathscr{G}$ on its $p$-center factors through the Frobenius twist $\mathscr{G}\fr$. 

To approach this, first note that given an $\Ring$-space $\sY$, for an $\Ring$-algebra $i: \Ring \rightarrow \Sing$, if we temporarily denote its $\Sing$ points by $\sY(\Sing, i)$, we have a canonical bijection 
$$\sY\fr(\Sing, i) \simeq \sY(\Sing, i \circ \mathfrak{Fr}).$$
From this, it easily follows that the Frobenius twist $\sM\fr$ of a module $\Ring$-space is again naturally a module $\Ring$-space, and similarly for restricted Lie algebras, Lie algebras, and algebras. 

\sss In general, there is no reason to expect an isomorphism $\on{Fr}^* \sM \simeq \sM\fr$. To address this, let us say a module $\Ring$-space $\sM$ is {\em quasicoherent} if for any map $\Sing \rightarrow \Ting$, the tautological map $$\Ting \underset{\Sing}{\otimes} \sM(\Sing) \rightarrow \sM(\Ting)$$is an isomorphism of $\Ting$-modules. Let us denote by $\on{QCoh}_\Ring$ the full subcategory of $\on{Mod}_\Ring$ consisting of quasicoherent module $\Ring$-spaces. 

\begin{Lem} \label{l:Fr*1} The category $\on{QCoh}_\Ring$ is preserved by the functors $\sM \mapsto \on{Fr}^* \sM$ and $\sM \mapsto \sM\fr$, and one has a canonical natural isomorphism
$$\on{Fr}^* \simeq (-)\fr: \QCoh_\Ring \rightarrow \QCoh_\Ring.$$   
\end{Lem}

\begin{proof} Given a quasicoherent $\Ring$-space $\sM$ and $\Ring$-algebra $i: \Ring \rightarrow \Sing$, we have 
$$(\on{Fr}^* \sM)(\Sing, i) = \Sing \underset{\Sing} \otimes \sM(\Sing) = \sM( \Sing, \mathfrak{Fr} \circ i) = \sM\fr(\Sing, i),$$and it is straightforward to check these identifications yield a natural isomorphism of functors. 
\end{proof}

Let us call a restricted Lie algebra $\Ring$-space quasicoherent if its underlying module $\Ring$-space is. 

\begin{Cor} \label{c:actiononpcenter}If a group $\Ring$-space $\mathscr{G}$ acts on a  quasicoherent restricted Lie algebra $\Ring$-space $\ff$, then the action of $\mathscr{G}$ on the $p$-center of $U(\ff)$ factors as the composition $\on{Fr}: \mathscr{G} \rightarrow \mathscr{G}\fr$ and the tautological action of the latter on $\ff\fr \simeq \on{Fr}^*(\ff).$
\end{Cor}

\begin{proof} This follows from combining Corollary \ref{c:funct} and Lemma \ref{l:Fr*1}. 
\end{proof}

\sss \label{sss_proqcoh} Finally, as our restricted Lie algebra $\Ring$-spaces of interest, namely Kac--Moody and Virasoro, are not quasicoherent restricted Lie algebras, we consider the following (minimal) modification. 

Let us call a functor $F: I \rightarrow \QCoh_\Ring, i \mapsto \sM_i$, a {\em pro-quasicoherent} module $\Ring$-space if $I$ is a countable codirected set, and for $i \rightarrow j$, the map $\sM_i(\Sing) \rightarrow \sM_j(\Sing)$ is surjective for all $\Ring$-algebras $\Sing$. In what follows, we denote such an object as a formal inverse limit $``\varprojlim_i" \sM_i$. 

Pro-quasicoherent $\Ring$-spaces naturally form a category $\on{ProQCoh}_\Ring$ enriched in $\Ring$-spaces, where we set
$$\uHom_{\on{ProQCoh}_\Ring}( ``\varprojlim_\alpha" \sM_\alpha, ``\varprojlim_\beta" \sN_\beta) := \varprojlim_\beta \varinjlim_\alpha \uHom_{\on{Mod}_\Ring}(\sM_\alpha, \sN_\beta).$$
We have a tautological conservative and faithful functor $$\on{ProQCoh}_\Ring \rightarrow \on{Mod}_\Ring, \quad \quad ``\varprojlim M_i" \mapsto \varprojlim M_i,$$
which sends a formal inverse limit to the actual limit computed in $\on{Mod}_\Ring$. Explicitly, we have $$(\varprojlim M_i)(\Sing) \simeq \varprojlim (M_i(\Sing)),$$i.e., the limit may be computed termwise. Note that the countability of the diagrams and surjectivity of the transition maps ensures this is indeed conservative and faithful. In particular, we may discuss the underlying $\Ring$-space of a pro-quasicoherent $\Ring$-space module. 

\sss We define the category of restricted Lie algebras in pro-quasicoherent $\Ring$-spaces as the fiber product
$$\on{ResLie}(\Pro\QCoh_\Ring) := \Pro\QCoh_\Ring \underset{\on{Mod}_\Ring} \times \on{ResLie}_\Ring.$$
In particular, given a pro-quasicoherent $\Ring$-space $``\varprojlim" \sM_i$, a restricted Lie algebra structure on it is  the datum of a restricted Lie algebra structure on $\varprojlim \sM_i$, such that the Lie bracket and restricted power maps lift to (unique) maps from $``\varprojlim" \sM_i$ to itself. In particular, this implies that both structures are continuous with respect to pro-topologies on $\varprojlim \sM_i(\Sing)$, for all $\Sing$. 

We define the category $\on{Lie}(\Pro\QCoh_\Ring)$ similarly, so in particular we have forgetful functors
$$\on{ResLie}(\Pro\QCoh_\Ring) \rightarrow \on{Lie}(\Pro\QCoh_\Ring) \rightarrow \Pro\QCoh_\Ring.$$

\sss \label{sss:completedenvelope} We have a natural functor of taking the completed enveloping algebra
$$\widetilde{U}(-): \on{Lie}(\Pro\QCoh_\Ring) \rightarrow \on{Alg}_\Ring, \quad \quad \ff \mapsto \widetilde{U}(\ff).$$
Explicitly, on objects, if the underlying pro-quasicoherent $\Ring$-space of $\ff$ is $``\varprojlim" \ff_\alpha, \alpha \in A$, for fixed $\beta \in A$ consider the `compact open' subspace $$\fk_\beta := ``\underset{\alpha \geqslant \beta}\varprojlim" \ker(\ff_\alpha \rightarrow \ff_\beta).$$Then on $\Sing$-points, we set
$$(\widetilde{U}(\ff))(\Sing) := \varprojlim U( \ff(\Sing)) / U(\ff(\Sing)) \cdot \fk_\alpha(\Sing), $$
with its natural structure as a $\Sing$-algebra, and with the natural functoriality along maps $\ff(\Sing) \rightarrow \ff(\Ting)$. The definition of what $\widetilde{U}(-)$ does to morphisms, i.e. the map 
$$\uHom_{\on{Lie}(\Pro\QCoh_\Ring)}(\ff_1, \ff_2) \rightarrow \uHom_{\on{Alg}_\Ring}(\widetilde{U}(\ff_1), \widetilde{U}(\ff_2))$$is similarly induced by applying the functoriality of taking the completed enveloping algebra of a given topological Lie algebra $\ff(\Sing)$.  

In particular, when discussing $p$-centers, we may again may consider the composition
\begin{equation} \label{route202} \on{ResLie}(\Pro\QCoh_\Ring) \xrightarrow{\on{Oblv}} \on{Lie}(\Pro\QCoh_\Ring) \xrightarrow{\widetilde{U}(-)} \on{Alg}_\Ring \xrightarrow{\Oblv} \on{Mod}_\Ring. \end{equation}

On the other hand, we may consider the functor
\begin{equation} \label{route201}\on{ResLie}(\Pro\QCoh_\Ring) \xrightarrow{\Oblv} \on{Mod}_\Ring \xrightarrow{(-)\fr} \on{Mod}_\Ring. \end{equation}

\begin{Prop} \label{p:pcenpro}There exists a canonical natural transformation from \eqref{route201} to \eqref{route202}, i.e., 
$$\xi_p: (-)\fr \circ \on{Oblv}^{\on{ResLie}}_{\on{Mod}} \rightarrow \on{Oblv}^{\on{Alg}}_{\on{Mod}} \circ \widetilde{U}(-) \circ \on{Oblv}^{\on{ResLie}}_{\on{Lie}}.$$  
\end{Prop}

\begin{proof} Suppose that $\ff \simeq ``\varprojlim" \ff_\alpha$ is a restricted Lie algebra in pro-quasicoherent $\Ring$-spaces. 

We first note that $(\on{Oblv}^{\on{ResLie}}_{\on{Mod}}(\ff))\fr$ is computed as follows. Using that Frobenius twist commutes with limits of $\Ring$-spaces, and that each $\ff_\alpha$ is quasicoherent, for an $\Ring$-algebra $\Sing$, we have
\begin{align*}
  (\on{Oblv}(\ff))\fr(\Sing)  = (\varprojlim \ff_\alpha)\fr(\Sing) & \simeq (\varprojlim \ff_\alpha\fr)(\Sing) \\ & \simeq \varprojlim (\ff_\alpha\fr(\Sing)) \\ & \simeq \varprojlim \Sing \underset{\Sing} \otimes \ff_\alpha(\Sing).
\end{align*}
To prove the proposition, it therefore suffices to check that the composition
$$\Sing \underset{\Sing} \otimes \ff(\Sing) \xrightarrow{\fZ_p} U(\ff(\Sing)) \rightarrow \widetilde{U}(\ff(\Sing))$$
is continuous, i.e., factors through $\varprojlim \Sing \underset{\Sing} \otimes \ff_\alpha(\Sing).$ I.e., for fixed $\beta$, we must show the composition
$$\Sing \underset{\Sing} \otimes \ff(\Sing) \xrightarrow{\fZ_p} U(\ff(\Sing)) \rightarrow U(\ff(\Sing)) / U(\ff(\Sing)) \cdot \fk_\beta,  \quad 1 \otimes X \mapsto X^p - X^{[p]} \text{ mod } U(\ff(\Sing)) \cdot \fk_\beta.$$
factors through $\Sing \underset{\Sing} \otimes \ff_\gamma(\Sing)$, for some $\gamma$, i.e., annihilates  $1 \otimes \fk_\gamma$. If $\gamma \geqslant \beta$, then $X^p \in U(\ff(\Sing)) \cdot \fk_\beta$, for any $X \in \fk_\gamma$. By the continuity of the restricted power map, for all sufficiently large $\gamma'$ we have $X^{[p]} \in \fk_{\beta}$ for all $X \in \fk_{\gamma'}$, and this $\gamma'$ may be chosen uniformly in $\Sing$. In particular, taking such a $\gamma'$ to be greater than $\beta$, we see the composition annihilates $1 \otimes \fk_{\gamma'}$, as desired. 
\end{proof}

\begin{Rem} Note that again, for a fixed $\ff$, the map $\xi_p: \ff\fr \rightarrow \widetilde{U}(\ff)$ factors through the center of the latter, i.e., $Z(\widetilde{U}(\ff))$. Indeed, this follows by continuity from the analogous assertion for uncompleted enveloping algebras, cf. Remark \ref{r:pcencen}. 
\end{Rem}

\sss Let us explicitly state the desired consequence of Proposition \ref{p:pcenpro} which will be used going forwards. 

\begin{Cor}\label{p:p-cent} Suppose $\ff$ is a pro-quasicoherent restricted Lie algebra in  $\Ring$-spaces. Then we have a canonical $p$-center map of $\Ring$-spaces
$$\xi_p: \ff\fr \rightarrow Z(\widetilde{U}(\ff)).$$
Moreover, if $\mathscr{G}$ is a group $\Ring$-space acting on $\ff$ by restricted Lie algebra automorphisms, then the map $\xi_p$ is $\mathscr{G}$-equivariant, where $\mathscr{G}$ acts on $\ff\fr$ via the relative Frobenius $\on{Fr}: \mathscr{G} \rightarrow \mathscr{G}\fr$, and the tautological action of the latter group $\Ring$-space on $\ff^{(1)}$. 
\end{Cor} 

In particular, one may take $\ff$ to be a Kac--Moody  or Virasoro Lie algebra, as in Section \ref{s:kmspace} and $\mathscr{G}$ to be the loop group $F(\!(t)\!)$ or $\on{Aut}(\D^\times)$, respectively.

\subsection{$p$-centers of $V_\kappa(\ff)$ and $\CDO_\kappa(F)$}
Let $F$ be a smooth algebraic group over $\Ring$, where $\Ring$ is a perfect field containing $\F_p$. Fix representations $W_1,\ldots,W_k$ of $F$ and form the corresponding group ind-scheme $F(\!(t)\!)^\flat$, see 
Section \ref{SSS_KM_groups}, so that we can talk about the $p$-center of $U(\ff(\!(t)\!)^\flat)$, see Section \ref{SSS:p_center}. In particular, this gives rise to the $p$-center in $U_\kappa(\hat{\ff})$.  

\sss\label{SSS_vertex_p_center}
Recall the $p$-center of $V_\kappa(\ff)$ following \cite[Section 2.3]{AW} and \cite[Section 3.3]{ATV}. By the construction of $V_\kappa(\ff)$ we can view it as a quotient of $U_\kappa(\hat{\ff})$, where $1$  goes to $|\varnothing\rangle$. By definition, the $p$-center of $V_\kappa(\ff)$, to be denoted by $\mathfrak{z}_{Fr}(V_\kappa(\ff))$, it is the image of the $p$-center in $V_\kappa(\ff)$. This is a central vertex subalgebra, see  \cite[Lemma 3.7]{ATV}.
The algebra  $\mathfrak{z}_{Fr}(V_\kappa(\ff))$ is identified with the algebra of regular functions on the space of $\AS(\kappa)$-connections $\AS(\kappa)\partial+ (\ff^*[[t]]dt)^{(1)}$. Here and below we write $\AS$ for the Artin-Schreier map $\Ring\rightarrow \Ring, \kappa\mapsto \kappa^p-\kappa$. 

\sss Now we proceed to defining and studying the $p$-center of $\CDO_\kappa(F)$.  
Consider the subalgebra $\Ring[(\Jet F)^{(1)}]\subset \Ring[\Jet F]$.


\begin{Lem}\label{Lem:CDO_p_center}
The following claims hold:
\begin{enumerate}
\item 
The subalgebras $\Ring[(\Jet F)^{(1)}],\iota_l(\mathfrak{z}_{Fr}(V_\kappa(\ff)))$
of $\CDO_\kappa(F)$ are central.
\item These subalgebras freely generate the center of $\CDO_\kappa(F)$. 
\end{enumerate}
\end{Lem}
\begin{proof}
(1): Note that, as a vertex algebra, $\CDO_\kappa(F)$ is generated by $\Ring[\Jet F]$
and $V_\kappa(\ff)$. Every endomorphism $xt^n$ for $x\in \ff$ and $n\geqslant 0$ of $\Ring[\Jet F]$ annihilates $\Ring[(\Jet F)^{(1)}]$. It follows that the elements $b\in \Ring[(\Jet F)^{(1)}]$ satisfy $Y(a,z)b\in V_\kappa(\ff)[[z]]$ for all $a\in V_\kappa(\ff)$
or $\Ring[\Jet F]$, so are central. 

Since $\mathfrak{z}_{Fr}(V_\kappa(\ff))$ is central in $V_\kappa(\ff)$, to check that $\mathfrak{z}_{Fr}(V_\kappa(\ff))$ is in the center of $\CDO_\kappa(F)$ it remains to show that $a_{(n)}f=0$ for all $a\in \mathfrak{z}_{Fr}(V_\kappa(\ff)),f\in \Ring[\Jet F]$ and $n\geqslant 0$. We note that $\mathfrak{z}_{Fr}(V_\kappa(\ff))$ is a commutative vertex algebra that is generated, as an associative algebra, by the elements $(xt^i)^p-x^{[p]}t^{ip}, i<0$. So we need to show that 
\begin{equation}\label{eq:p_center_vanishing}
[(xt^i)^p-x^{[p]}t^{ip}]_{(n)}f=0,  \forall n\geqslant 0.\end{equation}
This expression can be computed using \cite[Lemma 2.3]{AW}, which shows that for $j<0$ we have
\begin{equation}\label{eq:Y_pth_power}
Y((xt^{j})^p|\varnothing\rangle, z)=\sum_{n\in \Z}{-n-1\choose -j-1}(xt^n)^p z^{(j-n)p}.
\end{equation}
in $V_\kappa(\mathfrak{f})$. 
Note that the direct analog of (\ref{eq:Y_pth_power}) holds for any vertex algebra containing $V_\kappa(\g)$ including $\CDO_\kappa(G)$ with the same proof as in \cite{AW}. Since $\Ring[\Jet F]$ is a rational representation of the pro-algebraic group $\Jet F$, the elements $[xt^{in/p}]^p-x^{[p]}t^{in}=[xt^{in/p}]^p-[xt^{in/p}]^{[p]}$ act on $\Ring[\Jet F]$ by zero for all $n\geqslant 0$. Together with the generalization of
(\ref{eq:Y_pth_power}) for modules this implies (\ref{eq:p_center_vanishing}). This finishes the proof of (1).

(2): We note that since the naive filtration on $\CDO_\kappa(F)$ is almost commutative, the associated graded $\gr \CDO_\kappa(F)$ carries a natural Poisson vertex algebra structure, see Section \ref{SSS_Poisson_vertex}.  Recall the identification  $\gr\CDO_\kappa(F)\cong \Ring[\Jet(T^*F)]$,
a special case of (\ref{eq:gr_CDO}). We claim that the restriction of the bracket on $\gr\CDO_\kappa(F)$ to $\Ring[T^*F]\subset \Ring[\Jet(T^*F)]$ is the standard bracket. These easily follows from the identities $\iota_l(xt^{-1})_{(0)}\iota_l(yt^{-1})=\iota_l([x,y]t^{-1})$ (that is true because $\iota_l$ is an embedding of vertex algebras)
and $\iota_l(xt^{-1})_{(0)}f=x.f$ for $x\in \ff, f\in \Ring[F]$ that is a consequence of (\ref{eq:bracket_CDO}).

We will show that the Poisson center of $\Ring[\Jet(T^*F)]$ coincides with $\Ring[(\Jet(T^*F))^{(1)}]$, this will imply the claim of the lemma.
To show this, we will use the preliminary fact that if $X \rightarrow Y$ is an etale map of smooth varieties, then the natural map $\pi: X \xrightarrow{\sim} X^{(1)} \times_{Y^{(1)}} Y$ is an isomorphism. Indeed, by considering the maps to $X^{(1)}$, it follows that $\pi$ is finite, and by considering the maps to $Y$ it follows that $\pi$ is etale. By considering again the map to $X^{(1)}$, it follows $\pi$ is moreover finite etale of degree one, i.e., an isomorphism. In particular, we deduce that \begin{equation} \label{e:cacaw} T^* X \cong X \times_Y T^* Y \cong X^{(1)} \times_{Y^{(1)}} T^*Y.\end{equation}

We will apply the  isomorphism \eqref{e:cacaw} as follows. We can cover $F$ by open affines such that each of these open affines admits an etale morphism to $\mathbb{A}^n$, where $n$ is the dimension of $X$ over $\operatorname{Spec}(\Ring)$. Let $F^\diamond$ be one of these open affines.  We then have isomorphisms
\begin{equation}\label{eq:iso_etale_lift}
\Jet(T^*F^\diamond)\cong \Jet(F^{\diamond,(1)}\times_{\mathbb{A}^{n,(1)}} T^* \mathbb{A}^n ) \cong F^{\diamond,(1)}\times_{\mathbb{A}^{n,(1)}}\Jet(T^*\mathbb{A}^n),
\end{equation}
Note that $\Ring[\Jet(T^*\mathbb{A}^n)]$ carries a Poisson vertex algebra structure whose restriction to $\Ring[T^*\mathbb{A}^n]$ is the usual one (for example, as the associated graded of the chiral differential operators on the additive group $\mathbb{A}^n$). 
There is a unique extension of the Poisson structure from 
$\Ring[\Jet(T^*\mathbb{A}^n)]$ to $\Ring[F^{\diamond,(1)}\times_{\mathbb{A}^{n,(1)}}\Jet(T^*\mathbb{A}^n)]$ such that $\Ring[F^{\diamond,(1)}]$ is central.
Now (\ref{eq:iso_etale_lift}) gives two Poisson vertex algebra structures on $\Ring[\Jet(T^*F^\diamond)]$, whose restrictions to $\Ring[T^*F^\diamond]$ coincide with the usual bracket. By Section \ref{SSS_Poisson_structure_uniqueness}, they coincide. 

Let $X$ be a smooth affine finite type $\Ring$-scheme.  
Note that the $\mathbb{G}_m$-action on $T^*X$ by fiberwise dilations equips $\Ring[\Jet(T^*X)]$ with an energy grading. The graded pieces are finitely generated modules over $\Ring[X^{(1)}]$. 
The unique  Poisson vertex algebra structure on $\Ring[\Jet(T^*X)]$ such that the bracket on $\Ring[T^*X]$ is the standard one has the following properties: the Poisson structure is compatible with the energy grading, and the Poisson center is graded. 

Apply this to $X=T^*F^\diamond$.
Notice that, for each $m$, the equations giving the degree $m$ component of the Poisson center in $\Ring[F^{\diamond,(1)}\times_{\mathbb{A}^{n,(1)}}\Jet(T^*\mathbb{A}^n)]$ are obtained by pulling back the equations defining the degree $m$ component in the Poisson center of $\Ring[\Jet(T^*\mathbb{A}^n)]$. This reduces our task to showing that the Poisson center of  $\Ring[\Jet(T^*\mathbb{A}^n)]$ coincides with
$\Ring[\Jet(T^*\mathbb{A}^n)^{(1)}]$. This is verified by a direct check and is left as an exercise.
\end{proof}

\begin{Rem}\label{Rem:center_CDO_localization}
In fact, the proof implies that the subalgebras $\Ring[(\Jet F^\diamond)^{(1)}],\iota_l(\mathfrak{z}_{Fr}(V_\kappa(\ff)))$
freely generate the center of $\CDO_\kappa(F^\diamond)$. 
\end{Rem}

\sss Now we are going to give a description of the center  $\mathfrak{z}(\CDO_\kappa(G))$ as a commutative energy graded vertex algebra with an arc group action. 

Recall that $\AS(\kappa)$ stands for $\kappa^p-\kappa$. Recall the commutative vertex algebra $\gr \CDO_{\bone}(F)$. This algebra is a deformation 
of $\Ring[\Jet (T^*F)]$ over $\Ring[\mathbf{1}]$ described in Section \ref{SSS_arc_actions_CDO}. Consider the Frobenius twist $(\gr \CDO_{\bone}(F))^{(1)}$, an algebra over $\Ring[\mathbf{1}]^{(1)}$, and the quotient of this algebra by $\mathbf{1}-\AS(\kappa)$, to be denoted by $\Ring_{\AS(\kappa)}[\Jet(T^*F)^{(1)}]$.
The $\Jet F\times \Jet F$-action on $\gr \CDO_{\bone}(F)$,
see Section \ref{SSS_arc_actions_CDO}, gives rise to to an action of $(\Jet F)^{(1)}\times (\Jet F)^{(1)}$ on the quotient
$\Ring_{\AS(\kappa)}[\Jet(T^*F)^{(1)}]$.

It is an easy exercise to show that a vertex action of an arc group on a vertex algebra preserves the center. In particular, the group $\Jet F\times \Jet F$ acts on $\mathfrak{z}(\CDO_\kappa(F))$. 

Note that both $\Ring_{\AS(\kappa)}[\Jet(T^*F)^{(1)}]$ and $\mathfrak{z}(\CDO_\kappa(F))$ come with natural energy gradings. 

\begin{Lem}\label{Lem:CDO_center_description}
There is a $\Jet F\times \Jet F$-equivariant energy graded vertex algebra isomorphism
$$\Ring_{\AS(\kappa)}[\Jet(T^*F)^{(1)}]\xrightarrow{\sim} \mathfrak{z}(\CDO_\kappa(F)),$$
where the action of $\Jet F\times \Jet F$ on the source is pulled back from the action of $(\Jet F)^{(1)}\times (\Jet F)^{(1)}$.
\end{Lem}
\begin{proof}
By construction, $\Ring_{\AS(\kappa)}[\Jet(T^*F)^{(1)}]$ is freely generated by the subalgebras $$\Ring[(\Jet F)^{(1)}] \text{ and }
\Ring\left[\AS(\kappa)\partial+ (\ff^*[[t]]dt)^{(1)}\right].$$ These subalgebras are $(\Jet F\times \Jet F)^{(1)}$-stable, where $\Ring\left[\AS(\kappa)\partial+ (\ff^*[[t]]dt)^{(1)}\right]$ is invariant for the action from the right, while the action from the left comes from the adjoint action. Lemma \ref{Lem:CDO_p_center} says that $\mathfrak{z}(\CDO_\kappa(F))$ is freely generated by $\Ring[(\Jet F)^{(1)}]$ and $\iota_l(\mathfrak{z}_{Fr}(V_\kappa(\ff)))$. Recall, 
Section \ref{SSS_vertex_p_center}, that 
$\mathfrak{z}_{Fr}(V_\kappa(\ff))$ is identified with 
$\Ring\left[\AS(\kappa)\partial+ (\ff^*[[t]]dt)^{(1)}\right]$ and, by Corollary \ref{p:p-cent}, this identification is $\Jet F$-equivariant. And the claim that the identification of the subalgebras $\Ring[(\Jet F)^{(1)}]$ is equivariant is manifest. The claim of the lemma follows.
\end{proof}

\subsubsection{} We would like to relate the isomorphism in Lemma \ref{Lem:CDO_center_description} to its ''right handed" version. Namely, consider 
the embedding $\iota_r:V_{-\kappa}(\ff)\hookrightarrow \CDO_\kappa(F)$ from Section \ref{SSS_iota_r}. It restricts to an embedding $\mathfrak{z}_{Fr}(V_{-\kappa}(\ff))\hookrightarrow \mathfrak{z}(\CDO_\kappa(F))$. Using identifications of 
Section \ref{SSS_vertex_p_center} and 
Lemma \ref{Lem:CDO_center_description}, we get an embedding 
\begin{equation}\label{eq:Fr_center_right}
\Ring[-\AS(\kappa)\partial+(\ff^*[[t]]dt)^{(1)}]\hookrightarrow \Ring_{\AS(\kappa)}[\Jet(T^*F)^{(1)}]
\end{equation}
to be denoted by $\iota'_r$.

\begin{Lem}\label{Lem:Fr_center_right}
(\ref{eq:Fr_center_right}) coincides with the embedding given by the moment map for the action of $(\Jet F)^{(1)}$ from the right, to be denoted by $\mu_r^*$.
\end{Lem}
\begin{proof}
It is enough to prove the coincidence of these embeddings on the generating subspace
$t^{-1}\ff[t^{-1}]$. Let $\varphi: t^{-1}\ff[t^{-1}]\rightarrow \Ring_{\AS(\kappa)}[\Jet(T^*F)^{(1)}]$ denote the difference $\iota'_r-\mu_r^*$ restricted to $t^{-1}\ff[t^{-1}]$. The embeddings $\iota'_r,\mu_r^*$ coincide on the level of the associated graded algebras, hence the image of $\varphi$ lies in $\Ring[(\Jet F)^{(1)}]$.
Next, both $\iota_r',\mu_r^*$ are invariant for the action of $(\Jet F)^{(1)}$ from the left, hence so is $\varphi$. The only invariants in
$\Ring[(\Jet F)^{(1)}]$ are the scalars. Finally, both $\iota'_r,\mu_r^*$ are energy graded, hence so is $\varphi$. This shows $\varphi=0$.
\end{proof}



\subsection{Main result}
We are going to show that $\ffr_P: V_{\kappa}(\g)\rightarrow \CDO(N)\otimes V_\kappa(\lf)$ restricts to a map between the $p$-centers and describe the restriction. 

Consider the homomorphism (cf. Section \ref{SSS:ffr_deformed})
\begin{equation}\label{eq:ffr_p_center}
\ffr^{0,(1)}:\Ring\left[\AS(\kappa)\partial+ (\g^*[[t]]dt)^{(1)}\right]\rightarrow 
\Ring[\Jet (T^*N))^{(1)}]\otimes\Ring\left[\AS(\kappa)\partial+ (\lf^*[[t]]dt)^{(1)}\right]
\end{equation}
We identify $\frz_{Fr}(V_{\kappa}(\g))$ with $\Ring\left[\AS(\kappa)\partial+ (\g^*[[t]]dt)^{(1)}\right]$ and 
$\mathfrak{z}_{Fr}(\CDO(N))\otimes \mathfrak{z}_{Fr}(V_\kappa(\lf)))$ (which is a central subalgebra $\CDO(N)\otimes V_\kappa(\lf)$ in to be denoted by $\mathfrak{z}_{Fr}(\CDO(N)\otimes V_\kappa(\lf))$) with $\Ring[\Jet (T^*N))^{(1)}]\otimes\Ring\left[\AS(\kappa)\partial+ (\lf^*[[t]]dt)^{(1)}\right]$. So we can view 
$\ffr^{0,(1)}$ as a homomorphism 
$$\mathfrak{z}_{Fr}(V_\kappa(\g))\rightarrow 
\mathfrak{z}_{Fr}(\CDO(N)\otimes V_\kappa(\lf)).$$


\begin{Prop}\label{Prop:ffr_p_center}
The restriction of the homomorphism $\ffr_{P,\kappa}$ to 
$\frz_{Fr}(V_\kappa(\g))$ coincides with 
$\ffr^{0,(1)}_{P}$. 
\end{Prop}
\begin{proof}
Recall, Definition \ref{defi:ffrP}, that we realize $\ffr_{P,\kappa}$ as the composition of 
$$\theta: V_{\kappa}(\g) \hookrightarrow \CDO(N) \otimes \CDO_\kappa(L)\otimes V(\mathfrak{n}^-)$$
and the projection
$$\operatorname{id}_N\otimes\operatorname{id}_L\otimes\pi: \CDO(N)\otimes \CDO_\kappa(L)\otimes V(\mathfrak{n}^-)\twoheadrightarrow \CDO(N)\otimes \CDO_\kappa(L),$$
where $\pi$ is the projection $V(\mathfrak{n}^-)\twoheadrightarrow \Ring$.
The map $\ffr^{0,(1)}_P$ admits an analogous description: it is the composition of the embedding 
$$\theta^{0,(1)}:\Ring\left[\AS(\kappa)\partial+ (\g^*[[t]]dt)^{(1)}\right]\hookrightarrow \Ring[\Jet (T^*N))^{(1)}]\otimes\Ring\left[\AS(\kappa)\partial+ (\lf^*[[t]]dt)^{(1)}\right]\otimes \Ring[(\Jet\mathfrak{n}^{-,*})^{(1)}]$$
induced by the moment map and the projection that we are going to denote by $\operatorname{id}_N^{0,(1)}\otimes\operatorname{id}^{0,(1)}_L\otimes\pi^{0,(1)}$.

It is clear that the restriction 
of $\pi$ to $$\mathfrak{z}_{Fr}(V(\mathfrak{n}^-))\cong \Ring[(\Jet\mathfrak{n}^{-,*})^{(1)}]$$ coincides with its analog $\pi^{0,(1)}$. So it is sufficient to prove that 
\begin{itemize}
\item[(*)]
the restriction of $\theta$
to $\mathfrak{z}_{Fr}(V_\kappa(\g))$ coincides with the composition of $\theta^{0,(1)}$ and the inclusion 
\begin{align*}&\Ring[(\Jet T^*N)^{(1)}]\otimes \Ring_{\AS(\kappa)}[(\Jet T^*L)^{(1)}]\otimes \Ring[(\Jet\mathfrak{n}^{-,*})^{(1)}]\hookrightarrow \\&\CDO(N) \otimes \CDO_\kappa(L)\otimes V(\mathfrak{n}^-)\end{align*}
(as the p-center).
\end{itemize}

Thanks to Lemma \ref{Lem:CDO_center_description}, we have the following commutative diagram, where the top vertical maps are identifications and the bottom vertical maps are inclusions:
$$\xymatrix{\Ring\left[\AS(\kappa)\partial+ (\g^*[[t]]dt)^{(1)}\right] \ar[d] \ar[r] & \Ring_{\AS(\kappa)}[\Jet(T^*G)^{(1)}] \ar[d]\ar[r]& \Ring_{\AS(\kappa)}[\Jet(T^*G^\diamond)^{(1)}]\ar[d] \\
\mathfrak{z}_{Fr}(V_\kappa(\g)) \ar[d] \ar[r] & \mathfrak{z}(\CDO_\kappa(G)) \ar[d]\ar[r]& \mathfrak{z}(\CDO_\kappa(G^\diamond))\ar[d] \\ V_\kappa(\g) \ar[r]&  \CDO_\kappa(G)\ar[r] &\CDO_\kappa(G^\diamond)}$$
Next,  we claim that 
the following diagrams are commutative (the vertical arrows are still inclusions):
$$\xymatrix{\Ring_{\AS(\kappa)}[(\Jet T^*P)^{(1)}]\otimes \Ring[(\Jet T^*N^-)^{(1)}] \ar[d] \ar[r] & \Ring_{\AS(\kappa)}[(\Jet T^*G^\diamond)^{(1)}] \ar[d] \\ \CDO'_\kappa(P)\otimes \CDO(N^-) \ar[r]&  \CDO_\kappa(G^\diamond)}$$
$$\xymatrix{\Ring[(\Jet T^*N)^{(1)}]\otimes \Ring_{\AS(\kappa)}[(\Jet T^*L)^{(1)}]  \ar[d] \ar[r] & \Ring_{\AS(\kappa)}[(\Jet T^*P)^{(1)}] \ar[d] \\ \CDO(N)\otimes \CDO_\kappa(L) \ar[r]&  \CDO'_\kappa(P)}$$
(*) easily follows from these commutative diagrams. To see that they are indeed commutative, it is enough to show that the maps from the top left corner to the bottom right corner coincide on the generators, e.g., $$\Ring[(\Jet P)^{(1)}], \Ring[(\Jet N^-)^{(1)}],
\Ring[\AS(\kappa)\partial+ (\mathfrak{p}^*[[t]]dt)^{(1)}], 
\Ring[(\Jet\mathfrak{n}^{-,*})^{(1)}]$$
for the first diagram. For the first three subalgebras, the coincidence of the maps follows easily from the construction, while for the last one the coincidence is deduced from Lemma \ref{Lem:Fr_center_right}.
\end{proof}

\subsection{Map between $p$-centers of universal enveloping algebras}\label{SSS_univ_env_map}
Our goal in this section is to get an analog of Proposition \ref{Prop:ffr_p_center} for 
the universal enveloping algebras. 
Let $V=V_\kappa(\lf)$ or $\CDO(N)$ and $\frz_{Fr}(V)$ denote the $p$-center in the first case and the center in the 2nd case. 
The inclusion $\iota:\frz_{Fr}(V)\hookrightarrow V$ is a homomorphism of filtered vertex algebras leading to a homomorphism of filtered complete algebras $\widehat{U}(\iota):\widehat{U}(\frz_{Fr}(V))\hookrightarrow \widehat{U}(V)$.

\begin{Lem}\label{Lem:p_centers_universal_enveloping}
$V=V_\kappa(\lf)$ or $\CDO(N)$,
the homomorphism $\widehat{U}(\iota)$ is an embedding whose image is the $p$-center of $\widehat{U}(V)$.
\end{Lem}
\begin{proof}
We will consider the case when $V=V_\kappa(\lf)$, the other case is similar. 
It is enough to show that 
\begin{itemize}
\item[(*)]
for all $x\in \lf$ and $i\in \Z$, the topological generator $xt^i$ of $\widehat{U}(\mathfrak{z}(V))$ goes to $(xt^i)^p-x^{[p]}t^{ip}$. 
\end{itemize}
Recall (Section \ref{SSS:Utilde}) that for a vertex algebra $V'$ we defined the map $a\mapsto Y[a,z]:V'\rightarrow \widehat{U}(V')[[z^{\pm 1}]]$. By (\ref{eq:U_homom}), $\widehat{U}(\iota)$ sends $Y[a,z]$ to $Y[\iota(a),z]$ for all $a\in \mathfrak{z}_{Fr}(V)$ including $a=xt^j$ with $j<0$. So (*) reduces to 
\begin{equation}\label{eq:Y_sq_brack_id}
Y[\left((xt^j)^p-x^{[p]}t^{pj}\right)|\varnothing\rangle]=
\sum_{i\in \Z}{-i-1\choose -j-1}[(xt^{i})^p-x^{[p]}t^{ip}]z^{(j-i)p}.
\end{equation}
The proof repeats that of \cite[Lemma 2.3]{AW}.
\end{proof}

As a special case of Lemma \ref{Lem:p_centers_universal_enveloping}, the image of $\widehat{U}(\frz_{Fr}(V_\kappa(\g)))$
in $\widehat{U}_\kappa(\hat{\g})$ is nothing else but the completion of the $p$-center in $U_\kappa(\hat{\g})$.

Note also that $\widetilde{U}(\frz_{Fr}(V_\kappa(\g)))$ is identified 
with $\Ring_{\AS(\kappa)}[\Loop (\g^*)^{(1)}]$. Let
$\Ring_{\AS(\kappa)}[\Loop (\g^*)^{(1)}]_f$ denote the filtered part (=the union of filtered pieces) of $\Ring_{\AS(\kappa)}[\Loop (\g^*)^{(1)}]$. This filtered algebra is identified with the $p$-center in $\widehat{U}_\kappa(\hat{\g})$.

Using this lemma and Proposition \ref{Prop:ffr_p_center}, we deduce the following corollary.

\begin{Cor}\label{Cor:ffr_p_center}
The restriction of the homomorphism $$\Ffr_{P,\kappa}:\widehat{U}_\kappa(\hat{\g})\rightarrow \widehat{D}(\Loop N)\widehat{\otimes} \widehat{U}_\kappa(\hat{\lf})$$ to 
the $p$-center coincides with the composition
of 
$$\Ffr^{0,(1)}_{P,\AS(\kappa)}:\Ring_{\AS(\kappa)}[\Loop (\g^*)^{(1)}]_f\rightarrow \Ring[\Loop (T^*N)^{(1)}]_f\widehat{\otimes} \Ring_{\AS(\kappa)}[\Loop (\lf^*)^{(1)}]_f$$
and the embedding of the target to $\widehat{D}(\Loop N)\widehat{\otimes} \widehat{U}_\kappa(\hat{\lf})$ (the latter version of the tensor product was introduced in Section \ref{SSS_filt_compl_tensor} and the meaning of the target is similar).
\end{Cor}

\section{Affine Harish-Chandra homomorphism}\label{S_affine_HC}
\subsection{Construction}
The goal of this section is to prove the following 
\begin{Prop}\label{Prop:affine_HC}
The homomorphism 
$$\Ffr_P:\widehat{U}_\kappa(\hat{\g})\rightarrow 
\widehat{D}(\Loop N)\widehat{\otimes}\widehat{U}_\kappa(\hat{\lf})$$
restricts to $\widehat{U}_\kappa(\hat{\g})^{\Loop G}\rightarrow 
\widehat{U}_\kappa(\hat{\lf})^{\Loop L}$. 
\end{Prop}

For $L=H$, this map should be thought of as an affine analog of the Harish-Chandra homomorphism, hence the name of the section. The restriction of $\Ffr_P$ to $\widehat{U}_\kappa(\hat{\g})^{\Loop G}\rightarrow 
\widehat{U}_\kappa(\hat{\lf})^{\Loop L}$ will be denoted by $\HC_L$.

\begin{proof}
Recall, Section \ref{SSS_ffr_univ_enveloping}, that $\Ffr_P$ is $\Loop P$-equivariant, and so restricts to a map between $\Loop P$-invariants. Since 
$\widehat{U}_\kappa(\hat{\g})^{\Loop G}\hookrightarrow 
\widehat{U}_\kappa(\hat{\g})^{\Loop P}$, the claim of the proposition follows from the next lemma. 
\end{proof}

\begin{Lem}\label{Lem:invariant_containment}
The algebra of invariants
$[\widehat{D}(\Loop N)\widehat{\otimes}\widehat{U}_\kappa(\hat{\lf})]^{\Loop P}$ coincides with $\widehat{U}_\kappa(\hat{\lf})^{\Loop L}$ (embedded into the 2nd factor).
\end{Lem}

 In the proof and below we will use the following notation. Let $X$ be an affine variety equipped with a contracting action of $\mathbb{G}_m$. We will write $\Ring[\Loop X]_f$ for the ``finite degree part'' of $\Ring[\Loop X]$, in other words, $$\Ring[\Loop X]_f=\widehat{U}(\Ring[\Jet X])\subset \widetilde{U}(\Ring[\Jet X])=\Ring[\Loop X].$$

\begin{proof}
The proof is in several steps. 

{\it Step 1}. Here we prove that $\widehat{D}(\Loop N)^{\Loop N}$ coincides with the image of $\widehat{U}(\Loop \mathfrak{n})$ under the embedding  $\widehat{U}(\iota_r)$ (this image consists of invariants by the construction of $\iota_r$ in Section \ref{SSS_iota_r}). It is enough to show this claim 
after passing to the associated graded algebra. There 
$\Ring[\Loop(T^*N)]_f=\Ring[\Loop N]\widehat{\otimes}\Ring[\Loop \mathfrak{n}^*]_f$ with $\Loop N$ acting on the first factor by left translations. Take an element of degree $d$ for the naive grading
(by degree in $\Loop \mathfrak{n}^*$). It can be uniquely written as a converging sum $\sum_b f_b\otimes b$, where $b$ runs over the monomial topological basis in the degree $d$ component of  $\Ring[\Loop \mathfrak{n}^*]$ and $f_b\in \Ring[\Loop N]$. 
Let $\alpha$ denote the co-action map for the action of $\Loop N$ on $\Ring[\Loop N]$. Then the co-action map for the action of $\Loop N$
on $\Ring[\Loop(T^*N)]_f$ sends $\sum_b f_b\otimes b$ to 
$\sum_b \alpha(f_b)\otimes b$. This element equals
$\sum_b (1\otimes f_b)\otimes b$ if and only if $\alpha(f_b)=1\otimes f_b$. This is only possible if $f_b$ is in $\Ring$, yielding our claim. 

{\it Step 2}. Observe that there is a one-parameter subgroup
$\mathbb{G}_m\rightarrow Z(L)$ that acts on $\mathfrak{n}$
(and hence $\Loop \mathfrak{n}$) with positive weights. The invariants in $\widehat{U}(\Loop \mathfrak{n})$ for this action therefore coincide with $\Ring$. We conclude that  $\widehat{D}(\Loop N)^{\Loop N\rtimes \mathbb{G}_m}=\Ring$.

{\it Step 3}.
Note that the subgroup $\Loop N\rtimes \mathbb{G}_m$ acts trivially on 
$\widehat{U}_\kappa(\hat{\lf})$. By the previous step,  
$$[\widehat{D}(\Loop N)\widehat{\otimes}\widehat{U}_\kappa(\hat{\lf})]^{\Loop P}\subset \widehat{U}_\kappa(\hat{\lf}),$$
which easily implies the claim of the lemma.
\end{proof}

\sss We also note that, for the same reason, the homomorphism 
$$\Ffr_P^0: \widehat{S}_\kappa(\hat{\g})\rightarrow 
\Ring[\Loop(T^*N)]_f\widehat{\otimes}\widehat{S}_\kappa(\hat{\lf})$$ restricts to 
$$\HC^0_L: \widehat{S}_\kappa(\hat{\g})^{\Loop G}\rightarrow \widehat{S}_\kappa(\hat{\lf})^{\Loop L}.$$

\sss\label{SSS_HC_vertex1} 
We have similar (and easier) constructions on the level of vertex algebras, for example, we have $\HC_{L,\kappa}: V_\kappa(\g)^{\Jet G}\rightarrow V_\kappa(\lf)^{\Jet L}$. 

\subsection{Properties}
Here we study the properties of the maps 
$$\HC_L:\widehat{U}_\kappa(\hat{\g})^{\Loop G}\rightarrow 
\widehat{U}_\kappa(\hat{\lf})^{\Loop L}, \HC_L^0: \widehat{S}_\kappa(\hat{\g})^{\Loop G}\rightarrow \widehat{S}_\kappa(\hat{\lf})^{\Loop L}.$$ 

\sss Let $\HC^L_H$ denote the Harish-Chandra homomorphism
$\widehat{U}_\kappa(\hat{\lf})^{\Loop L}\rightarrow 
\widehat{U}_\kappa(\hat{\mathfrak{h}})^{\Loop H}$.

It follows from (\ref{eq:ffr_transitivity_algebra}) that 
\begin{equation}\label{eq:HC_transitivity}
\HC_{H,\kappa}=\HC_{H,\kappa}^L\circ \HC_{L,\kappa}.
\end{equation}
Similarly we have 
\begin{equation}\label{eq:HC_transitivity_cl}
\HC^0_{H,\kappa}=\HC_{H,\kappa}^{L,0}\circ \HC_{L,\kappa}^0.
\end{equation}

\sss

Below in this section we describe the map $\HC^0_{L,0}: \Ring[\Loop \g^*]_f^{\Loop G}\rightarrow \Ring[\Loop \lf^*]^{\Loop L}_f$. The embedding $\lf^*\hookrightarrow \g^*$
gives rise to an embedding $\Loop \lf^*\hookrightarrow \Loop \g^*$ of ind-affine ind-schemes. 

\begin{Lem}\label{Lem:HC0_interpretation}
The homomorphism $\HC^0_{L,0}$ coincides with the restriction to $\Loop \lf^*$.
\end{Lem}

\begin{proof}
$\HC^0_{L,0}$ coincides with the restriction to $\Ring[\Loop \g^*]_f^{\Loop G}$ of the composition  $$ \Ring[\Loop \g^*]_f\xrightarrow{\Ffr_P^0} \Ring[\Loop(T^*N)]_f\widehat{\otimes} \Ring[\Loop \lf^*]_f\xrightarrow{\pi\otimes \operatorname{id}}\Ring[\Loop \lf^*]_f,$$
where $\pi$ denotes the augmentation homomorphism $\Ring[\Loop(T^*N)]_f\rightarrow \Ring$. This composition is obtained by applying the functor $\widehat{U}$ to the analogous vertex algebra homomorphism $\Ring[\Jet \g^*]^{\Jet G}\rightarrow \Ring[\Jet \lf^*]$. It follows from Section \ref{SSS_ffrp_classical}
that the latter is the restriction for the embedding $\Jet \lf^*\hookrightarrow \Jet \g^*$. This implies the claim of the lemma. 
\end{proof}

\sss Our goal here is to prove the following lemma. 
\begin{Lem}\label{Lem:HC_injectivity}
The homomorphisms $\HC_{L,\kappa}$ and $\HC_{L,\kappa}^0$
are injective. 
\end{Lem}
\begin{proof}
Note that $\gr \widehat{U}_\kappa(\hat{\g})^{\Loop G}\hookrightarrow \Ring[\Loop \g^*]^{\Loop G}$ and similarly for $\lf$. Then $\gr \HC_{L,\kappa}$ is the restriction of $\HC^0_{L,0}$. And thanks to (\ref{eq:HC_transitivity_cl}), it is enough to show that $\HC^0_{H,0}$ is injective. 

By Lemma \ref{Lem:HC0_interpretation}, $\HC^0_H$ is the restriction to $\Loop \h^*$. A standard argument, cf. the proof of \cite[Proposition 4.3.4]{Frenkel_loop} reduces the proof of the injectivity of this map to showing the injectivity of the restriction map $\Ring[\Jet \g^*]^{\Jet G}\rightarrow \Ring[\Jet \mathfrak{h}^*]$. The injectivity of the latter is equivalent to the claim that the maps $\Ring[\Jet_n \g^*]^{\Jet_n G}\rightarrow \Ring[\Jet_n \mathfrak{h}^*]$ are injective for all $n$. This will follow once we show that the morphism $\varphi_n:\Jet_nG\times^{\Jet_n H}\Jet_n\mathfrak{h}^*\rightarrow \Jet_n\g^*$ induced by the action is dominant. Since both schemes are flat over $\Ring$, it is enough to assume that $\Ring$ is an algebraically closed field. In this case the claim that $\varphi_n$ is dominant follows by computation of the differential at a point $(1,x)$ with $x\in (\mathfrak{h}^*)^{reg}$: since the roots are nonzero on $\mathfrak{h}$, it is easy to see that $d_{1,x}\varphi_n$ is an isomorphism.  
\end{proof}

\sss\label{SSS_vertex_algebra_injectivity} For the same reason, $\HC_{L,\kappa}: V_\kappa(\g)^{\Jet G}\rightarrow V_\kappa(\lf)^{\Jet L}$ is injective. 


\subsection{The Harish-Chandra center in abelian case}
Here we describe the $\Ring$-algebra $\widehat{U}_\kappa(\hh)^{\Loop H}$, where $\Ring$ is an algebra over $\mathbb{F}_p$. 

\sss 
First, assume that $\Ring$ is a general commutative ring. Choose an energy graded basis $b_\alpha, \alpha\in A,$ in $\h[t^{\pm 1}]$ and
order the indexes so that, for $\alpha<\beta$, the energy of $b_\alpha$ does not exceed the energy of $b_\beta$. Then we have 
a filitered $\Ring[\mathbf{1}]$-module isomorphism $S(\hh)/(\mathbf{1}'-1)\xrightarrow{\sim} U(\hh)/(\mathbf{1}'-1)$ that sends an ordered monomial in $b_\alpha$'s to the same ordered monomial. 

\begin{Lem}\label{Lem:abel_iso_equiv}
This isomorphism extends by continuity to $\widehat{S}_\kappa(\hh)\xrightarrow{\sim} \widehat{U}_\kappa(\hh)$ and this extension is $\Loop H$-equivariant.   
\end{Lem}
\begin{proof}
The claim about the extension is manifest. To prove that the resulting isomorphism is $\Loop T$-equivariant, we argue as follows. Note that for an $\Ring'$-point $g$ of $\Loop T$ (where $\Ring'$ is an $\Ring$-algebra), we have $g.b_\alpha=b_\alpha+f_\alpha(g)\mathbf{1}$
for some $f_\alpha(g)\in \Ring$. Plugging the elements $b_\alpha+f_\alpha(g)\mathbf{1}$ into an ordered monomial instead of $b_\alpha$ and distributing, we get the sum of still ordered monomials. It follows that the isomorphism 
$\widehat{S}_\kappa(\hh)\xrightarrow{\sim} \widehat{U}_\kappa(\hh)$ is $\Loop H$-equivariant.
\end{proof}

\sss Now we are ready to describe $\widehat{U}_\kappa(\hh)^{\Loop H}$.
\begin{Prop}\label{Prop:HC_center_abelian}
Assume that $\Ring$ is an $\F_p$-algebra and $\kappa$ is invertible in $\Ring$. 
Then the elements of the form 
$(xt^j)^p-\kappa^{p-1}(xt^{jp})$ with $x\in \h_{\F_p}$ and $j\in \Z$
lie in $\widehat{U}_\kappa(\hh)^{\Loop H}$. Moreover, any element in  $\widehat{U}_\kappa(\hat{\h})^{\Loop T}_{\leqslant j}$
(the PBW filtration term) is uniquely written as the converging sum of degree $\leqslant j/p$ polynomial in the elements $(xt^j)^p-\kappa^{p-1}(xt^{jp})$.
\end{Prop}


\begin{proof}
The claim easily reduces to the case when $\dim H=1$, which is what we are going to assume. 
The proof is in several steps. 

{\it Step 1}. We use Lemma \ref{Lem:abel_iso_equiv} 
to see that the claim of the proposition is equivalent to its direct analog for $\widehat{S}_\kappa(\hh)$. Indeed, note that the image in  $\widehat{U}_\kappa(\hat{\h})^{\Loop H}_{\leqslant j}$ of an ordered monomial in elements  $(xt^j)^p-\kappa^{p-1}(xt^{jp})\in \widehat{S}_\kappa(\hh)$ is also the same ordered monomial but now 
$(xt^j)^p-\kappa^{p-1}(xt^{jp})$ are viewed as elements of $\widehat{U}_\kappa(\hh)$.

{\it Step 2}. In this and the subsequent steps we prove the complete analog of the claim of proposition for $\widehat{S}_\kappa(\hh)$, where we can assume that $\kappa=1$.
We write $b_n$ for $t^n$ viewed as an element of $\hh$. The action of $\Loop H$ on $b_n$ is given by $g.b_n=b_n +\chi_n(g)$, where $\chi_n(g)=\operatorname{Res}_{t=0}(t^n g^{-1}dg)$ thanks to (\ref{eq:loop_adjoint_action}). One can then check that $b_n^p-b_{np}$ is $\Loop H$-invariant, for example, using Step 1 and Corollary \ref{p:p-cent}.

{\it Step 3}. We claim that the elements $\chi_n$ (viewed as elements of $\Ring[\Loop H]$) are topologically linearly independent: the equality $\sum_{n\geqslant n_0}\alpha_n\chi_n=0$ for some $n_0$ and $\alpha_n\in \Ring$ implies $\alpha_n=0$.
Assume the contrary: there are elements $\alpha_n$ with 
$\sum_{n\geqslant n_0}\alpha_n\chi_n=0$. Let $f(t)=\sum_{n\geqslant N}\alpha_n t^n$.  Note that we have $\operatorname{Res}_{t=0}(f(t)g(t)^{-1}dg(t))=0$ for all $g(t)\in \Ring'(\!(t)\!)^\times$ for all possible $\Ring$-algebras $\Ring'$. 

{\it Step 4}. Take $\Ring':=\Ring[x]/(x^N)$ for a  sufficiently large integer $N$. 
Consider $g(t)=t-x$. Then $dg(t)=1$, while 
$g(t)^{-1}=t^{-1}(1-t^{-1}x)^{-1}=t^{-1}\sum_{i=0}^{N-1}t^{-i}x^i$. From $\operatorname{Res}_{t=0}(f(t)g(t)^{-1}dg(t))=0$
we deduce that $\sum_{i=0}^{N-1}\alpha_n x^n=0$ in $\Ring'$.
This shows that $\alpha_n=0$ for $0\leqslant n<N$. Similarly, considering $g(t)=t^{-1}-x$, we see that
$\alpha_n=0$ for $-N<n\leqslant 0$. We conclude that all 
$\alpha_n=0$. This finishes the proof of the claim in
Step 3 that the characters $\chi_n$ are topologically linearly independent. 



{\it Step 5}. Let $V$ be the (topological) span of monomials in the elements $b_i$, where each element has power strictly less than $p$. Note that $V$ is $\Loop H$-stable. We claim $V^{\Loop H}$ consists of scalars. Namely, let $f$ be an element in $V$ and let $f'$ be its top degree homogeneous component, of degree, say $d$. Then the top degree homogeneous component in 
$gf-f$ is 
\begin{equation}\label{eq:action_difference}
\sum_i \chi_i(g)\partial_i  f',
\end{equation}
where we write $\partial_i$ for the derivative with respect to $b_i$. Note that since $f$ is converging, for each monomial $F$ in the elements $b_i$ of degree $d-1$, there is $n_0(F)\in \Z$ such that $F$ does not appear in $\partial_i f'$ for $i<n_0(F)$.

Since the characters $\chi_i$ are topologically linearly independent, we see that $\partial_i f'=0$ for all $i$. Since
$f'\in V$, this implies $f'=0$.

{\it Step 6}. We note that, for each $k>0$, $\widehat{S}_1(\hh)_{\leqslant kp-1}$ is the completed tensor product of $V$ and the completion of $\F[b_i^p-b_{ip}]_{\leqslant k-1}$. The action of $\Loop H$ respects this decomposition and is trivial on the completion of 
$\F[b_i^p-b_{ip}]_{\leqslant k-1}$. It follows from Step 5
that the subspace of $\Loop H$-invariants in $\widehat{S}_1(\hh)_{\leqslant kp-1}$ equals to the completion of $\F[b_i^p-b_{ip}]_{\leqslant k-1}$. This finishes the proof.
\end{proof}

\sss\label{SSS_vertex_abelian} An analogous (and easier) argument shows that, for an invertible element $\kappa\in \Ring$, the algebra $V_\kappa(\h)^{\Jet H}$ is the algebra of polynomials in $(xt^j)^p-\kappa^{p-1}(x t^{jp})$ with $j<0$
and $x$ running through a basis in $\h_{\F_p}$. One just needs to modify Step 4 considering $g(t)=1-tx$.

\subsection{Case of integral $\kappa$}
In this section $\Ring$ is a perfect field of positive characteristic $p$. We assume that $p$ is bigger than the Coxeter number for any almost simple normal subgroup of $G$.

Recall the inclusion of the p-center $ \widehat{S}_{\AS(\kappa)}(\hg^{(1)})\hookrightarrow \widehat{U}_\kappa(\hat{\g})$ to be denoted by $\iota$. The map $\iota$ is $\Loop G$-equivariant, see Corollary \ref{p:p-cent}. 

Here is the main result. 

\begin{Thm}\label{Thm:integral_case}
The following claims are true:
\begin{enumerate}
\item $\gr \widehat{U}_\kappa(\hg)^{\Loop G}\subset  \widehat{S}\left((\Loop \g)^{(1)}\right)^{(\Loop G)^{(1)}}$, the containment of subalgebras of $\widehat{S}(\Loop \g)^{\Loop G}$.
\item If $\kappa\in \F_p\setminus \{0\}$, then
$$\widehat{U}_\kappa(\hg)^{\Loop G}=\iota\left(\widehat{S}\left((\Loop \g)^{(1)}\right)^{(\Loop G)^{(1)}}\right).$$
\end{enumerate}
\end{Thm}

The proof will be given below in this section.

\sss First, we need to describe the structures of $\Ring[\Jet \g^*]^{\Jet G}$
and $\Ring[\Loop \g^*]^{\Loop G}$.

Let $F_1,\ldots,F_r$ be free homogeneous generators of 
$\Ring[\g^*]^G$. We then can consider the elements $F_{i,n}\in \Ring[\Jet(\g^*)]^{\Jet G}$ for $i=1,\ldots,r$ and $n\leqslant -d_i$, where $d_i$ is the degree of $F_i$.

\begin{Lem}\label{Lem:Jet_Loop_invariants}
Recall that we assume that $p$ is larger than the Coxeter number of any of the almost simple normal subgroups of $G$.
The following claims hold:
\begin{enumerate}
\item The elements $F_{i,n}$ are free generators of $\Ring[\Jet(\g^*)]^{\Jet G}$.
Equivalently, $\Ring[\Jet\g^*]^{\Jet G}\xrightarrow{\sim}\Ring[\Jet(\g^*/\!/G)]$.
\item We have $\Ring[\Loop(\g^*)]^{\Loop G}=\varprojlim_{n\rightarrow +\infty}\Ring[(t^{-n}\g^*[[t]])]^{\Jet G}$
\end{enumerate}
\end{Lem}
\begin{proof}
(1) holds by By \cite[Theorem 4.4]{ATV}. (2) holds by the same argument as in the proof of 
\cite[Proposition 4.3.4]{Frenkel_loop}.
\end{proof}

\sss Recall the injective (by Lemma \ref{Lem:HC_injectivity}) homomorphism $$\HC_{H,0}^0: \widehat{S}(\Loop\g)^{\Loop G}\hookrightarrow 
\widehat{S}(\Loop\h),$$
it is given by restriction from $\Loop \g^*$ to $\Loop \mathfrak{h}^*$, see Lemma \ref{Lem:HC0_interpretation}.
Consider its Frobenius twist 
$$\HC_{H,0}^{0,(1)}:\widehat{S}([\Loop\g]^{(1)})^{[\Loop G]^{(1)}}\hookrightarrow 
\widehat{S}([\Loop\h]^{(1)})$$
Below we embed $\Ring[(\Loop \h^*)^{(1)}]_f$ into
$\Ring[(\Loop \h^*)]_f$ via $\operatorname{Fr}^*$
(the meaning of the subscript ``f'' was explained after the statemenet of Lemma \ref{Lem:invariant_containment}). We will need the following lemma. 

\begin{Lem}\label{Prop:inv_frob}
We have 
\begin{equation}\label{eq:intersection_coincide}
\operatorname{im}(\HC_{H,0}^{0,(1)})=\operatorname{im}\HC_{H,0}^0\cap \Ring[(\Loop \h^*)^{(1)}]_f.
\end{equation}
\end{Lem}
\begin{proof}
We have the inclusion of the left hand side into the right hand side, we need to show that it is the equality. The proof is in several steps.

{\it Step 1}.
Thanks to Lemma \ref{Lem:Jet_Loop_invariants}, we can deduce (\ref{eq:intersection_coincide}) from its jet analog, and the statements about jets follows from that for $n$th jets 
(for all $n$). In other words, let $\mathsf{res}_n$ denote the embedding $\Ring[\Jet_n(g^*)]^{\Jet_n G}\rightarrow \Ring[\Jet_n(\h^*)]$ (we have seen that this map is an embedding in the proof of Lemma \ref{Lem:HC_injectivity}). We need to show that 
\begin{equation}\label{eq:intersection_coincide_n}
\operatorname{im}(\mathsf{res}_n^{(1)})=\operatorname{im}(\mathsf{res}_n)\cap \Ring[(\Jet_n \h^*)^{(1)}].
\end{equation}

{\it Step 2}.
We note that since $\mathsf{res}_n^{(1)}$ is injective as well, the left hand side is a normal domain. The right hand side is integral over the left hand side because 
$\Ring[\Jet_n\g^*]^{\Jet_n G}$ is integral over $\Ring[(\Jet_n\g^*)^{(1)}]^{(\Jet_n G)^{(1)}}$
(the latter contains the $p$th powers of the elements of the former).
It remains to show that the right hand side has the same fraction field as the left hand side. Consider the element 
$\delta\in \Ring[\g^*]^G\subset \Ring[\Jet_n (\g^*)]^{\Jet G_n}$, the product of squares of all roots, and localize both sides of (\ref{eq:intersection_coincide_n}) at
$\mathsf{res}_n(\delta^p)$. Note that the action map 
$G\times^{N_G(H)}\h^{*,reg}\rightarrow \g^{*,reg}$
is an isomorphism, hence the same is true for 
$\Jet_nG\times^{N_G(H)\Jet_nH}\Jet_n\h^{*,reg}\rightarrow \Jet_n\g^{*,reg}$. So the localization of the left hand side of (\ref{eq:intersection_coincide_n}) is 
$\left(\Ring[(\Jet_n\h^*)^{(1)}][\mathsf{res}_n(\delta^{p})^{-1}]\right)^W$, while the localization of the right hand side is $\Ring[\Jet_n \h^*][\mathsf{res}_n(\delta^{p})^{-1}]^W\cap 
\Ring[(\Jet_n \h^*)^{(1)}][\mathsf{res}_n(\delta^{p})^{-1}]$.
It is easy to see that these two subalgebras are the same. 
\end{proof}

\sss Now we are ready to prove the theorem.
\begin{proof}[Proof of Theorem \ref{Thm:integral_case}]
Let $\mathscr{Z}$ denote $\widehat{U}_\kappa(\hg)^{\Loop G}$ and $\mathscr{Z}_0:=\iota(\widehat{S}\left((\Loop \g)^{(1)}\right)^{(\Loop G)^{(1)}})$. When $\kappa\in \F_p\setminus \{0\}$,  we have $\mathscr{Z}_0\subset \mathscr{Z}$. Both claims of the theorem will then follow if we show
\begin{equation}\label{eq:desired_inclusion_Thm_integral}
\gr \mathscr{Z}\subset \widehat{S}((\Loop \g)^{(1)})^{(\Loop G)^{(1)}}.
\end{equation}

Note that $\gr\mathscr{Z}\subset \widehat{S}(\Loop\g)^{\Loop G}$.
Thanks to Lemma \ref{Lem:HC_injectivity}, $\HC^0_{H,0}$ is injective. It follows that 
\begin{equation}\label{eq:gr_equality1}
\gr (\im \HC_{H,\kappa})
\subset \operatorname{im}\HC_{H,0}^0. \end{equation}
From Proposition \ref{Prop:HC_center_abelian} it follows that $\gr \widehat{U}_{\kappa}(\hh)^{\Loop H}=\widehat{S}((\Loop \h)^{(1)})$. So
\begin{equation}\label{eq:gr_containment}
\gr (\im \HC_{H,\kappa})\subset \widehat{S}((\Loop \h)^{(1)}).
\end{equation}

Combining 
 (\ref{eq:gr_equality1}) and  (\ref{eq:gr_containment}) we see that that $\HC_{H,0}^0(\gr \mathscr{Z})$ is contained in $$\HC_{H,0}^0\left(\widehat{S}(\Loop\g)^{\Loop G}\right)\cap \widehat{S}\left((\Loop\h)^{(1)}\right).$$ By Lemma \ref{Prop:inv_frob},
the latter coincides with $\HC_{H,0}^0\left(\widehat{S}((\Loop\g)^{(1)})^{(\Loop G)^{(1)}}\right)$.  Using that $\HC_{0,H}^0$ is injective again, we arrive at (\ref{eq:desired_inclusion_Thm_integral}).  
\end{proof}

\sss\label{SSS_center_vertex_algebra_reduction}
A complete analog of Theorem \ref{Thm:integral_case} holds for vertex algebras, in particular $$\gr V_\kappa(\g)^{\Jet G}\subset \Ring[(\Jet \g^*)^{(1)}]^{(\Jet G)^{(1)}}.$$ The proof of this claim is the same using an analog of Proposition \ref{Prop:HC_center_abelian} explained in Section \ref{SSS_vertex_abelian}. We also have the following result that will allow us to deduce a description of $\widehat{U}_\kappa(\hat{\g})^{\Loop G}$ from that of $V_\kappa(\g)^{\Jet G}$ in the case when $\kappa\not\in \F_p$.

\begin{Lem}\label{Lem:center_assoc_grad} Suppose $\gr [V_\kappa(\g)^{\Jet G}]= \Ring[(\Jet \g^*)^{(1)}]^{(\Jet G)^{(1)}}$. Consider the filtered algebra homomorphism
$\psi:\widehat{U}(V_\kappa(\g)^{\Jet G})\rightarrow  \widehat{U}_\kappa(\hg)^{\Loop G}$, cf. Lemma \ref{Lem:action_on_modes}. Then 
\begin{enumerate} 
\item this homomorphism is an isomorphism, 
\item and $\gr \widehat{U}_\kappa(\hg)^{\Loop G}=  \widehat{S}((\Loop \g)^{(1)})^{(\Loop G)^{(1)}}$.  
\end{enumerate}
\end{Lem}
\begin{proof}
Thanks to  Lemma \ref{Lem:Jet_Loop_invariants}, we have $\Ring[\Jet \g^*]^{\Jet G}\cong \Ring[\Jet (\g^*/\!/G)]$ and $\widehat{S}(\Loop \g)^{\Loop G}=\Ring[\Loop(\g^*/\!/G)]_f$. Thanks to 
$\gr (V_\kappa(\g)^{\Jet G})= \Ring[(\Jet (\g^{*,(1)}/\!/G^{(1)})]$, we have that the natural epimorphism
$$\varpi:\Ring[(\Loop (\g^*/\!/G))^{(1)}]_f\twoheadrightarrow \gr\widehat{U}\left(V_\kappa(\g)^{\Jet G}\right)$$
from Section \ref{SSS:filtered_enveloping} is an isomorphism.
The composition $\iota\circ [\gr\psi]\circ \varpi$, where $\iota$ stands for the inclusion   $\gr\widehat{U}_\kappa(\hg)^{\Loop G}\hookrightarrow \gr \widehat{U}_\kappa(\hg)$, is the inclusion 
$\Ring[(\Loop (\g^*/\!/G))^{(1)}]_f\hookrightarrow  \widehat{S}(\Loop \g)$, its image is  $\widehat{S}\left((\Loop \g)^{(1)}\right)^{(\Loop G)^{(1)}}$.

In particular, $\gr\psi$ is injective. And from (1) of Theorem \ref{Thm:integral_case} combined with the coincidence of $\iota\circ [\gr\psi]\circ \varpi$ with the composition $\Ring[(\Loop (\g^*/\!/G))^{(1)}]_f\xrightarrow{\sim} \widehat{S}\left((\Loop \g)^{(1)}\right)^{(\Loop G)^{(1)}}\hookrightarrow  \widehat{S}(\Loop \g)$ we conclude that $\gr\psi$ is an isomorphism.
\end{proof}

\section{The case of affine ${\slf}_2$}\label{S_sl2}
Here we suppose that $2$ is invertible in a commutative ring $\Ring$. Except in Section \ref{SS_pcenter_Poisson}, we assume that $G=\operatorname{SL}_2$.

\subsection{Sugawara elements}
Let $e,h,f$ be the standard basis elements of $\g$. 
Let $(\cdot,\cdot)$ be the trace form on $\g=\slf_2$.
This is the form we use in the construction of $\hat{\g}$. 

Pick a basis $x_1,x_2,x_3$ of (the $\Z[\frac{1}{2}]$-form of) $\g$ and the dual basis $x^1,x^2,x^3$ with respect to $(\cdot,\cdot)$.
For example, we can take $x_1=f, x_2=h, x_3=e$, then $x^1=f, x^2=h/2, x^3=e$. For $y\in \g$, we write $y_n$ for $yt^n\in \hg$.

\sss
Consider the Sugawara elements
\begin{equation}\label{eq:Sugawara_def}
L_n:=\begin{cases}&\frac{1}{2}\sum_{i=1}^3 \sum_{j\in \Z}x_{i,n+j}x^i_{-j}, \text{ for }n\neq 0, \\
&\frac{1}{2}\left(\sum_{i=1}^3 x_ix^i+2\sum_{i=1}^3\sum_{j>0}x_{i,-j}x^i_{j}\right),\text{ for }n=0.\end{cases}
\end{equation}

\sss
We will need the following well-known property of the elements $L_i$:
\begin{equation}\label{eq:Sugawara_derivation}
[L_i,x_n]=-\kappa nx_{i+n}, \forall x\in \g, i,n\in \Z. 
\end{equation}

Here is a consequence of (\ref{eq:Sugawara_derivation}). Let $\Ring'$ be an $\Ring$-algebra. Note that for an $\Ring'$-point $g$ of $\Loop G$ and an $\Ring'$-linear derivation $\partial$ of $\Ring'(\!(t)\!)$, we can view $g^{-1}\partial g$ as an $\Ring'$-point of $\Loop\g$. For example, for $g=\begin{pmatrix}1&0\\t&1\end{pmatrix}, $ and $\partial=t^n \frac{d}{dt}$, we have 
$g^{-1}\partial g=f_{n}$. Then

\begin{equation}\label{eq:Sugawara_adjoint}
\operatorname{Ad}(g^{-1})L_i=L_i-\kappa g^{-1}(t^{i+1}\frac{d}{dt})g. 
\end{equation}

\sss\label{SSS:Sugawara_classical}
We can consider the analogs of Sugawara elements $\underline{L}_i\in \widehat{S}_{\AS(\kappa)}(\hg)$.
They are given by 
\begin{equation}\label{eq:Sugawara_classical}
\underline{L}_n:=\frac{1}{2}\sum_{i=1}^3\sum_{j\in \Z} x_{i,n+j}x^i_{-j}.
\end{equation}
Further, we have the direct analog of (\ref{eq:Sugawara_derivation}), where we replace $L_i$ with $\underline{L}_i$.


\subsection{Construction of central elements}
Let $\Ring$ be an $\F_p$-algebra. We note that (\ref{eq:Sugawara_derivation}) implies the following claim. 

\begin{Lem}\label{Lem:central_1}
The elements 
$Y_i:=L_i^p$ for $i$ not divisible by $p$ and $Y_i:=L_i^p-\kappa^{p-1}L_{ip}$ for $i$ divisible by
$p$ are central in $\widehat{U}_{\kappa}(\hg)$.
\end{Lem}

Now we produce a sequence of $\Loop G$-invariant elements. Recall that $\iota$ denotes the inclusion $\widehat{S}_{\mathsf{AS}(\kappa)}(\hat{\g}^{(1)})\hookrightarrow \widehat{U}_\kappa(\hat{\g})$ of the $p$-center.
%


\begin{Prop}\label{Prop:HC_center_sl2_noninteger}
For all $n$,  $X_n:=\kappa^{p-1}\iota(\underline{L}_n)-(\kappa^{p-1}-1)Y_n\in\widehat{U}_\kappa(\hg)^{\Loop G}$.
\end{Prop}
\begin{proof}
The proof is in several steps. We note that it is enough to prove the claim for $\Ring=\F_p[\kappa]$, where $\kappa$ is an indeteminate. The claim for this ring will follow from the claim for its field of fractions, and then from the claim for the algebraic closure. So we can and will assume that $\kappa$ is invertible in an algebraically closed field $\Ring$. 

{\it Step 1}. Let us outline the strategy of the proof. 
Note that by the construction both $Y_n$ and $\underline{L}_n$
are $G$-invariant. In the next two steps we will see that
$\kappa^{p-1}\iota(\underline{L}_n)-(\kappa^{p-1}-1)Y_n$ is invariant with respect to the element  $g:=\begin{pmatrix}1&0\\t&1\end{pmatrix}$. Then we will see that this implies that it is $\Jet G$-invariant, and then that it is $\Loop G$-invariant. 

{\it Step 2}. Here we compute $\operatorname{Ad}(g^{-1})Y_n$ for 
$g=\begin{pmatrix}1&0\\t&1\end{pmatrix}$. By (\ref{eq:Sugawara_adjoint}) and the computation preceding this equation, we have $\operatorname{Ad}(g)^{-1}(L_n)=L_n-\kappa f_{n+1}$ for all $n$. 

Now we compute $\operatorname{Ad}(g^{-1})L_n^p=(\operatorname{Ad}g^{-1})L_n)^p=(L_n-\kappa f_{n+1})^p$. Recall that for elements $x,y$ of 
a Lie algebra over $\F_p$, the element $(x+y)^p-x^p-y^p$ in the universal 
enveloping algebra is actually a degree $p$ Lie polynomial,  $\sigma(x,y)=\sum_{i=1}^{p-1}\sigma_i(x,y)$, cf. Definition \ref{Defi:restricted}, where $\sigma_i$ has degree $i$ in $x$ and $p-i$ in $y$. 
If the elements 
\begin{itemize}
\item[(*)]
$\operatorname{ad}(x)^j y$ and $\operatorname{ad}(x)^k y$ commute 
for all $k,j$, 
\end{itemize}
then $\sigma_i=0$ for $i\neq p-1$. An easy check shows that this component coincides with
$\operatorname{ad}(x)^{p-1}y$. We see that 
\begin{align*}
(L_n-\kappa f_{n+1})^p&=L_n^p-\kappa\operatorname{ad}(L_n)^{p-1}f_{n+1}-\kappa^p f_{n+1}^p=\\&L_n^p-\kappa^{p}f_{n+1}^p-\kappa^{p}\left(\prod_{j=0}^{p-1}(jn+1)\right)f_{pn+1}.
\end{align*} We note that if $n$ is not divisible by $p$, then the last summand is equal to $0$, otherwise the coefficient of $f_{pn+1}$ is equal
to $-\kappa^p$. In both cases, we have 
\begin{equation}\label{eq:Y_conjugation}
\operatorname{Ad}(g)^{-1}Y_n=Y_n-\kappa^p f_{n+1}^p.    
\end{equation}

{\it Step 3}. By Section \ref{SSS:Sugawara_classical}, we have 
\begin{equation}\label{eq:Lunder_conjugation}
\operatorname{Ad}(g)^{-1}\underline{L}_n=\underline{L}_n-(\kappa^p-\kappa)f_{n+1}. 
\end{equation}
Note that $\iota(f_{n+1})=f_{n+1}^p$. Combining 
(\ref{eq:Y_conjugation}) and (\ref{eq:Lunder_conjugation}), we see that 
$\kappa^{p-1}(\iota(\underline{L}_n))-(\kappa^{p-1}-1)Y_n$ is $g$-invariant.

{\it Step 4}. We claim that any $G$- and $g$-invariant element  
$a\in\widehat{U}_\kappa(\hg)_{\leqslant i}$ (for some $i$) is $\Jet G$-invariant. Recall, Section \ref{SSS:filtered_enveloping}, that $V:=\widehat{U}_\kappa(\hg)_{\leqslant i}$ is a complete topological space whose base of neighborhoods of $0$ is formed by $U_N:=\widehat{U}_\kappa(\hg)_{\leqslant i-1}\g(\geqslant N), N>0$, where we write $\g(\geqslant N)$ for the span of $x_n$ with $n\geqslant N$. Each $U_N$ is $\Jet G$-stable. Moreover, $V/U_N$ is a rational representation of $\Jet G$, in particular, $a_N:=a+U_N$ lies in a subspace where the action of $\Jet G$ factors through $\Jet_{d(N)}G$ for some $d(N)>0$.
Note that $a$ is $\Jet G$-invariant if and only if $a_N$ is $\Jet_{d(N)}G$-invariant. It remains to notice that, as an abstract group, $\Jet_n G$ is generated by $G$ and $g$ for all $n$ (here we use  that $\Ring$ is an algebraically closed field). We conclude that $a$ is $\Jet G$-invariant.

{\it Step 5}. We now claim that any $\Jet G$-invariant element in $V$ is $\Loop G$-invariant. Consider the affine Grassmanian $\mathsf{Gr}$, the fpqc quotient  $\Loop G/\Jet G$. In Lemma \ref{Lem:aff_Gr_integral} we will see that the algebra of global functions on $\mathsf{Gr}$ is $\Ring$. Since $\mathsf{Gr}$ is the fpqc quotient, the subalgebra of $\Jet G$-invariants in $\Ring[\Loop G]$ is $\Ring$.

Now suppose that $v\in V$ is a $\Jet G$-invariant element. Let $\alpha:V\rightarrow V\widetilde{\otimes}\Ring[\Loop G]$ denote the coaction map. Let $\mathscr{B}$ be a topological basis of $V$ compatible with the energy grading. As any element 
of $V\widetilde{\otimes}\Ring[\Loop G]$, the element $\alpha(v)$ is uniquely written as a converging sum $\sum_{b\in \mathscr{B}}b\otimes f_b$ for $f_b\in \Ring[\Loop G]$. We see that all $f_b$ are $\Jet G$-invariant. So $f_b\in \Ring$ for all $b$, hence $v$ is $\Loop G$-invariant. 
\end{proof}

\begin{Lem}\label{Lem:aff_Gr_integral}
The algebra of global functions on $\mathsf{Gr}$ is $\Ring$.
\end{Lem}
\begin{proof}
The ind-scheme $\mathsf{Gr}$ is ind-projective: it is a countable union of finite type projective schemes: 
$\mathsf{Gr}=\bigcup_{i\geqslant 0}Z_i$. For $Z_i$ we can take the scheme of the following form. Take a vector space $V_i$ of dimension $2d_i$ equipped with a nilpotent linear operator. For $Z_i$ we take the subscheme of the Grassmannian of $d_i$-dimensional subspaces in $V_i$ that are stable with respect to the operator. 

We claim that $Z_i$ is connected. Indeed, $Z_i$ admits a surjective proper homomorphism from a suitable Springer fiber for $\operatorname{GL}(V_i)$. The nilpotent cone for $\operatorname{GL}(V_i)$ is normal (by the usual argument of Kostant applied to the slice consisting of companion matrices) and the Springer morphism is a resolution of singularities. So every Springer fiber is connected. Hence $Z_i$ is connected.

On the other hand, according to \cite[Theorem 6.1]{PR}, $\mathsf{Gr}$ is a reduced ind-scheme, i.e., is the countable union of reduced finite type schemes.  Now according to \cite[Lemma 6.3]{Beauville_Laszlo},   $\mathsf{Gr}=\bigcup_{i\geqslant 0}Z_{i,red}$, where the subscript ``red'' means the reduced scheme. We have $\Ring[\mathsf{Gr}]=\varprojlim_j \Ring[Z_{i,red}]$. Each $\Ring[Z_{i,red}]$ is $\Ring$ because $Z_{i,red}$ is proper and connected finishing the proof. 
\end{proof}

\subsection{Generators of the Harish-Chandra center}
Suppose that $\Ring$ is a perfect field of characteristic $p>2$.
The following is the main result of this section.
\begin{Thm}\label{Thm:sl2_HC_center}
The algebra $\widehat{U}_{\kappa}(\hg)^{\Loop G}$ is the filtered complete polynomial algebra in the elements $X_n$ from Proposition \ref{Prop:HC_center_sl2_noninteger}.
\end{Thm}
\begin{proof}
By (1) of Theorem \ref{Thm:integral_case}, we have a graded algebra embedding 
$\widehat{U}_{\kappa}(\hg)^{\Loop G}\hookrightarrow 
\widehat{S}\left((\Loop \g)^{(1)}\right)^{(\Loop G)^{(1)}}$. By Lemma \ref{Lem:Jet_Loop_invariants}, the algebra 
$\widehat{S}\left((\Loop \g)^{(1)}\right)^{(\Loop G)^{(1)}}$ is a graded complete polynomial algebra in homogeneous generators
\begin{equation}\label{eq:generators}
\frac{1}{2}\sum_{i=1}^3\sum_{j\in \Z}x_{i,n+j}^p(x^i_{-j})^p.
\end{equation}
Notice that (\ref{eq:generators}) is also the principal symbol of the element $X_n$. This is because
the elements $Y_n$ and $\iota(\underline{L}_n)$ are of PBW degree $2p$ and their principal symbols 
are (\ref{eq:generators}). The claim about the principal symbols of the $X_n$'s implies the claim of the theorem in the standard fashion. 
\end{proof}

Our next goal will be computing the images of the elements $X_n$ under the affine Harish-Chandra homomorphism $\operatorname{HC}_{\kappa, H}$. For this, we need a discussion of Poisson structures that is preceeded by a general result on the distribution algebras of algebraic groups and Frobenius homomorphisms.

\subsection{Distribution algebras and Frobenius homomorphisms}
In this section $\breve{F}$ is a smooth affine  group subscheme in some $\operatorname{GL}_n$ over $\mathbb{Z}_p$. Let $F$ denote its reduction mod $p$, a smooth group scheme over $\F_p$. To $\breve{F}$ (resp., $F$) we assign its distribution algebra $\operatorname{Dist}_1(\breve{F})$ (resp., $\operatorname{Dist}_1(F)$) over $\mathbb{Z}_p$ (resp., $\F_p$) so that $\operatorname{Dist}_1(F)=\operatorname{Dist}_1(\breve{F})/(p)$. The Frobenius homomorphism $F\rightarrow F^{(1)}$ gives rise to an algebra homomorphism $\operatorname{Fr}:\operatorname{Dist}_1(F)\rightarrow \operatorname{Dist}_1(F^{(1)})$. On the other hand, since $F$ is defined over $\F_p$, we have an isomorphism $F\xrightarrow{\sim}F^{(1)}$ that we fix throughout. 

Consider the Lie algebra $\breve{\ff}$ over $\Z_p$
and its quotient $\ff$ over $\F_p$. We have a commutative square
$$\xymatrix{U(\breve{\ff}) \ar[d] \ar[r] & \operatorname{Dist}_1(\breve{F}) \ar[d] \\ U(\ff) \ar[r]&  \operatorname{Dist}_1(F) }$$

Pick $x\in \ff_{\F_p}(=\ff_{\F_p}^{(1)})$. Our goal is to write down an element in $\operatorname{Fr}^{-1}(x)\subset  \operatorname{Dist}_1(F)$. Take  preimages $\breve{x},\breve{x}^{[p]}\in \breve{\ff}$ of $x$ and $x^{[p]}$. Then we have a well-defined element $\frac{1}{p!}(\breve{x}^p-\breve{x}^{[p]})\in \operatorname{Dist}_1(\breve{F})$. Let $y$ denote its image in $\operatorname{Dist}_1(F)$ (note that it depends on the choices of $\breve{x}$ and $\breve{x}^{[p]}$). 

Here is the main result of this section. 

\begin{Lem}\label{Lem:Fr_image}
We have $\operatorname{Fr}(y)=x$. 
\end{Lem}
\begin{proof}
The proof is in several steps. 

{\it Step 1}. First, we show that $\operatorname{Fr}(y)$ is independent of the choices: of $\breve{x}$ and $\breve{x}^{[p]}$. 
For this, observe that since the differential of $\operatorname{Fr}: F\rightarrow F$ is zero, it follows that $\operatorname{Fr}$ sends the image of $\ff$ in $\operatorname{Dist}_1(F)$ to $0$. If $\breve{x}'$ and $\breve{x}'^{[p]}$ are different lifts, then $\breve{x}^p-\breve{x}'^p, \breve{x}^{[p]}-\breve{x}'^{[p]}$ lie in the image of the ideal $p U(\breve{\ff})\breve{\ff}$ in $\operatorname{Dist}_1(\breve{F})$. So the difference between $\frac{1}{p!}(x^p-x^{[p]}), 
\frac{1}{p!}(x'^p-x'^{[p]})$ lies in the image of the ideal $U(\breve{\ff})\breve{\ff}$ in $\operatorname{Dist}_1(\breve{F})$. Hence its image in $\operatorname{Dist}_1(F)$ lies in the kernel of $\operatorname{Fr}$. 

{\it Step 2}. Now suppose $\breve{F}'$ is another smooth algebraic group scheme over $\Z_p$ equipped with an algebraic group homomorphism $\varphi: \breve{F}\rightarrow \breve{F}'$. We denote the induced homomorphisms $\operatorname{Dist}_1(\breve{F})\rightarrow \operatorname{Dist}_1(\breve{F}'), U(\ff)\rightarrow U(\ff')$ etc. also by $\varphi$. Observe that if the claim of the lemma holds for $(x,y)$, then it holds for $(\varphi(x), \varphi(y))$. Now suppose $\varphi:\breve{F}\rightarrow \breve{F}'$ is a closed embedding. 
Then $\varphi: \operatorname{Dist}_1(\breve{F})\rightarrow \operatorname{Dist}_1(\breve{F}')$ is injective, so if the claim of the lemma holds for $(\varphi(x),\varphi(y))$, then it also holds for $(x,y)$.

{\it Step 3}. Now we consider the following situation. Let $x_1,x_2\in \ff$ and $x:=x_1+x_2$.
Pick lifts $\breve{x}_1,\breve{x}_2, \breve{x}_1^{[p]},\breve{x}_2^{[p]}, \breve{x}^{[p]}$. Consider the corresponding elements $y_1,y_2$. For $y$ we take the image of 
$\frac{1}{p!}((\breve{x}_1+\breve{x}_2)^p-\breve{x}^{[p]})$. Assume that $\operatorname{Fr}(y_i)=x_i$ for $i=1,2$. We claim that $\operatorname{Fr}(y)=x$.

Recall that in $\ff$ we have 
$$(x_1+x_2)^p-x_1^p-x_2^p=(x_1+x_2)^{[p]}-x_1^{[p]}-x_2^{[p]}=\sigma(x_1,x_2)$$
for some Lie polynomial $\sigma$, see Section \ref{SSS_restricted_Lie}. From here we see that the difference 
$$\left((\breve{x}_1+\breve{x}_2)^p-(\breve{x}_1+\breve{x}_2)^{[p]}\right)-(\breve{x}_1^p-\breve{x}_1^{[p]})-(\breve{x}_2^p-\breve{x}_2^{[p]})$$
lies in $p U(\breve{\ff})$. As in Step 1, this implies the claim in the previous paragraph.

{\it Step 4}. In the next three steps we consider three special cases of $F$. First, assume that $F=\mathbb{G}_a$ so that $\Z_p[F]=\Z_p[t]$ and we can take $x=\partial_t$.  Then $x^{[p]}=0$ so we take $\breve{x}^{[p]}=0$. The claim that $\operatorname{Fr}(y)=x$ amounts to the following: for $f\in \Z_p[F]$, we have 
$$[\frac{1}{p!}\breve{x}^p].(f^p)-\breve{x}.f\in p\Z_p.$$
For constant $f$ or $f$ in the square of the ideal of $1$ in $\Z_p[F]$ both sides are $0$. 
So we only need to consider $f=t$ and $\breve{x}=\partial_t$, where the computation is easy. 

{\it Step 5}. Now consider the case when $F=\mathbb{G}_m$. Here $\Z_p[F]=\Z_p[t^{\pm 1}]$, and we can set $x=t\partial_t$. So $x^{[p]}=x$ and we take $\breve{x}^{[p]}=\breve{x}$. Again, we need to check that 
$$[\frac{1}{p!}(\breve{x}^p-\breve{x})].(f^p)-\breve{x}.f\in p\Z_p$$
and it is enough to do this for $f=t$. This reduces to the obvious congruence: $(p-1)!+1$ is divisible by $p$.

{\it Step 6}. Now consider the case when $G=\operatorname{GL}_n$. Thanks to Step 3, we can reduce to the case when $x$ is one of the matrix units $E_{ij}$. Using Step 2, we reduce the check for $x=E_{ij}$ to Step 4 if $i\neq j$ and to Step 5 if $i=j$.

{\it Step 7}. To handle the case of general $\breve{F}$, we recall that $\breve{F}\subset \operatorname{GL}_n$ and use Step 2 combined with Step 6.
\end{proof}

\subsection{Poisson structures}\label{SS_pcenter_Poisson}
We assume that $\Ring$ is an algebra over $\F_p$ with $p>2$. For $G$ we take a connected reductive group. 

The goal of this section is to show the equality between two Poisson brackets on the $p$-center of $\widehat{U}_\kappa(\hat{\g})$
(i.e., the image of the p-center of $\widehat{U}(\hat{\g})$ in $\widehat{U}_\kappa(\hat{\g})$). 

The first bracket is as follows. Consider the graded complete algebra $\widehat{S}(\hat{\g}^{(1)})$. The Lie bracket on $\hat{\g}^{(1)}$ induces a continuous degree $-1$ Poisson bracket on  $\widehat{S}(\hat{\g}^{(1)})$ and hence on its quotient $\widehat{S}_{\AS(\kappa)}(\hat{\g}^{(1)})$  that is the $p$-center of 
$\widehat{U}_\kappa(\hat{\g})$. This is the first Poisson bracket on the p-center that we need.

We proceed to the 2nd bracket. 
Assume for now that $\Ring=\F_p$.   We can form the algebra $\hat{\g}_{\Z_p}$ using the same $K_0$-class of representation $W$ as was used to form 
$\hat{\g}_{\F_p}$. In particular, $\widehat{U}(\hat{\g}_{\Z_p})/(p)$ is identified with $\widehat{U}(\hat{\g}_{\F_p})$ as an algebra over $\F_p[\mathbf{1},\mathbf{1}']$.

For $x,y$ in the center of $\widehat{U}(\hat{\g}_{\F_p})$, consider their lifts $\breve{x},\breve{y}$ in $U(\hat{\g}_{\Z_p})$. Then the image of $\{x,y\}:=\frac{1}{p!}[\tilde{x},\tilde{y}]$ is a well-defined element of the center of  $\widehat{U}(\hat{\g}_{\F_p})$. This gives a Poisson bracket on the entire center 
$Z(\widehat{U}(\hat{\g}_{\F_p}))$ of 
$\widehat{U}(\hat{\g}_{\F_p})$. We obtain 
the bracket on $\Ring\otimes_{\F_p}Z(\widehat{U}(\hat{\g}_{\F_p}))$. Note that the $p$-center of $\widehat{U}(\hat{\g})$ is a subalgebra of $\Ring\otimes_{\F_p}Z(\widehat{U}(\hat{\g}_{\F_p}))$.

\begin{Lem}\label{Lem:iota_Poisson}
The inclusion $\widehat{S}_{\AS(\kappa)}(\hat{\g}^{(1)})\hookrightarrow \Ring\otimes_{\F_p}Z(\widehat{U}(\hat{\g}_{\F_p}))$ is Poisson.
\end{Lem}
\begin{proof}
Let $\iota$ denote the inclusion map. 
It is enough to show that $\iota(\{x,y\})=\{\iota(x),\iota(y)\}$ when $x$ is an element of the form $x=e_{\alpha,(n)}(=e_\alpha t^n)$ or $x$ is of the form $h_{(n)}$ for a central element $h\in \mathfrak{g}$. Here $n\in \Z$. 
The reason is that these elements generate the Lie algebra $\hg$. 

{\it Case 1}.
First, assume $x=e_{\alpha,(n)}$. Note that $x$ lies in the image of an embedding $\mathfrak{sl}_2\hookrightarrow \hat{\g}$.

The claim that $\iota$ is 
equivariant from Proposition \ref{p:p-cent} implies that 
$x.\iota(y)=\iota(\{x,y\})$, where in the right hand side we view $x$ as an element of $\mathfrak{sl}_2^{(1)}$. Thanks to Lemma \ref{Lem:Fr_image}, $x$ acts on $\iota(y)$ as the divided power $x^{(p)}$ in the distribution algebra of $\operatorname{SL}_2$. 

It remains to show that $x^{(p)}.\iota(y)=\{\iota(x),\iota(y)\}$.
Let $\breve{y}$ be a lift of $\iota(y)$ to $\widehat{U}(\hg_{\Z_p})$ and 
$\breve{x}:=e_{\alpha,(n)}$. Then $x^{(p)}.\iota(y)$ is the image of 
$\breve{x}^{(p)}.\breve{y}=\frac{1}{p!}\operatorname{ad}(\breve{x})^p \breve{y}$. 
So we need to show that the image in $\widehat{U}(\hg_{\F_p})$ of 
$$\frac{1}{p!}\left(\operatorname{ad}(\breve{x})^p-\operatorname{ad}(\breve{x}^p)\right)\tilde{y}$$
is zero, equivalently, that the latter element is divisible by $p$. The element is equal to $\sum_{i=1}^{p-1}\frac{(-1)^i}{(p-i)!i!}\breve{x}^{p-i}\breve{y}x^i$
and its image is 
$$\sum_{i=1}^{p-1}\frac{(-1)^i}{(p-i)!i!}x^{p-i}\iota(y)x^i=\left(\sum_{i=1}^{p-1}\frac{(-1)^i}{(p-i)!i^!}\right)x^p\iota(y),$$
the last equality holds because $\iota(y)$ is central. It remains to note that 
$\sum_{i=1}^{p-1}\frac{(-1)^i}{(p-i)!i^!}=0$ in $\Z[\frac{1}{(p-1)!}]$: when we multiply this element by $p!$ we get $(1-1)^p$. 

{\it Case 2}. Let $x=h_{(n)}$ for central $h\in \g_{\F_p}$. By symmetry, we can also assume that $y=f_{(m)}$ for central $f\in \g_{\F_p}$. We can pick lifts $\breve{x}:=\breve{h}_{(n)}, \breve{y}:=\breve{f}_{(m)}$ for central $\breve{h},\breve{f}\in \g_{\Z_p}$. Set $\breve{c}:=\breve{\beta}_W(\breve{h},\breve{f})\mathbf{1}\in \hg_{\Z_p}$. The element $\{x,y\}$ is zero unless $m+n=0$, on the other hand $[\breve{h}^p_{(n)}-\breve{h}_{(np)}, \breve{f}^p_{(m)}-\breve{f}_{(mp)}]$ is easily seen to be divisible by $p^2$ if $m+n\neq 0$, so $\{\iota(x),\iota(y)\}=0$ in this case too. It remains to consider the case when $m+n=0$. Here $\{x,y\}=nc$ so $\iota(\{x,y\})=n(c^p-c)$. Now we examine $$[\breve{h}^p_{(n)}-\breve{h}_{(np)}, \breve{f}^p_{(m)}-\breve{f}_{(mp)}]=[\breve{h}^p_{(n)},\breve{f}^p_{(n)}]+np\breve{c}.$$
Similarly to Case 1, $[\breve{h}^p_{(n)},\breve{f}^p_{(n)}]$ is congruent 
to $\operatorname{ad}(h_{(n)})^p(\breve{f}^p_{(n)})=p!n^p\breve{c}^p$ modulo $p^2$.
And $\frac{1}{p!}(p!n^p\breve{c}^p+np\breve{c})$ is congruent to $n(\breve{c}^p-\breve{c})$ modulo $p$ finishing the proof.
 \end{proof}

\subsection{Poisson bracket on the centers of $\widehat{U}_\kappa(\hg), \widehat{U}_\kappa(\hh)$}
We continue to assume that $\Ring$ is a perfect field of characteristic $p>2$. And we return to the situation of $G=\operatorname{SL}_2$.

As was explained in Section \ref{SS_pcenter_Poisson}, the centers of $\widehat{U}_\kappa(\hg),\widehat{U}_\kappa(\hh)$ carry  natural Poisson structures. The purpose of this section is to compute the brackets between the elements $X_n,X_m$. This will be needed for the computation of $\HC(X_n)$ in  Section \ref{SS_HC_Xn_formula}.

\sss
We start with easier computations.

\begin{Lem}\label{Lem:Poisson_brackets1}
The following claims are true:
\begin{enumerate}
\item $\{\iota(\underline{L}_m),X_n\}=0$ for all $n,m\in \Z$
\item $\{\iota(\underline{L}_m),\iota(\underline{L}_n)\}=(\kappa^p-\kappa)(m-n)\iota(\underline{L}_{n+m})$.
\end{enumerate}
\end{Lem}
\begin{proof}
(1): Since $X_n$ is $\Loop G$-invariant,  see Proposition \ref{Prop:HC_center_sl2_noninteger}, the argument of the proof of 
Lemma \ref{Lem:iota_Poisson} implies that $\{\iota(y), X_n\}=0$
for $y$ of the form $e_{k}$ or $f_{k}, k\in \Z$. These two families topologically generate the Poisson algebra $\widehat{S}_{\AS(\kappa)}(\hg^{(1)})$  and so $\{\iota(y), X_n\}=0$ for all
$y\in \widehat{S}_{\AS(\kappa)}(\hg^{(1)})$, finishing the proof of (1).

(2): By Lemma \ref{Lem:iota_Poisson}, $\iota$ is a Poisson homomorphism.
So we need to prove that 
$$\{\underline{L}_m,\underline{L}_n\}=\AS(\kappa)(m-n)\underline{L}_{n+m}.$$
This formula is a straightforward consequence of its direct analog in $\widehat{S}(\hg)$.
That formula follows from the $\slf_2$-case of \cite[(12.8.8)]{Kac} by passing to the associated graded. This finishes the proof.
\end{proof}

\sss\label{SSS:Virasoro_computation}
It remains to compute the bracket between the elements $Y_n,Y_m$. Note that if $\kappa=0$, then the elements $L_i\in \widehat{U}_\kappa(\hg)$
are $\Loop G$-invariant for all $i$, see (\ref{eq:Sugawara_adjoint}). In particular, $\{Y_m,Y_n\}=0$. So, assume $\kappa\neq 0$. Then we can consider the elements $\kappa^{-1}L_i$ and $\kappa^{-p}Y_i$. 

Thanks to \cite[(12.8.8)]{Kac}, the elements $\kappa^{-1}L_i$ satisfy the following commutation relation
\begin{equation}\label{eq:Virasoro_commutation}
[\kappa^{-1}L_m,\kappa^{-1}L_n]=(m-n)\kappa^{-1}L_{n+m}+\delta_{m+n,0}\frac{3(\kappa-2)}{2\kappa}{m+1\choose 3}.
\end{equation}

Recall the Virasoro algebra $\mathfrak{Vir}$ from Section \ref{SSS_Virasoro}, it has a topological basis $L_i', 1$. In particular, $L_i'\mapsto \kappa^{-1}L_i, 1\mapsto -\frac{3(\kappa-2)}{2\kappa}$ defines a homomorphism $\psi:\widehat{U}(\mathfrak{Vir})\rightarrow \widehat{U}_\kappa(\hg)$. We note that this homomorphism also makes sense over $\Z_p$ (where we use a version of $\widehat{U}$, where $\kappa$ is an invertible indeterminate).

Both $\widehat{S}(\mathfrak{Vir}^{(1)})$ and the center $Z(\widehat{U}(\mathfrak{Vir}))$
are Poisson algebras similarly to Section \ref{SS_pcenter_Poisson}. 

\begin{Lem}\label{Lem:restricted_Vir}
The inclusion $\iota:\widehat{S}(\mathfrak{Vir}^{(1)})\hookrightarrow Z(\widehat{U}(\mathfrak{Vir}))$ is Poisson.
\end{Lem}
\begin{proof}
The proof is similar to that of Lemma \ref{Lem:iota_Poisson}. It reduces to showing $\iota(\{L_i',\cdot\})=\{\iota(L_i'),\cdot\}$. The proof is in several steps. 

Until the further notice assume that $i\geqslant 0$. 
Recall, Proposition \ref{p:p-cent}, that $\iota$ is $\operatorname{Aut}(\D^\times)$-equivariant, and, in particular, $\operatorname{Aut}^\bullet(\D)$-equivariant, { cf. Section \ref{sss:autd} for the definition of the latter.
} 
Note that $\operatorname{Aut}^\bullet(\D)$
is a pro-algebraic group, it is the inverse limit of the groups $\operatorname{Aut}^\bullet(\D)_i$: viewed as a functor this group sends a $\Ring$-algebra $S$ to the group of $S$-linear automorphisms of $S[t]/(t^{i+1})$.

We claim that the bracket with $\iota(L_i')$
in $Z(U(\mathfrak{Vir}))$ is the same as the action of $L_i'\in \operatorname{Lie}(\operatorname{Aut}^\bullet(\D)^{(1)})$. For this we will show that the latter action is induced from $\frac{1}{p!}(\operatorname{ad}(\breve{L}_i')^p-\operatorname{ad}(\breve{L}_i'^{[p]}))$. In order to do this, we will use Lemma \ref{Lem:Fr_image}. Note that the lemma deals with algebraic groups, while $\operatorname{Aut}^\bullet(\D)$ is a pro-algebraic group, but the claim generalizes to that setting. 


Consider now the restriction of $\iota$ 
to $S(\mathfrak{Vir}^{(1)})_{\leqslant j}\rightarrow Z(U(\mathfrak{Vir}))_{\leqslant pj}$. For $N\in \Z_{>0}$, let $I_N$ be the left ideal in $U(\mathfrak{Vir})$ generated by the elements $L'_i$ for $i\geqslant N$. Then $\iota$ induces a homomorphism
\begin{equation}\label{eq:iota_induced}
S(\mathfrak{Vir}^{(1)})_{\leqslant j}/(S(\mathfrak{Vir}^{(1)})_{\leqslant j}\cap I_N)\rightarrow Z(U(\mathfrak{Vir}))_{\leqslant pj}/
(Z(U(\mathfrak{Vir}))_{\leqslant pj}\cap I_N).
\end{equation}
Both sides are rational representations of the pro-algebraic group
$$\operatorname{Aut}^\bullet(\D)=\varprojlim_{i\rightarrow \infty}
\operatorname{Aut}^\bullet(\D)_i.$$ Apply the Lemma \ref{Lem:Fr_image} to $F:=\operatorname{Aut}(\D)_i$ (and its lift to $\Z_p$) for $i\gg 0$.
We see that the action of $L_i'$ on the target of
(\ref{eq:iota_induced}) is indeed induced by 
$\frac{1}{p!}(\operatorname{ad}(\breve{L}_i')^p-\operatorname{ad}(\breve{L}_i'^{[p]}))$. 

The action of $\frac{1}{p!}(\operatorname{ad}(\breve{L}_i')^p-\operatorname{ad}(\breve{L}_i'^{[p]}))$ on $S(\mathfrak{Vir}^{(1)})$
is by $\{L_i',\cdot\}$. Similarly to Step 1 in the proof of Lemma \ref{Lem:iota_Poisson}, the action of $\frac{1}{p!}(\operatorname{ad}(\breve{L}_i')^p-\operatorname{ad}(\breve{L}_i'^{[p]}))$ on $Z(U(\mathfrak{Vir}))$ is by $\{\iota(L_i'),\cdot\}$. This finishes the proof of  $\iota(\{L_i',\cdot\})=\{\iota(L_i'),\cdot\}$.

It remains to show that for all $i<0, a\in S(\mathfrak{Vir}^{(1)})$ we have 
$\{\iota(L_i'),a\}=\{\iota(L_i'),\iota(a)\}$. For this, consider the ``polynomial'' version of the Virasoro algebra, to be denoted by $\underline{\mathfrak{Vir}}$, one where the elements $L_i',1$
form a genuine, not a topological, basis. Note that this algebra has an automorphism sending $L_m'$ to $-L_{-m}', 1$ to $-1$. Now we use the case of 
$i\geqslant 0$ to deduce that $\iota(\{L_i',a\})=\{\iota(L_i'),\iota(a)\}$
for all $a\in S(\underline{\mathfrak{Vir}}^{(1)})$ and $i<0$. By continuity, we get the same equality for all $a\in \widehat{S}(\mathfrak{Vir}^{(1)})$.
\end{proof}

%

\sss Now we are ready to compute $\{X_m,X_n\}$.
\begin{Cor}\label{Cor:Y_brackets}
We have 
\begin{align*}
&\{Y_m,Y_n\}=\kappa^p\left((m-n) Y_{m+n}+\delta_{m+n,0}\frac{m^3-m}{2}(\kappa^{p-1}-1)\right),\\
&\{X_m,X_n\}=(1-\kappa^{p-1})\kappa^p\left((m-n)X_{m+n}-\delta_{m+n,0}\frac{m^3-m}{2}(1-\kappa^{p-1})^2\right).
\end{align*}
\end{Cor}
\begin{proof}
To prove (1) we use the notation from Section \ref{SSS:Virasoro_computation}.
(1) follows from the observation that $Y_i=\kappa^p\psi(\iota(L_i'))$ (recall from Section \ref{SSS:Virasoro_computation} that $\psi$ is the Sugawara homomorphism) for all $i$ 
and the formulas $\{\iota(L_i'),\iota(L_j')\}=\iota([L_i',L_j'])$ and $$\psi(\iota(1))=\psi(1^p-1)=-\frac{3(\kappa-2)^p}{2\kappa^p}+\frac{3(\kappa-2)}{2\kappa}=\frac{3(1-\kappa^{p-1})}{\kappa^p}.$$

Now we proceed to (2). By (1) of Lemma \ref{Lem:Poisson_brackets1}, 
we have$$\{X_m,X_n\}=\kappa^{2(p-1)}\iota(\{\underline{L}_m, \underline{L}_n\})-(\kappa^{p-1}-1)^2\{Y_m,Y_n\}.$$
The first summand is computed in (2) of Lemma \ref{Lem:Poisson_brackets1}, while the 2nd is computed using (1). Simplifying the resulting expression, we get (2).
\end{proof}

\sss Set $\tilde{h}_m:=h_m^p-\kappa^{p-1}h_{mp}$, this is an element in $\widehat{U}_\kappa(\hh)^{\Loop H}$. Moreover, the latter is a filtered complete algebra of polynomials in the variables $\tilde{h}_n$, this is a special case of Proposition \ref{Prop:HC_center_abelian}.
We proceed to computing $\{\tilde{h}_{m},\tilde{h}_{m}\}$ in the center of $\widehat{U}_\kappa(\hh)$. 
The following lemma is proved similarly to Case 2 in the proof of Lemma \ref{Lem:iota_Poisson}.

\begin{Lem}\label{Lem:tilde_h_brackets}
We have $\{\tilde{h}_{m}, \tilde{h}_{n}\}=\delta_{m+n,0}2m\kappa^p(1-\kappa^{p-1})$. 
\end{Lem}

\sss 
Set $X'_n:=\frac{1}{4}\left(\sum_{j\in \Z}\tilde{h}_{j+n}\tilde{h}_{-j}\right)$. We have the following standard corollary of Lemma \ref{Lem:tilde_h_brackets}.

\begin{Cor}\label{Cor:Xprime_brackets}
We have the following equalities
\begin{align*}
&\{X'_m, \tilde{h}_{j}\}=\kappa^p(\kappa^{p-1}-1)j \tilde{h}_{j+m},\\
&\{X_m',X_n'\}=\kappa^p(1-\kappa^{p-1})(m-n)X'_{m+n}.
\end{align*}
\end{Cor}

\subsection{Formula for $\HC(X_n)$}\label{SS_HC_Xn_formula}
\begin{Prop}\label{Prop:HC_quadratic}
We have $$\HC(X_n)=\frac{1}{4}\left(\sum_{j\in \Z}\tilde{h}_{j+n}\tilde{h}_{-j}\right)-\frac{1}{2}(1-\kappa^{p-1})(n+1)\tilde{h}_{n}.$$
\end{Prop}

\begin{proof}For ease of reading, we will break the proof into a series of steps. 

{\it Step 1.} By the compatibility between the maps of vertex algebras and enveloping algebras, (\ref{eq:U_homom}), note that the desired formula is equivalent to the assertion that 
$$\ffr_B(X_{-2} | \varnothing \rangle) \overset{?}= \frac{1}{4} \tilde{h}_{-1}^2 | \varnothing \rangle + \frac{1}{2}(1 - \kappa^{p-1}) \tilde{h}_{-2} | \varnothing \rangle. $$
To see this, we first note that $\HC(X_{-2})$ lies in the Harish--Chandra center of $V_{\kappa}(\fh)$, and is homogeneous of degree $-2p$ with respect to the energy grading. By Proposition \ref{Prop:HC_center_abelian}, we therefore must have 
$$\ffr_B(X_{-2} | \varnothing \rangle) = \beta \cdot \tilde{h}_{-2}^2 | \varnothing \rangle - \alpha \cdot \tilde{h}_{-2} | \varnothing \rangle,$$
for some scalars $\alpha,\beta\in \Ring.$

{\it Step 2.} By considering the associated graded with respect to the PBW filtration and using Lemma \ref{Lem:HC0_interpretation}, it is straightforward to see that $a = 1/4$, and hence on enveloping algebras we have 
$$\HC(X_{m}) = X'_m + \alpha \cdot (m+1) \cdot  \tilde{h}_m, \quad \quad m \in \mathbb{Z}. $$
To prove the proposition, it therefore remains to see that $\alpha = \frac{1}{2}(1 - \kappa^{p-1}).$ 

{\it Step 3.} To calculate $\alpha$, we will proceed via a calculation of Poisson brackets. Namely, we first note that $\HC$ is a Poisson algebra homomorphism because $\ffr_B$ restricts to the Poisson homomorphism from the center $\widehat{U}_\kappa(\hg)$ to the center of 
$\widehat{D}(\Loop N)\widehat{\otimes} \widehat{U}_\kappa(\hh)$ because it deforms to the algebras over $\Z_p$. Therefore, we have that 
$$\HC(\{ X_m, X_{-m} \} ) = \{ \HC(X_m), \HC(X_{-m}) \}  =  \{ X'_m + \alpha \cdot (m+1) \cdot \tilde{h}_m, X'_{-m} + \alpha \cdot (-m + 1) \cdot \tilde{h}_{-m} \}.$$
Using Corollaries \ref{Cor:Y_brackets}, 
\ref{Cor:Xprime_brackets}, and Lemma 
\ref{Lem:tilde_h_brackets} to compute the left and right hand sides separately, and acting on $| \varnothing \rangle$, we conclude that $$(m^3 - m) \cdot \alpha^2 =(m^3-m) \cdot \frac{1}{4}(1-\kappa^{p-1})^2,$$ 
and hence $\alpha = \pm \frac{1}{2} (1 - \kappa^{p-1}).$

{\it Step 4.} To calculate the remaining sign, we first note it is independent of the choice of $\kappa$, so we may assume that $\kappa=0$. In this case $X_n=Y_n=L_n^p$, and moreover $L_n$ is already in the Harish-Chandra center. It is well-known, see e.g. \cite[Proposition 6.2.2]{Frenkel_loop}, that 
$$\HC(L_n)=\Ffr_B(L_n)=\left(\frac{1}{4}\sum_{j\in \Z}h_{j+n}h_{-j}\right)-\frac{1}{2}(n+1) h_{n}.$$
Taking the $p$th power (note that $h_{i}$'s pairwise commute), we see that the coefficients of $\tilde{h}_n$ is $-\frac{1}{2}(n+1)$ and hence $\alpha=-\frac{1}{2}(1-\kappa^{p-1})$. This completes the proof.
\end{proof}

\section{Verma and Wakimoto modules}\label{S_Verma_Wak}
In this section we introduce two classes of modules over the (Rees versions of) affine Lie algebras and study homomorphisms between them. 
\subsection{Construction}
\sss First, we introduce algebras whose modules we are going to study. Set $\Ring:=\Z_p$ for $p>2$. Let $G$ be a split connected reductive group, $H\subset B$ be a maximal torus and Borel subgroup. We fix a parabolic subgroup $P\subset G$ containing $B$ and its Levi decomposition $P=L\ltimes N$, where $L$ contains $H$. Let $B=H\ltimes \tilde{N}$ be the Levi decomposition for $B$.
We write $B_L$ for $B\cap L$.

Consider the algebra $\widehat{U}^\hbar_{\bone}(\hg):=R_\hbar(\widehat{U}_{\bone}(\hg))$. This is an algebra over $\Ring[\mathbf{1},\hbar]$. Its specialization to $\bone=\kappa,\hbar=1$ is $\widehat{U}_\kappa(\hg)$, while its specialization to $\kappa=\hbar=0$ is $\Ring[\Loop(\g^*)]_f$.

The algebra $\widehat{U}^\hbar_{\bone}(\hg)$ carries the following 
compatible gradings. First, there is the naive grading, where 
the degrees of $\hbar,\bone$ and $\hg$ are equal to $1$. 

Second, we have a grading in a generalized sense to be defined now. Consider the root lattice $\Lambda$ of $\g$ and the imaginary root $\delta$. Form the affine root lattice $\Lambda^a:=\Lambda\oplus \Z\delta$.
Now let $V$ be a topological $\Ring$-module that is complete and separated with respect to its topology (here we view $\Ring$ as equipped with the discrete topology). By a {\it topological $\Lambda^a$-grading} on $V$
we mean the following data: 
\begin{itemize}
\item A system of neighborhoods $U_N, N>0,$ of $0$ with $U_N\supset U_{N+1}$ and $V=\varprojlim_{N\rightarrow \infty}V/U_N$,
\item and a collection of  $\Lambda^a$-gradings on $V/U_N$ for each $N$ such that $V/U_{N+1}\twoheadrightarrow V/U_N$ is graded. 
\end{itemize}

Each graded piece of the  algebra $\widehat{U}_{\bone}^\hbar(\hg)$ for the naive grading carries a topological $\Lambda^a$-grading induced from the $\Lambda^a$-grading on $\hg$ (with $\bone$ and $\hbar$ in degree $0$). In what follows we will abuse the terminology and say that the (graded complete) algebra $\widehat{U}_{\bone}^\hbar(\hg)$ is $\Lambda^a$-graded.  

Similarly, we can consider the $\Ring[\bone,\hbar]$-algebra $\widehat{U}^\hbar_{\bone}(\hl)$. It also carries compatible naive and $\Lambda^a$-gradings. 

Next, 
we consider the algebra $\widehat{D}^\hbar(\Loop N)$, the Rees algebra of $\widehat{D}(\Loop N)$. Let $\Phi_P^+,\partial_\alpha,y^\alpha$ for $\alpha\in \Phi^+_P$
have the same meaning as in Section \ref{SSS_Cartan_images_ffr}. The algebra $\widehat{D}_\hbar(\Loop N)$ has topological generators $\partial_{\alpha,n}, y^\beta_m$ with $\alpha,\beta\in \Phi^+_P$ and $n,m\in \Z$ with relations 
$[\partial_{\alpha,n},y^{\beta}_m]=\delta_{m+n,0}\delta_{\alpha,\beta}\hbar$
and all other brackets equal to $0$. This algebra also carries naive and $\Lambda^a$-gradings. For the naive grading, we have $\deg y^\beta_m=0, \deg \hbar=\deg \partial_{\alpha,n}=1$. The $\Lambda^a$-degrees are as follows:
$\deg \hbar=0, \deg y^\beta_m=-\beta+m\delta, \deg \partial_{\alpha,n}=\alpha+n\delta$. 

Consider the graded complete tensor product $\widehat{D}^\hbar(\Loop N)\widehat{\otimes}_{\Ring[\hbar]}\widehat{U}^\hbar_{\bone}(\hl)$. Recall, Section \ref{SSS_ffr_un}, that we have a naive graded $\Ring[\bone,\hbar]$-linear algebra homomorphism
$$\Ffr_{P,\bone}^\hbar:  \widehat{U}_{\bone}^\hbar(\hg)\rightarrow \widehat{D}^\hbar(\Loop N)\widehat{\otimes}_{\Ring[\hbar]}\widehat{U}^\hbar_{\bone}(\hl).$$
It is also $\Lambda^a$-graded. 

\sss\label{SSS_Verma_Wakimoto}
Now we define certain left ideals in the algebras $\widehat{U}_{\bone}^\hbar(\hg), \widehat{D}^\hbar(\Loop N)\widehat{\otimes}_{\Ring[\hbar]}\widehat{U}_{\bone}^{\hbar}(\hl)$ and consider the quotients by these left ideals. 

Consider the Iwahori subgroup $\Iw$, the preimage of $B^-$ under $\Jet G\twoheadrightarrow G$ (we emphasize that we need to consider $B^-$ instead of $B$ here). We choose the system of simple affine roots $I^a$
for this choice of Iwahori, and choose positive affine roots accordingly, explicitly, those are of the form $\beta+n\delta$, where either $n>0$ or $n=0$ but $\beta$ is a negative root for $B$. Let $I\subset I^a$ be the system of simple roots for $\g$ (=negative simple roots for $B$). 

Suppose $A$ is a (graded complete) $\Lambda^a$-graded algebra. Consider the left ideal $I_{>0}(A)$, the closure of the left ideal
generated by the graded components of degrees in the span of positive roots.

Consider the left ideal $I^{G,\hbar}_{\bone}=I_{>0}(\widehat{U}^\hbar_{\bone}(\hg))$. 
Recall that the loop group $\Loop G$ acts on $\widehat{U}_\bone^{\hbar}(\hg)$.
The ideal $I^{G,\hbar}_{\bone}$ is stable under $\mathsf{Iw}$. Form the quotient 
$\Delta^{G,\hbar}_{\bone}:=\widehat{U}_\bone^{\hbar}(\hg)/I^1_{\hbar,\bone}$, the universal Verma module. It is acted on by $\Iw$ and is a rational representation of this pro-algebraic group. 
It carries a naive grading and also a (genuine) $\Lambda^a$-grading. 
For $\lambda\in \Lambda^a$, we write $\Delta^{G,\hbar}_{\bone}[\lambda]$ for the $\lambda$-graded component. Then $\Delta^{G,\hbar}_{\bone}[\lambda]$  is a free $\Ring[\bone,\hbar,\h^*]$-module 
of rank equal to the affine Kostant partition function of $\lambda$
(=the number of ways to represent $\lambda$ as the sum of positive affine roots). We also note that we have a commuting action of $\Ring[\bone,\hbar,\h^*]$-action on $\Delta^{G}_{\bone,\hbar}$, it comes from the observation that $I^{G,\hbar}_{\bone}\Ring[\bone,\hbar,\h^*]\subset I^{G,\hbar}_{\bone}$. This action is compatible with the naive grading, preserves the $\Lambda^a$-grading, and is $\Iw$-invariant. We view $\Delta_{\bone,\hbar}^G$ as a $\widehat{U}_{\bone}^\hbar(\hat{\g})$-$\Ring[\bone,\hbar,\h^*]$-bimodule.

Similarly, we consider the left ideal $$I^{P,\hbar}_{\bone}=I_{>0}(\widehat{D}^\hbar(\Loop N)\widehat{\otimes}_{\Ring[\hbar]}\widehat{U}_\bone^{\hbar}(\hl)).$$
The group $\Loop P$ acts on $\widehat{D}^\hbar(\Loop N)\widehat{\otimes}_{\Ring[\hbar]}\widehat{U}_\bone^{\hbar}(\hl)$ by automorphisms 
and $I^{P,\hbar}_{\bone}$ is stable with respect to the intersection 
$\Iw_P:=\Iw\cap \Loop P$. Set
$$\WakDelta^{P,\hbar}_{\bone}:=(\widehat{D}^\hbar(\Loop N)\widehat{\otimes}_{\Ring[\hbar]}\widehat{U}_\bone^{\hbar}(\hl))/I^{P,\hbar}_{\bone}.$$
We can view $\WakDelta^{P,\hbar}_{\bone}$ as a $\widehat{U}^\hbar_{\bone}(\hg)$-module via the homomorphism $\Ffr^\hbar_{P,\bone}$.
This module carries compatible naive and $\Lambda^a$-grading.
Notice that 
$$\WakDelta^{P,\hbar}_{\bone}=[\widehat{D}_\hbar(\Loop N)/I_{>0}(\widehat{D}_\hbar(\Loop N))]\otimes_{\Ring[\hbar]}\Delta_{\bone}^{L,\hbar}.$$
From here we see that the weight $\lambda$ component $\WakDelta^P_{\bone,\hbar}[\lambda]$ is a free $\Ring[\bone,\hbar,\h^*]$-module of the same rank as 
$\Delta^{G}_{\bone,\hbar}[\lambda]$. Similarly to the previous paragraph,
$\Ring[\bone,\hbar,\h^*]$ acts on $\WakDelta^P_{\bone,\hbar}$
from the right turning $\WakDelta^P_{\bone,\hbar}$ into a $\widehat{U}^\hbar_{\bone}(\hat{\g})$-$\Ring[\bone,\hbar,\h^*]$-bimodule. 

When $P=B$, we write $\Wak_{\hbar,\bone}$ for $\WakDelta^B_{\hbar,\bone}$. This is a universal Wakimoto module. And the modules 
$\WakDelta^{P,\hbar}_{\bone}$ interpolate between Verma and Wakimoto modules.

\sss Let $F\subset \Iw_P\rtimes \mathbb{G}_m$ be an algebraic subgroup.
This group acts on $\widehat{U}^\hbar_{\bone}(\hg)$, where the $\mathbb{G}_m$ factor acts by loop rotations. 
A possible choice is $F:=H\times \mathbb{G}_m$.
We can consider the categories $\widehat{U}_\bone^{\hbar}(\hg)\operatorname{-mod}^F_{gr}$ of naive graded $\widehat{U}_{\bone,\hbar}(\hg)$-modules equipped with a rational action of $F$ that preserves the naive grading and makes the module structure map $F$-equivariant. Similarly, we can consider the category 
$\widehat{U}_{\bone,\hbar}(\hl)\operatorname{-mod}^F_{gr}$. We have the {\it Wakimotization} functor 
\begin{align*}
&\Wak^P: \widehat{U}_{\bone,\hbar}(\hl)\operatorname{-mod}^F_{gr}\rightarrow 
\widehat{U}_{\bone,\hbar}(\hl)\operatorname{-mod}^F_{gr},\\ 
&M\mapsto 
(\widehat{D}_\hbar(\Loop N)/I_{>0}(\widehat{D}_\hbar(\Loop N)))\otimes_{\Ring[\hbar]}M.
\end{align*}
Note that $\Wak^P(\Delta^{L,\hbar}_{\bone})\xrightarrow{\sim}\WakDelta^{P,\hbar}_\bone$.

\subsection{Main result and reductions}
\sss We start by introducing some notation.
Our base ring is $\Ring=\Z_p$.
Set $\h^{da}:=\operatorname{Spec}(\Ring[\bone,\hbar,\h^*])$, an affine space. For a field $\bk$ we write $\h^{da}(\bk)$ for the set of $\bk$-points of $\h^{da}$. Fix a point $\sx=(\kappa,\nu,\eta)\in \h^{da}$
(with $\kappa,\nu\in \bk, \eta\in \h^*(\bk)$). We write $\Delta_{\sx},
\WakDelta^P_{\sx}, \Wak_{\sx}$ for the specializations of the corresponding modules to $\sx$. 

Note that by the construction we have the unique  $\widehat{U}_\bone^{\hbar}(\hg)$-linear homomorphism $\varphi^P:\Delta^{G,\hbar}_{\bone}\rightarrow \WakDelta^{P,\hbar}_{\bone}$ sending the coset of $1$ to the coset of $1$. This homomorphism is $\Iw_P\rtimes \mathbb{G}_m$-equivariant and graded for the naive grading. We write $\varphi^P[\lambda]$ for the restriction of $\varphi^P$ to the components of degree $\lambda$, and $\varphi^P_{\sx}$ for the specialization of $\varphi$ to the point $\sx$. 

\sss\label{SSS_shift} Recall that the algebra $\Ring[\h^{da}]$ acts on  
$\Delta^{G,\hbar}_{\bone}, \WakDelta^{P,\hbar}_{\bone}$ from the right. The homomorphism $\varphi^P$ is not equivariant for this action, instead, it is semi-linear with respect to a certain automorphism of $\Ring[\h^{da}]$. Since $\hbar$ and $\bone$ are central, $\varphi^P$ is $\Ring[\bone,\hbar]$-linear.  
Set $\rho_P:=\frac{1}{2}\sum_{\alpha\in \Phi^+_P}\alpha$, where $\Phi_P^+$ stands for the set of roots of $\mathfrak{n}$.

\begin{Lem}\label{Lem:varphi_right_action}
For any $x\in \h,v\in \Delta^{G,\hbar}_\bone$, 
we have $\varphi^P(mx)=\varphi^P(m)(x+2\hbar\langle\rho_P,x\rangle)$.
\end{Lem}
\begin{proof}
It is enough assume that $m:=1+I^{G,\hbar}_\bone$ because $\Delta^G_{\bone,\hbar}$ is generated by this element as a $\widehat{U}^\hbar_{\bone}(\hat{\g})$-module. Note that, by Lemma \ref{Lem:ffr_Cartan_images}, 
$$\Ffr_P(x)=
(x-\sum_{\alpha\in \Phi^+_P}\langle\alpha,x\rangle y^\alpha \partial_\alpha)+?,$$
where $?$ lies in the left ideal generated by elements of positive energy. Note that $y^\alpha\partial_\alpha=\partial_\alpha y^\alpha-1$ in $\widehat{D}(\Loop N)\widehat{\otimes}_{\Ring}\widehat{U}_{\bone}(\hl)$, and $y^\alpha$ kills $\varphi^P(m)$. In the Rees algebras, we need to replace $1$ with $\hbar$. 
So 
\begin{align*}
\varphi^P(mx)=\varphi^P(xm)=\Ffr_P(x)\varphi^P(m)=\varphi^P(m)(x+2\hbar\langle\rho_P,x\rangle)
\end{align*}
finishing the proof.
\end{proof}

Define the $\Ring$-linear automorphism $\mathsf{t}_P$ of $\h^{da}=\Ring\bone\oplus \Ring\hbar\oplus \h^*$ by 
$\mathsf{t}_P(\kappa,\nu,\eta)=(\kappa,\nu,\eta+2\nu\rho_P)$. So, for any $x\in \Ring[\h^{da}]$, we have $\varphi^P(mx)=\varphi^P(m)\mathsf{t}^*_P(x)$.

In particular, $\varphi^P_{\sx}$ is a homomorphism $\Delta^G_{\sx}\rightarrow \WakDelta^P_{\mathsf{t}_P(\sx)}$.

\sss
We remark that thanks to the transitivity property (\ref{eq:ffr_transitivity_algebra}), we have 
\begin{equation}\label{eq:varphi_transitivity}
\varphi^B=\Wak^P(\varphi^{B_L})\circ \varphi^P.
\end{equation}

\sss\label{SSS_KK_hyperplanes}  Note that we can view the hyperplane $\hbar=1$ in $\h^{da}$
as the affine Cartan space (over $\Ring$). Kac and Kazhdan in \cite{Kac_Kazhdan} specified a collection of hyperplanes $H_{\beta,m}$ in the affine Cartan space. Let us now explain what we need to know about these hyperplanes. These hyperplanes are of the form $\{\eta|(\eta,\beta)+(\rho,\beta)-n(\beta,\beta)=0\}$, where $(\cdot,\cdot)$ is an affine Weyl group invariant pairing, $\beta$ is an indecomposable positive root, $n$ is a positive integer, $\rho$ is a suitable element. Note that with a suitable choice of $(\cdot,\cdot)$ (namely with the usual normalization, where all short coroots have length $2$) all these hyperplanes are defined over $\Ring$. Let $\tilde{H}_{\beta,m}$ denote the unique linear hyperplane in $\h^{da}$ whose intersection with $\hbar=1$ is $H_{\beta,m}$. 

Here is a property of this collection of hyperplanes established in \cite{Kac_Kazhdan} that we are going to need.
Suppose $\operatorname{char}\bk=0$ and $\sx\in \h^{da}(\bk)$
is of the form $(\kappa,1,\eta)$. Since $\tilde{H}_{\beta,m}$ is defined over $\Ring$, it makes sense to speak of $\sx$ lying in  $\tilde{H}_{\beta,m}$, we will just write $\sx\in \tilde{H}_{\beta,m}$ in this case. Then the following is true, see \cite[Proposition 3.1]{Kac_Kazhdan}:
\begin{itemize}
\item[(*)] Suppose $\sx$ does not lie on any of the hyperplanes $\tilde{H}_{\beta,m}$. Then the $\widehat{U}_{\kappa}(\hg)$-module $\Delta^G_{\sx}$ is irreducible. 
\end{itemize}


\sss
Now we state the main result of this section.

\begin{defi}\label{defi:P_generic}
We say that $\sx$ is {\it $P$-generic} if from $\sx\in \tilde{H}_{\beta,m}$
it follows that $\beta$ is a root of $\hl$. 
\end{defi}

\begin{Thm}\label{Thm:Verma_vs_Wakimoto}
Suppose that $\sx$ is $P$-generic. Then $\varphi^P_{\sx}:\Delta^\g_{\sx}\rightarrow \WakDelta^P_{\mathsf{t}_P(\sx)}$ is an isomorphism. 
\end{Thm}

\sss 
In the rest of the section we will reduce the proof to the following result (that will be proved in Section \ref{SS_classical_specialization}). 

\begin{Prop}\label{Prop:iso_classical}
Suppose  $\sx=(0,0,\eta)$. Further suppose that the only root hyperplanes in $\h^*(\bk)$ containing $\eta$ correspond to roots of $\lf$. 
Finally assume that $\operatorname{char}\bk=0$ or $P=B$.
Then $\varphi_{\sx}^P: \Delta^G_{\sx}\rightarrow \WakDelta^P_{\mathsf{t}_P(\sx)}$ is an isomorphism.  
\end{Prop}

\begin{proof}[Proof of Theorem \ref{Thm:Verma_vs_Wakimoto} modulo Proposition \ref{Prop:iso_classical}]
The proof is in several steps. 

{\it Step 1}. First, suppose that $\sx$ does not lie on $(\hbar)$, $(p)$, or 
any of $\tilde{H}_{\beta,m}$. Property (*) in Section \ref{SSS_KK_hyperplanes}
implies that $\Delta^{G}_x$ is simple. Moreover, by the construction in Section \ref{SSS_Verma_Wakimoto}, the $\Lambda^a$-graded characters of 
$\Delta^\g_{\sx}$ and $\WakDelta^P_{\mathsf{t}_P(\sx)}$ are the same. It follows that 
$\varphi^P_{\sx}$ is an isomorphism. 

{\it Step 2}. 
Note that $\varphi^P[\lambda]: \Delta^G[\lambda]\rightarrow 
\WakDelta^P[\lambda]$ is a $\Ring[\h^{da}]$-semi-linear homomorphism between $\Ring[\h^{da}]$-modules of the same rank. By Step 1, $\varphi^P[\lambda]$ is an isomorphism generically. 
It follows that the locus in $\h^{da}$, where this homomorphism fails to be an isomorphism is exactly a divisor, denote it by $D[\lambda]$. Applying Step 1 again, we see that the irreducible components of $D[\lambda]$ are of the form $(p),(\hbar)$ or $\tilde{H}_{\beta,m}$. Note that thanks to (\ref{eq:varphi_transitivity}), the map $\varphi^P_{\sx}$ is injective provided $\varphi^B_{\sx}$ is.
Proposition \ref{Prop:iso_classical} now shows that $(p)$ and $(\hbar)$ do not appear among the components. 

{\it Step 3}. It remains to show that $\tilde{H}_{\beta,m}$, where $\beta$ is not a root of $\hat{\lf}$, cannot a component in any $D[\lambda]$. The intersection of $\tilde{H}_{\beta,m}$
with $\h^*\subset \h^{da}$ is a root hyperplane that is given by a root of $\g$ but not of $\lf$. 
We arrive at a contradiction with Proposition \ref{Prop:iso_classical}.
\end{proof}

\subsection{Specializations with $\kappa=\nu=0$}\label{SS_classical_specialization}
\begin{proof}[Proof of Proposition \ref{Prop:iso_classical}]
The proof is in several steps. Note that it is enough to consider the case when $\bk$ is algebraically closed. Below we work over $\bk$. 


{\it Step 1}. 
Consider the finite analog of $\ffr_P^0$.
Let $\mathfrak{n}^{-,\perp}$ denote the annihilator of $\mathfrak{n}^-$ in $\mathfrak{g}^*$. 
Consider the action of $G$ on $T^*(G/N^-)=(G\times \mathfrak{n}^{-,\perp})/N^-$, where 
$N^-$ acts by $n.(g,x)=(gn^{-1}, \operatorname{Ad}^*(n)x)$. The moment map
$T^*(G/N^-)\rightarrow \g$ sends the orbit $[g,x]$ to $\operatorname{Ad}^*(g)x$. It is $L$-invariant and so descends to $G\times^{P^-}\mathfrak{n}^{-,\perp}$. Restrict it to the open subset $N\times \mathfrak{n}^{-,\perp}$ and identify $\mathfrak{n}^{-,\perp}$ with 
$\mathfrak{n}^*\oplus \mathfrak{l}^*$. The resulting map 
$N\times \mathfrak{n}^*\times \mathfrak{l}^*\rightarrow \g^*$ is therefore given by $(n,x,y)\rightarrow \operatorname{Ad}^*(n)(x+y), n\in N, x\in \mathfrak{n}^*, y\in \mathfrak{l}^*$.  Denote this map by $\mu$. 

{\it Step 2}. We need to understand $\mu^*:\g\rightarrow \bk[N]\otimes (\mathfrak{n}\oplus \mathfrak{l})$ more explicitly. 
Choose a one-parameter subgroup $\varrho:\mathbb{G}_m\rightarrow Z(L)$ such that its centralizer is $\lf$ and its eigen-characters on $\mathfrak{n}$ are positive. So we have decompositions
$\mathfrak{n}=\bigoplus_{i>0}\mathfrak{n}_i$ and $\mathfrak{n}^-=\bigoplus_{i<0}\mathfrak{n}^-_i$. Under the assumption of the proposition, the action of $L$ on $N$ linearizes leading to an $L$-equivariant isomorphism $N\xrightarrow{\sim} \mathfrak{n}$ that can be chosen so that its differential at $1$ is the identity, the construction is in Section \ref{SSS_Cartan_images_ffr}. Let $\mathfrak{m}$ denote the augmentation ideal in $\bk[N]$. We write $\mathfrak{m}_i$ for its degree $i$ component with respect to the $\mathbb{G}_m$-action, note that $\mathfrak{m}_i\neq \{0\}\Rightarrow i<0$.

The map $\mu^*$ sends $x\in \mathfrak{n}_i$ to the corresponding action vector field on $N$, we have
\begin{equation}\label{eq:moment_n}
\mu^*(x)\in x+\sum_{j>i}\mathfrak{m}_{i-j}\otimes \mathfrak{n}_j.
\end{equation}
Next, for $y\in \mathfrak{l}$, the image $\mu^*(y)$ is the sum of $y$ and the action vector field for the adjoint action of $L$ on $\mathfrak{n}$, so
\begin{equation}\label{eq:moment_l}
\mu^*(y)\in y+\sum_{i>0}\mathfrak{n}_i\otimes \mathfrak{n}_i^*.
\end{equation}
Finally, consider $w\in \mathfrak{n}^-_{-i}$. Then for the weight reasons, we have 
\begin{equation}\label{eq:moment_n-}
\mu^*(w)\in \theta(w)+(\mathfrak{m}^2)_{-i}\otimes \lf
\oplus \bigoplus_{j>0}\mathfrak{m}_{-i-j}\otimes \mathfrak{n}_j.  
\end{equation}
Here $\theta$ is an $L$-equivariant linear map $\mathfrak{n}^-_{-i}\rightarrow \mathfrak{n}^-_{-i}\otimes \lf$, where in the right hand side we embed 
$\mathfrak{n}^-_{-i}$ into $\mathfrak{m}$ via its identification with $\mathfrak{n}_i^*$ via the invariant form.

{\it Step 3}. For $y\in \lf^*$, let $\theta_y:\mathfrak{n}^-_i
\rightarrow \mathfrak{n}^-_i$ denote the pairing of $\theta$ with $y$. 
We can find $F\in \bk[\lf^*]^L$ such that the zeroes of $F$ in $\h^*\subset \lf^*$ is exactly the union of the root hyperplanes defined by roots that are not of $\lf$. We claim that if $F(y)\neq 0$, then $\theta_y$ is an isomorphism. 
In order to see this
identify $\mathfrak{n}^*$ with $\mathfrak{n}^-$
and $\lf^*$ with $\lf$ using the invariant form.
The differential $d_{(1,x,y)}\mu$ at the point 
$(\xi,x',y')$ for $\xi\in \mathfrak{n}, x,x'\in \mathfrak{n}^-, 
y,y'\in \lf$ is given by $x'+y'+[\xi,x+y]$. The condition that $F(y)\neq 0$ is equivalent to the operator $\xi\mapsto [\xi,x]$ being invertible on 
$\mathfrak{n}_i$. This is equivalent to $\theta_y$ being an isomorphism. 

{\it Step 4}. We need a description of $\Ffr^0_P: \Ring[\Loop (\g^*)]_f\rightarrow \Ring[\Loop(T^*N)]_f\widehat{\otimes} \Ring[\Loop(\lf^*)]_f$ that stems from the construction and 
the description of the jet analog in Section \ref{SSS_ffrp_classical}.
Namely, $\ffr^0_P$ is the pullback under the morphism $\Loop \mu:
\Loop(N\times \mathfrak{n}^*\times \mathfrak{l}^*)\rightarrow \Loop(\g^*)$.
To compute the image of $x_n:=xt^n\in \Loop \g$ under $\Ffr^0_P$ we form 
the infinite series $\sum_{n}x_nz^{-n-1}$, then expand $\mu^*(x)$ as a polynomial, say $Q$, in the elements $\partial_\alpha, y^\alpha$ and $y\in \lf$. Then $\sum_{n}\Ffr_P^0(x_n)z^{-n-1}$ is obtained from $Q$ by replacing each $y$ with $\sum_n (yt^n)z^{-n-1}$, each $\partial_{\alpha}$ with $\sum_{n}\partial_{\alpha,n}z^{-n-1}$ and each $y^\alpha$ with 
$\sum_{n}y_n^\alpha z^{-n}$. 

{\it Step 5}. In Steps 6 and 7 we will show that $\varphi^P_{\sx}:\Delta^G_{\sx}\rightarrow \WakDelta^P_{\sx}$ is surjective. Since these modules have the same graded characters, this will imply that $\varphi^P_{\sx}$ is an isomorphism. 

We will use the following version of the graded Nakayama lemma. Suppose that $A,B$ are two positively graded commutative $\bk$-algebras (``positive'' means that they are graded by $\Z_{\geqslant 0}$ and the degree $0$ components are $\bk$). Let $\varphi:A\rightarrow B$ be a graded algebra homomorphism and $y_i$, where $i$ is in some indexing set $I$, are homogeneous elements of $A$ of positive degree. If the induced homomorphism $$A/\operatorname{Span}_A(y_i|i\in I)\rightarrow B/\operatorname{Span}_B(\varphi(y_i)|i\in I)$$ 
is surjective, then $\varphi$ is surjective. 

{\it Step 6}. Note that both $\Delta^\g_{\sx}, \WakDelta^P_{\sx}$ are commutative $\bk$-algebras: we have $\Delta^{G}_{\sx}\cong S(t^{-1}\mathfrak{b}^-[t^{-1}]\oplus\tilde{\mathfrak{n}}[t^{-1}])$, where $\tilde{\mathfrak{n}}$ is the maximal nilpotent subalgebra of $\mathfrak{b}$. Further, 
$$\WakDelta^P_{\sx}=\bk[\partial_{\alpha,n}]_{n\leqslant 0}\otimes_{\bk}\bk[y^\alpha_n]_{n< 0}\otimes_{\bk} S(t^{-1}\mathfrak{b}^-_L[t^{-1}]\oplus\tilde{\mathfrak{n}}_L[t^{-1}]).$$
The map $\varphi^P_{\sx}:\Delta^\g_{\sx}\rightarrow \WakDelta^P_{\sx}$ is the homomorphism of quotients of the algebras $$\Ring[\Loop (\g^*)]_f\widehat{\otimes}_{\Ring[\h^{da}]}\bk,
\left(\Ring[\Loop(T^*N)]_f\widehat{\otimes} \Ring[\Loop(\lf^*)]_f\right)\widehat{\otimes}_{\Ring[\h^{da}]}\bk$$ induced by an algebra homomorphism $\Ffr_P^0$. So, $\varphi^P_{\sx}$ is an algebra homomorphism. Note that both algebras $\Delta^\g_{\sx}$ and $\WakDelta^P_{\sx}$ are graded by the span of positive affine roots, we turn this grading into a $\Z_{\geqslant 0}$-grading by using the height of affine roots. 

{\it Step 7}. We first apply the procedure in Step 5 to the elements 
$\mathfrak{n}[t^{-1}]\subset \Delta^\g_{\sx}$. From (\ref{eq:moment_n}) and the description of $\Ffr^0_P$ in Step 4, we easily conclude  that 
$$\operatorname{Span}_{\WakDelta^P_{\sx}}(\Ffr_P^0(\mathfrak{n}[t^{-1}]))=\operatorname{Span}_{\WakDelta^P_{\sx}}(\partial_{\alpha,i}| i\leqslant 0).$$ We reduce to proving that the induced homomorphism
$$S(t^{-1}\mathfrak{b}^-[t^{-1}]\oplus\tilde{\mathfrak{n}}_L[t^{-1}])
\rightarrow \bk[y^\alpha_n]_{n< 0}\otimes_{\bk} S(t^{-1}\mathfrak{b}^-_L[t^{-1}])\oplus\tilde{\mathfrak{n}}_L[t^{-1}])$$ is surjective. Now we similarly use (\ref{eq:moment_l}) combined with the description of $\Ffr^0_P$ from Step 4 to reduce the proof of surjectivty to showing that the following homomorphism of quotients is surjective:
\begin{equation}\label{eq:ffr_surjectivity_check}
S(t^{-1}\mathfrak{n}^-[t^{-1}])
\rightarrow \bk[y^\alpha_n]_{n< 0}= S(t^{-1}\mathfrak{n}^-[t^{-1}]).
\end{equation}
Thanks to (\ref{eq:moment_n-}), the map in question sends $w_{-n}$ for $w\in \mathfrak{n}^-_i$ and $n\geqslant 0$ to a polynomial of the form 
$\theta_\eta(w)_{-n}+?$. Here, recall $\sx=(0,0,\eta)$, and ``$?$'' stands for a polynomial in the elements $w'_{-n'}$ for $w'\in \mathfrak{n}^-$
and $n'<n$. Thanks to Step 3, $\theta_\eta$ is a bijection $\mathfrak{n}^-_{i}\rightarrow \mathfrak{n}^-_i$ for all $i$. Now an easy induction on $n$ shows that  (\ref{eq:ffr_surjectivity_check}) holds.
\end{proof}


\section{Endomorphisms of the universal Verma module}\label{S_Verma_endom}
\subsection{Main result}\label{SS:endo_main}
In this section $\Ring$ is a perfect characteristic $p$ field with $p>h$.
We assume $\kappa\neq 0$ for all simple factors of $\g$. Assume also that $G$ is semisimple.

Let $\h^a$ denote the affine subspace in $\h^{da}$ defined by the equation $\hbar=1$. Following Section \ref{SSS_Verma_Wakimoto}, consider the 
module $\Delta^{G}_{\kappa\hbar}:=\Delta^{G,\hbar}_{\bone}/(\bone-\kappa\hbar)\Delta^{G,\hbar}_{\bone}$ over $\Ring$. This is an object of the following category: $\Lambda^a$-graded and $\Iw$-equivariant $\widehat{U}^\hbar_{\kappa}(\hg)$-$\Ring[\h^*,\hbar]$-bimodules. Denote this category by $\sO(\h^*)$.
Our goal is to construct a homomorphism from a certain algebra, to be denoted 
by $\,^I\!\A^\hbar$, to $\End_{\sO(\h^*)}(\Delta^{G}_{\kappa\hbar})$. The construction is based on Theorem \ref{Thm:Verma_vs_Wakimoto}.

\sss Now we construct the algebra $\,^I\!\A$. Here $I$ is the set of simple roots for $\g$. Recall the elements $\tilde{h}_n:=h_n^p-\kappa^{p-1}h_{pn}, h\in \h, n\in \Z$.
Their significance is that for the simple coroot basis $h^i, i\in I,$ of $\h$, the algebra $\widehat{U}_\kappa(\hh)^{\Loop H}$ is the filtered complete algebra of polynomials in the elements $\tilde{h}^i_{n}$, see Proposition \ref{Prop:HC_center_abelian}. Define the algebra $\A$ as the algebra of polynomials in the variables $\tilde{h}^i_{n}$ with $n\leqslant 0$, a quotient of $\widehat{U}_\kappa(\hh)^{\Loop H}$. It carries a filtration, where the degree of $\tilde{h}^i_n$ is $p$. Let $\A^\hbar$ denote the Rees algebra of $\A$. 

We construct $\,^I\!\A^\hbar$ as the intersection $\bigcap_{i\in I}\,^i\!\A^\hbar$ of certain subalgebras $\,^i\!\A^\hbar$ of $\A^\hbar$. The subalgebra $\,^i\!\A\subset \A$ is defined as follows. Let $L_i$ be the minimal Levi in $G$ corresponding to $i\in I$. 

Let $X^i_n$ denote the elements from Proposition 
\ref{Prop:HC_center_sl2_noninteger} for the $\slf_2$-subalgebra $[\lf_i,\lf_i]$. 

{
\begin{Lem}\label{Lem:isogeny_Levi}
The algebra $\widehat{U}_\kappa(\hat{\lf}_i)^{\Loop L_i}$ consists of the converging sums of the elements $\tilde{h}_n$ for $h\in \mathfrak{z}(\mathfrak{l}_i)$ and $X^i_n$ for $n\in \Z$. 
\end{Lem}
\begin{proof}
It is enough to prove that these elements are $\Loop L_i$-invariant. Indeed, then their principal symbols topologically generate $\widehat{S}\left(\Loop \mathfrak{l}_i^{(1)}\right)^{(\Loop L_i)^{(1)}}$ and we are done by (1) of Theorem \ref{Thm:integral_case}.

Set $\tilde{L}_i:=Z(L_i)^\circ \operatorname{SL}_2$. The elements $\tilde{h}_n$
are $\Loop \tilde{L}_i$-invariant by Proposition \ref{Prop:HC_center_abelian} and the elements $X^i_n$ are $\Loop \tilde{L}_i$-invariant by Proposition \ref{Prop:HC_center_sl2_noninteger}. It remains to show any $\Loop \tilde{L}_i$-invariant element in $\widehat{U}(\Loop\mathfrak{l}_i)$ is $\Loop L_i$-invariant.

Notice that $\Jet\tilde{L}_i$ acts on $V_\kappa(\mathfrak{l}_i)$ via an epimorphism to $\Jet L_i$, in particular, $V_\kappa(\mathfrak{l}_i)^{\Jet L_i}=V_\kappa(\mathfrak{l}_i)^{\Jet \tilde{L}_i}$.
Thanks to Lemma \ref{Lem:center_assoc_grad}, it is enough to show that $\gr [V_\kappa(\mathfrak{l}_i)^{\Jet L_i}]=\Ring[(\Jet \mathfrak{l}_i^*)^{(1)}]^{(\Jet L_i)^{(1)}}$. This follows because the principal symbols of the elements $\tilde{h}_n|\varnothing\rangle, X^i_m|\varnothing\rangle$ for $n\leqslant -1,m\leqslant -2$ generate the right hand side finishing the proof.
\end{proof}
}

For $n\leqslant 0$, we set $\underline{X}^i_n$ to be the image of $\HC(X^i_n)$ in $\A$. Explicitly, thanks to Proposition 
\ref{Prop:HC_quadratic}, we have 
\begin{equation}\label{eq:underline_X_in}
\underline{X}^i_n=\frac{1}{4}\left(\sum_{j=n}^0\tilde{h}^i_{j+n}\tilde{h}^i_{-j}\right)-\frac{1}{2}(1-\kappa^{p-1})(n+1)\tilde{h}^i_{n}.
\end{equation}
For $\,^i\!\A$ we take the subalgebra of $\A$ generated by the elements $\tilde{h}_n$ for $n\leqslant 0$ and $h$ in the Lie algebra of the center of $L_i$ as well as the elements $\underline{X}^i_n$. 

Equivalently, by Lemma \ref{Lem:isogeny_Levi},  $\,^i\A$ is the image of $\widehat{U}_\kappa(\hat{\lf}_i)^{\Loop L_i}$ under $\HC_{H,\kappa}^{L_i}$.

And for $\,^i\!\A^\hbar$ we take the Rees algebra of $\,^i\!\A$ with its filtration restricted from $\A$. So $\,^i\!\A^\hbar$ is a graded subalgebra of $\A^\hbar$.

\sss
Note that $\A^\hbar$ naturally acts on $\Wak_{\kappa\hbar}$ by $\widehat{U}^\hbar_{\kappa}(\hg)$-$\Ring[\h^*,\hbar]$-linear endomorphisms. In particular, the subalgebra $\,^I\A^\hbar\subset \A^\hbar$ acts on $\Wak_{\kappa\hbar}$. 

We define the algebra endomorphism $\mathsf{t}_B$ of $\A^\hbar$ as follows: it is the identity on $\tilde{h}_n$ for $n<0$ and coincides with $\mathsf{t}_B$ from Section \ref{SSS_shift} on $\tilde{h}_0$.

The following is the main result of  Section \ref{S_Verma_endom}. 

\begin{Thm}\label{Thm:Verma_endomorphisms}
There is a unique algebra homomorphism $\,^I\!\A^\hbar\rightarrow \operatorname{End}_{\sO(\h^*)}(\Delta^G_{\kappa\hbar})$ such that the homomorphism $\varphi^B_{\kappa\hbar}: \Delta^G_{\kappa\hbar}\rightarrow \Wak_{\kappa\hbar}$ is $\,^I\!\A^\hbar$-semi-linear: for any $m\in \Delta^G_{\kappa\hbar}$ and $a\in \,^I\!\A^\hbar$, we have 
$\varphi^B_{\kappa\hbar}(ma)=\varphi^B_{\kappa\hbar}(m)\mathsf{t}_B(a)$. 
\end{Thm}

\sss
We now explain the steps of the proof. For this we need to introduce some notation. Let $P_i$ be the minimal parabolic containing $B$ with Levi $L_i$. We consider the hyperplane $\h^{da}_{\kappa\hbar}$ given by equation $\bone=\kappa\hbar$. Let $\Phi^+$ denote the set of positive roots for $B^-$. Let $\Sigma\subset \Phi^+$ be a subset. In Section \ref{SSS_KK_hyperplanes} we have introduced hyperplanes $\tilde{H}_{\beta,m}$.
We set
\begin{equation}
\h^{da}_{\kappa\hbar}(\Sigma)=\h^{da}_{\kappa\hbar}\setminus\bigcup_{\alpha\in \Phi^+\setminus\Sigma, m,n\in \Z}\widetilde{H}_{\alpha+n\delta,m}
\end{equation}
We will write $\h^{da}_{\kappa\hbar}(i)$ for 
$\h^{da}_{\kappa\hbar}(\{\alpha_i\})$.
We write $\Delta^{G,\Sigma}_{\kappa\hbar}$ for $\Delta^{G}_{\kappa\hbar}\otimes_{\Ring[\h^{da}_{\kappa\hbar}]}
\Ring[\h^{da}_{\kappa\hbar}(\Sigma)]$. The notation 
$\Wak^{P,\Sigma}_{\kappa\hbar}$ has a similar meaning. We write $\Delta_{\kappa\hbar}^{G,i}$ instead of $\Delta_{\kappa\hbar}^{G,\{i\}}$.

Thanks to Theorem \ref{Thm:Verma_vs_Wakimoto} we have the following isomorphisms
\begin{align}\label{eq:isom_Delta_Wak_gen}
& \varphi^B_{\kappa\hbar}: \Delta^{G,\varnothing}_{\kappa\hbar}\xrightarrow{\sim}  \Wak^{B,\varnothing}_{\kappa\hbar},\\\label{eq:isom_Delta_WakDelta_igen}
& \varphi^{P_i}_{\kappa\hbar}: \Delta^{G,i}_{\kappa\hbar}\xrightarrow{\sim}  \WakDelta^{P_i,i}_{\kappa\hbar}.
\end{align}

Using this, we will construct an action of $\,^I\!\A^\hbar$ on $\Delta^{G,I}_{\kappa\hbar}$
by $\widehat{U}^\hbar_\kappa(\hg)\otimes \Ring[\h^{da}_{\kappa\hbar}(I)]$-linear (this is relatively easy)
and $\Iw$-equivariant (this is harder) endomorphisms. Then we will show that this action preserves the $\Ring[\h^{da}_{\kappa\hbar}]$-lattice $\Delta^{G}_{\kappa\hbar}$ using intertwining functors from Section \ref{SS_intertwining}. 

\sss 
Note that $\,^i\!\A^\hbar$
acts on $\WakDelta_{\kappa\hbar}^{P_i}$ by endomorphisms, and $\,^I\!\A^\hbar\hookrightarrow \,^i\!\A^\hbar$. 
 
\begin{Cor}\label{Cor:Verma_endomorphisms}
The homomorphism $\varphi_{\kappa\hbar}^{P_i}:\Delta^{G}_{\kappa\hbar}\rightarrow 
\WakDelta^{P_i}_{\kappa\hbar}$ is $\,^I\!\A^\hbar$-semi-linear. 
\end{Cor}
\begin{proof}
Apply Theorem \ref{Thm:Verma_endomorphisms} to $L_i$ (instead of $G$) and $B\cap L_i$. We get that $\varphi_{\kappa\hbar}^{B\cap L_i}: \Delta^{L_i}_{\kappa\hbar}\rightarrow \Wak^{B\cap L_i}_{\kappa\hbar}$ is 
$\,^i\!\A^\hbar$-semi-linear. From here it follows that  $\Wak^P(\varphi_{\kappa\hbar}^{B\cap L_i}): \WakDelta^{P_i}_{\kappa\hbar}\rightarrow \Wak^B_{\kappa\hbar}$ is $\,^i\!\A^\hbar$-semi-linear. Note that the modules $\WakDelta^{P}_{\kappa\hbar}$ are flat over $\Ring[\h^{da}]/(\bone-\kappa \hbar)$.  So it is enough to show $\varphi_{\kappa\hbar}^{P_i}$ is $\,^I\!\A^\hbar$-semi-linear after localizing to the generic point of $\operatorname{Spec}(\Ring[\h^{da}]/(\bone-\kappa \hbar))$ (to be denoted by $\sx$) that is $B$-generic in the sense of Definition \ref{defi:P_generic}. Hence it is enough to prove the claim of the corollary after changing the base to $\sx$. By Theorem \ref{Thm:Verma_vs_Wakimoto}, the homomorphisms $\varphi_{\sx}^P$ are isomorphisms. Hence the semi-linearity of $\varphi_{\kappa\hbar}^{P_i}$ follows from (\ref{eq:varphi_transitivity}) -- note that the shifts match. 
\end{proof}

\subsection{Action on $\Delta^{G,I}$}
Recall, Lemma \ref{Lem:invariant_containment}, that 
$$[\widehat{D}(\Loop N)\widehat{\otimes}\widehat{U}_\kappa(\hat{\lf})]^{\Loop P}\xrightarrow{\sim}\widehat{U}_\kappa(\hat{\lf})^{\Loop L}.$$
This is a filtered algebra isomorphism. 
Applying this to $P=P_i$, we see that $\widehat{U}_{\kappa\hbar}(\hat{\lf}_i)^{\Loop L_i}$ acts on any $\Iw_{P_i}$-equivariant module over 
$\widehat{D}^\hbar(\Loop N_i)\widehat{\otimes}_{\Ring[\hbar]}\widehat{U}^\hbar_{\kappa}(\hat{\lf}_i)$, in particular, on $\WakDelta^{P_i}_{\kappa\hbar}$. The action is by $\Ring[\h^{da}_{\kappa\hbar}]$-linear automorphisms, so extends to $\WakDelta^{P_i,i}_{\kappa\hbar}$. Note that the resulting action factors through $\,^i\!\A^\hbar$. 

In particular, we can carry the action of $\,^i\!\A^\hbar$ to $\Delta^{G,i}_{\kappa\hbar}$ using isomorphism (\ref{eq:isom_Delta_WakDelta_igen}). 
Note that $\Delta^{G,I}_{\kappa\hbar}$
coincides with $\bigcap_{i\in I}\Delta_{\kappa\hbar}^{G,i}$, where the intersection is taken in $\Delta^{G,\varnothing}_{\kappa\hbar}$. It follows that $\,^I\!\A^\hbar$ acts on $\Delta^{G,I}_{\kappa\hbar}$. The action is by 
$\Ring[\h^{da}_{\kappa\hbar}(I)]$-linear endomorphisms that are equivariant for the action of $\Iw_{P_i}$ for all $i\in I$. Note that, thanks to (\ref{eq:varphi_transitivity}), the condition that $\varphi^B:\Delta^{G,\varnothing}_{\kappa\hbar}\rightarrow \Wak^{B,\varnothing}_{\kappa\hbar}$ is $\,^I\!\A^\hbar$-semi-linear and the flatness of $\Delta^G_{\kappa\hbar}$ over $\Ring[\h^{da}]/(\bone-\kappa \hbar)$ recovers the action of $\,^I\!\A^\hbar$ on $\Delta^{G,I}_{\kappa\hbar}$ uniquely. 

\begin{Lem}\label{Lem:Iw_equivariance}
The action of $\,^I\!\A^\hbar$ on  $\Delta^{G,I}_{\kappa\hbar}$ is by $\Iw$-equivariant endomorphisms.
\end{Lem}
\begin{proof}

Note that $\Iw$ is a pro-algebraic group, and its representation in 
$\Delta^{G,I}_{\kappa\hbar}$ is rational. For $m\geqslant 0$, let $\Iw_m$
denote the preimage of $B^-$ in $\Jet_m G$. Define $\Iw_{P_i,m}$ as the intersection of $\Iw_m$ with $\Jet_m P_i$.
So $\Iw_{P_i,m}\subset \Iw_m$ are algebraic groups, they are connected, and $\Iw=\varprojlim_{m}\Iw_m, \Iw_{P_i}=\varprojlim_{m}\Iw_{P_i,m}$. Note that $\Delta_{\kappa\hbar}^{G,I}$ is a cyclic module over $\widehat{U}^\hbar_\kappa(\hg)\otimes \Ring[\h^{da}_{\kappa\hbar}(I)]$ generated by an $\mathsf{Iw}$-invariant vector. So, the claim of the lemma boils down to showing that any element in a rational representation of $\Iw_m$ that is invariant under $\Iw_{P_i,m}$ for all $i$, is also $\Iw_m$-invariant. This will follow if we show that the subalgebras $\operatorname{Dist}_1(\Iw_{P_i,m})$ for $i\in I$ generate the algebra $\operatorname{Dist}_1(\Iw_m)$.

Note that $\Jet_m \tilde{N}^-\subset \Iw_m$. Note that $\Iw_m\cap \Jet_m B$ is the semidirect product of $\Jet_m T$ and the first congruence subgroup in $\Jet_m \tilde{N}$. Since $\Iw_m$ lies in the open Bruhat cell of $\Jet_m G$, we see that
the product map $\Jet_m \tilde{N}^-\times (\Iw_m\cap \Jet_m B)\rightarrow \Iw_m$ is an isomorphism. 

Since $\Iw_m\cap \Jet_m B\subset \Iw_{P_i,m}$
for all $i$, it suffices to show that the subalgebras $\operatorname{Dist}_1(\Jet_m \tilde{N}_i^-)$ for $i\in I$ generate 
$\operatorname{Dist}_1(\Jet_m \tilde{N}^-)$, where $\widetilde{N}_i:=L_i\cap \widetilde{N}$. First of all, the claim on the level of Lie algebras is clear. We reduce to showing that the subalgebras 
$\operatorname{Dist}_1([\Jet_m \tilde{N}_i^-]^{(1)})$ generate $\operatorname{Dist}_1([\Jet_m \tilde{N}^-]^{(1)})$. 
This follows once we observe that the algebraic groups $\Jet_m \tilde{N}_i^-, \Jet_m\tilde{N}^-$ are defined over $\F_p$ as are the embeddings  
$\Jet_m \tilde{N}_i^-\hookrightarrow \Jet_m \tilde{N}^-$, while the Frobenius homomorphisms $\Jet_m \tilde{N}_i^-\rightarrow (\Jet_m \tilde{N}_i^-)^{(1)}$
and $\Jet_m \tilde{N}^-\rightarrow (\Jet_m \tilde{N}^-)^{(1)}$ are epimorphisms (because these algebraic groups are smooth over $\operatorname{Spec}(\Ring)$).
\end{proof}

Our next and final task in proving Theorem \ref{Thm:Verma_endomorphisms} is to show that the action of $\,^I\!\A^\hbar$ on $\Delta^{G,I}_{\kappa\hbar}$ preserves $\Delta^G_{\kappa\hbar}$. In order to show this (by using intertwining functors) we need to digress and discuss intertwining functors for suitable categories of equivariant modules. 

\subsection{Equivariant modules}
Here we recall some basics on weak and strong actions of algebraic groups and the corresponding categories of equivariant modules. In this section, all algebraic groups are taken to be smooth.

\sss  Let $K$ be an algebraic group over a field $\bk$, and denote its Lie algebra by $\fk$.
Let $U^{\hbar}(\fk)$ denote the Rees algebra of the universal enveloping algebra $U(\fk)$
so that for $x,y\in \fk\subset U^\hbar(\fk)$ we have
$$\xi \cdot \eta - \eta \cdot \xi = \hbar \cdot [\xi,\eta], \forall \xi,\eta \in \fk.$$

\sss For later use, let us observe that if $\operatorname{char}\bk=p$, and $\fk$ is equipped with a restricted Lie algebra structure, which we denote by $\xi \mapsto \xi^{[p]}$, then this equips $U^\hbar(\fk)$ with a $p$-center. Namely, we have a  $\bk$-linear map  
$$\fk^{(1)} \rightarrow U^\hbar(\fk),  \xi \mapsto \xi^p - \hbar^{p-1} \xi^{[p]},$$
whose image consists of central elements. The image freely generates a central  $\bk[\hbar]$-subalgebra $$Z_p^\hbar(\fk) \hookrightarrow U^\hbar(\fk)$$
whose spectrum, if we write $\mathbb{A}^1_\hbar := \Spec \bk[\hbar]$, is therefore given by $$(\fk^*)^{(1)} \times \mathbb{A}^1_\hbar.$$

\sss\label{SSS_quantum_comoment} 
{ We write $\OO_K$ for the algebra of regular functions on $K$, and $\on{Rep}(K)$ for the category of $\OO_K$-comodules in $\bk[\hbar]\lmod$. Consider a  graded complete topological algebra $$A = \underset{i \geqslant 0} \oplus \hspace{.5mm} A_{i}$$over $\bk[\hbar]$, where $\hbar$ has degree 1, cf. Definition \ref{defi:filtered_complete}. 

\begin{defi} The datum of a {\em weak action} of $K$ on $A$ consists of 
\begin{enumerate}
\item A continuous $\bk[\hbar]$-linear action $\alpha$ of $K$ on each $A_{ i}$, i.e., a lift of each $A_i$ from $\Pro(\bk[\hbar]\lmod)$ to $\Pro(\Rep(K))$, so that $K$ acts trivially on $\bk[\hbar] \subset A$.

In particular each $A_{ i}$ admits a neighborhood basis of open subspaces $A_{ i}^\alpha \subset A_{ i}$ which are $K$-invariant, and the induced action of $K$ on the discrete $\bk$-module $A_{ i}/A_{ i}^\alpha$ is rational, compatibly as one varies $\alpha$.

\item For all $i,j \geqslant 0$, the multiplication map $A_{ i} \otimes A_{ j} \rightarrow A_{ i + j}$ is $K$-equivariant. I.e., for any $A_{ i + j}^\alpha$ as above, consider the obtained map 
$$A_{ i} \otimes A_{ j}/A_{ j}^\beta \rightarrow A_{ i + j}/A_{ i + j}^\alpha$$for some $\beta$, and for any  $K$-invariant subspace $S^\beta \subset A_{ j}/A_{ j}^\beta$ which is finite dimensional as a $\bk$-module consider the induced map 
\begin{equation} \label{gumgum}A_{ i}/A_{ i}^\gamma \otimes S^\beta \rightarrow A_{ i + j}/A_{ i + j}^\alpha\end{equation}
for some $\gamma$. Then we ask that the maps \eqref{gumgum} are maps in $\Rep(K)$, for all $\alpha, \beta, S^\beta$, and $\gamma$ as above.  
\end{enumerate}
\end{defi}

\begin{defi} The datum of a {\em strong action} of $K$ on $A$ is a weak action of $K$ on $A$ along with a homomorphism of graded $\bk[\hbar]$-algebras (known as the {\it quantum comoment map})
$$\phi: U^\hbar(\fk) \rightarrow A$$
satisfying the following condition: for any $\xi \in \fk$, if we write $d\alpha_\xi: A \rightarrow A$ for the $\bk[\hbar]$-linear endomorphism of $A$ obtained by differentiating the action of $K$, then we have the equality$$\hbar \cdot d\alpha_\xi(-) = [\phi(\xi), -]: A \rightarrow A.$$
\end{defi}

Graded complete algebras equipped with strong actions $(A, \alpha, \phi)$ form a category, where the morphisms are graded continuous $K$-equivariant homomorphisms intertwining the quantum comoment maps. }

\sss An example of primary interest for us, concerning a discrete case of the above, is the following.
\begin{Ex}  \label{ex:variety} If $X$ is a smooth affine variety over $\bk$ equipped with an action of $K$,  write $D^\hbar(X)$ for its Rees algebra of differential operators. Then an action of $K$ on $X$ equips     $D^\hbar(X)$ with a strong action of $K$: for $\xi\in \fk$ we set $\phi(\xi)=\xi_X$, where $\xi_X$ is the  vector field defined by $\xi\in \fk$.

Now assume that $\operatorname{char}\bk=p$.
We recall that $D^\hbar(X)$ has a $p$-center $Z^\hbar_p(X) \hookrightarrow D^\hbar(X)$. Namely, for each vector field $\zeta$ on $X$, with restricted power $\zeta^{[p]}$ (again, a vector field on $X$), one has the corresponding central element$$\zeta^p - \hbar^{p-1} \zeta^{[p]} \in D^\hbar(X),$$and these, along with $p^{th}$ powers of functions, generate the $p$-center, so that canonically one has
$$\Spec Z^\hbar_p(X) \simeq (T^* X)^{(1)} \times \mathbb{A}^1_\hbar.$$
We recall that in this case, as the map from $\fk$ to vector fields on $X$ is a morphism of restricted Lie algebras, the quantum comoment map $\phi: U^\hbar(\fk) \rightarrow D^\hbar(\fk)$ restricts to a map $\underline{\phi}: Z^\hbar_p(\fk) \rightarrow Z^\hbar_p(X)$. Moreover, on spectra, if we write $\mu: T^* X \rightarrow \fk^*$ for the moment map, $\underline{\phi}$ corresponds to 
$$\mu^{(1)} \times \on{id}: (T^*X)^{(1)} \times \mathbb{A}^1_\hbar \rightarrow (\fk^*)^{(1)} \times \mathbb{A}^1_\hbar.$$\end{Ex}

{ \sss Given an algebra $A$ as in Section \ref{SSS_quantum_comoment}, we recall that a module $M$ is {\em discrete} if, for any finitely generated $\bk[\hbar]$-submodule $S \subset M$ and any $i \geqslant 0$, the action map $A_{ i} \otimes S \rightarrow M$ factors through 
\begin{equation} \label{discrete} A_{ i}/A_{ i}^\alpha \otimes S \rightarrow M \end{equation}
for some $\alpha$. Let us denote by $A\lmod$ the category of discrete $A$-modules, which is an abelian category.

\begin{defi} A datum of {\em weak equivariance} on an $A$-module $M$ is a $\bk[\hbar]$-linear action of $K$ on $M$
$$\beta: M \rightarrow \OO_K \underset {\bk} \otimes M$$
which is compatible with the action of $A$ in the following sense: for any $K$-invariant subspace $S \subset M$ finitely generated over $\bk[\hbar]$, and any $i \geqslant 0$, any factorization of the form \eqref{discrete}, where $A^\alpha_{ i}$ is $K$-invariant, is a map of $K$-representations.    
\end{defi}

 We denote the usual category of weakly equivariant $A$-modules by $A\lmod^{(K,w)}$; this is again an abelian category. Note that if we write $Z(A)$ for the center of $A$, then $A\lmod^{(K,w)}$ is naturally a $Z(A)^K$-linear category.

\begin{defi} A datum of {\em strong equivariance} on $M$ is a datum of weak equivariance satisfying the further condition that, for any $\xi \in \fk$, if we write $d\beta_\xi: M \rightarrow M$ for the endomorphism of $M$ obtained by differentiating the action of $K$ on $M$, then we have an equality of endomorphisms
$$\hbar \cdot d\beta_\xi(-) = \phi(\xi) \cdot -: M \rightarrow M.$$
\end{defi}
The category of strongly equivariant modules is denoted by $A\lmod^K$, and is again abelian. }
We have the forgetful functor $$\on{Oblv}: A\lmod^{K} \rightarrow A\lmod^{(K,w)}.$$
and note that both $A\lmod^K$ and the functor are again naturally $Z(A)^K$-linear. 

\sss The obstruction to a weakly equivariant object being strongly equivariant implies a further linearity for weakly equivariant objects, as follows. 

{ First, note that if $A^1$ and $A^2$ are graded complete topological algebras,  the graded complete tensor product $$A^1 \widehat{\otimes}_{\bk[\hbar]} A^2$$
(defined as $\bigoplus_{i\geqslant 0}(\bigoplus_{j=0}^i A^1_j\widetilde{\otimes}_{\bk[\hbar]}A^2_{j-i})$, where $\widetilde{\otimes}$ stands for the completed tensor product of topological $\bk[\hbar]$-modules)
is naturally a graded complete algebra. Moreover, if $A^1$ and $A^2$ both have weak or strong actions of $K$, their tensor product inherits a natural weak or strong action, respectively.

Note that e.g. if $A^2$ moreover is discrete, and each $A^2_i$ is a free module over finite rank over $\bk[\hbar]$, the graded complete tensor product agrees with the naive tensor product, i.e., 
$$ A^1 \widehat{\otimes}_{\bk[\hbar]} A^2 \simeq A^1 \otimes_{\bk[\hbar]} A^2.$$ 
}

In particular, consider the tautological map $$A \rightarrow A  \otimes_{\bk[\hbar]} U^{\hbar}(\fk), \quad \quad a \mapsto a \otimes 1,$$
and the associated restriction functor $A \otimes_{\bk[\hbar]} U^\hbar(\fk) \lmod \rightarrow A\lmod.$

\begin{Lem}\label{l:bigbird} Restriction to the factor $A$ defines an equivalence
$$A \otimes_{\bk[\hbar]} U^\hbar(\fk) \lmod^K \xrightarrow{\sim} A\lmod^{(K,w)}.$$
In particular, $A\lmod^{(K,w)}$ is naturally $U^{\hbar}(\fk)^K$-linear. 
\end{Lem}

\begin{proof} A quasi-inverse functor is constructed as follows. For $(M, \beta)$ an object of $A\lmod^{(K,w)}$ the formula 
$$\xi \mapsto  \hbar \cdot d\beta_\xi - \phi(\xi)$$
defines an action of $U^{\hbar}(\fk)$ on $M$ commuting with the action of $A$, and tautologically $(M, \beta)$ is strongly equivariant with respect to this extended action. To check that this is indeed a quasi-inverse functor is left as an exercise.     
\end{proof}

\sss\label{SSS_subgroup_equivariance} 
Note that if we have a homomorphism of groups $S \rightarrow K$, then the strong action of $K$ on $A$ induces a strong action of $S$. Moreover, let us fix $\underline{S}$ normal in $S$, with quotient $Q = S/\underline{S}$, and write the Lie algebra of the latter as $\mathfrak{q}$. Then, in evident notation, we can consider the intermediate category of modules which are only strongly equivariant for $\underline{S}$, 
$$A\lmod^S \hookrightarrow A\lmod^{\underline{S}, (Q,w)} \hookrightarrow A\lmod^{(S,w)},$$
and this intermediate category is naturally linear over $U^{\hbar}(\mathfrak{q})^{Q}$. 

{ Explicitly, let us write $\underline{\mathfrak{s}} \hookrightarrow \mathfrak{s}$ for the Lie algebras of $\underline{S} \hookrightarrow S$. Then given an object $(M, \beta)$ of $A\lmod^{\underline{S}, (Q,w)}$, consider the action of $U^\hbar(\mathfrak{s})$ on $M$ obtained as in the proof of Lemma \ref{l:bigbird}. By the strong equivariance of $M$ with respect to $\underline{S}$, this action is trivial on the two sided ideal generated by $\underline{\mathfrak{s}}$, which we denote by $(\underline{\mathfrak{s}}) \subset U^\hbar(\mathfrak{s})$, hence factors through an action of $$U^\hbar(\mathfrak{q}) \simeq U^\hbar(\mathfrak{s}) / (\underline{\mathfrak{s}});$$by construction the action of the subalgebra $U^\hbar(\mathfrak{q})^Q = U^\hbar(\mathfrak{q})^S$ commutes with both the action of $A$ and the coaction of $S$, as desired. 
}

\sss The following standard fact will be useful to us. The natural action of $K \times K$ on $K$ by left and right translations yields a strong action of $K \times K$ on $D^\hbar(K)$. Let us equip $$A \widehat{\otimes}_{\bk[\hbar]} D^\hbar(K)$$with a strong $K^\ell \times K^r := K \times K$ action, where the action of $K^\ell$ is the tensor product of the action on $A$ and the left translation action on $D^\hbar(K)$, and the action of $K^r$ is the right translation action on $D^\hbar(K)$. { I.e., if we write $\mu^\ell, \mu^r$ for the left and right translation actions of $K$ on $D^\hbar(K)$,  $\alpha$ for the action of $K$ on $A$, and $\beta^\ell$ and $\beta^r$ for the presently defined actions of $K^\ell$ and $K^r$, we have 
$$\beta^\ell = \alpha \otimes \mu^\ell \quad \text{and} \quad \beta^r = \on{id}_A \otimes \mu^r.$$
}

\begin{Lem} \label{l:dividebyK}There is a canonical equivalence, functorial in $A$
$$A\lmod \simeq (A \widehat{\otimes}_{\bk[\hbar]}D^\hbar(K))\lmod^{K^\ell}.$$    
Moreover, for any $\underline{S}$ and $Q$ as in Section \ref{SSS_subgroup_equivariance}, this induces an equivalence 
$$A\lmod^{\underline{S}, (Q,w)} \simeq (A \widehat{\otimes}_{\bk[\hbar]} D^\hbar(K))\lmod^{K^\ell \times \underline{S}^r, (K^\ell\times Q^r, w)}.$$
\end{Lem}

{
\begin{Rem} Let us make an orienting comment for the reader. The proof of the lemma, in the ostensibly special case of Example \ref{ex:variety}, amounts to considering the change of variables 
$$X \times K \simeq X \times K, \quad \quad (x, k) \mapsto (k^{-1}x, k);$$
the general case is a straightforward generalization of this.  \end{Rem}

\begin{proof} We break the proof into several steps for ease of reading. 

{\it Step 1.} We will construct below an automorphism of graded complete algebras
$$\Xi: A \wotimes D^\hbar(K) \xrightarrow{\sim} A \wotimes D^\hbar(K)$$
satisfying the following properties. 

\begin{enumerate}
\item The strong action of $K^\ell$ is exchanged with the tensor product of the trivial action of $K$ on $A$ and the left translation action of $K$ on $D^\hbar(K)$, i.e., in evident notation
$$(\alpha \otimes \mu^\ell) \circ \Xi   =   \Xi \circ (\on{id}_A \otimes \mu^\ell).$$

\item The strong action of $K^r$ is exchanged with the tensor product of the action $\alpha$ of $K$ on $A$ and the right translation action of $K$ on $D^\hbar(K)$, i.e., 
$$(\on{id}_A \otimes \mu^r) \circ \Xi  =  \Xi \circ (\alpha \otimes \mu^r).$$
    
\end{enumerate}
To construct $\Xi$, first note that as a graded $\bk[\hbar]$-module we have $D^\hbar(K) \simeq \bk[K] \otimes U^\hbar(\fk)$, where $U^\hbar(\fk) \subset D^\hbar(K)$ is the subalgebra of left translation invariants. As a consequence, as graded complete $\bk[\hbar]$-modules we have 
$$A \wotimes D^\hbar(K) \simeq (A \widehat{\otimes} \Ring[K]) \wotimes U^\hbar(\fk);$$
in particular, $\Xi$ is determined by its restrictions to $A \widehat{\otimes} \bk[K]$ and $U^\hbar(\fk)$, which we denote by $\Xi_1$ and $\Xi_2$. Consider the composition 
$$A \widehat{\otimes} \bk[K] \xrightarrow{\on{coact} \widehat{\otimes} \on{id}} \bk[K] \widehat{\otimes} A \widehat{\otimes} \bk[K]  \simeq A \widehat{\otimes} \bk[K] \widehat{\otimes} \bk[K] \simeq A \widehat{\otimes} (\bk[K] \otimes \bk[K]) \xrightarrow{ \on{id} \widehat{\otimes} \Delta} A \cotimes \bk[K] \hookrightarrow A \wotimes D^\hbar(K)$$
and the quantum comoment map $U^\hbar(\fk) \rightarrow A \wotimes D^\hbar(K)$ associated to  $\alpha \otimes \mu^\ell$. It is then a standard exercise to verify that there is a unique $\Xi$ for which these are $\Xi_1$ and $\Xi_2$, respectively, and that this satisfies the claimed properties (1) and (2) - to verify $\Xi$ is an isomorphism one may write down its inverse by similar means.

{\em Step 2.} The automorphism $\Xi$ constructed in the previous step yields an equivalence 
$$(A \wotimes D^\hbar(K))\lmod^{K^\ell} \simeq (A \wotimes D^\hbar(K))\lmod^{K} $$
where the $K$-equivariance on the right-hand side is imposed with respect to the strong action $\on{id}_A \otimes \mu^\ell$. In particular, an object of the latter category consists of a $D^\hbar$-module on $K$, strongly equivariant with respect to left translations, equipped with a commuting action of $A$, i.e., 
$$(A \wotimes D^\hbar(K))\lmod^K \simeq A\lmod(D^\hbar(K)\lmod^K).$$
However, the category $D^\hbar(K)\lmod^K$ is well known to be simply $\bk[\hbar]\lmod$, with generating object $\bk[K] \otimes \bk[\hbar]$, so we obtain 
$$A\lmod(D^\hbar(K)\lmod^K) \simeq A\lmod(\bk[\hbar]\lmod) \simeq A\lmod,$$
as desired. 
\end{proof}
}

\subsection{Localization}

\sss We now specialize the preceding discussion to $K = M$, a split connected reductive group, and $\bk$ a field of characteristic $p > h$, where $h$ denotes the maximum of the Coxeter numbers of all the almost simple normal subgroups of $M$. We will need a $\bk[\hbar]$-linear, families version of the localization theorem in positive characteristic, cf. \cite{BMR}, albeit only for suitably generic infinitesimal characters. We were unable to locate the needed formulation in the literature, so we now set it up.

\sss Fix a Borel subgroup with Levi quotient
$$T \leftarrow B \rightarrow M,$$
i.e., if we write $N$ for the unipotent radical of $B$, $T$ is the universal Cartan $B/N$. Write  {$\mathfrak{t}, \mathfrak{b}, \mathfrak{n}$ for the Lie algebras of $T$, $B$, and $N$, respectively,} and $\check{\Phi} \subset \mathfrak{t}$ for the coroots. 

\sss For an algebra $A$ with a strong action of $M$, we can consider the category of modules
$$A\lmod^{N, (T,w)},$$
which is naturally linear over $U^\hbar(\mathfrak{t})^T = U^\hbar(\mathfrak{t}).$ Below, we will refer to objects of this category as `weak highest weight' modules, as we relax the strong $T$-equivariance in the definition of highest weight modules to merely weak $T$-equivariance. 

\begin{defi} \label{d:genmodule}Let us call a $U^\hbar(\ft)$-module $\cM$ {\em generic} if the elements $$\halpha - n\hbar, \quad \text{for } \halpha \in \check{\Phi}, n \in \mathbb{F}_p,$$act invertibly on $\cM$.
\end{defi} 

Denote the localization of $U^\hbar(\ft)$ by all the elements $\check{\alpha}-n\hbar$ by $U^\hbar(\ft)^\circ$, and the corresponding full subcategory of generic objects by  
$$j_*: A\lmod^{N,(T,w),\circ} \hookrightarrow A\lmod^{N, (T,w) };$$ 
note the inclusion $j_*$ admits a left adjoint $j^*$, which explicitly is given by base change to $U_\hbar(\ft)^\circ$ . 

\sss For later use, let us note that one can detect the genericity of a module $\cM$ using only the action of the $p$-center of $U^\hbar(\ft)$. Namely, recall the canonical identification $$\Spec Z_p^\hbar(\ft) \simeq (\ft^{*})^{(1)} \times \mathbb{A}^1_\hbar.$$In particular, we may consider the open subset
$$(\ft^*_{rs})^{(1)} \times \mathbb{A}^1_\hbar \hookrightarrow (\ft^*)^{(1)} \times \mathbb{A}^1_\hbar,$$
where $\ft^*_{rs}$ is the regular semisimple locus, i.e., the complement of the root hyperplanes. Denote the corresponding localization of $Z^\hbar_p(\ft)$ by $Z^\hbar_p(\ft)^\circ$. 

\begin{Lem}\label{l:genpcent} A module over $U^\hbar(\ft)$ is generic if and only if the action of $Z_p^\hbar(\ft)$ factors through $Z_p^\hbar(\ft)^\circ$. Equivalently, we have a canonical isomorphism of algebras
$$U^\hbar(\ft)^\circ \simeq U^\hbar(\ft) \underset{Z^\hbar_p(\ft)} \otimes Z^\hbar_p(\ft)^\circ.$$
\end{Lem}

\begin{proof} Note the identity, for each $\halpha \in \check{\Phi}$, 
$$\underset{n \in \mathbb{F}_p}\Pi \hspace{.5mm} (\halpha - \hbar \cdot n) = \halpha^p - \hbar^{p-1} \cdot \halpha;$$
the same identity holds with $\halpha$ and $\hbar$ replaced by any two commuting elements in an $\mathbb{F}_p$-algebra. In particular, inverting the product on the left-hand side is equivalent to inverting one $p$-central element, namely the right-hand side. 
\end{proof}

\sss Note that, as for any torus, if we write $\Lambda$ for the character lattice of $T$, we have an equivalence of cocomplete abelian categories$$U^\hbar(\ft)\lmod^{(T,w)}  \simeq \underset{\lambda \in \Lambda}\oplus \hspace{.5mm} U^\hbar(\ft)\lmod.$$
Explicitly, this corresponds to the set of projective generators $U^\hbar(\ft)_\lambda$, $\lambda \in \Lambda$, where $U^\hbar(\ft)_\lambda$ denotes the Rees enveloping algebra with an action of $T$ purely by the character $\lambda$. In particular, we have
$$\on{Hom}(U^\hbar(\ft)_\lambda, U^\hbar(\ft)_\lambda) \simeq U^\hbar(\ft), \quad \quad \lambda \in \Lambda.$$
Set $U^\hbar(\ft)^\circ_\lambda := j^* U^\hbar(\ft)_\lambda$. The objects $U^\hbar(\ft)^\circ_\lambda$ compactly generate the generic category, and we again have 
$$U^\hbar(\ft)\lmod^{(T,w), \circ} \simeq  \underset{\lambda \in \Lambda}\oplus \hspace{.5mm} U^\hbar(\ft)^\circ\lmod;$$
we recall in passing that the analogous claim holds for the category of weakly equivariant modules over any algebra with trivial $T$-action.

\sss \label{s:hch} Consider the parabolic induction functor
$$\on{pind}: U^\hbar(\ft)\lmod^{(T,w)} \rightarrow U^\hbar(\fm)\lmod^{N, (T,w)}, \quad \quad \cM \mapsto U^\hbar(\fm) \underset{U^\hbar(\fb)} \otimes \cM.$$
We note this functor is $U^\hbar(\ft)$-linear, and in particular restricts to a functor
\begin{equation}
    \label{e:pindg} \on{pind}:  U^\hbar(\ft)\lmod^{(T,w), \circ} \rightarrow U^\hbar(\fm)\lmod^{N, (T,w), \circ} 
\end{equation}
In particular, we have the generic universal Verma module with weak highest weight $\lambda$:
\begin{equation}\label{eq:universal_Verma_defn}
\Delta_\lambda^\circ := \on{pind}(U^\hbar(\ft)_\lambda^\circ)\in U^\hbar(\fm)\lmod^{N,(T,w), \circ},  \quad \lambda \in \Lambda.
\end{equation}

\begin{Rem} In fact, one may show that \eqref{e:pindg} is equivalence of categories; we do not need this, so omit the straightforward proof. 
\end{Rem}

{ \sss  Let us specialize to $\lambda = 0$; similar results hold for any $\lambda \in \Lambda$.} Note that the negation map $X \mapsto -X, X \in \ft$, is tautologically an isomorphism of Lie algebras and hence prolongs to an automorphism $$\on{neg}: U^\hbar(\ft) \rightarrow U^\hbar(\ft).$$
Consider the composition \begin{equation} \label{signsigns} U^\hbar(\ft) \xrightarrow{\on{neg}} U^\hbar(\ft) \rightarrow \End(\Delta_0) \end{equation} {where the second arrow arising from the $N, (T, w)$-equivariance of $\Delta_0$. This composition is an isomorphism of algebras, as follows by noting that any endomorphism of $\Delta_0$ must in particular preserve its subspace of $B$-invariants, but this is precisely the image of the tautological unit map $$U^\hbar(\ft) \hookrightarrow \Delta_0 = \on{ind}_{\mathfrak{b}}^{\mathfrak{m}}(U^\hbar(\ft)),$$i.e., the highest weight space.

Let us unwind an explicit description of how $U^\hbar(\ft)$ acts on the highest weight space $(\Delta_0)^B$ via \eqref{signsigns}. On the one hand, the action of $U^\hbar(\fm)$ on $\Delta_0$ induces, via restriction, an action of $U^\hbar(\fb)$ on $\Delta_0$. This by construction preserves the highest weight space, where it factors through an action of $U^\hbar(\fb/\mathfrak{n}) \simeq U^\hbar(\ft)$. For $\check{\lambda} \in \ft$, let us write 
$$\check{\lambda} \overset{\ell} \cdot -:  (\Delta_0)^B \rightarrow (\Delta_0)^B$$
for this induced action. On the other hand, via \eqref{signsigns}, to $\check{\lambda}$ we also obtain a self-map of the highest weight space, which we denote by 
$$\check{\lambda} \overset{r} \cdot -: (\Delta_0)^B \rightarrow (\Delta_0)^B.$$As the difference of these actions integrates to the trivial action of $B$, where the difference rather than the sum is due to the presence of `neg' in \eqref{signsigns}, these two actions agree, i.e. 
$$\check{\lambda} \overset{\ell} \cdot - = \check{\lambda} \overset{r} \cdot -.$$

\sss \label{sss:hcstar} As a consequence, the induced map 
\begin{equation} \label{hcsign} U^\hbar(\fm)^M \rightarrow \End(\Delta_0) \overset{\eqref{signsigns}}\simeq  U^\hbar(\ft)\end{equation}
agrees with the homogenized Harish--Chandra homomorphism.  

Explicitly, if we denote the usual linear action of $W$ on $\ft$ by $(y, h) \mapsto y(h)$, and we write $\rho \in \ft^*$ for one half of the image in $\ft^*$ of the sum of the positive roots, the dot action of $W$ on generators $\ft\subset U^\hbar(\ft)$  is given by 
$$\check{\lambda} \mapsto y \overset{\hbar}\cdot \check{\lambda} := y(\check{\lambda}) + \hbar \cdot \langle y(\check{\lambda}) - \check{\lambda}, \rho \rangle, \quad \quad \check{\lambda} \in \ft.$$
On the spectrum $\Spec U^\hbar(\ft) \simeq \ft^* \times \mathbb{A}^1_\hbar$, the induced action takes the form 
$$y \overset{\hbar} \cdot (\lambda, \hbar) = ( y(\lambda + \hbar \rho) - \hbar \rho, \hbar).$$
Then the map \eqref{hcsign} factors through an isomorphism 
$$U^\hbar(\fm)^M \simeq U^\hbar(\ft)^W,$$
where on the right hand side we are taking invariants with respect to the dot action
of $W$ on $U^\hbar(\ft)$. 
\iffalse Namely, if for $\check{\lambda} \in \ft$ we write $\check{\lambda} \overset{\ell} \cdot$ for the induced action of $\ft \simeq \fb/\fn$ on the highest weight space $(\Delta^\circ_0)_0$}
\fi 

}

\sssRecall that the action of left and right invariant vector fields on $M$ give rise to a map of $\bk[\hbar]$-algebras
$$U^\hbar(\fm) \underset{U^\hbar(\fm)^M} \otimes U^\hbar(\fm)\lmod \rightarrow  D^\hbar(M).$$
Consider the associated induction and restriction functors
$$\on{Ind}: U^\hbar(\fm) \underset{U^\hbar(\fm)^M} \otimes U^\hbar(\fm)\lmod \rightleftarrows  D^\hbar(M)\lmod: \on{Oblv}.$$
Passing to weak highest weight objects with respect to right translations, these adjunctions give rise to an adjunction 
$$\on{Ind}: U^\hbar(\fm) \underset{U^\hbar(\fm)^M} \otimes U^\hbar(\fm)\lmod^{N^r, (T^r, w), \circ} \rightleftarrows  D^\hbar(M)\lmod^{N^r, (T^r, w), \circ}: \on{Oblv}.$$
Precomposing with the tautological adjunction 
$$U^{\hbar}(\fm) \underset{U^\hbar(\fm)^M} \otimes U^\hbar(\ft)^\circ\lmod \rightleftarrows U^{\hbar}(\fm) \underset{U^\hbar(\fm)^M} \otimes U^{\hbar}(\fm)\lmod^{N^r, (T^r, w), \circ}$$afforded by $\Delta_0^\circ$, where we use the identification $\End(\Delta_0^\circ) \simeq U^\hbar(\ft)^\circ$ induced by \eqref{signsigns}, we obtain the adjunction relevant for localization, namely
\begin{equation} \label{e:loc}  \on{Loc}: U^\hbar(\fm) \underset{U^\hbar(\fm)^M} \otimes U^\hbar(\ft)^\circ\lmod \rightleftarrows D^\hbar(M)\lmod^{N^r, (T^r, w), \circ}: \on{Loc}^R.\end{equation}

\sss The desired localization in families, over the generic locus, now reads as follows. 

\begin{Prop} \label{t:loc} The functors in \eqref{e:loc} are mutually inverse equivalences of abelian categories. 
\end{Prop}

\begin{proof} For ease of reading, we denote the enhanced and usual flag varieties by $$\wsB  := M/N \xrightarrow{\pi} M/B =: \sB,$$ 
and break the argument into several steps.  We will first show that $\on{Loc}^R$ is exact and faithful.

{\noindent \em Step 1.} Note that the exactness and faithfulness of $\on{Loc}^R$ may be checked at the level of the further composition, which we denote by $\Gamma$, 
$$\Gamma: D^\hbar(M)\lmod^{N^r, (T^r, w), \circ} \xrightarrow{\on{Loc}^R} U^\hbar(\fm) \underset{U^\hbar(\fm)^M} \otimes U^\hbar(\ft)^\circ\lmod \xrightarrow{\on{Oblv}} \on{Vect}.$$
Let us begin by recalling an alternative description of $\Gamma$. Note that we have a natural forgetful functor \begin{equation} \label{e:2qcoh}\on{Oblv}: D^\hbar(M)\lmod^{N^r, (T^r, w), \circ} 
\rightarrow \QCoh(\sB),\end{equation}which passes from an $(N^r, (T^r, w), \circ)$-equivariant $D^\hbar(M)$-module to its underlying $B^r$-equivariant quasicoherent sheaf on $M$. With this, $\Gamma$ then canonically identifies with the composition 
$$D^\hbar(M)\lmod^{N^r, (T^r, w), \circ} \xrightarrow{\on{Oblv}} \QCoh(\sB) \xrightarrow{\Pi_*} \on{Vect},$$
where $\Pi: \sB \rightarrow \on{pt}$ is the projection to a point.  That is, $\Gamma$ simply calculates the global sections on $\sB$ of the underlying quasicoherent sheaf. 

To proceed, recall the projection $\pi: \wsB \rightarrow \sB$. Write $\Dh$ for the sheaf of Rees differential operators on $\wsB$, and consider the associated sheaf of algebras $\pi_*(\Dh)$ on $\sB$, along with its subsheaf of right invariant operators $$\mathscr{R} := \pi_* (\Dh)^{T^r}.$$ 
By construction, we may factor the functor \eqref{e:2qcoh} as 
\begin{equation} \label{e:2qcoh2} D^\hbar(M)\lmod^{N^r, (T^r, w), \circ} 
\xrightarrow{\Oblv^{enh}} \mathscr{R}\lmod(\QCoh(\sB)) \rightarrow \QCoh(\sB),\end{equation}
where $\on{Oblv}^{enh}$ simply remembers the natural action of $\mathscr{R}$ on the sheaf of $B^r$-invariant sections of an $(N^r, T^r, w, \circ)$-equivariant $D^\hbar(M)$-module. In particular, we  obtain a similar factorization for $\Gamma$. 

{\em \noindent Step 2.} Consider the natural maps  $$\mathscr{O}_{\sB} \rightarrow \mathscr{R} \quad \text{and} \quad \phi: U^\hbar(\fm) \underset{U^\hbar(\fm)^M}\otimes  U^\hbar(\ft)  \rightarrow  \Pi_* \mathscr{R}.$$
Let us write $\sT(\wsB)$ for the tangent sheaf of $\wsB$, and let us consider the corresponding restricted Lie algebroid on $\sB$, namely $$\widetilde{\sT} := \pi_*\sT(\wsB)^{T^r}.$$ Consider the tautological short exact sequence of restricted Lie algebroids
$$0 \rightarrow \fb_{univ} \rightarrow \OO_{\sB} \otimes (\fm \oplus \ft)  \rightarrow \widetilde{\sT} \rightarrow 0.$$
If we pass to their Rees enveloping sheaves of algebras, we recover $\phi$ from the global sections of the surjective map
\begin{equation}\label{eq:mapofsheavesofalgebras} \mathscr{U}^\hbar(\OO_{\sB} \otimes (\fm \oplus \ft)) \twoheadrightarrow  \mathscr{U}^\hbar(\widetilde{\sT}) \simeq \mathscr{R},\end{equation}
and moreover deduce that the kernel of (\ref{eq:mapofsheavesofalgebras}) is the two-sided ideal generated by the augmentation ideal of $\mathscr{U}^\hbar(\fb_{univ})$. 

Similarly, passing further to the $p$-centers of the Rees enveloping sheaves of algebras, we obtain a surjective map 
\begin{equation} \label{e:groth}\mathscr{Z}^\hbar_p(\OO_{\sB} \otimes (\fm \oplus \ft)) \twoheadrightarrow \mathscr{Z}^\hbar_p(\widetilde{\sT})\end{equation}
and deduce that its kernel is the ideal generated by the augmentation ideal of $\mathscr{Z}^\hbar_p(\fb_{univ})$. We also note that explicitly, on spectra, if we write $\widetilde{\fm}^* \rightarrow \fm^*$ for the Grothendieck alteration, then the map \eqref{e:groth} corresponds to the tautological closed embedding 
\begin{equation} \label{e:groth2}  \iota^{(1)} \times \on{id}: (\widetilde{\fm}^*)^{(1)} \times \mathbb{A}^1_\hbar \hookrightarrow  (\sB \times \fm^* \times \ft^*)^{(1)} \times \mathbb{A}^1_\hbar,\end{equation}
where $\iota: \widetilde{\fm}^* \rightarrow \sB \times \fm^* \times \ft^*$ is the usual closed embedding. Namely, recall that $\widetilde{\fm}^*$ consists of pairs $(\xi, \fb')$, where $\xi \in \fm^*$ and $\fb'$ is a Borel subalgebra whose nilradical $\mathfrak{n}'$ is annihilated by $\xi$; $\iota$ sends such a pair $(\xi, \fb')$ to the triple $(\fb', \xi, \xi')$ where $\xi'$ denotes the induced map $\ft \simeq \fb'/\mathfrak{n}' \xrightarrow{\xi} \bk.$ 

{\em \noindent Step 3.} Let us now show the claimed exactness and faithfulness of $\on{Loc}^R$. Namely, note that $\Gamma$ factors as 
\begin{align*}D^\hbar(M)\lmod^{N^r, (T^r, w), \circ} \xrightarrow{\on{Oblv}^{enh}} \mathscr{R}\lmod(\QCoh(\sB)) \\ \xrightarrow{\on{Fr}_*^{enh}} \mathscr{Z}^{\hbar}_p(\widetilde{\sT})\lmod(\QCoh(\sB^{(1)}) \xrightarrow{\Pi_*} \on{Vect};\end{align*}
here $\on{Fr}_*$ denotes pushforward along the relative Frobenius, the superscript `$enh$' on $\on{Fr}_*$ indicates that we remember the residual action of the $p$-center of $\mathscr{R}$, and $\Pi_*$ again denotes pushforward to a point. 

To finish, we note that $\on{Oblv}^{enh}$ and $\on{Fr}_*^{enh}$ are tautologically exact and faithful. Moreover, by our assumption of genericity on our $D^\hbar(M)$-module, as a quasicoherent sheaf on the relative spectrum of $\mathscr{Z}^\hbar_p(\widetilde{\sT}),$ it is pushed forward from the open set 
$$(\widetilde{\fm}^*)^{(1)} \times \mathbb{A}^1_{\hbar} \underset{(\sB \times \fm^* \times \ft^*)^{(1)} \times \mathbb{A}^1_\hbar} \times (\sB \times \fm^* \times \ft_{rs}^*)^{(1)} \times \mathbb{A}^1_\hbar = (\widetilde{\fm}^*_{rs})^{(1)} \times \mathbb{A}^1_\hbar.$$
But it is standard that over the regular semisimple locus the Grothendieck alteration $\widetilde{\fm}^*_{rs} \rightarrow \fm_{rs}$ is a finite Weyl group torsor, and in particular finite over affine, whence affine, which shows the exactness and faithfulness of the final pushforward to a point.

{\em \noindent Step 4.} We next claim that the unit map $\on{id} \rightarrow \on{Loc}^R \circ \on{Loc}$, when applied to the generator $$U^\circ := U^{\hbar}(\fm) \underset{U^\hbar(\fm)^M} \otimes U^\hbar(\ft)^\circ,$$is an isomorphism. However, we note that the adjunction \eqref{e:loc} is obtained by base changing to the generic locus a similar adjunction 
$$\on{Loc}: U^{\hbar}(\fm) \underset{U^\hbar(\fm)^M} \otimes U^\hbar(\ft) \rightleftarrows D^\hbar(M)\lmod^{N, (T^r,w)}: \on{Loc}^R.$$By the commutation of both appearing functors with 
filtered colimits, it is therefore enough to verify the analogous claim for 
$$U :=  U^{\hbar}(\fm) \underset{U^\hbar(\fm)^M} \otimes U^\hbar(\ft).$$
Note that $U$ is $\mathbb{Z}_{\geqslant 0}$-graded and torsion free $\bk[\hbar]$-algebra, as is $\on{Loc}^R \circ \on{Loc}(U)$. Therefore, it is enough to check the assertion holds at the central fiber, i.e., at $\hbar = 0$, where it is the standard fact that the projection $$\widetilde{\fm}^* \rightarrow \fm^* \underset{\fm^* /\!/ M} \times \ft^*$$exhibits the target as the affinization of the source.

{\em \noindent Step 5.} We next claim that the unit map $\on{id} \rightarrow \on{Loc}^R \circ \on{Loc}$ is an isomorphism for any $U^\circ$-module. However, by writing it as the cokernel of a map between two free $U^\circ$-modules, this follows from the previous step, the exactness of $\on{Loc}^R$, and its commutation with arbitrary direct sums. 

{\em \noindent Step 6.} We now claim we are done on general grounds. Namely, we are considering an adjunction between categories
$$F: \sC \leftrightarrows \sD: G$$
where we know (i) $F$ is fully faithful, and (ii) $G$ reflects isomorphisms, i.e., a map $d \rightarrow d'$ is an isomorphism if and only if $G(d) \rightarrow G(d')$ is an isomorphism; for functors between abelian categories this follows from faithfulness, which we know holds for 
 $\on{Loc}^R$ by Step 3. Let us see that this implies $F$ and $G$ are mutually inverse equivalences. 

Indeed, it remains to see that $G$ is fully faithful, i.e., that for object $d$ of $\sD$ we have the counit map $FG(d) \rightarrow d$ is an isomorphism. It is enough to see that the map $GFG(d) \rightarrow G(d)$ is an isomorphism. But for any adjunction the composition $$G(d) \rightarrow GFG(d) \rightarrow G(d)$$is an isomorphism, and the first arrow $G(d) \rightarrow GFG(d)$ is an isomorphism by the fully faithfulness of $F$, whence the second is an isomorphism as well, as desired. \end{proof}

\subsection{Intertwining functors}\label{SS_intertwining}

\sss In the previous section, we set up a localization theorem in families over the generic locus of the universal Cartan. We will now apply it to construct intertwining functors which naturally act on categories of generic weak highest weight modules, and establish some of their basic properties.

\sss We retain the notation of the previous section. Let us first record the following useful consequence of Proposition \ref{t:loc}.

\begin{Prop}\label{p:loc} For any graded complete topological $\bk[\hbar]$-algebra $A$ equipped with a strong action of $M$, there is a canonical equivalence, functorial in $A$ and linear over $Z(A)^M$
$$A\lmod^{N, (T,w), \circ} \simeq (A  \otimes_{U^\hbar(\fm)^M} U^\hbar(\ft)^\circ)\lmod^{(M, w)}.$$
\end{Prop}

\begin{proof} Recall from Lemma \ref{l:dividebyK} the canonical identification 
$$A\lmod \simeq (A \wotimes D^\hbar(M))\lmod^{M^\ell},$$
and, in particular, the induced identification 
\begin{align*}A\lmod^{N, (T, w),\circ} &\simeq (A \wotimes D^\hbar(M))\lmod^{M^\ell\times N^r, (T^r, w), \circ}. \intertext{Note next that Proposition \ref{t:loc} identifies the latter with} & (A  \otimes_{\bk[\hbar]} U^\hbar(\fm)  \otimes_{U^\hbar(\fm)^M} U^\hbar(\ft)^\circ)\lmod^{M^\ell}, \end{align*}which we may by Lemma \ref{l:bigbird} rewrite in the desired form, namely  $(A  \otimes_{U^\hbar(\fm)^M} U^\hbar(\ft)^\circ)\lmod^{(M,w)}.$
\end{proof}

\sss {Let us now use the previous proposition to construct the desired intertwining functors. As preparation, denote the abstract Weyl group of $M$ by $W$, and recall from Section \ref{sss:hcstar} the dot action of $W$ on $U^\hbar(\ft)$. 

\begin{Cor} \label{c:intops}For any graded complete topological $\bk[\hbar]$-algebra $A$ equipped with a strong action of $M$, there is a canonical action of $W$ on  $A\lmod^{N, (T, w), \circ}$  $$\sM \mapsto y \star \sM, \quad \sM \in A\lmod^{N, (T, w), \circ}, y \in W.$$
enjoying the following properties. 

\begin{enumerate}
\item This action is $Z(A)^M$-linear.

\item This action is $U^\hbar(\ft)$-semilinear. That is, for any $\sM$ as above, and $\xi \in U^\hbar(\ft)$ , we have that applying $y \star -$ to $\xi: \sM \rightarrow \sM$ gives $$y \overset{\hbar} \cdot \xi: y \star \sM \rightarrow y \star \sM.$$
{ \item This action is functorial in the algebra $A$. That is, for any $A_1, A_2$ as above, and a (discrete)  $A_1$-$A_2$ bimodule $B$ strongly equivariant for the diagonal action of $M$, the associated functor
$$B \underset{A_2} \otimes -: A_2\lmod^{N, (T, w), \circ} \rightarrow A_1\lmod^{N, (T, w), \circ}$$
carries a canonical datum of $W$-equivariance. }
\end{enumerate}
\end{Cor}

\begin{proof}By  Proposition \ref{p:loc}, it is enough to produce an action of $W$ on $$A  \otimes_{U^\hbar(\fm)^M} U^\hbar(\ft)^\circ,$$viewed as an algebra with a strong $M$-action, functorially in the argument $A$. However, as the dot action of $W$ on $U^\hbar(\ft)$ fixes $U^\hbar(\fm)^M$ and permutes the affine hyperplanes $\halpha - \hbar \cdot n 
= 0$, $\halpha \in \check{\Phi}, n \in \mathbb{F}_p$, we have an induced action of $W$ on the above tensor product of the desired form. The claimed properties follow from the construction. 
\end{proof}

\sss We finish with one relevant calculation of the intertwining functors.

{
\begin{Prop} \label{p:univverm} 
Consider the generic universal Verma module $$\Delta^\circ_{0} \in U^\hbar(\fm)\lmod^{N, (T, w), \circ}$$of weak highest weight $0$, cf. (\ref{eq:universal_Verma_defn}). Then we have isomorphisms of objects of  $U^\hbar(\fm)\lmod^{N, (T, w), \circ}$
$$y \star \Delta^\circ_{0} \simeq \Delta^\circ_{0}, \quad \quad y \in W.$$
\end{Prop}

\begin{proof} { For ease of reading, we will break the argument into steps. } 

{{\em Step 1.}} It is enough to verify that, under the equivalence 
$$U^\hbar(\fm)\lmod^{N, (T,w), \circ} \simeq U^{\hbar}(\fm)  \otimes_{U^\hbar(\fm)^M} U^\hbar(\ft)^\circ \lmod^{(M,w)}$$
afforded by Proposition \ref{p:loc}, the object $\Delta^\circ_0$ is exchanged with the algebra $$A := U^{\hbar}(\fm)  \otimes_{U^\hbar(\fm)^M} U^\hbar(\ft)^\circ$$itself, equipped with its tautological datum of $(M,w)$-equivariance, as the latter object carries a tautological datum of $W$-equivariance.  

To see this, we will push the object $A$ through the equivalences used in the proof of Proposition \ref{p:loc}. 

{{\em Step 2.}} Set $$\tilde{A}:=U^{\hbar}(\fm)  \otimes_{\bk[\hbar]} A.$$
Under the equivalence 
$$A \lmod^{(M,w)} \simeq  \tilde{A}\lmod^{M^\ell},$$by construction $A$ is exchanged with the object 
$$B := \on{ind}_{U^\hbar(\fm^\ell)}^{\tilde{A}} \bk[\hbar],$$
i.e., the induction from the diagonal copy { $\fm^\ell \subset \tilde{A}$} of the trivial representation. 

{ {\em Step 3.} Next, We would like to calculate the image of $B$ under the equivalence 
$$\Phi: \tilde{A}\lmod^{M^\ell} \simeq U^\hbar(\fm)  \otimes_{\bk[\hbar]} D^\hbar(M)\lmod^{M^\ell \times N^r, (T^r, w), \circ}.$$
To approach this, we first collect some preliminary facts in this step, and return to our calculation of the image of $B$ in its sequel. 

Note the {  algebra isomorphism \begin{equation} \label{e:2down1across}(D^\hbar(M)/D^\hbar(M) \cdot \mathfrak{n}^r)^{B^r} \simeq  U^{\hbar}(\fm) \otimes_{{U^\hbar(\fm)^M}} U^\hbar(\ft),\end{equation} which follows from Step 4 of the proof of Proposition \ref{t:loc}. By construction, if we denote the left translation quantum comoment map for $D^\hbar(M)$ by $\mu^l: U^\hbar(\fm^l) \rightarrow D^\hbar(M),$ the composition 
\begin{equation} \label{e:comoment1} U^\hbar(\fm) = U^\hbar(\fm^l) \xrightarrow{\mu^l} D^\hbar(M)^{B^r} \rightarrow (D^\hbar(M)/D^\hbar(M) \cdot \mathfrak{n}^r)^{B^r} \overset{\eqref{e:2down1across}}{\simeq}  U^{\hbar}(\fm)  \otimes_{{U^\hbar(\fm)^M}} U^\hbar(\ft)\end{equation}agrees with the tautological insertion 
\begin{equation} \label{e:comoment2} \on{id} \otimes 1: U^\hbar(\fm) \rightarrow  U^{\hbar}(\fm)  \otimes_{{U^\hbar(\fm)^M}} U^\hbar(\ft) .\end{equation}In other words, the isomorphism \eqref{e:2down1across} is one of algebras with strong actions of $M$. By further base changing \eqref{e:2down1across}, we   deduce that \begin{equation} \label{e:basechange2down} U^{\hbar}(\fm) \otimes_{{U^\hbar(\fm)^M}} U^\hbar(\ft)^\circ \simeq \left([D^\hbar(M)/D^\hbar(M) \cdot \mathfrak{n}^r]  \otimes_{{U^\hbar(\fm^l)^{M^l}}} U^\hbar(\ft^l)^\circ\right)^{B^r}.\end{equation}

{ As a final preliminary fact,  we claim that in the right hand side of \eqref{e:basechange2down} we may replace our base change  with respect to  $U^\hbar(\fm^l)^{M^l}$ with the analogous base change with respect to $U^\hbar(\fm^r)^{M^r}$, i.e., that we have a further identification 
\begin{equation} \label{e:leftvsright} ([D^\hbar(M)/D^\hbar(M) \cdot \mathfrak{n}^r]  \otimes_{{U^\hbar(\fm^l)^{M^l}}} U^\hbar(\ft^l)^\circ)^{B^r} \simeq ([D^\hbar(M)/D^\hbar(M) \cdot \mathfrak{n}^r]  \otimes_{{U^\hbar(\fm^r)^{M^r}}} U^\hbar(\ft^r)^\circ)^{B^r}.\end{equation} Indeed, to see this, recall that we have the equality of subalgebras $$U^\hbar(\fm^l)^{M^l} \subset D^\hbar(M) \supset U^\hbar(\fm^r)^{M^r}.$$Let us also consider the identification $\ft^l$ and $\ft^r$ given by $X \mapsto -w_\circ(X)$, for $X \in \ft$, where $w_\circ$ denotes the longest element of $W$; we similarly denote the corresponding isomorphism of Rees enveloping algebras by 
$$-w_\circ: U^\hbar(\ft^l) \simeq U^\hbar(\ft^r).$$Then by a standard argument involving restriction to the big cell $M^\diamond \subset M$, i.e., the open $B^l \times B^r$ orbit, and considering its Hamiltonian reduction by $N^l \times N^r$, we have a commutative diagram 
$$\xymatrix{ U^\hbar(\ft^l) \ar[r]^{-w_\circ} & U^\hbar(\ft^r) \\ U^\hbar(\fm^l)^{M^l} \ar[u] \ar@{=}[r] & U^\hbar(\fm^r)^{M^r} \ar[u],}$$ 
where the vertical arrows are the Harish--Chandra homomorphisms for $\fm^l$ and $\fm^r$, respectively. As $-w_\circ$ extends uniquely to an isomorphism 
$$-w_\circ: U^\hbar(\ft^l)^\circ \simeq U^\hbar(\ft^r)^\circ,$$
i.e., permutes the relevant affine root hyperplanes in Definition \ref{d:genmodule}, the claimed identity \eqref{e:leftvsright} follows.}

{\em Step 4.} We are now ready to calculate the image of $B$ under the equivalence 
$$\Phi: \tilde{A}\lmod^{M^\ell} \simeq U^\hbar(\fm)  \otimes_{\bk[\hbar]} D^\hbar(M)\lmod^{M^\ell \times N^r, (T^r, w), \circ}.$$
To compute this, we first note that $\Phi$ evaluated on $\tilde{A}$ itself yields
$$\Phi: \tilde{A} = U^\hbar(\fm) \otimes_{\bk[\hbar]} U^\hbar(\fm) \otimes_{U^\hbar(\fm)^M} U^\hbar(\ft)^\circ \overset{\on{id} \otimes  (\eqref{e:leftvsright} \circ  \eqref{e:2down1across})}\simeq U^\hbar(\fm) \otimes_{\bk[\hbar]}  ([D^\hbar(M)/D^\hbar(M) \cdot \mathfrak{n}^r]  \otimes_{{U^\hbar(\fm^r)^{M^r}}} U^\hbar(\ft^r)^\circ)^{B^r}.$$
By the discussion near \eqref{e:comoment1} and \eqref{e:comoment2}, $\Phi$ exchanges the diagonal comoment embedding $\fm^{\ell, \tilde{A}} \subset \tilde{A}$ and the similarly named diagonal comoment embedding $$\fm^\ell \subset U^\hbar(\fm) \otimes_{\bk[\hbar]}  ([D^\hbar(M)/D^\hbar(M) \cdot \mathfrak{n}^r)]  \otimes_{{U^\hbar(\fm^r)^{M^r}}} U^\hbar(\ft^r)^\circ)^{B^r}.$$ 
It therefore follows from this, the very definition of $\Phi$, and its $t$-exactness (in particular, the right exactness) that
\begin{equation}\label{eq:PhiB}\Phi(B) := \Phi( \tilde{A} \otimes_{U^\hbar(\fm^{\ell, \tilde{A}})} \bk[\hbar])\simeq \Phi(\tilde{A})   \otimes_{U^\hbar(\fm^\ell)} \bk[\hbar].\end{equation}

Therefore, we see that the right hand side of (\ref{eq:PhiB}) is naturally isomorphic to
$$(U^\hbar(\fm)  \otimes_{\bk[\hbar]} [D^\hbar(M)/ D^\hbar(M) \cdot \mathfrak{n}^r]  \otimes_{U^\hbar(\fm^r)^{M^r}} U^\hbar(\ft^r)^\circ)  \otimes_{U^\hbar(\fm^\ell)} \bk[\hbar].$$
and therefore deduce the presentation of $\Phi(B)$ by generators and relations: $$[U^\hbar(\fm)  \otimes_{\bk[\hbar]} D^\hbar(M)]/ [U^\hbar(\fm)  \otimes_{\bk[\hbar]} D^\hbar(M) \cdot (\fm^\ell + \mathfrak{n}^r)]  \otimes_{U^\hbar(\fm^r)^{M^r}} U^\hbar(\ft^r)^\circ =: C.$$

{\em Step 5.} Finally, it remains to compute the image of $C$ under the equivalence 
$$\Psi: U^\hbar(\fm)  \otimes_{\bk[\hbar]} D^\hbar(M)\lmod^{M^\ell \times N^r, (T^r, w), \circ} \simeq U^\hbar(\fm)\lmod^{N, (T, w), \circ}. $$By definition of $C$, we have 
\begin{align*}
\Psi(C) &= \Psi( [U^\hbar(\fm)  \otimes_{\bk[\hbar]} D^\hbar(M)]/ [U^\hbar(\fm) \otimes_{\bk[\hbar]}  D^\hbar(M) \cdot (\fm^\ell + \mathfrak{n}^r)]\otimes_{U^\hbar(\fm^r)^{M^r}}  U^\hbar(\ft^r)^\circ). \intertext{Let us tautologically rewrite the appearing relation corresponding to $\mathfrak{n}^r$ as} &\simeq \Psi(\left([U^\hbar(\fm)  \otimes_{\bk[\hbar]} D^\hbar(M)]/ [U^\hbar(\fm)  \otimes_{\bk[\hbar]} D^\hbar(M) \cdot \fm^\ell] \otimes_{U^\hbar(\mathfrak{n}^r)} \bk[\hbar]\right) \otimes_{U^\hbar(\fm^r)^{M^r}} U^\hbar(\ft^r)^\circ).
\intertext{Recall that for a general algebra $A$ equipped with a strong action of $M$, the canonical equivalence 
$$A \otimes_{\bk[\hbar]} D^\hbar(M)\lmod^{M^\ell} \simeq A\lmod$$exchanges the object $$\xi := [A \otimes_{\bk[\hbar]} D^\hbar(M)]/[A \otimes_{\bk[\hbar]} D^\hbar(M) \cdot \fm^\ell]$$ with $A$ itself. Moreover if for a ring $R$ we write $R^{\on{op}}$ for the ring with its reversed multiplication, the natural map $$U^\hbar(\fm^r)^{\on{op}} \rightarrow \End_{A \otimes_{\bk[\hbar]} D^\hbar(M)\lmod^{M^\ell}}(\xi)$$is exchanged with the quantum comoment map $$U^\hbar(\fm)^{\on{op}} \rightarrow A^{\on{op}} \simeq \End_{A\lmod}(A),$$and relatedly the tautological $M^r$-weakly equivariant structure on $\xi$ is exchanged with the tautological $M$-weakly equivariant structure on $A$. Applying this to $A = U^\hbar(\fm)$, we may continue as}
&\simeq \left(U^\hbar(\fm)   \otimes_{U^\hbar(\mathfrak{n})} \bk[\hbar]\right)  \otimes_{U^\hbar(\fm)^M} U^\hbar(\ft)^\circ 
\\ & \simeq \Delta_0  \otimes_{U^\hbar(\fm)^M} U^\hbar(\ft)^\circ \\ & \simeq \Delta_0^\circ,
\end{align*}
as desired. }}\end{proof}

}

\subsection{Action of $\,^I\!\A^\hbar$ on $\Delta^G_{\kappa\hbar}$}
We have established an action of $\,^I\!\A^\hbar$ on $\Delta^{G,I}_{\kappa\hbar}$. It is by $\Iw$-equivariant endomorphisms and makes 
$\varphi^B_{\kappa\hbar}: \Delta^{G,I}_{\kappa\hbar}\rightarrow \Wak^I_{\kappa\hbar}$ semilinear with respect to $\,^I\!\A^\hbar$. To finish the proof of Theorem \ref{Thm:Verma_endomorphisms} it suffices to prove the following lemma.

\begin{Lem} The action of $\,^I\!\A^\hbar$ preserves $\Delta^{G}_{\kappa\hbar}\subset \Delta^{G,I}_{\kappa\hbar}$.
\end{Lem}
\begin{proof} We write $\Delta^{G,\overline{i}}_{\kappa\hbar}$ for 
$\Delta^{G,\Phi^+\setminus \{\alpha_i\}}_{\kappa\hbar}$. We note that $\Delta^{G,\overline{i}}_{\kappa\hbar}, 
\Delta^{L_i,\varnothing}_{\kappa,\hbar}$ are parabolically induced from the Verma module
$\Delta^\circ_{0; i}$ for $L_i$ defined by (\ref{eq:universal_Verma_defn}). From Proposition \ref{p:univverm} we get an automorphism $s_i$ of $\Delta^\circ_{0; i}$ that acts as the operator $s_i\overset{\hbar}\cdot-$ on $\h^{da,*}_{\kappa\hbar}$ by (2) of 
Corollary \ref{c:intops}. 

{ Note that $\widehat{U}^\hbar_{\kappa}(\hat{\g})$ and $\widehat{D}^\hbar(\Loop N_i)\widehat{\otimes}_{\Ring[\hbar]}\widehat{U}^\hbar_{\kappa}(\hat{\lf}_i)$ are graded complete topological $\Ring[\hbar]$-algebras equipped with strong actions of $L_i$.} Applying (1) of Corollary \ref{c:intops} to the action of $L_i$ on 
$A:=\widehat{D}^\hbar(\Loop N_i)\widehat{\otimes}_{\Ring[\hbar]}\widehat{U}^\hbar_{\kappa}(\hat{\lf}_i)$ we see that the action of $\,^I\!\A^\hbar$ on $\WakDelta^{P_i,\overline{i}}_{\kappa\hbar}$ intertwines the action of the endomorphism $s_i$. And applying (3) of Corollary \ref{c:intops} to the homomorphism $\Ffr_{P_i}$ we see that $\varphi^{P_i}_{\kappa\hbar}$ also intertwines the endomorphisms $s_i$. It follows that the action of $\,^I\!\A^\hbar$ on
$\Delta^{G,\varnothing}_{\kappa\hbar}$ is invariant under all $s_i$, hence under $W$. 

Now observe that the action of $s_i$ on 
$\Delta^{G,\overline{i}}_{\kappa\hbar}$ intertwines the localizations 
$$\Delta^{G,\Phi^+\setminus\Sigma}_{\kappa\hbar}, 
\Delta^{G,\Phi^+\setminus s_i\Sigma}_{\kappa\hbar}$$
for all subsets $\Sigma$ of $\Phi^+$ containing $\alpha_i$, where we write $s_i\Sigma$ for $s_i(\Sigma\setminus \{\alpha_i\})\sqcup \{\alpha_i\}$. Since 
$\,^I \A^\hbar$ preserves the localization 
$\Delta^{G,i}_{\kappa\hbar}$ for all $i\in I$, it follows that $\,^I \A^\hbar$ preserves the localization $\Delta^{G,\{\alpha\}}_{\kappa\hbar}$ for any positive root $\alpha$. Hence $\,^I\!\A^\hbar$ preserves the intersection of these localizations, i.e., $\Delta^G_{\kappa\hbar}$.
\end{proof}

\section{Description of $\,^I\!\A$}\label{S_IA_description}
Let $\Ring$ be a perfect characteristic $p$ field. Assume $\kappa\in \Ring\setminus \F_p$.

\subsection{Main result}

\sss
Let $\A^{\hbar,\natural}$ denote the algebra of polynomials in the variables $h^{i,\natural}_n, n\leqslant 0$, and $\hbar^\natural$.
We have $\Ring[\hbar]\otimes_{\Ring[\hbar^\natural]}\A^{\natural,\hbar}\hookrightarrow \A^\hbar$ via 
$$h^{i,\natural}_n\mapsto \frac{(h^i_n)^p-(\kappa\hbar)^{p-1}h^i_{pn}}{1-\kappa^{p-1}},\quad \hbar^\natural\mapsto \hbar^p.$$ 
Note that $\A^{\hbar,\natural}$ is (naively) graded with 
$\deg h^{i,\natural}_n=\deg \hbar^\natural=1$ and the isomorphism 
with $\A^\hbar$ rescales the degrees $p$ times. 

Define a homogeneous degree $2$ element $\underline{X}^{i,\natural}_n\in \A^{\hbar,\natural}$ (with $n\leqslant 0$) by
\begin{equation}\label{eq:underline_X_in_natural}
\underline{X}^{i,\natural}_n:=\frac{1}{4}\left(\sum_{j=n}^0h^{i,\natural}_{j+n}h^{i,\natural}_{-j}\right)-\frac{1}{2}(n+1)\hbar^\natural h^{i,\natural}_{n}.
\end{equation}
Note that under the homomorphism $\A^{\hbar,\natural}\rightarrow \A^\hbar$, the element $X^{i,\natural}_n$ is sent to the homogenous analog of $\HC(X^i_n)$, see Section \ref{SS_HC_Xn_formula}.

Consider the $\Ring[\hbar^\natural]$-subalgebra $\,^i\!\A^{\hbar,\natural}$ generated by the elements $X^{i,\natural}_n$ and $h^{j,\natural}_n$ for $j\neq i$.
Then $\Ring[\hbar]\otimes_{\Ring[\hbar^\natural]}\,^i\!\A^{\hbar,\natural}\xrightarrow{\sim} \,^i\!\A^{\hbar}$. 



\sss
Set $\,^i\!\A^\natural=\,^i\!\A^{\hbar,\natural}/(\hbar^\natural-1)$.
We notice that, over a characteristic $0$ field, by results of \cite{Frenkel_loop}, the analog of $\bigcap_{i\in I}\,^i\!\A^\natural$ is the algebra on the space of $\check{G}$-opers on $\D^\times$ with regular singularities. Here $\check{G}$ is the Langlands dual group of adjoint type. We will see that this is the case in our situation as well (up to twists). 

\sss 
Choose a maximal torus and a Borel subgroup $\check{T}\subset\check{B}\subset \check{G}$.  Let $\rho:\mathbb{G}_m\rightarrow \check{T}$ be the one-parameter subgroup corresponding to the sum of fundamental weights. Let $\check{\g}=\bigoplus_{i\in \Z}\check{\g}_i$ be the corresponding grading. 
Let $\check{f}\in \check{\g}_{-1}$ be the sum of roots vectors corresponding to negative simple roots. 

For $k\in \Z$, define the subscheme $\widetilde{\Op}^\hbar_k\subset \operatorname{Spec}(\Ring[\hbar^\natural])\times \Loop \check{\g}$ as 
$$\{\lambda\partial+\check{f}+\sum_{i\geqslant 0}F_i| F_i\in t^{-k(i+1)}\check{\g}_i[[t]]\},$$
here $\hbar^\natural$ is the function sending $\lambda\partial+\check{f}+\sum_{i\geqslant 0}F_i$ to $\lambda$. We will write $\widetilde{\Op}^\hbar_k(\check{G})$ when we want to indicate the dependence on $\check{G}$. The fiber of $\widetilde{\Op}^\hbar_k$ at $\hbar^\natural=1$ will be denoted by 
$\widetilde{\Op}_k$.

Define the subgroup scheme $\Jet_{(k)} \check{N}\subset \Loop \check{N}$ as 
$$\{\exp(\sum_{i>0}x_i)| x_i\in t^{-ki}\check{\g}_{i}[[t]]\}.$$
Note that since $p>h$, the morphism $\exp:\mathfrak{n}\rightarrow N$ is well-defined.
We also remark that 
$$\widetilde{\Op}^\hbar_k\subset 
\widetilde{\Op}^\hbar_{k+1}, \Jet_{(k)} \check{N}\subset 
\Jet_{(k+1)}\check{N}, \forall k.$$
Define the actions of $\mathbb{G}_m$ and $\Jet_{(k)} \check{N}$ on  $\widetilde{\Op}^\hbar_k(\check{\g})$ by 
\begin{equation}\label{eq:Kostant_actions}
\begin{split}
&n.(\lambda\partial+F)=\lambda\partial+\operatorname{Ad}(n)F-\lambda(\partial. n)n^{-1}, n\in \Jet_{(k)} N\\
&z.(\lambda\partial+F)=\lambda z^{-1}\partial+ z^{-1}\check{\rho}(z)^{-1}F.
\end{split}
\end{equation}
In particular, note that $f$ is $\mathbb{G}_m$-stable. 
These two actions combine into an action of $\mathbb{G}_m\ltimes \Jet_{(k)} \check{N}$, where $\mathbb{G}_m$ acts on $\Jet_{(k)} \check{N}$ by $\operatorname{Ad}(\rho)^{-1}$.

\sss\label{SSS_opers_coord}
Now we introduce a ``deformed jet'' version of the Kostant-Slodowy slice, a space of opers. 

Consider the $\slf_2$-triple $(\check{e},2\rho,\check{f})$ and 
the affine subspace 
$$\Op^\hbar_k:=\{\lambda\partial+f+\sum_{i\geqslant 1}F_i| F_i\in t^{-(i+1)k}\mathfrak{z}_{\check{\g}}(\check{e})_i[[t]]\},$$
where $\mathfrak{z}_{\check{\g}}(\check{e})$ stands for the centralizer of $\check{e}$ in $\check{\g}$.
Then we have the following result. 

\begin{Lem}\label{Lem:Kostant_section}
Suppose $p>h$. Then the multiplication map
\begin{equation}\label{eq:jet_product_action}
\Jet_{(k)} \check{N}\times \Op^\hbar_k \rightarrow \widetilde{\Op}^\hbar_k 
\end{equation}
is an isomorphism.
\end{Lem}
\begin{proof}
We note that, under our conditions on $p$, the multiplication map $\check{N}\times \check{S}\rightarrow \check{f}+\check{\mathfrak{b}}$ is an isomorphism, \cite[Proposition 3.2.1]{Riche_Kostant}. By passing to arc spaces we see that $\Jet \check{N}\times \Jet \check{S}\rightarrow \Jet(\check{f}+\check{\mathfrak{b}})$ is an isomorphism. 
From here we see that $\Jet_{(k)} \check{N}\times \Jet_{(k)} \check{S}\rightarrow \Jet_{(k)}(\check{f}+\check{\mathfrak{b}})$ is an isomorphism for all $k$. Here we write 
$\Jet_{(k)} \check{S}, \Jet_{(k)}(\check{f}+\check{\mathfrak{b}})$ for the fibers of $\Op^\hbar_{k},\widetilde{\Op}^\hbar_{k}$ over $\hbar=0$.
Namely, we have an automorphism of $t^\rho$ of $ \Loop \check{\g}$ that on $\Loop \check{\g}_i$ is given by the multiplication by $t^i$. Applying $t^{-k\rho}$ to $\Jet_{(k)}\check{S}, 
\Jet_{(k)}(\check{f}+\check{\mathfrak{b}}),\Jet_{(k)} \check{N}$ we get 
$t^{-k}\Jet \check{S}, t^{-k}\Jet(\check{f}+\check{\mathfrak{b}})\subset \Loop \check{\g}$ and $\Jet \check{N}\subset \Loop \check{N}$. The claim that the multiplication morphism is an isomorphism follows.

Now note that the $\mathbb{G}_m$-actions on $\Jet_{(k)} \check{N},\Op^\hbar_k, \widetilde{\Op}^\hbar_k$ give positive gradings on the algebras of functions. So, a version of the graded Nakayama lemma tells us that once the multiplication map (\ref{eq:jet_product_action}) is an isomorphism  over $0\in \operatorname{Spec}(\Ring[\hbar])$, it is an isomorphism. 
\end{proof}

We remark that the fibers at $\hbar=1$ of the affine schemes $\Op^\hbar_0$
and $\Op^\hbar_1$ can be interpreted as the space of $\check{G}$-opers on
$\D=\operatorname{Spec}(\Ring[[t])$ and the space of $\check{G}$-opers on $\D^\times$ with regular singularities, respectively.

\sss 

To make the notation compatible, below we will write $\Op^{\hbar^\natural}_k$ instead of
$\Op^{\hbar}_k$: the element $\hbar^\natural$ still sends $\lambda+\check{f}+\sum F_i$ to $\lambda$. 
It turns out that there is a (naively) graded algebra homomorphism 
$\mathsf{M}_k:\Ring[\Op^{\hbar^\natural}_k]\rightarrow \Ring[t^{-k}\Jet\check{\h}][\hbar^\natural]$ called the {\it Miura map}. This will be recalled in Section \ref{SS_Miura}. Our main result, to be proved in Section \ref{SS_Miura_main}, is as follows.

\begin{Prop}\label{Prop:Miura_isomorphism}
$\mathsf{M}_1$ is injective and its image coincides with $\,^I\!\A^{\hbar,\natural}$.
\end{Prop}

\subsection{Miura map}\label{SS_Miura}
\sss We start by recalling how $\mathsf{M}_k$ is defined for arbitrary $k$.  Consider the embedding
\begin{equation}\label{eq:jet_Cartan_embedding}
\lambda\partial+F\mapsto \lambda\partial+\check{f}+F: \operatorname{Spec}(\Ring[\hbar^\natural])\times t^{-k}\Jet(\check{\h})\hookrightarrow 
\widetilde{\Op}^{\hbar^\natural}_k.
\end{equation}
The pullback map $\Ring[\widetilde{\Op}^\hbar_k]\rightarrow \Ring[t^{-k}\Jet\check{\h}]$ is a (naively) graded algebra homomorphism. Let $\mathsf{M}_k$ be the composition of this pullback map with the embedding
$$\Ring[\Op^{\hbar^\natural}_k]\xrightarrow{\sim} \Ring[\widetilde{\Op}^{\hbar^\natural}_k]^{\Jet \check{N}}\hookrightarrow 
\Ring[\widetilde{\Op}^{\hbar^\natural}_k].$$

Note that by the construction the following identity holds:
\begin{equation}\label{eq:Miura_relation}
\operatorname{Ad}(t^\rho)^*\circ \mathsf{M}_k=\mathsf{M}_{k-1}\circ \operatorname{Ad}(t^\rho)^*, \forall k\in \Z. 
\end{equation}

\sss\label{SSS_opers_SL_2_example} A crucial example to understand is for $\check{G}=\operatorname{PGL}_2$. We can write an element of $t^{-k}\Jet(\check{h})$ as $\operatorname{diag}(a(t),-a(t))$, where $a(t)$ is a formal Laurent series $\sum_{j=-k}^\infty a_j t^j$. The function $h_n$ with $n\in \Z$ sends $a(t)$ to $2a_{-n}$. 
Under the map to $\Op^{\hbar^\natural}_k$, the element $\lambda\partial+\operatorname{diag}(a(t),-a(t))$ goes to 
$$\lambda\partial+\begin{pmatrix}0&b(t)\\1&0\end{pmatrix}, b(t):=a(t)^2+\lambda \partial a(t).$$
The free generator $x_n$ of $\Ring[\zeta+\partial+\Jet\check{S}]$ sends 
$b(t)=\sum_{j=0}^\infty b_j t^j$ to $b_{-n}$. So, for $k=1$, 
$\mathsf{M}_1(x_n)$ is exactly the element $\underline{X}_n^\natural$. In particular,  $\operatorname{im}\mathsf{M}_1=\,^I\!\A^{\hbar,\natural}$.

\sss We finish this section with the following lemma.

\begin{Lem}\label{Lem:Miura_image_contained}
For general $\check{\g}$ we have $\operatorname{im}\mathsf{M}_1\subset \,^I\!\A^{\hbar,\natural}$.
\end{Lem}
\begin{proof}
Consider the minimal Levi subgroup $\check{L}_i$ and its maximal unipotent subgroup $\check{N}_i$. Let $f_i,\check{S}_i$ be the analogs of $\check{f}, \check{S}$ for $L_i$. We have the analog of $\mathsf{M}$ for $L_i$, denote it by $\mathsf{M}^i$. Thanks to 
Section \ref{SSS_opers_SL_2_example}, for $k=1$, we have 
\begin{equation}\label{eq:Miura_i_image}
\operatorname{im}\mathsf{M}^i=\,^i\!\A^{\hbar,\natural}.
\end{equation}

Note that we have 
a $\Jet_{(k)} N_i$-equivariant embedding 
\begin{equation}\label{eq:jet_Levi_embedding}
\widetilde{\Op}^{\hbar^\natural}_k(\check{L}_i)\hookrightarrow \widetilde{\Op}^{\hbar^\natural}_k(\check{G}),
\lambda \partial+\check{f}_i+F\mapsto \lambda\partial+\check{f}+F.
\end{equation}
The embedding (\ref{eq:jet_Cartan_embedding}) factors into the composition of its analog for $H\subset L_i$ and (\ref{eq:jet_Levi_embedding}). It follows that $\mathsf{M}$ factors as $\mathsf{M}^i\circ \varphi$, where $\varphi$ is an algebra homomorphism $\Ring[\Op^\hbar_1(\check{G})]\rightarrow 
\Ring[\Op^\hbar_1(\check{L}_i)]$.  The claim of the lemma now follows from (\ref{eq:Miura_i_image}).
\end{proof}


\subsection{Proof of the main result}\label{SS_Miura_main}
\begin{proof}[Proof of Proposition \ref{Prop:Miura_isomorphism}]
The proof is in several steps. 

{\it Step 1}. Thanks to Lemma \ref{Lem:Miura_image_contained}, we have $\operatorname{im}\mathsf{M}_1\subset \,^I\!\A^{\hbar,\natural}$. 
Recall that $\mathsf{M}_1$ is $\Ring[\hbar^\natural]$-linear and graded. Its specialization at $\hbar^\natural=0$ is a homomorphism 
\begin{equation}\label{eq:jet_slice_Cartan}
\Ring[\Jet_{(1)}\check{S}]\rightarrow \Ring[t^{-1}\Jet\check{\h}].
\end{equation}
This homomorphism can be understood as follows. We note that the restriction of the quotient morphism $\check{\g}\rightarrow \check{\g}/\!/\check{G}$ to
$\check{S}$ is an isomorphism under our restrictions on $p$, see \cite[Theorem 3.2.2]{Riche_Kostant}. From here we conclude that the pullback 
under the composition 
$$\check{\h}\hookrightarrow \check{f}+\check{\mathfrak{b}}\twoheadrightarrow 
(\check{f}+\check{\mathfrak{b}})/\check{N}\xrightarrow{\sim} \check{S}\xrightarrow{\sim} \check{\g}/\!/\check{G}$$
is the Chevalley restriction map to be denoted by $\mathsf{res}$.

As in the proof of Lemma \ref{Lem:Kostant_section}, we can identify 
$\Jet_{(1)}\check{S}$ with $\Jet\check{S}$ and $t^{-1}\Jet\check{\h}$ with 
$\Jet\check{\h}$.
Under this identification, the homomorphism (\ref{eq:jet_slice_Cartan}) is nothing else but $\Jet(\mathsf{res})$. As we have seen in the proof of Lemma 
\ref{Lem:HC_injectivity}, the map $\Jet(\mathsf{res})$ is injective. 
It follows that $\mathsf{M}_1$ is injective. It remains to show that 
\begin{equation}\label{eq:image_coincidence}
\operatorname{im}\mathsf{M}_1=\,^I\!\A^{\natural,\hbar}. 
\end{equation}

{\it Step 2}.  
We claim that the cokernel of $\,^i\!\A^{\hbar,\natural}\hookrightarrow \Ring[t^{-1}\Jet\check{\h}][\hbar^\natural]$, equivalently of $\Ring[\Op_1^{\hbar^\natural}(L_i)]\hookrightarrow \Ring[t^{-1}\Jet\check{\h}][\hbar^\natural]$, has no $\hbar^\natural$-torsion. Using (\ref{eq:Miura_relation}) we reduce to showing that the cokernel of 
\begin{equation}\label{eq:fun_oper_embedding}
\Ring[\Op_0^{\hbar^\natural}(L_i)]\hookrightarrow \Ring[\Jet\check{\h}][\hbar^\natural]\end{equation}
Note that (\ref{eq:fun_oper_embedding}) is energy graded.
So, it is enough to show that each energy graded component of the cokernel of (\ref{eq:fun_oper_embedding}) is flat over $\Ring[\hbar^\natural]$. The energy graded components  of $\Ring[\Op_0^{\hbar^\natural}(L_i)], \Ring[\Jet\check{\h}][\hbar^\natural]$ are free finite rank modules over $\Ring[\hbar^\natural]$. Since the specialization of $\Ring[\Op_0^{\hbar^\natural}(L_i)]\hookrightarrow \Ring[\Jet\check{\h}][\hbar^\natural]$ to $\hbar^\natural=0$ is injective (by the straightforward modification of Step 1), the claim in the beginning of the step follows.

{\it Step 3}. Let $W_i\cong S_2$ denote the Weyl group of $\check{L}_i$
and $W$ denote the Weyl group of $\check{G}$. The quotient morphisms $\check{\h}\rightarrow \check{\h}/W_i, \check{\h}\rightarrow \check{\h}/W$ give rise to the embeddings $\Ring[\Jet(\check{\h}/W_i)],\Ring[\Jet(\check{\h}/W)]\hookrightarrow \Ring[\Jet\check{\h}]$. Note that $\Ring[\Jet(\check{\h}/W)]\subset \Ring[\Jet(\check{\h}/W_i)]$.
Note also that the image of $\mathsf{M}_1$ in $\Ring[\Jet \check{\h}]$   
contains
$\operatorname{im}\Jet(\mathsf{res})$, while the image of $\,^i\!\A^{\hbar,\natural}$ in $\Ring[\Jet\check{\h}]$
coincides with
$\Ring[\Jet(\check{\h}/W_i)]$.
Thanks to Step 2, we reduce (\ref{eq:image_coincidence}) to checking that 
\begin{equation}\label{eq:jet_intersection}
\Ring[\Jet(\check{\h}/W)]=\bigcap_{i\in I}\Ring[\Jet(\check{\h}/W_i)].
\end{equation}

{\it Step 4}. Let $(\check{\h}/W_i)^{reg}$ denote the locus in $\check{\h}/W_i$, where the natural morphism between quotients 
$\check{\h}/W_i\rightarrow \check{\h}/W$ is etale. This is a principal open subset. Note that for a smooth finite type scheme $X$ over $\operatorname{Spec}(\Ring)$, the morphism $\Jet X\rightarrow X$ is flat. 
Therefore we see that the natural morphisms 
$\Ring[\Jet(\check{\h}/W_i)]\rightarrow \Ring[\Jet\left((\check{\h}/W_i)^{reg}\right)]\rightarrow 
\Ring[\Jet(\check{\h}^{reg})]$ are embeddings. So (\ref{eq:jet_intersection}) 
will follow from
\begin{equation}\label{eq:jet_intersection_reg}
\Ring[\Jet(\check{\h}/W)]=\bigcap_{i\in I}\Ring[\Jet\left((\check{\h}/W_i)^{reg}\right)].
\end{equation}

{\it Step 5}. Note that if $Y\rightarrow X$ is an etale morphism of finite type smooth schemes over $\operatorname{Spec}(\Ring)$, then 
$\Jet Y\xrightarrow{\sim} Y\times_X \Jet X$. Applying this to $X=\check{\h}/W, Y=(\check{\h}/W_i)^{reg}$ or $Y=\check{\h}^{reg}$,  
we see that 
\begin{align*}
&\Ring[\Jet((\check{\h}/W_i)^{reg})]=
\Ring[(\check{\h}/W_i)^{reg}]
\otimes_{\Ring[\check{\h}/W]}\Ring[\Jet(\check{\h}/W)], \\
&\Ring[\Jet(\check{\h}^{reg})]=\Ring[\check{\h}^{reg}]
\otimes_{\Ring[\check{\h}/W]}\Ring[\Jet(\check{\h}/W)].
\end{align*}
Since $\Ring[\Jet(\check{\h}/W)]$ is a flat $\Ring[\check{\h}/W]$-module, tensoring with this module preserves the intersections of submodules, so
(\ref{eq:jet_intersection_reg}) will follow from 
\begin{equation}\label{eq:intersection_reg}
\Ring[\check{\h}/W]=\bigcap_{i\in I}\Ring[(\check{\h}/W_i)^{reg}].
\end{equation}
For the latter, note that the right hand side is contained in 
$\Ring[\check{\h}^{reg}/W]$ and consists of all elements there that do not have poles on $(\check{\h}/W)\setminus (\check{\h}^{reg}/W)$. Now 
(\ref{eq:intersection_reg}) follows.
\end{proof}

\begin{Rem}\label{rem:Miura_general_k}
The same argument shows that the homomorphism $\mathsf{M}_0:\Ring[\Op^{\hbar^\natural}_0]\rightarrow \Ring[\Jet \check{h}][\hbar^\natural]$ is injective and its image is the intersection of the images of the homomorphisms $\mathsf{M}_0^i$.
\end{Rem}

\begin{Rem}
Over a characteristic $0$ field, one can interpret $\mathsf{M}_1$ (specialized to $\hbar^\natural=1$) as follows. One identifies $\A^\natural$ with the algebra of functions on the space of generic Miura opers with regular singularities. This endows $\operatorname{Spec}(\A^\natural)$ with a free action of $\check{N}$ 
and the natural map $\operatorname{Spec}(\A^\natural)\rightarrow \operatorname{Spec}(\Op^1_1)$ is a quotient for this action, see \cite[Section 8.2]{Frenkel_loop}. This construction fails in characteristic $p$: for example, 
the map from the space of generic Miura opers to the space of opers can be shown to be invariant for $\check{N}[[t^p]]$.
\end{Rem}

\section{Computation of $V_\kappa(\g)^{\Jet G}$}\label{S_V_inv}
In this section we assume that $\Ring$ is a perfect field of characteristic $p>h$. We further assume that $\kappa\not\in \F_p$.
The goal of this section is to compute the commutative vertex algebra 
$V_\kappa(\g)^{\Jet G}$. Recall the embedding $\HC_H:
V_\kappa(\g)^{\Jet G}\hookrightarrow V_\kappa(\h)^{\Jet H}$ of commutative vertex algebras, see Section \ref{SSS_center_vertex_algebra_reduction}. We are going to describe the image. 
For this observe that we have a surjective algebra homomorphism
$\A\twoheadrightarrow V_\kappa(\h)^{\Jet H}$ sending $\tilde{h}_n$ to $\tilde{h}_n$ for $n<0$ and to zero else. Under the isomorphism 
$\A^\natural\xrightarrow{\sim}\A$, the quotient $V_\kappa(\h)^{\Jet H}$ gets identified with the polynomial algebra $\A^\natural_0$ in the variables $h^{i,\natural}_n$ with $n<0$. 

\subsection{Main result}
The following is the main result of this section:

\begin{Thm}\label{Thm:VA_invariants} \label{t:tungus}
Under the isomorphism $\A^\natural_0\cong V_\kappa(\h)^{\Jet H}$, the subalgebra $V_\kappa(\g)^{\Jet G}$ is sent to $\operatorname{im}\mathsf{M}_0$. 
\end{Thm}

Here is a scheme of proof of this theorem that will be given below in this section. 

{\it Step 1}: Recall, Theorem \ref{Thm:Verma_endomorphisms}, that $\,^I\!\A^\hbar$ acts on $\Delta^G_{\kappa\hbar}$ by $\Iw$-equivariant and $\widehat{U}^{\hbar}_\kappa(\hat{\g})$-linear endomorphisms. 
The vertex algebra $V^\hbar_\kappa(\g)$ is a quotient of $\Delta^G_{\kappa\hbar}$. We will see that the action of $\,^I\!\A^\hbar$ descends to $V^\hbar_\kappa(\g)$. We will further see that the action is in fact by $\Jet G$-equivariant endomorphisms. This will give rise to a graded $\Ring[\hbar]$-algebra homomorphism $\,^I\!\A^\hbar\rightarrow V^\hbar_\kappa(\g)^{\Jet G}$.

{\it Step 2}. We will analyze the specialization of the homomorphism 
$\,^I\!\A^\hbar\rightarrow V^\hbar_\kappa(\g)^{\Jet G}$ at $\hbar=0$.
We will see that it is surjective. That will allow us to determine the energy graded character of $V_\kappa(\g)^{\Jet G}$.

{\it Step 3}. Using the previous steps we will finish the proof of the theorem.

\subsection{Homomorphism $\,^I\!\A^\hbar\rightarrow V^\hbar_\kappa(\g)^{\Jet G}$}
We have an associative algebra isomorphism 
$V^\hbar_\kappa(\g)^{\Jet G}\xrightarrow{\sim} \operatorname{End}(V^\hbar_\kappa(\g))$, where by $\operatorname{End}$ we mean $\widehat{U}^\hbar_\kappa(\hg)$-linear and $\Jet G$-equivariant endomorphisms. See \cite[(3.3-3)]{Frenkel_loop} for a completely analogous statement.

Note that $V^\hbar_\kappa(\g)$ is a quotient of  $\Delta^G_{\kappa\hbar}$ viewed as an object in the category $\sO(\h^*)$ introduced in Section \ref{SS:endo_main}. Recall, Theorem \ref{Thm:Verma_endomorphisms}, that there is a unique algebra homomorphism $\theta:\,^I\!\A^\hbar\rightarrow \operatorname{End}_{\sO(\h^*)}(\Delta^G_{\kappa\hbar})$ such that the homomorphism $\varphi^B_\kappa: \Delta^G_{\kappa\hbar}\rightarrow \Wak^B_{\kappa\hbar}$ is $\,^I\!\A^\hbar$-semi-linear
(with respect to the automorphism $\mathsf{t}_B$ of $\A^\hbar$). 

The following is the main result of this section. 

\begin{Prop}\label{Prop:from_Verma_to_vacuum}
For every $a\in \,^I\!\A^\hbar$, the following claims hold:
\begin{enumerate}
\item The endomorphism $\theta(a)$ of $\Delta^G_{\kappa\hbar}$ descends to a $\widehat{U}^\hbar_\kappa(\hg)$-linear endomorphism of $V^\hbar_\kappa(\g)$, to be denoted by $\underline{\theta}(a)$.
\item The endomorphism $\underline{\theta}(a)$ is $\Jet G$-equivariant.
\end{enumerate}
\end{Prop}
\begin{proof}
We start by proving (1). The proof is in several steps.

{\it Step 1}. Let $v$ denote a generator in the highest weight space $\Delta^G_{\kappa\hbar}[0]$. The quotient $V^\hbar_\kappa(\g)$ of $\Delta^G_{\kappa\hbar}[0]$ is given by relations $xv=0$ for all $x\in \h$
and $e_iv=0$ for all $i\in I$, where $e_i$ denote the simple root vector (where the positive roots are those for $B^-$). Let $\check{\varpi}_i$ denote the fundamental coweight labelled $i$. We will show that the action of 
$\,^I\!\A^\hbar$ on $\Delta^G_{\kappa\hbar}$ preserves the submodules generated by $\check{\varpi}_iv$ and $e_i v$. 

{\it Step 2}. 
Let $P$ be the minimal parabolic $P_i$,
and $L$ and $N$ be its Levi subgroup and the maximal unipotent subgroup. 
Let $v'$ and $v''$ be generators of 
$\widehat{D}_\hbar(\Loop N)/I_{>0}(\widehat{D}_\hbar(\Loop N))][0]$ and $\Delta_{\kappa\hbar}^{L}[0]$ so that $v'\otimes v''$ is a generator 
of $\Wak^P_{\kappa\hbar}$. We can assume that $\varphi^P(v)=v'\otimes v''$. 
One can get a formula for $\ffr_P(e_it^{-1})$ using Lemma 
\ref{Lem:ffr_Cartan_images}: $\ffr_P(e_it^{-1})$ is represented as the sum of two summands, one in $\CDO(N)$, the other in $V_\kappa(\lf)$, and the latter is $|\varnothing\rangle\otimes (e_it^{-1})|\varnothing\rangle$. Consequently from looking at the $\Lambda^a$-grading, one sees that $e_i(v'\otimes v'')=v'\otimes e_i v''$. 

{\it Step 3}. 
For $z\in \Ring$, we write $\Delta^G_{\kappa\hbar}(z)$ and $\Wak^P_{\kappa\hbar}(z)$ for the quotients of $\Delta^G_{\kappa\hbar}, \Wak^P_{\kappa\hbar}$ by  
$\check{\varpi}_i-z\hbar$ under the right action (and we write 
$\Delta^G_{\kappa\hbar,i}(z)$ if we want to indicate the dependence on $i$). Then $\varphi^P$ induces 
$\varphi^P(0):\Delta^G_{\kappa\hbar}(0)\rightarrow \Wak^P_{\kappa\hbar}(0)$. 
Under the action of $\widehat{U}^\hbar_{\kappa}(\hg)$, the vector $e_iv\in \Delta^G_{\kappa\hbar}(0)$ generates a copy of $\Delta^G_{\kappa\hbar}(2)$.
Similarly, under the action of $\widehat{D}^\hbar(\Loop N)\widehat{\otimes}_{\Ring[\hbar]}\widehat{U}^\hbar_\kappa(\hat{\lf})$, 
the vector $v'\otimes e_iv''\in \Wak^P_{\kappa\hbar}(0)$ generates a copy of $\Wak^P_{\kappa\hbar}(2)$. 

{\it Step 4}. We claim that 
\begin{itemize}
\item[(*)] $\,^I\!\A^\hbar$ preserves $\Delta^G_{\kappa\hbar}(2)\subset \Delta^G_{\kappa\hbar}(0)$. 
\end{itemize}
Once we know (*) for all $i$, the claim of (1) follows. Indeed, the action of $\,^I\A^\hbar$ preserves the submodule $\Delta^G_{\kappa\hbar}\h$. And $V_{\kappa}^\hbar(\g)$ is the quotient of $\Delta^G_{\kappa\hbar}/ \Delta^G_{\kappa\hbar}\h$ by the sum of images of $\Delta^G_{\kappa\hbar,i}(2)$ over all $i$.

{\it Step 5}. We proceed to proving (*).
Note that $\varphi^P(0)$ restricts to a homomorphism 
$\Delta^G_{\kappa\hbar}(2)\rightarrow \Wak^P_{\kappa\hbar}(2)$ that sends a generator to a generator. So this restriction, to be denoted by $\varphi^P(2)$ coincides with the specialization of $\varphi^P$ (possibly after multiplying by an invertible element of $\Ring$). The action of $\,^i\A^\hbar$ on $\Wak^P_{\kappa\hbar}(0)$, by the construction, comes from a homomorphism to the center of $\widehat{D}^\hbar(\Loop N)\widehat{\otimes}_{\Ring[\hbar]}\widehat{U}^\hbar_\kappa(\hat{\lf})$.
So it preserves the submodule $\Wak^P_{\kappa\hbar}(2)$. We note that $\,^i\!\A^\hbar$ is stable under the automorphism $\mathsf{t}_P$. It follows that  $\Wak^P_{\kappa\hbar}(2)$ is also stable under the action of $\,^I\!\A^\hbar$. 

{\it Step 6}. Let $\h^{da}_{\kappa\hbar}(2)$ be the hyperplane 
in $\h^{da}_{\kappa\hbar}$ defined by $\check{\varpi}_i=2\hbar$. 
Let $\sx$ be the generic point of this hyperplane (i.e., 
$\bk=\operatorname{Frac}(\h^{da}_{\kappa\hbar}(2))$ and the homomorphism $\Ring[\h^{da}_{\kappa\hbar}(2)]\rightarrow 
\bk$ is the canonical one). We make two observations. 
First 
\begin{itemize}
\item[(**)]
$\Delta^G_{\kappa\hbar}(0)/\Delta^G_{\kappa\hbar}(2)$ is free over $\Ring[\h^{da}_{\kappa\hbar}(2)]$. 
\end{itemize}
This is because this module is parabolically induced from the $U^\hbar(\lf)$-module $\Ring[\h^{da}_{\kappa\hbar}(2)]$ (where $\mathbf{1}$ acts by 
$\kappa\hbar$). Second, the point $\sx$ is $P$-generic in the sense of 
Definition \ref{defi:P_generic}.  Theorem \ref{Thm:Verma_vs_Wakimoto} implies that 
\begin{itemize}
\item $\varphi^P(2):\Delta^G_{\sx}\rightarrow 
\Wak^P_{\sx}$ is an isomorphism. 
\end{itemize}
It follows from (**) that $\,^I\!\A^\hbar$ preserves $\Delta^G_{\sx}$ (viewed as a submodule in a suitable localization of $\Delta^G_{\kappa\hbar}(0)$). Thanks to (**) we see that $\Delta^{G}_{\kappa\hbar}(2)=\Delta^G_{\sx}\cap 
\Delta^{G}_{\kappa\hbar}(0)$. Hence we see that $\Delta^{G}_{\kappa\hbar}(2)$
is fixed by $\,^I\!\A^\hbar$ finishing the proof of (1).

Now we proceed to (2). Note that $\underline{\theta}(a)$ is an $\Iw$-equivariant endomorphism of $V^\hbar_\kappa(\g)$. We need to show that it is $\Jet G$-equivariant. Note that since $V^\hbar_\kappa(\g)$ is a cyclic 
$\widehat{U}_{\kappa\hbar}(\hg)$-module, the space of its $\Iw$- (resp., $\Jet G$-) equivariant endomorphisms is identified with 
$V^\hbar_\kappa(\g)^{\Iw}$ (resp., $V^\hbar_\kappa(\g)^{\Jet G}$). We need to show that $V_\kappa(\g)^{\Jet G}=V_\kappa(\g)^{\Iw}$. Since $V_\kappa(\g)$ is a rational representation of $\Jet G$, our claim amounts to the classical claim that any $B^-$-invariant vector in a rational representation of $G$ is $G$-invariant.
\end{proof}

\subsection{The graded character of $V_\kappa(\g)^{\Jet G}$}
The goal of this section is to prove the following result.

\begin{Prop}\label{Prop:graded_characters}
The energy graded characters of $\Ring[(\Jet \g^*)^{(1)}]^{(\Jet G)^{(1)}}$
and $V_\kappa(\g)^{\Jet G}$ are the same.
\end{Prop}
\begin{proof}
The proof is in several steps. 

{\it Step 1}.
Set $\,^I\!\A^0:=\,^I\!\A^\hbar/(\hbar)$. This algebra acts by endomorphisms on $\Delta^G_0:=S(\mathfrak{b}^-\oplus t^{-1}\g[t^{-1}])$. Moreover, the generic point of $\h^*\subset \h^{da}$ is also generic in the sense of Definition \ref{defi:P_generic}, so the homomorphism 
$\varphi^B:\Delta^G_0\rightarrow \Wak^B_0$ is generically injective. 
The action of $\,^I\A^0$ on $\Delta^G_0$ is uniquely characterized  
by the condition that the homomorphism $\Delta^G_0\rightarrow \Wak^B_0$
is $\,^I\!\A^0$-equivariant.

{\it Step 2}. 
Note that $S(\g[t^{-1}]^{(1)})^{(\Jet G)^{(1)}}$ acts on $\Delta^G_0$ by $\Iw$-equivariant endomorphisms. Consider the Chevalley restriction map
$S(\g[t^{-1}]^{(1)})^{(\Jet G)^{(1)}}\rightarrow S(\h[t^{-1}]^{(1)})$. This gives an action of $S(\g[t^{-1}]^{(1)})^{(\Jet G)^{(1)}}$ on $\Wak^B_0$. The homomorphism $\Delta^G_0\rightarrow \Wak^B_0$ is equivariant for the action of $S(\g[t^{-1}]^{(1)})^{(\Jet G)^{(1)}}$.
Proposition \ref{Prop:Miura_isomorphism} and its proof (Steps 2 and 3) show that $\,^I\!\A^0$
embeds into $S(\h[t^{-1}]^{(1)})$ and the image coincides with that of the Chevalley restriction map. This shows that $\,^I\!\A^0$ acts on $\Delta^G_0$ via its isomorphism with $S(\g[t^{-1}]^{(1)})^{(\Jet G)^{(1)}}$.

{\it Step 3}. It follows from Step 2 that $\,^I\!\A^0$ acts on $\gr V_\kappa(\g)=S(t^{-1}\g[t^{-1}])$ via its natural homomorphism to
$S\left((t^{-1}\g[t^{-1}])^{(1)}\right)^{(\Jet G)^{(1)}}$. This homomorphism is surjective. On the other hand, the homomorphism $\,^I\!\A^0\rightarrow 
S(t^{-1}\g[t^{-1}])$ factors through $\gr V_\kappa(\g)$. This realizes $S\left((t^{-1}\g[t^{-1}])^{(1)}\right)^{(\Jet G)^{(1)}}$ as an energy graded subquotient of $\gr V_\kappa(\g)$. On the other hand, Section \ref{SSS_center_vertex_algebra_reduction} yields an energy graded inclusion in the opposite direction. The claim of the lemma follows.
\end{proof}

Now we are ready to complete the proof of Theorem \ref{Thm:VA_invariants}.
\begin{proof}[Proof of Theorem \ref{Thm:VA_invariants}]
By Lemma \ref{Lem:HC_injectivity} and (\ref{eq:HC_transitivity}), we see that 
$\HC_H$ embeds $V_\kappa(\g)^{\Jet G}$ into $\bigcap_{i\in I}\operatorname{im}\HC_{L_i}$, where $L_i$ is the minimal Levi associated to $i$. Under the isomorphism of $V(\h)^{\Jet H}$ and $\A^\natural_0$, the subalgebra $\operatorname{im} \HC_{L_i}$ is identified with $\operatorname{im}\mathsf{M}^i_0$. By Remark \ref{rem:Miura_general_k},
$\bigcap_{i\in I}\operatorname{im}\mathsf{M}^i_0=\operatorname{im}\mathsf{M}_0$. So under the isomorphism $\A^\natural_0\cong V_\kappa(\h)^{\Jet H}$, the subalgebra $V_\kappa(\g)^{\Jet G}$ is injectively mapped to $\operatorname{im}\mathsf{M}_0$.

The image of
$\operatorname{im}\mathsf{M}_0$ under the isomorphism $\A^\natural_0\cong V_\kappa(\h)^{\Jet H}$
has the same energy graded character as $S((t^{-1}\g[t^{-1}])^{(1)})^{(\Jet G)^{(1)}}$. By Proposition \ref{Prop:graded_characters}, $V_\kappa(\g)^{\Jet G}$ has the same energy graded character. The claim of the theorem follows. 
 \end{proof}

\section{Coordinate invariant description of the Harish--Chandra center}\label{S_coord_indep}

In this section, we let $\Ring$ be a perfect field of characteristic $p>h$, and {assume $\kappa$ is nondegenerate}.

In performing the relevant calculations for Theorem \ref{Thm:VA_invariants}, we have throughout fixed a coordinate $t$ on our formal disc $\D$. Our goal in this section is to record the equivariance of Theorem \ref{Thm:VA_invariants} with respect to changes of coordinates, i.e., the action of $\Aut(\sD)$. {We refer the reader to Section \ref{sss:autd} for a reminder of the definition of the latter group ind-scheme.} 

{ In addition, to speak of the Langlands dual group $\check{G}$ of $G$, we fix pinnings on $G$ and $\check{G}$.}

{\iffalse Without loss of generality we may assume the derived subgroup of $\check{G}$ of adjoint type. {\bf Justify this}} \fi 

\subsection{Coordinate invariance for the free field realization}  

\sss{} Let us first recall from Section \ref{ss:diglets} the free field homomorphism associated to the Borel subgroup, and state its construction in a manner which is visibly $\Aut(\sD)$-equivariant. 

{To identify the target of the free field realization in Langlands dual terms, we will need to account for a certain $\rho$-shift. For this reason, we will assume that the derived subgroup of $\check{G}$ is of adjoint type from now until Section \ref{ss:hccenD}. }

\sss For a scheme $\mathscr{X} \rightarrow \D$ we will denote by $\Jet \mathscr{X}$ its scheme of sections, i.e., the functor which sends an $\Ring$-algebra $\Sing$ to the set of sections of { $\mathscr{X}_\Sing \rightarrow \D_\Sing$} \iffalse is represented by the scheme $\mathscr{Y}_O$?}\fi. In particular, for a scheme $\mathscr{X}$ over $\Ring$, note that this definition applied to its base change $\mathscr{X} \underset{\Spec \Ring} \times \sD$  agrees with the previous definition of its jet scheme.

\sss{} If $\mathscr{X}$ is equipped with an action of $\Aut(\sD)$ compatible with the projection map to $\sD$, its jet scheme $\Jet \mathscr{X}$ carries a tautological action of $\Aut(\sD)$.  Similarly, for any affine scheme $\mathscr{Y} \rightarrow \pD$ with a compatible action of $\Aut(\pD)$, its loop ind-scheme $\Loop \mathscr{Y}$ of sections carries a tautological action of $\Aut(\pD)$ and hence by restriction $\Aut(\sD)$. Moreover, if $\mathscr{Y}$ is the base change of $\mathscr{X}$, compatibly with the actions of $\Aut(\sD)$, the natural restriction map $\Jet \mathscr{X} \rightarrow \Loop \mathscr{Y}$ is $\Aut(\sD)$-equivariant.

In particular, for an algebraic group $F$ with Lie algebra $\ff$, we obtain that $\Aut(\sD)$ acts compatibly on the arc Lie algebra $\Jet \ff$ and the loop Lie algebra $\Loop \ff$ by Lie algebra automorphisms, where we view $\Jet \ff$ and $\Loop \ff$ as Lie algebra $\Ring$-spaces. 

\sss For any Kac--Moody level $\kappa$, recall that the underlying vector space of the central extension is canonically split 
$$\hat{\ff}_\kappa \simeq \Loop \ff \oplus \Ring.$$
As the Kac--Moody two cocycle is manifestly coordinate invariant, if we have $\Aut(\pD)$ act on $\Loop \ff$ in its natural way and trivially on $\Ring$, the induced action on $\hat{\ff}_\kappa$ is again an action by Lie algebra automorphisms.

\sss \label{sss_jetscdoaction} In particular, the vacuum algebra 
$$V_\kappa(\ff) \simeq \widetilde{U}_\kappa(\hat{\ff}) \underset{\widetilde{U}( \Jet \ff \oplus \Ring )} \otimes \Ring $$
inherits by its construction a canonical action of $\Aut(\sD)$; we recall the induced action of $\Aut(\sD)$ on the appearing completed enveloping algebras follows from the functoriality established in Section \ref{sss:completedenvelope}. By the same reasoning, the level $\kappa$ chiral differential operators 
$$\CDO_\kappa(F) \simeq \widetilde{U}_\kappa(\hat{\ff}) \underset{\widetilde{U}( \Jet \ff \oplus \Ring )} \otimes \Ring[\Jet F] $$
carries a canonical action of $\Aut(\sD)$, compatible with the embedding of left invariant vector fields $\iota_\ell: V_\kappa(\ff) \rightarrow \CDO_\kappa(F)$. 

Moreover, $\CDO_\kappa(F)$ carries compatible left and right translation actions of $\Jet F$, i.e., an action of $$\Aut(\sD) \ltimes (\Jet F^\ell \times \Jet F^r).$$Namely, the left action of $\Jet F^\ell \simeq \Jet F$ is induced by the adjoint action on $\widetilde{U}_\kappa(\hat{\ff})$ and the left translation action on $\Ring[\Jet F]$, and the right action of $\Jet F^r \simeq \Jet F$ is induced by the trivial action on $\widetilde{U}_\kappa(\hat{\ff})$ and the right translation action on $\Ring[\Jet F]$.

\sss{} \label{re:rightCDO} We remind the reader of the basic complication that, for general $F$, the subalgebra of $\CDO_\kappa(F)$ corresponding to right invariant vector fields, as constructed in Section \ref{SSS_iota_r}, is {\em not} $\Aut(\sD)$-equivariantly isomorphic to $V_{-\kappa}(\ff)$. However, we have the following result, which suffices for our purposes. 

\begin{Lem} \label{l:rightvectfields} Suppose $F$ is unimodular, i.e., the adjoint action of $\ff$ on $\det \ff^*$ is trivial. Then the map 
$$ V_{-\kappa}(\ff) \rightarrow \CDO_\kappa(F), $$
    as constructed in Section \ref{SSS_iota_r} with respect to a chosen coordinate $t$, is $\Aut(\D)$-equivariant.  
\end{Lem}

\begin{proof} Recall that $\Aut(\sD)$ is the product of (i) the formal neighborhood of $0$ in $\mathbb{A}^1$, (ii) one copy of $\mathbb{G}_m$, and (iii) a countable product of copies of $\mathbb{A}^1$, cf. Section \ref{sss:autd}. 

Therefore, when defined over $\mathbb{Z}$, $\Aut(\D)$ is an ind-affine group ind-scheme expressible as a filtered colimit of affine schemes flat over $\mathbb{Z}$. Similarly,  its comodules  $V_{-\kappa}(\ff)$ and $\CDO_\kappa(F)$, for variable $\kappa$, are similarly flat over $\mathbb{Z}[\frac{1}{2}]$, cf. Lemma \ref{Lem:jet_smoothness_flatness} and Section \ref{ss:defcdo}. 

By base change, it is therefore enough to address the claim when the relevant vertex algebras are defined over a field of characteristic zero, where it is well known, cf.  Section 3.6 of \cite{Arkhipov_Gaitsgory}.
\end{proof}

\sss As in Section \ref{ss:diglets}, let us consider the free field realization associated to the Borel. Namely, exactly as in the construction for $\CDO_\kappa(G)$, $\CDO_\kappa(G^\diamond)$ carries by construction a canonical action of $$\Aut(\sD) \ltimes (\Jet B^\ell \times \Jet B^{-, r}).$$ 
Moreover, by construction, the embedding$$\CDO_\kappa'(B) \rightarrow \CDO_\kappa(G^\diamond)$$is equivariant for the action of $$\Aut(\sD) \ltimes (\Jet B^\ell \times \Jet H^r);$$we recall the superscript `$'$' in $\CDO_\kappa'(B)$ is simply notation to record a change in conventions regarding critical shifts.

\sss In particular, the construction of Section \ref{ss:diglets} produced the target of the free field realization as a subquotient of $\CDO_\kappa(G^\diamond)$; let us give a manifestly $\Aut(\sD)$-equivariant presentation of this subquotient.  Namely, we may first pass to the $\Jet N^{-, r}$-invariants, i.e., we consider the subalgebra  
$$\CDO_\kappa(G^\diamond)^{\Jet N^{-, r}} \hookrightarrow \CDO_\kappa(G^\diamond).$$
Thanks to Lemma \ref{l:rightvectfields}, we may consider the $\Aut(\D)$-equivariant composition 
$$V(\mathfrak{n}^-) \hookrightarrow V_{-\kappa}(\fg) \rightarrow \CDO_\kappa(G)^{\Jet N^{-, r}} \rightarrow \CDO_\kappa(G^\diamond)^{\Jet N^{-, r}},$$
and in particular may view $\CDO_\kappa(G^\diamond)^{\Jet N^{-, r}}$ as a module over $\Loop \mathfrak{n}^-$. With this, the constructed subquotient was simply the $\Jet H^r$-invariants in  the $\Loop \mathfrak{n}^-$-coinvariants, i.e., 
$$[(\CDO_\kappa(G^\diamond)^{\Jet N^{-, r}})_{\Loop \mathfrak{n}^-}]^{\Jet H^r} \hookrightarrow  (\CDO_\kappa(G^\diamond)^{\Jet N^{-, r}})_{\Loop \mathfrak{n}^-} \twoheadleftarrow \CDO_\kappa(G^\diamond)^{\Jet N^{-, r}} \hookrightarrow \CDO_\kappa(G^\diamond).$$
We recall this subquotient was mapped isomorphically onto by $\CDO_\kappa'(B)^{\Jet H^r}$. In particular, we may regard the free field realization as an $\Aut(\sD) \ltimes \Jet B$-equivariant map 
$$\ffr_B: V_\kappa(\fg) \rightarrow \CDO_\kappa'(B)^{\Jet H^r}.$$
For ease of notation in what follows, we set 
$$\diglet := \CDO_\kappa'(B)^{\Jet H^r},$$and refer to the $\Jet B^\ell$ action on $\diglet$ simply as $\Jet B$.






\sss Note that the subgroup $(\Jet \tilde{N}) \rtimes H\subset \Jet B$} is $\on{Aut}(\D)$-stable. In particular, we obtain an $\Aut(\sD)$-equivariant map 
$$\ffr_B: V_\kappa(\fg)^{\Jet G} \rightarrow \diglet^{(\Jet \tilde{N}) \rtimes H}.$$
As in Proposition \ref{Prop:affine_HC}, a choice of a coordinate $t$ yields an isomorphism 
\begin{equation}\iota_t: \diglet^{(\Jet \tilde{N})\rtimes H} \xrightarrow{\sim} V_\kappa(\fh). \label{e:autdeq}\end{equation}
As alluded to in the discussion of Section  \ref{re:rightCDO}, this isomorphism is {\em not} $\on{Aut}(\D)$-equivariant, for its natural action on both sides, but instead must be corrected.  We now turn to this correction. 


\sss{} Set $\check{\mathfrak{h}}:=\mathfrak{h}^*$.
As a preparatory step, note that the isomorphism $-\kappa: \fh \simeq \check{\fh}$ prolongs to an $\on{Aut}(\D)$-equivariant isomorphism of vertex algebras 
$$\iota: V_\kappa(\fh) \xrightarrow{\sim} V_{\check{\kappa}}(\check{\fh}),$$
where $\check{\kappa}$ is the symmetric bilinear form on $\check{\mathfrak{h}}$ given by $\kappa^{-1}$.

Via the left presentation, this vertex algebra carries an action of $\Jet H$, and via the right presentation, an action of $\Jet \check{H}$, and we have the following basic assertion, which will be convenient to use below. 

\begin{Lem}The actions of $\Jet H$ and $\Jet\check{H}$ on $V_\kappa(\fh) \simeq V_{\ckappa}(\cft)$ commute.  \label{l:twoHeiscommute}
\end{Lem}
\begin{proof} This follows immediately from Lemma \ref{Lem:KM_cocycle_group}. \end{proof}

From the actions of $\Jet H$ and $\Jet\check{H}$, we obtain two Harish--Chandra centers
$$V_\kappa(\fh)^{\Jet H} \hookrightarrow V_\kappa(\fh) \simeq V_{\check{\kappa}}(\check{\fh}) \hookleftarrow V_{\check{\kappa}}(\check{\fh})^{\Jet\check{H}}.$$
Similarly, via either presentation, this vertex algebra acquires a corresponding $p$-center
$$\mathfrak{z}_{Fr}(V_\kappa(\fh)) \hookrightarrow V_\kappa(\fh) \simeq V_{\ckappa}(\check{\mathfrak{h}}) \hookleftarrow \mathfrak{z}_{Fr}(V_{\ckappa}(\check{\mathfrak{h}})).$$
A basic instance of quantum geometric Langlands duality in the current setting is that, under $\iota$, the Harish--Chandra and $p$-centers are not preserved, but rather are swapped. 

\begin{Prop} \label{p:pvsHCswap}The isomorphism $\iota$ restricts to isomorphisms 
$$\iota: V_\kappa(\fh)^{\Jet H} \simeq \frz_{Fr}(V_{\ckappa}(\check{\mathfrak{h}})) \quad \text{and} \quad \iota: \frz_{Fr}(V_\kappa(\fh)) \simeq V_{\ckappa}(\check{\mathfrak{h}})^{\Jet \check{H}}.$$
\end{Prop}

\begin{proof} By symmetry, it is enough to establish the first equality of subalgebras, and to do so we may moreover choose a coordinate $t$ on $\D$.

For an element $\lambda$ in the cocharacter lattice $\Lambda$ of $\check{H}$, we, by abuse of notation, denote again by $\lambda$ the corresponding element of $\check{\mathfrak{h}}$. The vertex subalgebra $\frz_{Fr}(V_{\ckappa}(\check{\mathfrak{h}}))$ (defined in Section \ref{SSS_vertex_p_center} in a more general situation) is generated by the elements 
$$(\lambda_{-1}^p - \lambda_{-p})|\varnothing\rangle, \quad \quad \lambda \in \Lambda.$$
If we apply $\iota^{-1}$ we then obtain the elements
$$(-\ckappa(\lambda)_{-1}^p + \ckappa(\lambda)_{-p})|\varnothing\rangle, \quad \quad \lambda \in \Lambda,$$
but these are precisely the generators of $V_{\kappa}(\fh)^{\Jet H}$, cf. Proposition \ref{Prop:HC_center_abelian}. \end{proof}

\sss{} To explain the necessary correction to the failure of equivariance of \eqref{e:autdeq}, let us write $\omega$ for the canonical bundle of $\D$. If we write $\rho$ for the half sum of the positive roots of $G$, by our assumption that the derived subgroup of $\check{G}$ is of adjoint type, we may view this as a cocharacter of $\check{H}$. We may therefore consider the $\check{H}$-bundle $\omega^\rho$ over $\D$. 

Consider the scheme of sections $\Jet \omega^\rho$, and note that, by the trivializability of $\omega^\rho$ as an $\check{H}$-bundle, $\Jet \omega^\rho$ is a trivializable $\Jet\check{H}$-torsor. 

Note that any choice of coordinate $t$ yields a trivialization of $\omega$ on $\D$, namely the section $dt$. This choice in turn induces by transport of structure a trivialization of the $\Jet\check{H}$-torsor $\Jet \omega^\rho$. 

\sss{} In particular,  consider the twist 
$$V_{\ckappa'}(\check{\mathfrak{h}}) := V_{\ckappa}(\check{\mathfrak{h}}) \overset{\Jet\check{H}} \times\Jet \omega^\rho.$$
{ This vertex algebra is isomorphic to $V_{\check{\kappa}}(\check{\mathfrak{h}})$ but not $\on{Aut}(\D)$-equivariantly.}
Any choice of coordinate yields a composite equivalence 
\begin{equation} \label{e:cheat} \diglet^{\Jet \tilde{N}\rtimes H} \simeq V_\kappa(\fh) \simeq V_{\ckappa}(\check{\mathfrak{h}}) \simeq V_{\ckappa'}(\check{\mathfrak{h}}).\end{equation}
While the intermediate isomorphisms are not coordinate invariant, we have the following.

\begin{Prop} \label{p:coordinvarianceofcheat} The composite isomorphism \eqref{e:cheat} is $\on{Aut}(\D)$-equivariant, and in particular independent of the coordinate $t$. 
\end{Prop}

\begin{proof} It is enough to verify the case of $\Ring = \Z[\check{\kappa}]$, the general case is obtained from here by base change. 

As $\on{Aut}(\D)$ is a countable  inductive limit of flat affine $\mathbb{Z}$-schemes, and similarly $V_{\ckappa'}(\check{\mathfrak{h}})$ is an inductive limit of free $\mathbb{Z}$-modules stable under the action of $\on{Aut}(\D)$, by faithful flatness we may  reduce to the case of $\Ring=\mathbb{Q}(\check{\kappa})$, where it is known. More precisely, in the latter case, it is standard that the canonical insertion map into the BRST complex
$$\diglet^{\Jet \tilde{N}\rtimes H} \rightarrow \diglet \rightarrow \on{C}^{\frac{\infty}{2} + *}(\Loop \tilde{\mathfrak{n}}, \Jet\tilde{\mathfrak{n}}, \diglet)  $$ 
is a quasi-isomorphism, and the identification of the latter with $V_{\ckappa'}(\check{\mathfrak{h}})$ is well known, see e.g.    \cite[Section 10.2]{FG}. 
\end{proof}

In particular, we deduce the following consequence for the original free field realization.

\begin{Cor} \label{c:freefieldcoordinv} We have an $\Aut(\sD)$-equivariant equivalence 
$$\CDO'_\kappa(B)^{\Jet H^r} \simeq \CDO(\tilde{N}) \otimes V_{\ckappa'}(\check{\fh}).$$
In particular, the free field realization associated to the Borel is an $\Aut(\sD)$-equivariant map 
$$\on{ffr}_B: V_\kappa(\fg) \rightarrow \CDO(\tilde{N}) \otimes V_{\ckappa'}(\check{\fh}).$$
\end{Cor}

\subsection{The scheme of $\kappa$-connections} 
\label{ss:kconnections}
\sss Recall, Section \ref{SSS_HC_vertex1}, that the free field realization restricts to a map $$\ffr: V_\kappa(\fg)^{\Jet G} \hookrightarrow \diglet^{\Jet B} = (\diglet^{\Jet \tilde{N} \rtimes H})^{\Jet H}.$$
Our next task, which is accomplished in Proposition \ref{p:digjet} below, is to identify the target more explicitly. To do so, we first need some preparation. 

\sss Consider a smooth algebraic group $K$ over $\Ring$, and write $\fk$ for its Lie algebra.  Fix a $K$-bundle $\pi: \sP \rightarrow \sD$. As notation, given a scheme $X$ with an action of $K$, let us write $$X_\sP := \sP \overset{K} \times X$$for the associated $X$-bundle over $\sD$. Let us write $\on{T}_\sD \rightarrow \sD$ for the vector bundle of (continuous) $\Ring$-linear derivations on $\sD$, and similarly for $\on{T}_\sP \rightarrow \sP$, so that we have the Atiyah exact sequence 
\begin{equation} \label{e:atiyah} 0 \rightarrow \fk_\sP \rightarrow \pi_*(\on{T}_\sP)^K \rightarrow \on{T}_\sD \rightarrow 0.\end{equation}
Recall that a connection on $\sP$ is exactly a splitting of \eqref{e:atiyah}.

Let us write $\Conn(\sD, \sP)$ for the scheme of connections on $\sP$. I.e., given a test $\Ring$-algebra $\Sing$, an $\Sing$-point of $\Conn(\sD, \sP)$ is a connection on the base change $\sP_\Sing \rightarrow \sD_\Sing$. 

\sss Let us recall the basic structures this space of connections carries. If we consider $K$, with its adjoint action, and write $\on{Aut}(\sP) \rightarrow \sD$ for the group scheme of automorphisms of $\sP$, we have a canonical isomorphism $$\on{Aut}(\sP) \simeq K_{\sP},$$
and in particular an $\Sing$-point of the jet scheme $\Jet K_\sP$ is exactly an automorphism of $\sP_\Sing \rightarrow \sD_\Sing$. 

Then $\Conn(\sD, \sP)$ is canonically a torsor over $\Jet (\fk_\sP \otimes \omega)$, and carries a natural action by $\Jet K_\sP$. Moreover, if $\sP$ is given an $\Aut(\sD)$-equivariant structure, then $\Conn(\sD, \sP)$ similarly acquires an action of $\Aut(\sD) \ltimes \Jet K_\sP$. 

\sss Consider the case of the trivial bundle $\sP^\circ \rightarrow \sD$. In this case, it comes with a canonical connection, which we denote by $d$, and so we have a canonical identification $$\Conn(\sD, \sP^\circ) \simeq \Jet (\fk \otimes \omega); \quad \quad d + \alpha \mapsto \alpha, \quad \quad \alpha \in \Jet(\fk \otimes \omega)(\Sing).$$
We recall this is {\em not} $\Jet K$-equivariant. Indeed, if we denote the standard action of $\Jet K$ on $\Jet (\fk \otimes \omega)$ on $\Sing$-points by 
$$ k \cdot \alpha =: \on{Ad}_{k}(\alpha), \quad \quad k \in \Jet K (\Sing), \quad \alpha \in \Jet (\fk \otimes \omega) (\Sing),$$
we instead on $\Conn(\sD, \sP^\circ)$ have the formula 
$$ k \cdot (d + \alpha) = d + \on{Ad}_k(\alpha) + k^{-1} \cdot dk.$$

In particular, if we write $\Ring[\Conn(\sD, \sP^\circ)]^{\leqslant 1}$ for the  space of functions on this affine space of degree at most one, we have a canonically split, and in particular $\Aut(\sD)$-equivariantly split, short exact sequence 
$$0 \rightarrow \Ring \cdot 1 \rightarrow \Ring[\Conn(\sD, \sP^\circ)]^{\leqslant 1} \rightarrow  \Loop \fk^* / \Jet \fk^*  \rightarrow 0,$$
where we have used the residue perfect pairing 
$$\langle -, - \rangle: \Loop \fk^* / \Jet \fk^* \otimes \Jet(\fk \otimes \omega) \rightarrow \Ring.$$
In particular the natural action of $\Jet K$ on this split extension is determined by the formula 
\begin{equation} \label{e:actconn} k \cdot X = \on{Ad}_k(X) - \langle X, k^{-1} \cdot dk \rangle \cdot 1, \quad \quad k \in \Jet K(\Sing), \quad X \in \Loop \fk^* / \Jet \fk^* (\Sing) .\end{equation}

\sss To proceed, we need a variation on the above moduli spaces in the presence of a Kac--Moody level
$$\kappa \in \Sym^2(\fk^*)^K.$$
Namely, if we view $\kappa$ as a $K$-equivariant map $\kappa: \fk \rightarrow \fk^*$, let us push out the Atiyah extension \eqref{e:atiyah} along {\em minus} $\kappa$ to obtain an exact sequence of vector bundles on $\sD$
\begin{equation}  \xymatrix{0 \ar[r] & \fk_\sP \ar[d]_{-\kappa} \ar[r] & \pi_*(\on{T}_\sP)^K \ar[r] \ar[d] & \on{T}_\sD \ar[d]^{\on{id}} \ar[r] & 0 \\ 0 \ar[r] & \fk^*_\sP \ar[r] & \mathscr{E}_{\sP, \kappa} \ar[r] & \on{T}_\sD \ar[r] & 0.  } \end{equation}

\begin{defi} We define a {\em $\kappa$-connection} on $\sP$ to be a splitting of the sequence
    \begin{equation} \label{e:atiyahkappa} 0 \longrightarrow \fk_\sP^* \longrightarrow \mathscr{E}_{\sP, \kappa} \longrightarrow \on{T}_\sD \longrightarrow 0,\end{equation}
    and denote the scheme of $\kappa$-connections on $\sP$ by $\Conn_\kappa(\sD, \sP).$
\end{defi}

Explicitly, an $\Sing$-point of $\Conn_\kappa(\sD, \sP)$ consists of a splitting of \eqref{e:atiyahkappa} over $\sD_\Sing$. 

\sss Note that $\Conn_\kappa(\sD, \sP)$ is now an affine space over $\Jet (\fk^*_\sP \otimes \omega)$, and in particular in the case of the trivial bundle $\sP^\circ$ we have an exact sequence of $\Aut(\sD) \ltimes \Jet K$ modules, canonically split at the level of underlying $\Aut(\sD)$-modules 
$$0 \rightarrow \Ring \cdot 1 \rightarrow \Ring[\Conn_\kappa(\sD, \sP^\circ)]^{\leqslant 1} \rightarrow \Loop \fk / \Jet \fk \rightarrow 0.$$
In particular, by applying $-\kappa$ to \eqref{e:actconn}, it follows the natural action of $\Jet K$ on this extension is determined by the formula 
\begin{equation} \label{E:actkappaconn} k \cdot X = \on{Ad}_k(X) + \langle \kappa(X) , k^{-1} \cdot dk \rangle \cdot 1, \quad \quad k \in \Jet K(\Sing), \quad X \in \Loop \fk / \Jet \fk (\Sing) .\end{equation}

\sss Let us now suppose $\kappa$ is an $\Ring$-linear combination of extensions as in Remark \ref{r:traceforms}. In this case, the first step of the PBW filtration on affine vertex algebra also yields an exact sequence of $\Aut(\sD) \ltimes \Jet K$-modules, canonically split at the level of underlying $\Aut(\sD)$-modules
$$0 \rightarrow \Ring \cdot |\varnothing\rangle \rightarrow V_\kappa(\fk)^{\leqslant 1} \rightarrow \Loop \fk / \Jet \fk \rightarrow 0.$$

\begin{Lem} \label{l:KMvskappaconnections} The composition 
$$\Ring[\Conn_\kappa(\sD, \sP^\circ)]^{\leqslant 1} \simeq \Ring \cdot 1 \oplus \Loop \fk / \Jet \fk \simeq \Ring \cdot | \varnothing \rangle \oplus \Loop \fk / \Jet \fk \simeq V_\kappa(\fk)^{\leqslant 1},$$where the outer equivalences are the canonical $\Aut(\sD)$-equivariant splittings, and the middle equivalence sends $1$ to $| \varnothing \rangle$ and acts as the identity on $\Loop \fk / \Jet \fk$, 
is an $\Aut(\sD) \ltimes \Jet K$-equivariant isomorphism.  
\end{Lem}

\begin{proof} This follows from comparing formula \eqref{E:actkappaconn} and Lemma \ref{Lem:KM_cocycle_group}.     
\end{proof}

\sss Let us deduce from the above a coordinate invariant description of the Frobenius center of the Kac--Moody vertex algebra. To do so, let us for simplicity suppose that $K$ is defined over $\mathbb{F}_p$, so that its space of invariant forms $\Sym^2(\fk^*)^K$, viewed as a group scheme over $\Ring$, comes with a relative Frobenius endomorphism $\on{Fr}$, and in particular an Artin--Schreier map 
$$\on{AS}: \Sym^2(\fk^*)^K \rightarrow \Sym^2(\fk^*)^K, \quad \quad \kappa \mapsto \on{Fr}(\kappa) - \kappa =: \kappa^p - \kappa.$$

In addition, let us make the following basic observation regarding mapping spaces. Namely, given $\Ring$-spaces $X$ and $Y$, recall that we can form the inner Hom space 
$$\uHom(X,Y)(\Sing) := \Hom_{\Sing\on{-spaces}}(X_\Sing, Y_\Sing),$$
which satisfies the universal property 
\begin{equation*} \Hom_{\Ring\on{-spaces}}( -, \uHom(X,Y)) \simeq \Hom_{\Ring\on{-spaces}}( - \times X, Y). \end{equation*}
As the Frobenius twist \begin{equation} \label{e:twists}(-)\fr: \Ring\on{-spaces} \rightarrow \Ring\on{-spaces}\end{equation}is (symmetric) monoidal, for any $X$ and $Y$ we have a canonical map 
\begin{equation} \label{e:twists2} \uHom(X,Y)\fr \rightarrow \uHom(X\fr, Y\fr),\end{equation}
and as $\Ring$ is perfect, \eqref{e:twists} is an equivalence hence so is \eqref{e:twists2}.  Moreover, by construction, the composition 
$$\uHom(X,Y) \xrightarrow{\on{Fr}} \uHom(X,Y)\fr \simeq \uHom(X\fr, Y\fr)$$
sends on $\Sing$-points a map $\phi: X_\Sing \rightarrow Y_\Sing$ to its Frobenius twist 
$$(X\fr)_\Sing \simeq (X_\Sing)\fr \xrightarrow{\phi\fr} (Y_\Sing)\fr \simeq (Y\fr)_\Sing.$$

Let us highlight two relevant special cases. First, if we write $\underline{\on{Isom}}(X,Y) \subset \uHom(X,Y)$ for the subfunctor of isomorphisms, the above specializes to an equivalence $$\underline{\on{Isom}}(X,Y)\fr \simeq \underline{\on{Isom}}(X\fr, Y\fr).$$
In particular, by specializing to the automorphism space of the ind-scheme $\sD$, i.e.,  $\Aut(\sD) = \underline{\on{Isom}}(\sD, \sD)$, we have a canonical map 
$$\Aut(\sD) \xrightarrow{\on{Fr}} \Aut(\sD)\fr \simeq \Aut(\sD\fr).$$
Second, given a map $\pi: X \rightarrow Y$ of $\Ring$-spaces, we may form its space of sections 
$$\underline{\Gamma}(Y,X) := \uHom(Y,X) \underset{\uHom(Y,Y)}\times  \on{id}_Y  \subset \uHom(Y,X).$$
and the above specializes to an equivalence 
$$\underline{\Gamma}(Y,X)\fr \simeq \underline{\Gamma}(Y\fr, X\fr).$$
In particular, by specializing to the case of a space $\pi: X \rightarrow \sD$, we obtain a canonical map 
$$\Jet X \simeq \underline{\Gamma}(\sD, X) \xrightarrow{\on{Fr}} \underline{\Gamma}(\sD, X)\fr \simeq \underline{\Gamma}(\sD\fr, X\fr) =: \Jet\fr X\fr.$$
By combining the above two points, for $K \simeq K\fr$ as above, we have a canonical map of group schemes 
$$\on{Fr}: \Aut(\sD) \ltimes \Jet K \rightarrow \Aut(\sD\fr) \ltimes \Jet\fr K.$$

\begin{Lem}  \label{l:coordinvpcenter}Let $\kappa$ be as in Lemma \ref{l:KMvskappaconnections}. Then we have a canonical $\Aut(\sD) \ltimes \Jet K$-equivariant equivalence 
$$\mathfrak{z}_{Fr}(V_\kappa(\fk)) \simeq \Ring[\Conn_{\kappa^p - \kappa}(\sD\fr, \sP^\circ)]. $$
\end{Lem}

\begin{proof} Given a map $\pi: X \rightarrow \sD^\times$, consider its Frobenius twist $\pi\fr: X\fr \rightarrow (\sD^\times)\fr$, and denote the resulting space of sections by $\Loop\fr X\fr$. For $X$ the constant scheme $\fk \times \mathscr{D}^\times$, we may use the $\mathbb{F}_p$-form of $K$ to further identify $$\Loop\fr \fk\fr \simeq \Loop\fr \fk.$$

Within $\mathfrak{z}_{Fr}(\widetilde{U}_\kappa(\fk))$, the first step of its PBW filtration yields by Proposition \ref{p:p-cent} an $\Aut(\sD) \ltimes \Jet K$-equivariant short exact sequence 
$$0 \rightarrow \Ring \rightarrow \mathscr{E} \rightarrow \Loop\fr \fk \rightarrow 0.$$Moreover, recalling the Kac--Moody extension comes from a central extension of the loop group by a torus, it follows that $\mathscr{E}$ is $\Aut(\sD) \ltimes \Jet K$-equivariantly identified with the extension 
$$0 \rightarrow \Ring \rightarrow \widehat{\Loop\fr \fk}_{\kappa^p - \kappa} \rightarrow \Loop\fr \fk \rightarrow 0,$$i.e., the Kac--Moody extension of $\Loop\fr \fk$ at level $\kappa^p - \kappa$. 
Passing to the image of $\mathscr{E}$ in $\mathfrak{z}_{Fr}(V_\kappa)$, we are therefore done by Lemma \ref{l:KMvskappaconnections}. 
\end{proof}

\sss To proceed, we need two final pieces of general preparation regarding trivializable, but not trivialized, bundles $\sP$.

The first piece of general preparation concerns Frobenius twists of bundles and the corresponding jet schemes. 
\begin{Lem}\label{l:twistandjet}If $\sP \rightarrow \sD$ is a trivializable bundle, consider its Frobenius twist $\sP\fr \rightarrow \sD\fr$. Then we have a canonical $\Jet K_\sP$-equivariant isomorphism 
$$ \Jet\fr \sP\fr \simeq \Jet\fr K \overset{\Jet K} \times \Jet \sP,$$where the contracted product is along $\on{Fr}: \Jet K \rightarrow \Jet\fr K.$
Moreover, if $\sP$ is equipped with a datum of $\Aut(\sD)$-equivariance, this isomorphism is $\Aut(\sD) \ltimes \Jet K_\sP$-equivariant. 
\end{Lem}

\begin{proof} Note that an $\Sing$-point of $\Jet\fr \sP\fr$ is an isomorphism $$\iota': \sP\fr_\Sing \simeq \sP^\circ_\Sing$$over $\sD_\Sing\fr$, an $\Sing$-point of $\Jet\fr K$ is an automorphism $\phi$ of $\sP^\circ_\Sing$ over $\sD\fr_\Sing$, an $\Sing$-point of $\Jet K$ is an automorphism $\alpha$ of $\sP^\circ_\Sing$ over $\sD_\Sing$, and an $\Sing$-point of $\Jet \sP$ is an isomorphism $$\iota: \sP_\Sing \simeq \sP^\circ_\Sing$$over $\sD_\Sing$. Then the desired isomorphism is induced by the map 
$$\Jet \fr \sP\fr \leftarrow \Jet\fr K \times \Jet \sP$$sending a pair $(\phi, \iota)$ to the composition of isomorphisms over $\sD\fr_\Sing$
$$\sP\fr_\Sing \xrightarrow{\iota\fr}  \sP^\circ_\Sing \xrightarrow{\phi}\sP^\circ_\Sing.$$
\end{proof}

The second piece of general preparation concerns $\kappa$-connections on trivializable bundles. 
\begin{Lem} \label{l:connectandjet}If $\sP \rightarrow \sD$ is a trivializable bundle, then we have a canonical $\Jet K_\sP$-equivariant isomorphism
$$\Conn_\kappa(\sD, \sP) \simeq \Conn_\kappa(\sD, \sP^\circ) \overset{ \Jet K} \times \Jet \sP.$$
Moreover, if $\sP$ is equipped with a datum of $\Aut(\sD)$-equivariance, this isomorphism is $\Aut(\sD) \ltimes \Jet K_\sP$-equivariant. 
\end{Lem}

\begin{proof} Note that an $\Sing$-point of $\Conn_\kappa(\sD, \sP^\circ)$ is a $\kappa$-connection $\nabla$ on $\sP^\circ_\Sing \rightarrow \sD_\Sing$, an $\Sing$-point of $\Jet K$ is an automorphism $\alpha$ of $\sP^\circ_\Sing \rightarrow \sD_\Sing$, and an $\Sing$-point of $\Jet \sP$ is an isomorphism $\iota: \sP_\Sing \simeq \sP^\circ_\Sing$. Then the desired equivalence $$\Conn_\kappa(\sD, \sP) \simeq \Conn_\kappa(\sD, \sP^\circ) \overset{ \Jet K} \times \Jet \sP$$ 
sends a pair $(\nabla, \iota)$ to the $\kappa$-connection on $\sP$ interchanged with $\nabla$ by $\iota$, and by construction satisfies the mentioned properties.
\end{proof}

\subsection{The Harish--Chandra center of the free field algebra} 

\label{ss_hcff}
 As promised in Section \ref{ss:kconnections}, we are now ready to describe the subalgebra  $\diglet^{\Jet B} \subset \diglet$ coordinate invariantly. 

\begin{Prop} \label{p:digjet} There is a canonical $\Aut(\D)$-equivariant isomorphism 
$$\diglet^{\Jet B} \simeq \Ring[\Conn_{\ckappa^p - \ckappa}(\sD\fr, \omega^\rho)].$$
\end{Prop}

\begin{proof}We begin our computation by noting that an (easier) variant of Lemma \ref{Lem:invariant_containment} yields an equality
\begin{align*}
    \diglet^{\Jet B} &= V_{\ckappa'}(\check{\mathfrak{h}})^{\Jet H}. \intertext{Using Lemma  \ref{l:twoHeiscommute} and Proposition \ref{p:coordinvarianceofcheat} we may rewrite this as} &\simeq V_{\ckappa}(\check{\mathfrak{h}})^{\Jet H} \overset{\Jet\check{H}} \times \Jet \omega^\rho  \intertext{and using Proposition \ref{p:pvsHCswap} we may further rewrite this as} & \simeq \frz_{Fr}(V_{\ckappa}(\check{\mathfrak{h}})) \overset{\Jet\check{H}} \times\Jet \omega^\rho. \intertext{Using Lemma \ref{l:coordinvpcenter} we may rewrite the appearing Frobenius center as} &\simeq \Ring[\on{Conn}_{\ckappa^p - \ckappa}(\sD\fr, \sP^\circ)] \overset{\Jet \check{H}} \times \Jet \omega^\rho. \label{align1} \intertext{Recalling from Lemma \ref{l:coordinvpcenter} that the action of $\Jet \check{H}$ on the appearing space of $\ckappa$-connections factors through the relative Frobenius map $\on{Fr}: \Jet \check{H} \rightarrow \Jet\fr \check{H}$, we may rewrite this as} &\simeq \Ring[\on{Conn}_{\ckappa^p - \ckappa}(\sD\fr, \sP^\circ)] \overset{ \Jet\fr \check{H}} \times \Jet\fr \check{H} \overset{\Jet \check{H}} \times \Jet \omega^\rho. \intertext{Using Lemma \ref{l:twistandjet} and noting that applying $(-)\fr$ to $\omega^\rho \rightarrow \sD$ yields $\omega^\rho \rightarrow \sD\fr$, we obtain} &\simeq \Ring[\on{Conn}_{\ckappa^p - \ckappa}(\sD\fr, \sP^\circ)] \overset{ \Jet\fr \check{H}} \times \Jet\fr \omega^\rho. \intertext{Finally, using Lemma \ref{l:connectandjet}, we obtain } &\simeq \Ring[\on{Conn}_{\ckappa^p - \ckappa}(\sD\fr, \omega^\rho)],\end{align*}as desired. 
\end{proof}

\subsection{The scheme of $\ckappa$-opers}
 \label{ss:kappaopers}
\sss In the previous section, we identified the target of the free field realization 
$$V_\kappa(\fg)^{\Jet G} \hookrightarrow \diglet^{\Jet B}$$
as the algebra of functions on a space of $\check{H}$-connections. We would like to use this to deduce a coordinate invariant description of the Harish--Chandra center itself. In this section, we introduce the relevant moduli space of $\check{G}$-connections, and in its sequel provide the desired identification with the Harish--Chandra center. 

\sss Consider the dual Lie algebra $\check{\fg}^*$ with its coadjoint action of $\check{G}_{\on{ad}}$, where $\check{G}_{\on{ad}} \simeq [\check{G}, \check{G}]$ denotes the adjoint form of $\check{G}$. Consider the Borel $\check{B} \subset \check{G}$ afforded by our pinning, and write $\check{B}_{\on{ad}} \subset \check{G}_{\on{ad}}$ for its image in $\check{G}_{\on{ad}}$. In particular, if we write $\check{\mathfrak{n}} \subset \check{\fg}$ for its Lie algebra of its nilradical, we may pass to its orthogonal subspace and obtain a short exact sequence of $\check{B}$-modules 
$$0 \rightarrow \check{\mathfrak{n}}^{{}^\perp} \rightarrow \check{\fg}^* \rightarrow \check{\fg}^*/\check{\mathfrak{n}}^{{}^\perp} \rightarrow 0.$$
The socle of $\check{\fg}^*/\check{\mathfrak{n}}^{{}^\perp}$ decomposes as a direct sum of lines, with weights with respect to the abstract Cartan $\check{H}_{\on{ad}}$ given by the negative simple roots of $\check{G}$, i.e., we obtain a correspondence of $\check{B}_{\on{ad}}$-representations 
$$\underset{i} \oplus \hspace{.5mm} \check{\fg}^*_{-\check{\alpha}_i} \hookrightarrow \check{\fg}^*/\check{\mathfrak{n}}^{{}^\perp} \twoheadleftarrow \check{\fg}^*.$$
In particular, given a $\check{B}_{\on{ad}}$-bundle $\sP \rightarrow \sD$, we obtain a correspondence of vector bundles over $\sD$
$$\underset{i} \oplus \hspace{.5mm} (\check{\fg}^*_{-\check{\alpha}_i})_{\sP} \hookrightarrow (\check{\fg}^*/\check{\mathfrak{n}}^{{}^\perp})_{\sP} \twoheadleftarrow \check{\fg}^*_{\sP}.$$

\sss{} Fix a Kac--Moody level ${\ckappa} \in \Sym^2( \check{\fg}^*)^{\check{G}}$ and a $\check{G}_\on{ad}$-bundle $\sP' \rightarrow \sD$. We may again as in \eqref{e:atiyahkappa} form the Atiyah exact sequence 
\begin{equation} \label{e:atiyahad}0 \rightarrow \check{\fg}^*_{\sP'} \rightarrow \mathscr{E}_{\sP', \ckappa} \rightarrow \on{T}_\sD \rightarrow 0,\end{equation}
and define a ${\ckappa}$-connection on $\sP'$ to be a splitting $\nabla_{\ckappa}: \on{T}_\sD \rightarrow \mathscr{E}_{\sP', \ckappa}$ of \eqref{e:atiyahad}.

Now suppose in addition $\sP'$ is equipped with a $\check{B}_{\on{ad}}$-reduction $\sP \subset \sP'$. If we form the pushout of \eqref{e:atiyahad} along $\check{\fg}_{\sP}^* \twoheadrightarrow (\check{\fg}^*/\check{\fb}^{{}^\perp})_{\sP}$, i.e.,  
$$0 \rightarrow (\check{\fg}^*/\check{\fb}^{{}^\perp})_{\sP} \rightarrow \mathscr{F}_{\sP, \ckappa} \rightarrow \on{T}_{\sD} \rightarrow 0$$
this is canonically split, so that from a $\ckappa$-connection on $\sP'$, the induced splitting 
$\on{T}_\sD \xrightarrow{\nabla_{\ckappa}} \mathscr{E}_{\sP', \ckappa} \twoheadrightarrow \mathscr{F}_{\sP, \ckappa}$ 
is equivalent to a map which we denote by
$$\overline{\nabla}_{\ckappa}: \on{T}_{\sD} \rightarrow (\check{\fg}^*/\check{\fb}^{{}^\perp})_{\sP}.$$

\begin{defi} \label{d:kappaoper} A {\em ${\ckappa}$-oper} is a pair $(\sP, \nabla_{\ckappa})$ where:

\begin{enumerate}
\item $\sP$ is a $\check{B}_{\on{ad}}$-bundle on $\sD$, 
\item $\nabla_{\ckappa}$ is a ${\ckappa}$-connection on $\sP \overset{\check{B}_{\on{ad}}} \times \check{G}_{\on{ad}}$,
\end{enumerate}
such that (i) the induced map $\overline{\nabla}_{\ckappa}$ factors as $$\overline{\nabla}_{\ckappa}: \on{T}_\sD \rightarrow \underset{i} \oplus \hspace{.5mm} (\check{\fg}^*_{-\check{\alpha}_i})_{\sP} \hookrightarrow (\check{\fg}^*/\check{\mathfrak{n}}^{{}^\perp})_{\sP},$$and (ii) with respect to this factorization, for each $i \in I$ the composite map 
$$\on{T}_\sD \xrightarrow{\overline{\nabla}_{\ckappa}}  \underset{i} \oplus \hspace{.5mm} (\check{\fg}^*_{-\check{\alpha}_i})_{\sP} \twoheadrightarrow (\check{\fg}^*_{-\check{\alpha}_i})_{\sP}$$
is an isomorphism. 
\end{defi}

The collection of all ${\ckappa}$-opers on $\sD$ forms a groupoid in a natural way, and a standard argument shows that any individual ${\ckappa}$-oper has trivial automorphism group, cf. the proof of Lemma \ref{l:kosslice} below. In particular, ${\ckappa}$-opers form a discrete groupoid. 

\begin{Lem} \label{l:kosslice} The $\Ring$-space $\Op_{\ckappa}(\sD)$ defined by sending a test $\Ring$-algebra $\Sing$ to the set of isomorphism classes of $\ckappa$-opers over $\sD_\Sing$ is representable by a scheme. \end{Lem}

\begin{proof} Let us fix a coordinate $t$ on $\sD$ and a $\check{G}$-invariant isomorphism $\epsilon: \check{\fg}^* \simeq \check{\fg}$. Given an $\Sing$-point $$(\sP, \nabla_{\ckappa}) \in [\Op_{\ckappa}(\sD)](\Sing),$$we may \'{e}tale locally on $\Sing$ trivialize $\sP$ on $\sD_\Sing$. Given $\sP$ which is trivializable, if we choose a principal $\mathfrak{sl}_2$-triple $f,h,e \in \check{\fg}$ such that $h,e \in \check{\fb}$, the space of oper structures on $\sP$ is identified via $\epsilon$ with 
$$(f + \check{\mathfrak{b}}_\Sing[[t]]) / \check{N}_\Sing[[t]],$$
where $\check{N} \subset \check{B}$ denotes the unipotent radical and the action of $\check{N}_\Sing[[t]]$ is induced by \eqref{E:actkappaconn}. A standard argument shows that each $\check{N}_\Sing[[t]]$-orbit intersects the Kostant slice $f + \mathfrak{z}_{\check{\fg}}(e)_\Sing[[t]]$ in a unique point, and each has trivial stabilizer in $\check{N}_\Sing[[t]]$. Using the uniqueness of this canonical form, and its compatibility with base change in $\Sing$, it follows in particular that any $\sP$ is globally trivializable on $\Sing$, and moreover that 
$$\on{Op}_{\ckappa}(\sD) \simeq f + \mathfrak{z}_{\check{\fg}}(e)[[t]],$$yielding the claimed representability. 
\end{proof}

{In particular, we see that a choice of the coordinate $t$ on $\D$ identifies $\on{Op}_{\check{\kappa}}(\D)$ with the fiber of $\operatorname{Op}^\hbar_0$ from Section \ref{SSS_opers_coord} at $\hbar=\check{\kappa}$.}

\sss{} The last property of $\ckappa$-opers we need is the following. Recall that $\check{G}$ comes equipped with a pinning. In particular, we have a splitting $\check{H}_{\on{ad}} \hookrightarrow \check{B}_{\on{ad}}$, and trivializations $$(\check{\fg}_{\check{\alpha}_i})^* \simeq \Ring, \quad \quad i \in I.$$

\begin{Prop} There exists a canonical {\em Miura transform} map 
$$\mu: \Conn_{\ckappa}(\sD, \omega^\rho) \rightarrow \Op_{\ckappa}(\sD)$$
which is moreover dominant. \label{p:miuradisc}
\end{Prop}

\begin{proof} Note that, by our assumption on $\check{G}$ { that $[\check{G},\check{G}]$ is of adjoint type}, we have $\check{H} \simeq Z_{\check{G}} \times \check{H}_{\on{ad}},$ and with respect to this decomposition $\omega^\rho$ is induced from a $\check{H}_{\on{ad}} \subset \check{H}$-bundle, which we denote by the same letter. 

We will define our desired Miura transform on $\Sing$-points. Thus, given an $\Sing$-point of $\Conn_{\ckappa}(\sD, \omega^\rho)$, i.e., a $\ckappa$-connection on $\omega^\rho$:$$\nabla_{\ckappa}: \on{T}_\sD \rightarrow \mathscr{E}_{\omega^\rho, \ckappa},$$on $\sD_\Sing$, we need to produce a $\ckappa$-oper $(\sP, \nabla'_{\ckappa})$ on $\sD_\Sing$. Using the splitting $\check{H}_{\on{ad}} \hookrightarrow \check{B}_{\on{ad}}$, we take $\sP$ to be the induced bundle 
$$\sP := \omega^\rho \overset{\check{H}_{\on{ad}}} \times \check{B}_{\on{ad}}.$$
Note that for any $\check{G}_{\on{ad}}$-bundle $\sP_{\check{G}}$ induced from an $\check{H}_{\on{ad}}$-bundle $\sP_{\check{H}}$, if we let $\halpha$ run over the set of roots $\check{\Phi}$ of $\check{G}$, i.e., the nonzero weights of $\check{\fg}^*$ as a $\check{H}_{\on{ad}}$-module, we have canonically 
$$\sE_{\sP_{\check{G}}, \ckappa} \simeq   [\underset{\check{\alpha}} \oplus \hspace{.5mm} (\check{\fg}^*_{\check{\alpha}})_{\sP_{\check{H}}}] \oplus \sE_{\sP_{\check{H}}, \ckappa}, $$
compatibly with the natural projections to $\on{T}_{\sD}$. In particular, the data of a $\ckappa$-connection $\nabla'_{\ckappa}$ on $\sP_{\check{G}}$ is the data of a tuple $(\nabla_{\ckappa}, (\xi_{\halpha})_{\halpha \in \check{\Phi}})$, where $\nabla_{\ckappa}$ is a $\ckappa$-connection on $\sP_{\check{H}}$, and each $\xi_{\halpha}$, $\halpha \in \check{\Phi}$, is a map of bundles
$$\xi_{\halpha}: \on{T}_\sD \rightarrow (\check{\fg}^*_{\check{\alpha}})_{\sP_{\check{H}}}.$$

Specializing to the case of $\sP_{\check{H}} = \omega^\rho$, we obtain 
$$\sE_{\sP, \ckappa} \simeq [\underset{\check{\alpha}} \oplus \hspace{.5mm} (\check{\fg}^*_{\check{\alpha}}) \otimes_{\Sing} \omega^{\langle \rho, \halpha \rangle}] \oplus \sE_{\omega^\rho, \ckappa}.$$
We define the desired oper structure $\nabla'_{\ckappa}$ by the tuple $(\nabla_{\ckappa}, (\xi_{\halpha})_{\halpha \in \check{\Phi}})$, where $\nabla_{\ckappa}$ is the original point of $\Conn_{\ckappa}(\sD, \omega^\rho)$,  $\xi_{\halpha} = 0$ unless $\halpha$ is a negative simple root, and for each negative simple root $-\halpha_i$ we set using our pinning
$$\xi_{-\halpha_i}: \on{T}_{\sD} \rightarrow  \fg^*_{-\halpha_i} \otimes_\Sing \omega^{\langle \rho, - \halpha_i \rangle} \simeq \Sing \otimes_\Sing \on{T}_{\sD} \simeq \on{T}_{\sD}$$
to be the identity map.

It is clear the above construction yields an oper structure on $\sP$ and is moreover compatible with base change in $\Sing$, which defines the desired map $\mu$. Moreover, its dominance may be checked after fixing a coordinate $t$ on $\sD$, whence it follows from Step 1 of the proof of Proposition \ref{Prop:Miura_isomorphism}. 
\end{proof}

\subsection{Identification of the Harish--Chandra center}

\label{ss:hccenD}
Finally, we are ready to deduce the main result of this section. { In particular, we no longer assume that the derived subgroup of $\check{G}$ is of adjoint type.}

\begin{Prop} \label{p:autDdisc}There exists a canonical $\Aut(\sD)$-equivariant isomorphism 
$$V_\kappa(\fg)^{\Jet G} \simeq \Ring[\on{Op}_{\ckappa^p - \ckappa}(\sD\fr)].$$In particular, the natural action of $\Aut(\sD)$ on $V_\kappa(\fg)^{\Jet G}$ factors through the relative Frobenius $$\on{Fr}: \Aut(\sD) \rightarrow \Aut(\sD\fr).$$
\end{Prop}

\begin{proof} { We first claim the case of general $\check{G}$ reduces to the case of $\check{G}$ with $[\check{G}, \check{G}]$ of adjoint type. Indeed, on the spectral side, we note that the moduli space of $\ckappa^p - \ckappa$-opers for $\check{G}$ by definition only depends on $\check{G}_{\on{ad}}$. On the representation theoretic side, we note that $V_\kappa(\fg)^{\Jet G}$ is unchanged by a finite  isogeny $G' \rightarrow G$ of degree coprime to the characteristic. By our assumption that $p > h$, we may always find such a $G'$ such that its Langlands dual has a derived subgroup of adjoint type, which completes the reduction. Therefore, in what follows we assume that $\check{G}$ has a derived subgroup of adjoint type.}

It is enough to verify that the isomorphism of Proposition \ref{p:digjet}$$\diglet^{\Jet B} \simeq \Ring[\Conn_{\ckappa^p - \ckappa}(\sD\fr, \omega^\rho)]$$intertwines the inclusions of subalgebras 
$$V_\kappa(\fg)^{\Jet G} \hookrightarrow \diglet^{\Jet B} \quad \text{and} \quad \mu^*: \Ring[\on{Op}_{\ckappa^p - \ckappa}(\sD\fr)] \hookrightarrow \Ring[\Conn_{\ckappa^p - \ckappa}(\sD\fr, \omega^\rho)].$$
Moreover, it is straightforward to reduce this claim to the analogous claim for $G$ either a torus or an almost simple, simply connected group. In the case of $G$ a torus, this follows from Lemma \ref{l:twoHeiscommute}, Proposition \ref{p:pvsHCswap}, and Lemma \ref{l:coordinvpcenter}, exactly as in the proof of Proposition \ref{p:digjet}. 

In what remains, we will therefore assume $G$ to be an almost simple, simply connected group. We further break into two cases, depending on whether $\kappa$ is integral or non-integral. 

If $\kappa$ is an integral level, so is $\ckappa$, and in particular $\kappa^p - \kappa = 0$ and $\ckappa^p - \ckappa = 0$. Note that with respect to the isomorphism 
$$V_\kappa(\fh) \simeq V_{\ckappa}(\check{\fh})$$
Proposition \ref{p:pvsHCswap} simplifies to the chain of equalities
$$V_\kappa(\fh)^{\Jet H} = \frz_{Fr}(V_\kappa(\fh)) = \frz_{Fr}(V_{\ckappa}(\check{\fh})) = V_{\ckappa}(\check{\fh})^{\Jet \check{H}}.$$Explicitly, the middle identification is given via Lemma \ref{l:coordinvpcenter} by 
$$\frz_{Fr}(V_{\kappa}(\fh)) \simeq \Ring[\Conn_{\kappa^p - \kappa}(\sD\fr, \sP_H^\circ)] \simeq \Ring[\Jet\fr(\check{\fh} \otimes \omega)] \overset{\phi}\simeq \Ring[\Jet\fr(\fh \otimes \omega)] \simeq \Ring[\Conn_{\ckappa^p - \ckappa}(\sD\fr, \sP_{\check{H}}^\circ)] \simeq \frz_{Fr}(V_{\ckappa}(\check{\fh})),$$where $\phi$ is induced by \begin{equation} \label{e:mk} -\kappa: \fh \simeq \check{\fh}.\end{equation}Let us think of $\fh \otimes \omega$ as the associated bundle $\fh_\omega$, where by a mild abuse of notation we again use $\omega \rightarrow \sD\fr$ to denote the $\mathbb{G}_m$-bundle over $\sD\fr$ corresponding to the canonical line bundle, and are considering $\mathbb{G}_m$ acting on $\fh$ via dilations. As \eqref{e:mk} is $W$-equivariant, we deduce a commutative square over $\sD\fr$
$$\xymatrix{ \fh_\omega \ar[d] \ar@{-}[rr]^\sim_{\eqref{e:mk}} & &\check{\fh}_\omega \ar[d] \\ (\fh /\!/ W)_\omega \ar@{-}[rr]^\sim_{\eqref{e:mk}} && (\check{\fh}/\!/W)_\omega},$$
where the vertical maps are the tautological projections. In particular, we deduce, after passing  to functions on the corresponding schemes of sections a commutative diagram 
$$\xymatrix{ \Ring[\Jet\fr(\fh \otimes \omega)]  \ar@{-}[rr]^\sim_{\eqref{e:mk}} & & \Ring[\Jet\fr (\check{\fh} \otimes \omega)] \\ \Ring[\Jet\fr(\fh /\!/ W)_\omega] \ar[u] \ar@{-}[rr]^\sim_{\eqref{e:mk}} && \Ring[\Jet\fr(\check{\fh}/\!/W)_\omega], \ar[u]}$$
so the claim now follows from (the proof of) Theorem \ref{Thm:integral_case}. Namely, the latter theorem showed that, using a coordinate $t$ to produce an embedding $V_\kappa(\fg)^{\Jet G} \subset V_\kappa(\fh)^{\Jet H}$, the image was precisely $$\Ring[\Jet\fr(\check{\fh} /\!/ W)_\omega] \subset \Ring[\Jet\fr(\check{\fh} \otimes \omega)]$$  
but the latter is also easily seen to be (using the identifications afforded by the coordinate $t$) the image of the level zero Miura map
$$\Ring[\Op_0(\sD\fr)] \subset \Ring[\Conn_0(\sD\fr, \omega^\rho)].$$

It remains to address the case of $G$ an almost simple, simply connected group, and $\kappa$ a non-integral level. It is enough to check the desired equality of subalgebras after a choice of coordinate $t$, but this is precisely the assertion of Theorem \ref{t:tungus}. 
\end{proof}

\subsection{Passing to the punctured disc}

\sss{} Having proven the coordinate invariant description of $V_\kappa(\fg)^{\Jet G}$, i.e., the Harish--Chandra center of the vertex algebra, let us turn to the Harish--Chandra center of the enveloping algebra. In particular, the relevant group of coordinate changes now is the ind-affine group ind-scheme $\Aut(\pD)$, cf. Section \ref{SSS_Virasoro}.  

\sss On the representation-theoretic side, we note the following. 

\begin{Lem} The filtered complete Harish--Chandra center $\widehat{U}_{\kappa}(\hat{\fg})^{\Loop G}$ carries a natural action of $\Aut(\pD)$. \label{lem:autpdactonskm}
\end{Lem}

\begin{proof} We first note that the Kac--Moody algebra $\hat{\fg}_\kappa$ carries a canonical action of $\Aut(\pD) \ltimes \Loop G$. Indeed, this follows straightforwardly from the construction reviewed in Section \ref{SSS_KM_definition}. Namely, following the terminology of {\em loc. cit.}, with $F$ replaced by $G$, we first note the map 
$$\Loop G \rightarrow GL(\Loop W)$$is $\Aut(\pD) \ltimes \Loop G$ equivariant. 
Moreover, as the action of $\Aut(\pD) \ltimes \Loop G$ on $\Loop W$ is linear, i.e., factors through a { group ind-scheme homomorphism} $$\Aut(\pD) \ltimes \Loop G \rightarrow GL(\Loop W)$$it follows that $\Aut(\pD) \ltimes \Loop G$ acts on the Tate extension of $GL(\Loop W)$ and hence by restriction on the Kac--Moody extension of $\Loop G$. By Corollary \ref{c:ernie}, it follows that the $\Sing$-points of $\Aut(\pD) \ltimes \Loop G$ act on the $\Sing$-points of $\hat{\fg}_\kappa$ for every $\Ring$-algebra $\Sing$, cf. the argument of Section \ref{sss:adjointactions}, which yields the claim. 

By the functoriality of completed enveloping algebras, cf. Section \ref{sss:completedenvelope}, it follows that $\Aut(\pD) \ltimes \Loop G$ acts on $\widetilde{U}_\kappa(\widehat{\fg})$. Moreover, the construction of {\em loc. cit.} implies this action restricts to the filtered complete enveloping algebra $\widehat{U}_{\kappa}(\hat{\fg})$, and in particular the $\Loop G$ invariants in the latter inherit an action of $\Aut(\pD)$. 
\end{proof}

\sss{} \label{sss:operspd} On the spectral side, consider the moduli space of $\ckappa$-opers on the punctured disc, which we denote by $\Op_{\ckappa}(\pD)$. Namely, its groupoid of $\Sing$-points is defined exactly as in Definition \ref{d:kappaoper}, after replacing $\sD$ with $\pD$, and insisting the $\check{B}_{\on{ad}}$-bundle be trivializable. The resulting moduli space is now an ind-scheme; namely, as in Lemma \ref{l:kosslice}, any choice of coordinate $t$ on $\pD$ yields an isomorphism 
$$\Op_{\ckappa}(\pD) \simeq f + \mathfrak{z}_{\check{\fg}}(e)(\!(t)\!).$$
Via the choice of coordinate $t$, we in particular obtain an identification of $\Op_{\ckappa}(\pD)$ with an affine space over $\check{\fg}^e(\!(t)\!),$ and in particular a filtration on its algebra of functions. It is straightforward to see that this filtration is in fact independent of the coordinate $t$, cf. Equations (4.2-6)-(4.2-8) of \cite{Frenkel_loop}. As in Section \ref{SSS_univ_env_map}, let us denote the corresponding filtered complete subalgebra by 
$$\Ring[\Op_{\ckappa}(\pD)]_f \hookrightarrow \Ring[\Op_{\ckappa}(\pD)].$$

It is clear that $\Aut(\pD)$ acts naturally on $\Op_{\ckappa}(\pD)$, and hence on $\Ring[\Op_{\ckappa}(\pD)]$, as e.g. one has a manifest action on $\Sing$-points. Moreover, by the aforementioned coordinate invariance of the filtration, this action preserves the filtered complete subalgebra.

\sss{} Having reviewed the natural actions of $\Aut(\pD)$ on the representation-theoretic and spectral sides, we may state the main result of this subsection.

\begin{Prop} \label{p:coordinvariancepD} There is a canonical $\Aut(\pD)$-equivariant isomorphism of filtered complete algebras
\begin{equation} \label{e:coordinvhc} \widehat{U}_\kappa(\hat{\fg})^{\Loop G} \simeq \Ring[\Op_{\ckappa^p - \ckappa}({\pD}\fr)]_f.\end{equation}
\end{Prop}

\sss{} Before giving the proof of Proposition \ref{p:coordinvariancepD}, it will be helpful to establish a group theoretic property of $\Aut(\pD)$. This will be insensitive to the choice of coefficient ring $\Ring$, which we now take to be a general commutative ring. 

As motivation, recall that the closed embedding $\Aut(\D) \hookrightarrow \Aut(\pD)$ yields on Lie algebras the embedding $\Ring[[t]]\partial_t \hookrightarrow \Ring (\!(t)\!)\partial_t$. This admits a complementary Lie subalgebra, namely $t^{-1}\Ring[t^{-1}]\partial_t$. We would like to now establish the counterpart of this fact at the level of group ind-schemes. Recall from Section \ref{SSS_Virasoro} that the $\Sing$-points of $\Aut(\pD)$ consist of Laurent series \begin{equation} \label{e:phi} \phi(t) = \underset{i} \Sigma \hspace{.5mm} a_i t^i \in \Sing(\!(t)\!)\end{equation}with (i) $a_i, i \leqslant 0$, nilpotent, and (ii) $a_1$ a unit. 

\begin{defi} We define $$\Am \rightarrow \Aut(\pD)$$
to be the subfunctor with $\Sing$-points Laurent series $\phi(t)$ as in Equation \eqref{e:phi} satisfying the two conditions that (i) the $a_i, i < 0$, are nilpotent, (ii) we have  $a_i - \delta_{i, 1} = 0$, for all $i \geqslant 0$. 
\end{defi}
That is, the $\Sing$-points of $\Am$ consist of series of the form 
\begin{equation} \label{e:autminus} \phi(t) = t + a_{-1} t^{-1} + \cdots + a_{-n} t^{-n}\end{equation}
with $a_i, i < 0$, nilpotent. We will provide a slightly alternative presentation of $\Am$ in Remark \ref{r:rightmovers} below, after we establish the following basic assertions. 

\begin{Lem}\label{l:autpd-} The following hold. 

\begin{enumerate} \item $\Am$ carries a unique structure of group ind-scheme compatible with its inclusion into $\Aut(\pD)$. 

\item The multiplication map  
$$\Aut(\sD) \times \Am \hookrightarrow \Aut(\pD) \times \Aut(\pD)  \xrightarrow{\on{mult}} \Aut(\pD)$$
is an isomorphism of ind-schemes, i.e., we have a factorization of underlying ind-schemes $\Aut(\sD) \times \Am \simeq \Aut(\pD)$.  

\item The Lie algebra of $\Am$ is given by $t^{-1} \Ring[t^{-1}]\partial_t$. 

\end{enumerate}
\end{Lem}

\begin{proof} For (1), we need to show that the $\Sing$-points of $\Am$ are a subgroup of the $\Sing$-points of $\Aut(\pD)$, i.e., are closed under multiplication, and inversion.  

For the closure under multiplication, given two such series 
$$ \phi(t) = t + a_{-1} t^{-1} + \cdots + a_{-n}t^{-n} \quad \text{and} \quad \psi(t) = t + b_{-1} t^{-1} + \cdots + b_{-m} t^{-m},$$for some integers $n,m > 0$, let us consider their composition
\begin{align*}
   (\phi \circ \psi)(t) &= \psi(t) + \underset{i > 0} \Sigma \hspace{.5mm} a_{-i} \frac{1}{(t + b_{-1} t^{-1} + \cdots + b_{-m}t^{-m})^i}. \intertext{For fixed $i$, let us divide the numerator and denominator of the appearing fraction by $t^{i}$ to obtain } &=\psi(t) + \underset{i > 0} \Sigma \hspace{.5mm} a_{-i} \frac{t^{-i}}{(1 + b_{-1} t^{-2} + \cdots + b_{-m}t^{-m-1})^i}  \intertext{Recalling that in any ring an element of the form $1 + \eta$, for $\eta$ nilpotent, is invertible by expansion into a geometric series, we may continue as} &= \psi(t) + \underset{i > 0} \Sigma \hspace{.5mm} a_{-i}t^{-i} [\underset{j \geqslant 0} \Sigma (-1)^j(b_{-1} t^{-2} + \cdots + b_{-m}t^{-m-1})^j]^i, 
\end{align*}
and it is now manifest the composition again lies in $\Am$.

We will address the closure of $\Am$ under inversion simultaneously with proving (2). 
To do so, note that by definition, the intersection of the $\Sing$-points of $\Aut(\sD)$ and $\Am$ 
within the $\Sing$-points of $\Aut(\pD)$ consists only of the identity element $\phi(t) = t$. Consider an arbitrary $\Sing$-point of $\Aut(\pD)$
$$\rho(t) = \underset{i} \Sigma \hspace{.5mm} \eta_i t^i.$$
We claim that 
\begin{itemize}
\item[(*)]
there exists an $\Sing$-point $\phi(t)$ of $\Am$ such that the composition $\rho \circ \phi$ is an $\Sing$-point of  $\Aut(\sD)$. 
\end{itemize}
Assuming (*) momentarily, we note that if $\rho(t)$ is an $\Sing$-point of $\Am$, then $\rho \circ \phi$ is simultaneously an $\Sing$-point of $\Am$ and $\Aut(\sD)$, i.e., the identity, which establishes the closure of $\Am$ under inversion. Having shown $\Am$ is a group, the claim combined with the above intersection property similarly easily implies (2), i.e., that multiplication produces a bijection on $\Sing$-points. 

It therefore remains to show (*). Let us write $I_{\rho}$ for the ideal of $\Sing$ generated by the $\eta_i$, for $i < 0$. As this is an ideal generated by finitely many nilpotent elements, we have $I_\rho^{N_\rho} = 0$ for some minimal positive integer $N_\rho \geqslant 0$.  We will prove the claim by induction on $N_\rho$. 

For the base case, note that $N_\rho = 1$ if and only if $\rho(t) \in \Aut(\sD)$, so we may take $\phi$ to be the identity. 

For the inductive step, consider a transformation of the form $\phi(t) = t + a_{-1}t^{-1} + \cdots a_{-n}t^{-n} =: t + q(t^{-1})$ with $a_i \in I_\rho$ for all $i < 0$.
We claim that 
\begin{equation} \label{e:desiredidentity} \rho \circ \phi(t) \equiv \rho(t) + q(t^{-1})\partial_t \rho(t) \quad \quad \on{mod} \quad I_\rho^2(\!(t)\!),\end{equation}
i.e., the two elements agree after the projection $S(\!(t)\!) \twoheadrightarrow (S/I^2_\rho)(\!(t)\!).$ To see this, one may first note that for any $k \geqslant 0$, we have 
\begin{equation} \label{e:helpfulidentity1} \phi(t)^k \equiv t^k + q(t^{-1}) \partial_t t^k \quad \on{mod} \quad I^2_\rho(\!(t)\!),  \end{equation}
as follows directly from binomial expansion. Similarly, for any $k < 0$ and $\iota \in I_\rho$, we have 
\begin{equation} \label{e:helpfulidentity2} \iota \phi(t)^k \equiv \iota t^k  \equiv \iota t^k + \iota q(t^{-1}) \partial_t t^k \quad \on{mod} \quad I^2_\rho(\!(t)\!),\end{equation}
as follows from a geometric series manipulation similar to the proof of the closure of $\Am$ under composition above. The desired identity \eqref{e:desiredidentity} now follows simply by expanding $\rho \circ \phi(t)$ into powers of $\phi(t)$ and applying \eqref{e:helpfulidentity1} and \eqref{e:helpfulidentity2}. 

To finish the inductive step, let us write $n > 0$ for the order of the pole of $\rho(t)$, i.e., we have  $$\rho(t) = \eta_{-n}t^{-n} + \eta_{-n+1}t^{-n+1} + \cdots$$with $\eta_{-n}$ nonzero. Recalling that $\eta_1$ is a unit, it follows from \eqref{e:desiredidentity}, and in particular the intermediate identity in \eqref{e:helpfulidentity2}, that 
$\rho( t - \eta_1^{-1} \eta_{-n}t^{-n})$ has a pole of order at most $n - 1$ mod $I^2_\rho(\!(t)\!)$; iterating this, we may clear out all the strictly negative powers of $t$ appearing in $\rho(t)$ mod $I^2_\rho(\!(t)\!)$. As the resulting element $\tilde{\rho}(t)$ has all its coefficients of negative powers of $t$ in $I^2_\rho$, it follows that we have a strict inequality $N_{\tilde{\rho}} < N_\rho$, which completes the inductive step. 

Having proven (1) and (2), claim (3) follows from inspecting the $\Ring[\epsilon]/(\epsilon^2)$ points of $\Am$.
\end{proof}

\begin{Rem}\label{r:rightmovers} Here is a slightly alternative way to think about $\Am$. For an automorphism $\alpha$ of $\pD = \Spec \Ring(\!(t)\!),$ the pushforward on functions applied to the uniformizer, i.e.,  $$\alpha_*(t) = \underset{i}\Sigma \hspace{.5mm} a_i t^i$$is the Laurent series $\phi(t)$ we typically work with, i.e., the corresponding $\Sing$-point of $\Aut(\pD)$. However, the automorphism is equally well determined by the pushforward of the inverse of the uniformizer, i.e., $\alpha_*(t^{-1})$. With this, $\Am$ consists of those automorphisms $\alpha$ for which 
$$\alpha_*(t^{-1}) = t^{-1} + b_{-3}t^{-3} + \cdots + b_{-n}t^{-n},$$
for $b_i, i < -1$, nilpotent. This may be proven straightforwardly using the geometric series manipulations in the proof of Lemma \ref{l:autpd-}(1) above. 

In particular, this presentation makes manifest the sense in which $\Am$ is a formal polynomial analogue of the derived subgroup of the prounipotent radical of the automorphisms of the pointed formal disc.
\end{Rem}

\begin{Rem}Although we do not need it, let us mention that by similar arguments, one can prove a triangular decomposition for $\Aut(\pD)$, i.e., that multiplication yields an isomorphism of ind-schemes $$\Aut(\pD) \simeq \Aut_{>0}(\pD) \times \Aut_0(\pD) \times \Aut_{<0}(\pD),$$where the $\Sing$-points of the three appearing subgroups are defined as follows. \begin{enumerate}  \item $\Aut_{> 0}(\pD)$ consists of Laurent series $\phi(t) = \underset{i}\Sigma \hspace{.5mm} a_it^i$ for which $a_i - \delta_{i,1} = 0$ for $i \leqslant 1$.

\item $\Aut_0(\pD)$ consists of Laurent series with $a_i = 0$ for $i \neq 1$, and $a_1$ a unit, i.e., $\Aut^0(\pD)$ is the loop rotation $\mathbb{G}_m$.

\item $\Aut_{<0}(\pD)$ consists of Laurent series with $a_i  - \delta_{i,1} = 0$, for $i \geqslant 1$, and $a_i, i < 1$, nilpotent. 

\end{enumerate}
This moreover lifts to a triangular decomposition of the Virasoro group ind-scheme, as in Section \ref{SSS_Virasoro}, which on Lie algebras recovers the standard triangular decomposition of the Virasoro Lie algebra; we leave the details to the interested reader. 
\end{Rem}

\sss{} We will need one more result concerning $\Am$ in the special case when $\Ring$ is a $\mathbb{Q}$-algebra. 

  Write $\mathbb{G}^\wedge_{a}$ for the formal group of $\mathbb{G}_a$, i.e., its formal completion at the identity. Explicitly, the $\Sing$-points of $\mathbb{G}^\wedge_{a}$ are nilpotent elements $\nu$ of $\Sing$, viewed as a group under addition. 

To state it, we note the following. 
\begin{Lem} \label{l:formexp} For any element $X \in t^{-1}\Ring[t^{-1}]\partial_t$, we have an associated exponential map of group ind-schemes 
$\mathbb{G}_a^\wedge \rightarrow \Am$ 
which on $\Sing$-points sends a nilpotent element $\nu$ of $\Sing$ to the element
\begin{equation} \label{e:formexp} e^{\nu X}(t) := \underset{i \geqslant 0} \Sigma \hspace{.5mm} \frac{1}{i!} \nu^i X^i(t).\end{equation}
\end{Lem}
Explicitly, in forming the above exponential, we are viewing $X$ as an $\Sing$-linear endomorphism of $\Sing(\!(t)\!)$.  
\begin{proof} We first check that \eqref{e:formexp} defines a point of $\Am$ as in \eqref{e:autminus}. Indeed, using the fact that $\nu$ is nilpotent, it is clear the sum in \eqref{e:formexp} is finite. Moreover, noting that with respect to the loop rotation grading on $\Ring[t^{\pm 1}]$, $X$ is a sum of operators of degree $\leqslant -2$, it is clear that $e^{\nu X}(t)$ is of the form \eqref{e:autminus}, as claimed.  

To check that \eqref{e:formexp} is a map of group ind-schemes, on $\Sing$-points we must verify the equality, for any pair of nilpotent elements $\nu_1, \nu_2$ of $\Sing$, 
$$e^{(\nu_1 + \nu_2)X}(t) = e^{\nu_1 X}(e^{\nu_2 X}(t)).$$
This is standard, and follows from applying a binomial expansion to $(\nu_1 + \nu_2)^i$ and using that $X$ is a derivation. 
\end{proof}

Let us write $\ell_i := -t^{i+1} \partial_t$. We will apply Lemma \ref{l:formexp} to the sequence of elements $\ell_i$, for $i \leqslant -2$, as follows. For any $j \leqslant -2$,  write $\underline{j}$ for the set of integers $j \leqslant i \leqslant -2$, and consider the map of ind-schemes
$$\underset{\underline{j}}\Pi \hspace{.5mm} \mathbb{G}_a^\wedge \rightarrow \Am \quad \quad (\nu_{j}, \nu_{j+1}, \ldots, \nu_{-2}) \mapsto e^{\nu_j \ell_j} \circ e^{\nu_{j+1} \ell_{j+1}} \circ \cdots \circ e^{\nu_{-2} \ell_{-2}} (t).$$For ease of notation, we will denote this map below on $\Sing$-points simply by $\nu \mapsto e^{\nu \cdot \ell}(t)$. 
These maps form an inductive system as one decreases $j$ in an evident sense, so that we may form the colimit
\begin{equation} \label{e:cells} \varinjlim \underset{\underline{j}}\Pi \hspace{.5mm} \mathbb{G}_a^\wedge \rightarrow \Am. \end{equation}

\begin{Lem} \label{l:cells} The map \eqref{e:cells} is an isomorphism of ind-schemes. \end{Lem}

\begin{proof} By base change, it is enough to address the case of $\Ring = \mathbb{Q}$. If we consider the induced map on $\mathbb{Q}[\epsilon]/(\epsilon^2)$-points, this is a bijection, as follows from the fact that the $\ell_i, i \leqslant -2$, are a basis for $t^{-1}\mathbb{Q}[t^{-1}]\partial_t$.   

We now claim we are done on general grounds. Namely, consider an ind-affine ind-scheme $X = \varinjlim X_\alpha$ over $\mathbb{Q}$, and write $A = \varprojlim A_\alpha$ for its topological ring of functions.  

Let us say $A$ is {\em convenient} if the following two conditions hold. First, suppose
\begin{itemize}
\item[(i)] the corresponding reduced ind-scheme is $\Spec \mathbb{Q}$, i.e., for each $\alpha$ the composite map $\mathbb{Q} \rightarrow A_\alpha \rightarrow A_\alpha/I_\alpha$ is an isomorphism, where $I_\alpha$ denotes the nilradical of $A_\alpha$. 
\end{itemize}
Let us write $I = \varprojlim I_\alpha$, and note that we have an isomorphism of topological algebras 
\begin{equation}\label{e:convenient} A \simeq \varprojlim A/I^n,\end{equation}
where we emphasize that each $A/I^n$ is a topological algebra. Consider $I/I^2 \simeq \varprojlim I_\alpha/I_\alpha^2$, and suppose that 
\begin{itemize}
\item[(ii)] for each $n$, the natural map of pro-vector spaces 
$$\Sym^n(I/I^2) \rightarrow I^n/I^{n+1}$$
is an isomorphism. 
\end{itemize}
Then given convenient topological algebras $(A,I)$ and $(B,J)$, note that by (i) any map of topological algebras $A \rightarrow B$ sends $I$ into $J$. We claim that if the induced map $I/I^2 \rightarrow J/J^2$ is an isomorphism of pro-vector spaces, then the map $A \rightarrow B$ is an isomorphism. Indeed, by \eqref{e:convenient}, it is enough to verify the induced maps of topological algebras $A/I^n \rightarrow B/J^n$ are isomorphisms. For $n = 1$, this follows from (i). For each $n \geqslant 2$, we have the commutative diagram with exact rows in the abelian category of pro-vector spaces
$$\xymatrix{0 \ar[r]&  I^{n-1}/I^{n} \ar[d] \ar[r] & A/I^{n} \ar[d] \ar[r] & A/I^{n-1} \ar[d] \ar[r] & 0 \\ 0 \ar[r] & J^{n-1}/J^{n} \ar[r] & B/J^{n} \ar[r] & B/J^{n-1} \ar[r] & 0.}$$
Using the assumption on $I/I^2 \simeq J/J^2$ and condition (ii) of convenience, the result follows by induction and the five lemma. 

It is clear that both ind-schemes in the statement of the Lemma are convenient, so we are done by the result on $\mathbb{Q}[\epsilon]/(\epsilon^2)$-points. 
\end{proof}

\sss{} Finally, we also need the following basic fact concerning actions of $\mathbb{G}_a^\wedge$ when $\Ring$ is a $\mathbb{Q}$-algebra; in particular, we continue to assume the latter holds. 

Namely, let us view $\Ring[\mathbb{G}_a^\wedge]$ as a coalgebra in $\Pro(\Ring\mod)$. Explicitly, we have $\Ring[\mathbb{G}_a^\wedge] \simeq \Ring[[x]] = \varprojlim \Ring[x]/(x^n)$, and in particular for any object $V$ of $\Pro(\Ring\mod)$ we have 
$$V \widetilde{\otimes}_\Ring \Ring[[x]] \simeq \varprojlim V \widetilde{\otimes}_\Ring \Ring[x]/(x^n) \simeq \varprojlim V[x]/(x^n) =: V[[x]];$$
explicitly, as an object of $\Pro(\Ring\mod)$ $V[[x]]$ is simply a countable product of copies of $V$. Moreover, the coproduct on $\Ring[\mathbb{G}_a^\wedge]$ is given by 
$$\Ring[[x]] \rightarrow \Ring[[x_1]] \widetilde{\otimes}_\Ring \Ring[[x_2]] \simeq \Ring[[x_1, x_2]], \quad \quad x \mapsto x_1 + x_2.$$
We note in addition that the Lie algebra $\mathfrak{g}_a$ of $\mathbb{G}_a^\wedge$ is canonically $\Ring$; i.e., a continuous derivation $$d: \Ring[[x]] \rightarrow e_* \Ring$$is the same data as an element $d(x) \in \Ring$.

As in Section \ref{ss:restrI}, we may form the category of comodules 
$$\on{Rep}(\mathbb{G}_a^\wedge) := \Ring[\mathbb{G}_a^\wedge]\comod(\Pro(\Ring\mod)).$$
In addition, we may consider the category of modules 
$$\Rep(\mathfrak{g}_a) := \mathfrak{g}_a\mod(\Pro(\Ring\mod));$$
explicitly, this is the category of pairs $(V, \Phi)$, where $V$ is an object of $\Pro(\Ring\mod)$ and $\Phi: V \rightarrow V$ is an endomorphism. Moreover, by Lemma \ref{l:derivaction}, we have a natural forgetful functor 
\begin{equation} \label{e:formalgroup} \on{Oblv}: \Rep(\mathbb{G}_a^\wedge) \rightarrow \Rep(\fg_a). \end{equation}

\begin{Lem} \label{l:formalga} The functor \eqref{e:formalgroup} is an equivalence, i.e., a coaction of $\mathbb{G}_a^\wedge$ on a pro-$\Ring$-module can be uniquely reconstructed from the action of its Lie algebra.
\end{Lem}

\begin{proof} An action of $\mathbb{G}_a^\wedge$ on $V \in \Pro(\Ring\mod)$ is the data of a map 
$$\on{coact}: V \rightarrow V[[x]]$$
satisfying the coassociativity constraint identifying two maps $V \rightarrow V[[x_1, x_2]]$. Recalling that $V[[x]]$ is a countable product of copies of $V$, we may equivalently regard the datum of the coaction as a formal sum 
$$\underset{k \geqslant 0} \Sigma \hspace{.5mm} \Phi_k x^k,$$
for some uniquely determined endomorphisms $\Phi_k: V \rightarrow V, k \geqslant 0$. I.e.,  $\Phi_k$ is given by the composition 
$$V \xrightarrow{\on{coact}} V[[x]] \rightarrow V \otimes_\Ring \Ring[x]/(x^{k+1}) \rightarrow V \otimes_\Ring \Ring \simeq V,$$
where the final map is induced by the map $\Ring[x]/(x^{k+1}) \rightarrow \Ring$ sending a polynomial to the coefficient of $x^k$. 

Coassociativity is then the identity of formal sums of endomorphisms
$$\underset{k \geqslant 0} \Sigma \hspace{.5mm}   \Phi_k (x_1 + x_2)^k = \underset{i,j \geqslant 0} \Sigma  \hspace{.5mm} (\Phi_i \circ \Phi_j) x_1^i x_2^j.$$
As $\Ring$ is a $\mathbb{Q}$-algebra, it follows that $\Phi_k = \frac{1}{k!} \Phi_1^k$, i.e., that the coaction is just the formal exponential of $\Phi_1$, and in particular is equivalent to the datum of $\Phi_1$. The claims of the lemma now easily follow. 
\end{proof}

\sss{} We are now ready to prove Proposition \ref{p:coordinvariancepD}. In particular, $\Ring$ again is a perfect field of characteristic greater than $h$.

\begin{proof}[Proof of Proposition \ref{p:coordinvariancepD}] The main idea of the proof is that we will pass to the filtered complete enveloping algebras of the vertex algebras appearing on both sides of Proposition \ref{p:autDdisc}. For ease of reading, we break the proof into steps. 


{\em Step 1.} Our first task is to address $\Aut(\pD)$-equivariance in the free field realization. For this, consider an algebra of chiral differential operators $\CDO_\kappa(F)$ as in Section \ref{SSS_CDO_construction}. 

We first claim that its enveloping algebra $$\widetilde{D}_\kappa(\Loop F) := \widetilde{U}(\CDO_\kappa(F))$$carries a natural action of $\Aut(\pD)$. Indeed, recall the universal property of the completed enveloping algebra, and the generators and relations for $\CDO_\kappa(F)$. Note in addition that $\hat{\ff}_\kappa \ltimes \Ring[\Loop F]$ is a pro-quasicoherent $\Ring$-space, in the sense of Section \ref{sss_proqcoh}, so that we may form its completed enveloping algebra. Combining these, it is straightforward to deduce $\widetilde{D}_\kappa(\Loop F)$ is the quotient of the enveloping algebra $\widetilde{U}( \hat{\ff}_\kappa \ltimes \Ring[\Loop F] )$ modulo the closed two-sided ideal generated by the relations (i) identifying the central element in $\hat{\ff}_\kappa$ with the multiplicative unit and (ii) identifying the product of two elements in $\widetilde{U}(\hat{\ff}_\kappa \ltimes \Ring[\Loop F])$ given by elements of $\Ring[\Loop F]$ with the corresponding element in $\Ring[\Loop F]$. 

From this presentation, which is compatibly defined on $\Sing$-points for every $\Ring$-algebra $\Sing$, it is clear that $\widetilde{D}_\kappa(\Loop F)$ inherits an action of $\Aut(\pD)$ from its natural action on $\hat{\ff}_\kappa \ltimes \Ring[\Loop F]$. Moreover, as in Section \ref{sss_jetscdoaction}, this action extends naturally to an action of 
$$\Aut(\pD) \ltimes (\Loop F^\ell \times \Loop F^r),$$
induced from the natural action of this group on $\hat{\ff}_\kappa \ltimes \Ring[\Loop F].$

{\em Step 2.} Let us specialize our discussion to the free field realization. Namely, if we pass to the completed enveloping algebras in Corollary \ref{c:freefieldcoordinv}, we obtain an $\Aut(\sD)$-equivariant map
$$\widetilde{U}_\kappa(\hat{\fg}) \rightarrow \widetilde{U}(\CDO(\tilde{N}) \otimes_\Ring V_{\ckappa'}(\check{\fh})) \simeq \widetilde{D}(\Loop \tilde{N}) \widetilde{\otimes}_\Ring \widetilde{U}( V_{\ckappa'}(\check{\fh})) =: \widetilde{D}(\Loop \tilde{N}) \widetilde{\otimes}_\Ring \widetilde{U}_{\ckappa'}(\check{\fh}),$$
where the tensor product factorization is by Lemma \ref{Lem:univ_env_tensor}. Let us denote the composition by
\begin{equation} \label{e:ffnocoord} \on{Ffr}_B: \widetilde{U}_\kappa(\hat{\fg}) \rightarrow \widetilde{D}(\Loop \tilde{N}) \widetilde{\otimes}_\Ring \widetilde{U}_{\ckappa'}(\check{\fh}). \end{equation}

We next note that both sides of \eqref{e:ffnocoord} in fact carry natural actions of $\Aut(\pD)$. For the left hand side, this was constructed in the proof of Lemma \ref{lem:autpdactonskm}. For the right hand side, the action of $\Aut(\pD)$ on $\widetilde{D}(\Loop \tilde{N})$ was constructed in the previous step, and as we have that 
$$\widetilde{U}_{\ckappa'}(\check{\fh}) \simeq \widetilde{U}_{\ckappa}(\check{\fh}) \overset{\Loop \check{H}} \times \Loop \omega^{\rho}$$
 this again carries a natural action of $\Aut(\pD)$. In particular, the right hand side of \eqref{e:ffnocoord} inherits a natural tensor product action of $\Aut(\pD)$. 

{\em Step 3.} We moreover claim that the map \eqref{e:ffnocoord} is $\Aut(\pD)$-equivariant. 
Indeed, as in the proof of Proposition \ref{p:coordinvarianceofcheat}, we first pass to integral forms, i.e., the corresponding objects defined over $\Spec \Z$, rather than $\Ring$, and with a variable level.  We then note that $\Aut(\pD)$ is again an inductive limit of schemes flat over $\Z$ under closed embeddings, and the relevant filtered complete algebras both admit a topological basis of left ideals for which the quotients are flat over $\Z$. Therefore, by  flatness we are reduced to the analogous assertion over $\mathbb{Q}$.

If we replace $\Aut(\pD)$ with its Lie algebra $\mathfrak{Witt}$, the fact that \eqref{e:ffnocoord} is equivariant for the action of $\mathfrak{Witt}$ follows from the fact that the map of vertex algebras $\on{ffr}_B$ exchanges the standard conformal vectors on both sides, cf. \cite[Proposition 6.2.2]{Frenkel_loop}.

As we already know the map is $\Aut(\D)$-equivariant, it suffices by Lemma \ref{l:autpd-}(2) to show it is $\Am$-equivariant. By Lemma \ref{l:cells}, it is enough check the equivariance when restricted to one-parameter formal subgroups as in Lemma \ref{l:formexp}. However, for those the equivariance follows from the assertion on Lie algebras by Lemma \ref{l:formalga}, so we are done.

{\em Step 4.} Note that as in Lemma \ref{l:twoHeiscommute}, the actions of $\Loop H$ and $\Loop \check{H}$ on $\widetilde{U}_\kappa(\hat{\fh}) \simeq \widetilde{U}_{\ckappa}({\check{\fh}})$ commute; in particular the twist $\widetilde{U}_{\ckappa'}(\check{\fh})$ retains an action of $\Aut(\pD) \ltimes \Loop H$ preserving its filtered complete subalgebra $\widehat{U}_{\ckappa'}(\check{\fh})$.

If we consider the $\Aut(\pD)$-equivariant composition 
$$\widehat{U}_\kappa(\hat{\fg})^{\Loop G} \rightarrow \widetilde{U}_\kappa(\hat{\fg}) \rightarrow \widetilde{D}(\Loop \tilde{N}) \widetilde{\otimes}_\Ring \widetilde{U}_{\ckappa'}(\check{\fh})$$
by Proposition \ref{Prop:affine_HC} this factors through the natural map 
$$\widehat{U}_\kappa(\hat{\fg})^{\Loop G} \rightarrow \widehat{U}_{\ckappa'}(\check{\fh})^{\Loop H} \hookrightarrow \widetilde{D}(\Loop \tilde{N}) \widetilde{\otimes}_\Ring \widetilde{U}_{\ckappa'}(\check{\fh}).$$As the second embedding is visibly $\Aut(\pD)$-equivariant, we deduce the first map is $\Aut(\pD)$-equivariant.

From here, the arguments of Sections \ref{ss:kconnections}-\ref{ss_hcff} apply {\em mutatis mutandis} and imply that we have a canonical $\Aut(\pD)$-equivariant isomorphism 
$$\widehat{U}_{\ckappa'}(\check{\fh})^{\Loop H} \simeq \Ring[\Conn_{\ckappa^p - \ckappa}({\pD}\fr, \omega^\rho)]_f.$$

As in Proposition \ref{p:miuradisc}, we similarly have the $\Aut(\pD)$-equivariant  Miura map on the punctured disc
$$\Ring[\Op_{\ckappa^p - \ckappa}({\pD}\fr)]_f \rightarrow \Ring[\Conn_{\ckappa^p - \ckappa}({\pD}\fr, \omega^\rho)]_f.$$

Putting things together, we obtain a diagram of $\Aut(\pD)$-equivariant maps
$$\widehat{U}_\kappa(\hat{\fg})^{\Loop G} \xrightarrow{\alpha}  \Ring[\Conn_{\ckappa^p - \ckappa}({\pD}\fr, \omega^\rho)]_f \xleftarrow{\beta} \Ring[\Op_{\ckappa^p - \ckappa}({\pD}\fr)]_f.$$
By Lemma \ref{Lem:HC_injectivity}, $\alpha$ is an embedding. By passing to filtered complete enveloping algebras in Theorem \ref{t:tungus}, it follows that $\beta$ is as well, and the two maps have the same image, which provides the desired $\Aut(\pD)$-equivariant identification. \end{proof}

\appendix

%
%

\end{document}